\documentclass[11pt]{article}
\usepackage{amsfonts,amsmath,amssymb}
\usepackage{hyperref}
\usepackage{graphicx}
\usepackage{color}
\usepackage{makeidx}

\makeindex

\newcommand{\detail}[1]{\par\noi{\bf [Proof detail\ }{#1}
\hfill{\bf ]}\par\noi\hspace{-4pt}}
\renewcommand{\detail}[1]{}

\newcommand{\dis}{\displaystyle}
\newcommand{\txt}{\textstyle}

\newcommand{\noi}{\noindent}

\newcommand{\halmos}{\rule{1ex}{1.4ex}}
\def \qed {\nopagebreak{\hspace*{\fill}$\halmos$\medskip}}
\newcommand{\med}{\medskip}

\newtheorem{theorem}{Theorem}[section]
\newtheorem{proposition}[theorem]{Proposition}
\newtheorem{corollary}[theorem]{Corollary}
\newtheorem{conjecture}[theorem]{Conjecture}
\newtheorem{lemma}[theorem]{Lemma}

\newtheorem{remark}[theorem]{Remark}
\newtheorem{defi}[theorem]{Definition}
\newcommand{\bt}{\begin{theorem}}
\newcommand{\et}{\end{theorem}}
\newcommand{\bl}{\begin{lemma}}
\newcommand{\el}{\end{lemma}}
\newcommand{\bp}{\begin{proposition}}
\newcommand{\ep}{\end{proposition}}
\newcommand{\bcor}{\begin{corollary}}
\newcommand{\ecor}{\end{corollary}}
\newcommand{\br}{\begin{remark}\rm}
\newcommand{\er}{\end{remark}}
\newcommand{\bcon}{\begin{conjecture}}
\newcommand{\econ}{\end{conjecture}}
\newcommand{\bd}{\begin{defi}}
\newcommand{\ed}{\end{defi}}

\newcommand{\be}{\begin{equation}}
\newcommand{\ee}{\end{equation}}
\newcommand{\ba}{\begin{array}}
\newcommand{\ea}{\end{array}}
\newcommand{\bc}{\be\begin{array}{r@{\,}c@{\,}l}}
\newcommand{\ec}{\end{array}\ee}

\newcommand{\al}{\alpha}
\newcommand{\bet}{\beta}
\newcommand{\ga}{\gamma}
\newcommand{\Ga}{\Gamma}
\newcommand{\de}{\delta}
\newcommand{\De}{\Delta}
\newcommand{\eps}{\varepsilon}

\newcommand{\la}{\lambda}
\newcommand{\La}{\Lambda}
\newcommand{\sig}{\sigma}
\newcommand{\Sig}{\Sigma}
\newcommand{\tet}{\theta}
\newcommand{\Tet}{\Theta}
\newcommand{\om}{\omega}
\newcommand{\Om}{\Omega}

\newcommand{\si}{\ensuremath{\sigma}}

\newcommand{\Ai}{{\cal A}}
\newcommand{\Bi}{{\cal B}}
\newcommand{\Ci}{{\cal C}}
\newcommand{\Di}{{\cal D}}

\newcommand{\Fi}{{\cal F}}
\newcommand{\Gi}{{\cal G}}
\newcommand{\Hi}{{\cal H}}
\newcommand{\Ii}{{\cal I}}
\newcommand{\Ki}{{\cal K}}
\newcommand{\Li}{{\cal L}}
\newcommand{\Mi}{{\cal M}}
\newcommand{\Ni}{{\cal N}}
\newcommand{\Oi}{{\cal O}}

\newcommand{\Ti}{{\cal T}}
\newcommand{\Ui}{{\cal U}}
\newcommand{\Vi}{{\cal V}}
\newcommand{\Wi}{{\cal W}}

\newcommand{\R}{{\mathbb R}}
\newcommand{\N}{{\mathbb N}}
\newcommand{\Z}{{\mathbb Z}}

\newcommand{\Q}{{\mathbb Q}}

\renewcommand{\P}{{\mathbb P}}
\newcommand{\E}{{\mathbb E}}

\newcommand{\li}{\langle}
\newcommand{\re}{\rangle}

\newcommand{\volgt}{\ensuremath{\Rightarrow}}
\newcommand{\up}{\uparrow}
\newcommand{\down}{\downarrow}
\newcommand{\sub}{\subset}
\newcommand{\beh}{\backslash}

\newcommand{\asto}[1]{\underset{{#1}\to\infty}{\longrightarrow}}
\newcommand{\astoo}[2]{\underset{{#1}\to{#2}}{\longrightarrow}}
\newcommand{\Asto}[1]{\underset{{#1}\to\infty}{\Longrightarrow}}
\newcommand{\symto}[1]{\underset{{#1}\to\infty}{\sim}}

\newcommand{\ti}{\tilde}

\newcommand{\ov}{\overline}

\newcommand{\subb}[2]{_{\ba{c}\scriptstyle{#1}\\[-.15cm]\scriptstyle{#2}\ea}}

\newcommand{\ffrac}[2]{{\textstyle\frac{{#1}}{{#2}}}}

\newcommand{\nab}{\nabla}

\newcommand{\di}{\mathrm{d}}
\newcommand{\half}{{[0,\infty)}}
\newcommand{\expo}{\mbox{\large\it e}}
\newcommand{\ex}[1]{\expo^{\,\textstyle{#1}}}

\newcommand{\Zev}{{\Z^2_{\rm even}}}
\newcommand{\Zod}{{\Z^2_{\rm odd}}}

\newcommand{\Rc}{R^2_{\rm c}}
\newcommand{\Wl}{\Wi^{\rm l}}
\newcommand{\Wr}{\Wi^{\rm r}}
\newcommand{\sign}{{\rm sign}}
\newcommand{\switch}{{\rm switch}}
\newcommand{\hop}{{\rm hop}}
\newcommand{\supp}{{\rm supp}}
\newcommand{\Qdis}{{\bf Q}}
\newcommand{\pad}{p}
\newcommand{\Kl}{{\it Closed}}
\newcommand{\Count}{{\it Count}}
\newcommand{\Yr}{{\rm Y}}

\newcommand{\varal}{\mathring\alpha}

\begin{document}

\makeatletter\@addtoreset{equation}{section}
\makeatother\def\theequation{\thesection.\arabic{equation}}

\renewcommand{\labelenumi}{{(\roman{enumi})}}

\title{\vspace{-3cm}Stochastic flows in the Brownian web and net}
\author{Emmanuel~Schertzer \and Rongfeng~Sun \and Jan~M.~Swart}
\date{May 28, 2013}

\maketitle\vspace{-30pt}

\begin{abstract}\noi
It is known that certain one-dimensional nearest-neighbor random walks in
i.i.d.\ random space-time environments have diffusive scaling limits. Here, in
the continuum limit, the random environment is represented by a `stochastic
flow of kernels', which is a collection of random kernels that can be loosely
interpreted as the transition probabilities of a Markov process in a random
environment. The theory of stochastic flows of kernels was first developed by
Le Jan and Raimond, who showed that each such flow is characterized by its
$n$-point motions. Our work focuses on a class of stochastic flows of kernels
with Brownian $n$-point motions which, after their inventors, will be called
Howitt-Warren flows.

Our main result gives a graphical construction of general Howitt-Warren flows,
where the underlying random environment takes on the form of a suitably marked
Brownian web. This extends earlier work of Howitt and Warren who showed that a
special case, the so-called `erosion flow', can be constructed from two
coupled `sticky Brownian webs'. Our construction for general Howitt-Warren
flows is based on a Poisson marking procedure developed by Newman, Ravishankar
and Schertzer for the Brownian web. Alternatively, we show that a special
subclass of the Howitt-Warren flows can be constructed as random flows of mass
in a Brownian net, introduced by Sun and Swart.

Using these constructions, we prove some new results for the Howitt-Warren
flows. In particular, we show that the kernels spread with a finite speed and
have a locally finite support at deterministic times if and only if the flow
is embeddable in a Brownian net. We show that the kernels are always purely
atomic at deterministic times, but, with the exception of the erosion flows,
exhibit random times when the kernels are purely non-atomic. We moreover prove
ergodic statements for a class of measure-valued processes induced by the
Howitt-Warren flows.

Our work also yields some new results in the theory of the Brownian web and
net. In particular, we prove several new results about coupled sticky
Brownian webs and about a natural coupling of a Brownian web with a Brownian
net. We also introduce a `finite graph representation' which gives a precise
description of how paths in the Brownian net move between deterministic times.
\end{abstract}

\noi
{\it Acknowledgement.} R.~Sun is supported by grants R-146-000-119-133 and
R-146-000-148-112 from the National University of Singapore. J.M.~Swart is
sponsored by GA\v CR grants 201/07/0237 and 201/09/1931.

{\small\setlength{\parskip}{-2pt}\tableofcontents}

\newpage

\section{Introduction}

\subsection{Overview}\label{S:intro}

In \cite{LR04AOP}, Le Jan and Raimond introduced the notion of a {\it
  stochastic flow of kernels}, which is a collection of random probability
kernels that can be loosely viewed as the transition kernels of a Markov
process in a random space-time environment, where restrictions of the
environment to disjoint time intervals are independent and the environment is
stationary in time. For suitable versions of such a stochastic flow of
kernels (when they exist), this loose interpretation is exact, see Definition~\ref{D:stochflow}
below and the remark following it. Given the environment, one
can sample $n$ independent copies of the Markov process and then average over
the environment. This defines the {\em $n$-point motion} for the flow, which
satisfies a natural consistency condition: namely, the marginal distribution
of any $k$ components of an $n$-point motion is necessarily a $k$-point
motion. A fundamental result of Le Jan and Raimond \cite{LR04AOP} shows that
conversely, any family of Feller processes that is consistent in this way
gives rise to an (essentially) unique stochastic flow of kernels.

As an example, in \cite{LR04PTRF}, the authors used Dirichlet forms to
construct a consistent family of reversible $n$-point motions on the circle,
which are $\al$-stable L\'evy processes with some form of sticky interaction
characterized by a real parameter $\tet$. In particular, for $\al=2$, these
are sticky Brownian motions. Subsequently, Howitt and Warren \cite{HW09a} used
a martingale problem approach to construct a much larger class of consistent
Feller processes on $\R$, which are Brownian motions with some form of sticky
interaction characterized by a finite measure $\nu$ on $[0,1]$. In particular,
if $\nu$ is a multiple of the Lebesgue measure, these are the sticky Brownian
motions of Le Jan and Raimond. From now on, and throughout this paper, we
specialize to the case of Browian underlying motions. By the general result of
Le Jan and Raimond mentioned above, the sticky Brownian motions of Le Jan and
Raimond, resp.\ Howitt and Warren, are the $n$-point motions of an
(essentially) unique stochastic flow of kernels on $\R$, which we call a {\em
  Le Jan-Raimond flow}, resp.\ {\em Howitt-Warren flow} (the former being a
special case of the latter). It has been shown in \cite{LL04,HW09a} that these
objects can be obtained as diffusive scaling limits of one-dimensional random
walks in i.i.d.\ random space-time environments.

The main goal of the present paper is to give a graphical construction of
Howitt-Warren flows that follows as closely as possible the discrete
construction of random walks in an i.i.d.\ random environment. In particular,
we want to make explicit what represents the random environment in the
continuum setting. The original construction of Howitt-Warren flows using
$n$-point motions does not tell us much about this.  In \cite{HW09b}, it was
shown that the Howitt-Warren flow with $\nu=\de_0+\de_1$, known as the {\em
  erosion flow}, can be constructed using two coupled Brownian webs, where one
Brownian web serves as the random space-time environment, while the
conditional law of the second Brownian web determines the stochastic flow of
kernels.

We will extend this construction to general Howitt-Warren flows, where in the
general case, the random environment consists of a Brownian web together with
a marked Poisson point process which is concentrated on the so-called points
of type $(1,2)$ of the Brownian web. A central tool in this construction is a
Poisson marking procedure invented by Newman, Ravishankar and Schertzer in
\cite{NRS10}. Of course, we also make extensive use of the theory of the
Brownian web developed in \cite{TW98,FINR04}. For a special subclass of the
Howitt-Warren flows, we will show that alternatively the random space-time
environment can be represented as a Brownian net, plus a countable collection
of i.i.d.\ marks attached to its so-called separation points. Here, we use the
theory of the Brownian net, which was developed in \cite{SS08} and
\cite{SSS09}.

Using our graphical construction, we prove a number of new properties for the
Howitt-Warren flows. In particular, we give necessary and sufficient
conditions in terms of the measure $\nu$ for the random kernels to spread with
finite speed, for their support to consist of isolated points at deterministic
times, and for the existence of random times when the kernels are non-atomic
(Theorems~\ref{T:speed}, \ref{T:supp} and \ref{T:atom} below). We moreover use
our construction to prove the existence of versions of Howitt-Warren flows
with nice regularity properties (Proposition~\ref{P:regul} below), in
particular, versions which can be interpreted as bona fide transition kernels
in a random space-time environment. Lastly, we study the invariant laws for
measure-valued processes associated with the Howitt-Warren flows
(Theorem~\ref{T:HIL}).

Our graphical construction of the Howitt-Warren flows is to a large extent
motivated by its discrete space-time counterpart, i.e., random walks in
i.i.d.\ random space-time environments on $\Z$. Many of our proofs will also
be based on discrete approximation. Therefore, in the rest of the
introduction, we will introduce a class of random walks in i.i.d.\ random
space-time environments and some related objects of interest, and sketch
heuristically how the Brownian web and the Brownian net will arise in the
representation of the random space-time environment for the Howitt-Warren
flows. An outline of the rest of the paper will be given at the end of the
introduction.

Incidentally, we note that random walks in i.i.d.\ random space-time environments
have been used in the physics literature to model the flow of stress in a
granular medium, called the $q$ model, see e.g.~\cite{LMY01, JM11} and the references therein.
The Howitt-Warren flows we consider are effectively scaling limits of so-called near-critical
$q$ models.

\subsection{Discrete Howitt-Warren flows}\label{S:dis}
\begin{figure}[htb]
\begin{center}
\includegraphics[width=13cm]{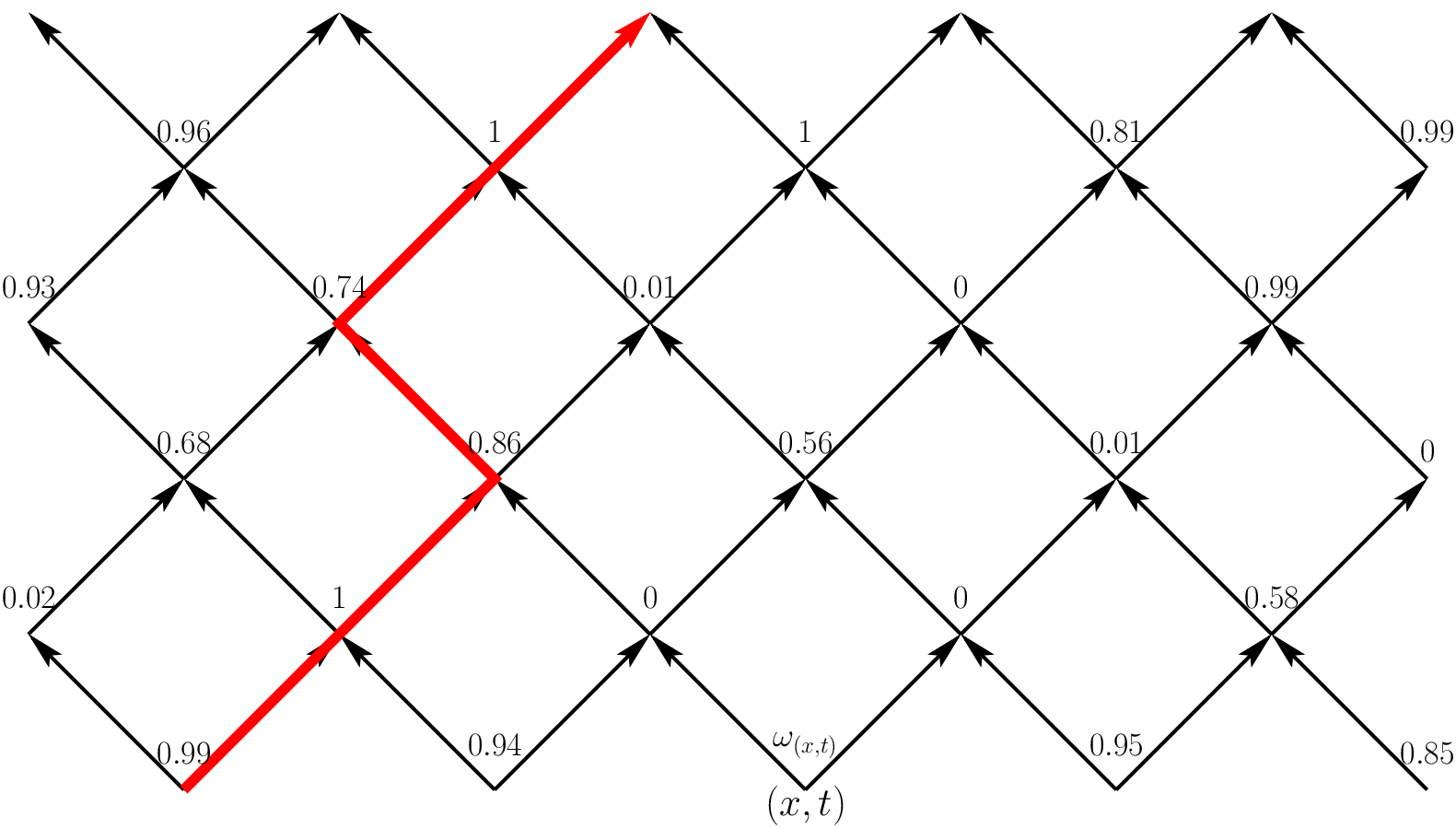}
\caption{Random walk on $\Zev$ in a random environment $\om$.}
\label{fig:disflow}
\end{center}
\end{figure}

\index{discrete!Howitt-Warren flow}
Let $\Zev:=\{(x,t):x,t\in\Z,\ x+t\mbox{ is even}\}$ be the even sublattice of
$\Z^2$. We interpret the first coordinate $x$ as space and the second
coordinate $t$ as time, which is plotted vertically in figures. Let
$\om:=(\om_z)_{z\in \Zev}$ be i.i.d.\ $[0,1]$-valued random variables
with common distribution $\mu$. We view $\om$ as a random space-time
environment for a random walk, such that conditional on
the environment $\om$, if the random walk is at time $t$ at the position $x$,
then in the next unit time step the walk jumps to $x+1$ with probability
$\om_{(x,t)}$ and to $x-1$ with the remaining probability $1-\om_{(x,t)}$ (see
Figure~\ref{fig:disflow}).

To formalize this, let $\P$ denote the law of the environment $\om$ and for
each $(x,s)\in\Zev$, let $\Qdis^\om_{(x,s)}$ denote the conditional law, given
the random environment $\om$, of the random walk in random environment
$X=(X(t))_{t\geq s}$ we have just described, started at time $s$ at position
$X(s)=x$. Since parts of the random environment belonging to different times
are independent, it is not hard to see that under the averaged (or `annealed')
law $\int\P(\di \om)\Qdis^\om_{(x,s)}$, the process $X$ is still a Markov chain,
which in each time step jumps to the right with probability $\int\mu(\di q)q$
and to the left with the remaining probability $\int\mu(\di q)(1-q)$. Note
that this is quite different from the usual random walk in random environment
(RWRE) where the randomness is fixed for all time, and the averaged motion no
longer has the Markov property.

We will be interested in three objects associated with the random walks in the
i.i.d.\ random space-time environment $\om$, namely: random transition
kernels, $n$-point motions, and a measure-valued process. The law of each of
these objects is uniquely characterized by $\mu$ and, conversely, uniquely
determines $\mu$.

First of all, the random environment $\om$ determines a family of random
transition probability kernels,
\be\label{disK}
K^\om_{s,t}(x,y):=\Qdis^\om_{(x,s)}\big[X(t)=y\big]
\qquad\big(s\leq t,\ (x,s),(y,t)\in\Zev\big),
\ee
which satisfy
\begin{enumerate}
\item $\dis\sum_{y:\,(y,t)\in\Zev}K^\om_{s,t}(x,y)K^\om_{t,u}(y,z)
  =K^\om_{s,u}(x,z)\qquad\big(s\leq t\leq u,\ (x,s),(z,u)\in\Zev)$.
\item For each $t_0<\cdots<t_n$, the random variables
  $(K^\om_{t_{i-1},t_i})_{i=1,\ldots,n}$ are independent.
\item $K^\om_{s,t}$ and $K^\om_{s+u,t+u}$ are equal in law for each
  $u\in\Z_{\rm even}:=\{2x:x\in\Z\}$.
\end{enumerate}
We call the collection of random probability kernels $(K^\om_{s,t})_{s\leq
  t}$ the {\em discrete Howitt-Warren flow}\index{flow!discrete} with {\em characteristic measure}
$\mu$. \index{characteristic measure!for discrete flow} Such a collection is a discrete time analogue of a stochastic flow of
kernels as introduced by Le Jan and Raimond in \cite{LR04AOP} (see
Definition~\ref{D:stochflow} below).

Next, given the environment $\om$, we can sample a collection of independent
random walks
\be\label{disnpoint}
(\vec X(t))_{t\geq 0}=\big(X_1(t),\ldots,X_n(t)\big)_{t\geq 0}
\ee
in the random environment $\om$, started at time zero from deterministic sites
$x_1,\ldots,x_n\in\Z_{\rm even}$, respectively. It is easy to see that under
the averaged law
\be\label{disaver}
\int\P(\di \om)\bigotimes_{i=1}^n\Qdis^\om_{(x_i,0)},
\ee
the process $\vec X=(\vec X(t))_{t\geq 0}$ is still a Markov chain, which we
call the {\em discrete $n$-point motion}. \index{n-point motion!discrete} Its transition probabilities are
given by
\be
\P^{(n)}_{s,t}(\vec x,\vec y)=\int\P(\di\om)\prod_{i=1}^nK^\om_{s,t}(x_i,y_i)
\qquad\big(s\leq t,\ (x_i,s),(y_i,t)\in\Zev,\ i=1,\ldots,n\big).
\ee
Note that these discrete $n$-point motions are consistent in the sense that
any $k$ coordinates of $\vec X$ are distributed as a discrete $k$-point
motion. Each coordinate $X_i$ is distributed as a nearest-neighbor random walk
thats makes jumps to the right with probability $\int\mu(\di q)q$. Because of
the spatial independence of the random environment, the coordinates evolve
independently when they are at different positions. To see that there is some
nontrivial interaction when they are at the same position, note that if $k+l$
coordinates are at position $x$ at time $t$, then the probability that in the
next time step the first $k$ coordinates jump to $x+1$ while the last $l$
coordinates jump to $x-1$ equals $\int\mu(\di q)q^k(1-q)^l$, which in general
does not factor into $(\int\mu(\di q)q)^k(\int\mu(\di q)(1-q))^l$. Note that
the law of $\om_{(0,0)}$ is uniquely determined by its moments, which are
in turn determined by the transition probabilities of the discrete $n$-point
motions (for each $n$).

Finally, based on the family of kernels $(K^\om_{s,t})_{s\leq t}$, we can
define a measure-valued process
\be\label{disrhoK}
\rho_t(x)=\sum_{y\in\Z_{\rm even}}\rho_0(y)K^\om_{0,t}(y,x)
\qquad\big(t\geq 0,\ (x,t)\in\Zev\big),
\ee
where $\rho_0$ is any locally finite initial measure on $\Z_{\rm
  even}$. Note that conditional on $\om$, the process
$\rho=(\rho_t)_{t\geq 0}$ evolves deterministically, with
\be\label{disrho}
\rho_{t+1}(x):=\om_{(x-1,t)}\rho_t(x-1)+(1-\om_{(x+1,t)})\rho_t(x+1)
\qquad\big((x,t+1)\in\Zev,\ t\geq 0\big).
\ee
Under the law $\P$, the process $\rho$ is a Markov chain, taking values
alternatively in the spaces of finite measures on $\Z_{\rm even}$ and $\Z_{\rm
  odd}:=\{2x+1:x\in\Z\}$. Note that (\ref{disrho}) says that in the time step
from $t$ to $t+1$, an $\om_{(x,t)}$-fraction of the mass at $x$ is sent to
$x+1$ and the rest is sent to $x-1$. Obviously, this dynamics preserves the
total mass. In particular, if $\rho_0$ is a probability measure, then $\rho_t$
is a probability measure for all $t\geq 0$. We call $\rho$ the {\em discrete
  Howitt-Warren process}. \index{Howitt-Warren!process!discrete}

We will be interested in the diffusive scaling limits of all these objects,
which will be (continuum) Howitt-Warren flows and their associated $n$-point
motions and measure-valued processes, respectively. Note that the discrete
Howitt-Warren flow $(K^\om_{s,t})_{s\leq t}$ determines the random environment
$\om$ a.s.\ uniquely. The law of $(K^\om_{s,t})_{s\leq t}$ is uniquely
determined by either the law of its $n$-point motions or the law of its
associated measure-valued process.

\subsection{Scaling limits of discrete Howitt-Warren flows}\label{S:heur}

We now recall from \cite{HW09a} the conditions under which the $n$-point
motions of a sequence of discrete Howitt-Warren flows converge to the
$n$-point motions of a (continuum) stochastic flow of kernels, which we call a
Howitt-Warren flow. We will then use discrete approximation to sketch
heuristically how such a Howitt-Warren flow can be constructed from a
Brownian web or net.

Let $(\eps_k)_{k\in\N}$ be positive constants tending to zero, and let
$(\mu_k)_{k\in\N}$ be probability laws on $[0,1]$ satisfying\footnote{We
  follow \cite{HW09a} in our definition of $\nu$. Many of our formulas,
  however, such as (\ref{betplusm}), (\ref{stickyparam}) or (\ref{Poisint})
  are more easily expressed in terms of $2\nu$ than in $\nu$. Loosely
  speaking, the reason for this is that in (\ref{mucon})~(ii), the weight
  function $q(1-q)$ arises from the fact that if $\al^1,\al^2$ are independent
  $\{-1,+1\}$-valued random variables with $\P[\al^i=+1]=q$ $(i=1,2)$, then
  $\P[\al^1\neq\al^2]=2q(1-q)$.}
\be\ba{rr@{\,}c@{\,}l}\label{mucon}
{\rm(i)}&\dis\eps_k^{-1}\int(2q-1)\mu_k(\di q)&\dis\asto{k}&\dis\bet,\\[10pt]
{\rm(ii)}&\dis\eps_k^{-1}q(1-q)\mu_k(\di q)&\dis\Asto{k}&\dis\nu(\di q)
\ec
for some $\bet\in\R$ and finite measure $\nu$ on $[0,1]$, where $\Rightarrow$
denotes weak convergence. Howitt and Warren \cite{HW09a}\footnote{Actually,
  the paper \cite{HW09a} considers a continuous-time analogue of the discrete
  $n$-point motions defined in Section~\ref{S:dis}, but their proof, with
  minor modifications, also works in the discrete time setting. In
  Appendix~\ref{A:HWMP} we present a similar, but somewhat simplified
  convergence proof.} proved that
under condition (\ref{mucon}), if we scale space by $\eps_k$ and time by
$\eps_k^2$, then the discrete $n$-point motions with characteristic measure
$\mu_k$ converge to a collection of Brownian motions with drift $\bet$ and
some form of sticky interaction characterized by the measure $\nu$. These
Brownian motions form a consistent family of Feller processes, hence by the
general result of Le Jan and Raimond mentioned in Section~\ref{S:intro}, they
are the $n$-point motions of some stochastic flow of kernels, which we call
the {\em Howitt-Warren flow}
with {\em drift} $\bet$ and {\em characteristic measure}
 $\nu$. The definition of Howitt-Warren flows and their $n$-point
motions will be given more precisely in Section~\ref{S:HWclass}.

Now let us use discrete approximation to explain heuristically how to
construct a Howitt-Warren flow based on a Brownian web or net. The
construction based on the Brownian net is conceptually easier, so we consider
this case first.

Let $\bet\in\R$ and let $\nu$ be a finite measure on $[0,1]$. Assuming, as we
must in this case, that $\int\frac{\nu(\di q)}{q(1-q)}<\infty$, we may define
a sequence of probability measures $\mu_k$ on $[0,1]$ by
\be\ba{l}
\dis\mu_k:=b\eps_k\bar\nu+\ffrac{1}{2}(1-(b+c)\eps_k)\de_0
+\ffrac{1}{2}(1-(b-c)\eps_k)\de_1\\[5pt]
\dis\quad\mbox{where}\quad
b:=\int\frac{\nu(\di q)}{q(1-q)},\quad
c:=\bet-\int(2q-1)\frac{\nu(\di q)}{q(1-q)},\quad
\bar\nu(\di q):=\frac{\nu(\di q)}{bq(1-q)}.
\ec
Then $\mu_k$ is a probability measure on $[0,1]$ for $k$ sufficiently large
(such that $1-(b+|c|)\eps_k\geq 0$), and the $\mu_k$ satisfy (\ref{mucon}).
Thus, when space is rescaled by $\eps_k$ and time by $\eps_k^2$, the
discrete Howitt-Warren flow with characteristic measure $\mu_k$ approximates a
Howitt-Warren flow with drift $\bet$ and characteristic measure $\nu$.

\begin{figure}[htb] 
\centering
\includegraphics[width=13cm]{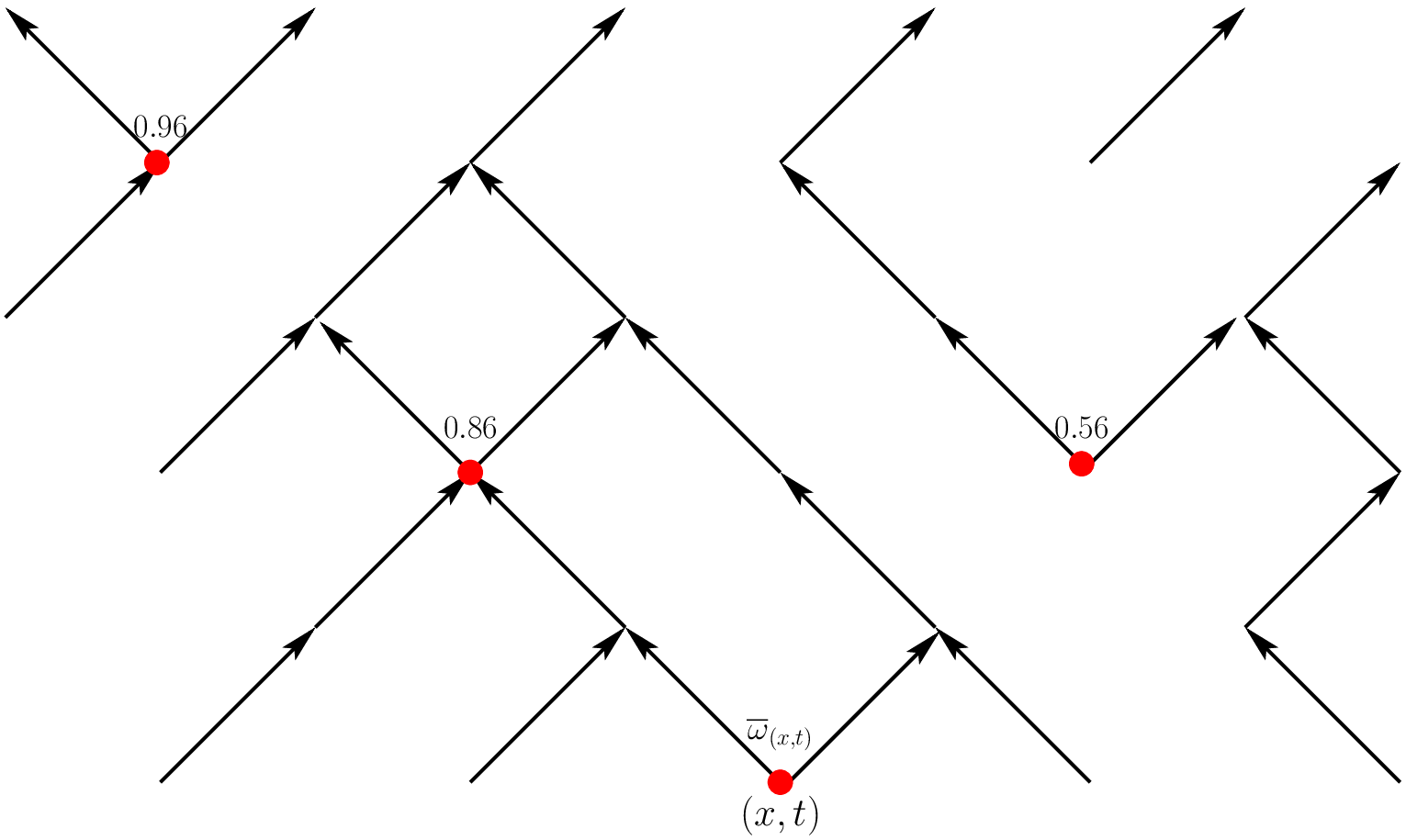}
\caption{Representation of the random environment $(\om^{\langle
    k\rangle }_z)_{z\in\Zev}$ in terms of a marked discrete net $(N^{\langle
    k\rangle},\bar\om^{\langle k\rangle})$.}
\label{fig:netenviron}
\end{figure}

Let $\om^{\langle k\rangle}:=(\om^{\langle k\rangle}_z)_{z\in \Zev}$ be
i.i.d.\ with common law $\mu_k$, which serves as the random environment for a
discrete Howitt-Warren flow with characteristic measure $\mu_k$. We observe
that for large $k$, most of the $\om^{\langle k\rangle}_z$ are either zero or
one. In view of this, it is convenient to alternatively encode $\om^{\langle
  k\rangle}$ as follows. For each $z=(x,t)\in \Zev$, if $\om^{\langle
  k\rangle}_z\in (0,1)$, then we call $z$ a separation point, set
$\bar\om^{\langle k\rangle}_z=\om^{\langle k\rangle}_z$, and we draw two
arrows from $z$, leading respectively to $(x\pm 1, t+1)$. When $\om^{\langle
  k\rangle}_z=0$, resp.~$1$, we draw a single arrow from $z$ to $(x-1, t+1)$,
resp.~$(x+1, t+1)$. Note that the collection of arrows $N^{\langle k\rangle}$
generates a branching-coalescing structure, called discrete net, on $\Zev$ (see
Figure~\ref{fig:netenviron}) and conditional on $N^{\langle k\rangle}$, the
$\bar\om^{\langle k\rangle}_z$ at separation points $z$ of $N^{\langle
  k\rangle}$ are independent with common law $\bar\nu$. Therefore the random
environment $\om^{\langle k\rangle}$ can be represented by the pair
$(N^{\langle k\rangle}, \bar\om^{\langle k\rangle})$, where a walk in such an
environment must navigate along $N^{\langle k\rangle}$, and when it encounters
a separation point $z$, it jumps either left or right with probability
$1-\bar\om^{\langle k\rangle}$, resp.~$\bar\om^{\langle k\rangle}$.

It turns out that the pair $(N^{\langle k\rangle}, \bar\om^{\langle
  k\rangle})$ has a meaningful diffusive scaling limit. In particular, if
space is scaled by $\eps_k$ and time by $\eps_k^2$, then $N^{\langle
  k\rangle}$ converges to a limiting branching-coalescing structure $\Ni$
called the {\em Brownian net}, the theory of which was developed in
\cite{SS08, SSS09}. In particular, the separation points of $N^{\langle
  k\rangle}$ have a continuum analogue, the so-called {\em separation points}
of $\Ni$, where incoming trajectories can continue along two groups of
outgoing trajectories. These separation points are dense in space and time,
but countable. Conditional on $\Ni$, we can then assign i.i.d.\ random
variables $\bar\om_z$ with common law $\bar\nu$ to the separation points of
$\Ni$. The pair $(\Ni, \bar\om)$ provides a representation for the random
space-time environment underlying the Howitt-Warren flow with drift $\beta$
and characteristic measure $\nu$. A random motion in such a random environment
must navigate along $\Ni$, and whenever it comes to a separation point $z$,
with probability $1-\bar\om_z$ resp.\ $\bar\om_z$, it continues along the left
resp.\ right of the two groups of outgoing trajectories in $\Ni$ at $z$. We
will recall the formal definition of the Brownian net and give a rigorous
construction of a random motion navigating in $\Ni$ in Section~\ref{S:net}.

\begin{figure}[htb] 
\centering
\includegraphics[width=13cm]{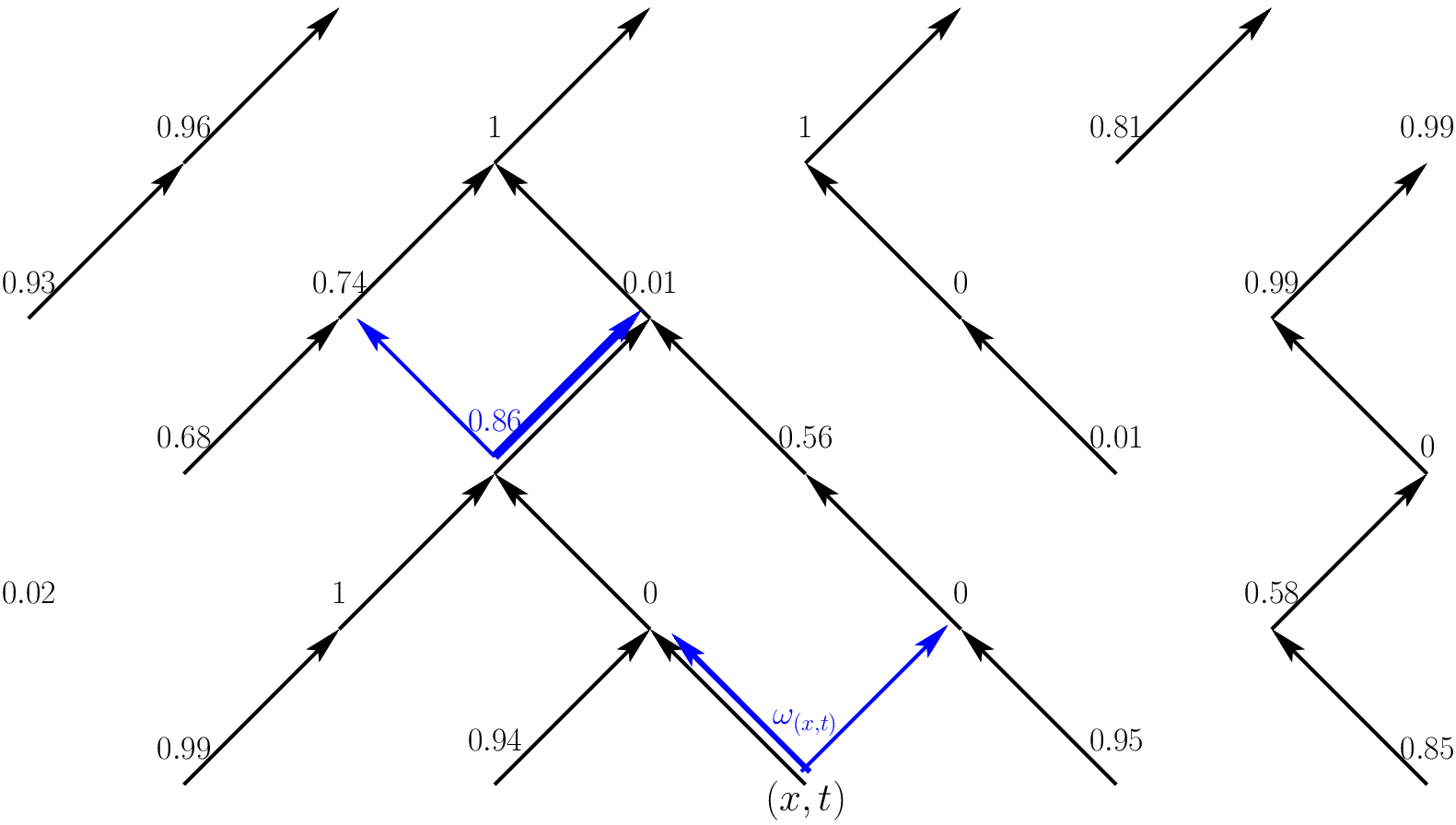}
\caption{Representation of the random environment $(\om^{\langle
    k\rangle }_z)_{z\in\Zev}$ in terms of a marked discrete web $(W^{\langle
    k\rangle}_0,\om^{\langle k\rangle})$.}
\label{fig:webenviron}
\end{figure}

We now consider Howitt-Warren flows whose characteristic measure is a general
finite measure $\nu$. Let $(\mu_k)_{k\in\N}$ satisfy (\ref{mucon}) and let
$\om^{\langle k\rangle}:=(\om^{\langle k\rangle}_z)_{z\in \Zev}$ be an
i.i.d.\ random space-time environment with common law $\mu_k$. Contrary to the
previous situation, it will now in general not be true that the most of the
$\om^{\langle k\rangle}_z$'s are either zero or one. Nevertheless, it is still
true that for large $k$, most of the $\om^{\langle k\rangle}_z$'s are either
close to zero or to one. To take advantage of this fact, conditional on
$\om^{\langle k\rangle}$, we sample independent $\{-1,+1\}$-valued random
variables $(\varal^{\li k\re}_z)_{z\in\Zev}$ such that $\varal^{\li k\re}_z=+1$
with probability $\om^{\li k\re}_z$. For each $z=(x,t)\in\Zev$, we draw an
arrow from $(x,t)$ to $(x+\varal^{\li k\re}_z,t+1)$. These arrows define a
coalescing structure $W^{\li k\re}_0$, called discrete web, on $\Zev$ (see
Figure~\ref{fig:webenviron}). Think of these arrows as assigning to each point
$z$ a preferred direction, which, in most cases, will be $+1$ if $\om^{\li
  k\re}_z$ is close to one and $-1$ if $\om^{\li k\re}_z$ is close to zero.

Now let us describe the joint law of $(\om^{\li k\re},\varal^{\li k\re})$
differently. First of all, if we forget about $\om^{\li k\re}$, then the
$(\varal^{\li k\re}_z)_{z\in\Zev}$ are just i.i.d.\ $\{-1,+1\}$-valued random
variables which take the value $+1$ with probability $\int q\mu_k(\di
q)$. Second, conditional on $\varal^{\li k\re}$, the random variables $(\om^{\li
  k\re}_z)_{z\in\Zev}$ are independent with distribution
\be\label{mulr}
\mu_k^{\rm l}:=\frac{(1-q)\mu_k({\rm d}q)}{\int (1-q)\mu_k({\rm d}q)},
\qquad\mbox{resp.}\qquad
\mu_k^{\rm r}:= \frac{q\mu_k({\rm d}q)}{\int q\mu_k({\rm d}q)}
\ee
depending on whether $\varal^{\li k\re}_z=-1$ resp.\ $+1$. Therefore, we can
alternatively construct our random space-time environment $\om^{\langle
  k\rangle}$ in such a way, that first we construct an i.i.d.\ collection
$\varal^{\li k\re}$ as above, and then conditional on $\varal^{\li k\re}$,
independently for each $z\in\Zev$, we choose $\om^{\langle k\rangle}_z$ with
law $\mu_k^{\rm l}$ if $\varal^{\li k\re}_z=-1$ and law $\mu_k^{\rm r}$ if
$\varal^{\li k\re}_z=+1$.

Let $W^{\langle k\rangle}_0$ denote the coalescing structure on $\Zev$ generated
by the arrows associated with $(\varal^{\li k\re}_z)_{z\in\Zev}$ (see
Figure~\ref{fig:webenviron}). Then $(W^{\langle k\rangle}_0,\om^{\langle
  k\rangle})$ gives an alternative representation of the random environment
$\om^{\langle k\rangle}$. A random walk in such an environment navigates in
such a way that whenever it comes to a point $z\in\Zev$, the walk jumps
to the right with probability $\om^{\langle k\rangle}_z$ and to the left
with the remaining probability. The important thing to observe is that if $k$
is large, then $\om^{\langle k\rangle}_z$ is with large probability close
to zero if $\varal^{\li k\re}=-1$ and close to one if $\varal^{\li k\re}=+1$. In
view of this, the random walk in the random environment $(W^{\langle
  k\rangle}_0,\om^{\langle k\rangle})$ will most of its time walk along paths
in $W^{\langle k\rangle}_0$.

It turns out that $(W^{\langle k\rangle}_0,\om^{\langle k\rangle})$ has a
meaningful diffusive scaling limit. In particular, if space is scaled by
$\eps_k$ and time by $\eps_k^2$, then the coalescing structure $W^{\langle
  k\rangle}_0$ converges to a limit $\Wi_0$ called the {\em Brownian web} (with
drift $\bet$), which loosely speaking is a collection of coalescing Brownian
motions starting from every point in space and time. These provide the default
paths a motion in the limiting random environment must follow. The
i.i.d.\ random variables $\om^{\langle k\rangle}_z$ turn out to converge to a
marked Poisson point process which is concentrated on so-called points of type
$(1,2)$ in $\Wi_0$, which are points where there is one incoming path and two
outgoing paths. These points are divided into points of type $(1,2)_{\rm l}$
and $(1,2)_{\rm r}$, depending on whether the incoming path continues on the
left or right. A random motion in such an environment follows paths in $\Wi_0$
by default, but whenever it comes to a marked point $z$ of type $(1,2)$, it
continues along the left resp.\ right outgoing path with probability $1-\om_z$
resp.\ $\om_z$.\footnote{In fact, this is not the full story, but
  describes only what happens if the measure $\nu$ from (\ref{mucon}) is
  concentrated on $(0,1)$. If $\nu$ puts mass on the boundary of $[0,1]$, then
  a random motion in $\Wi_0$ will in addition, with a certain Poisson rate,
  decide to follow the non-default outgoing path at some unmarked points of
  type $(1,2)$. In particular, this is the only mechanism if $\nu$ is
  concentrated on $\{0,1\}$, i.e., for so-called {\em erosion flows}.}
We will give the rigorous construction in Section~\ref{S:web}. The procedure
of marking a Poisson set of points of type $(1,2)$ that we need here was first
developed by Newman, Ravishankar and Schertzer in \cite{NRS10}, who used it
(among other things) to give an alternative construction of the Brownian net.

\subsection{Outline and discussion}

The rest of the paper is organized as follows. Sections
\ref{S:HWclass}--\ref{S:net} provide an extended introduction where we
rigorously state our results. In Section~\ref{S:HWclass}, we recall the notion
of a stochastic flow of kernels, first introduced in \cite{LR04AOP}, and
Howitt and Warren's \cite{HW09a} sticky Brownian motions, to give a rigorous
definition of Howitt-Warren flows. We then state out main results for these
Howitt-Warren flows, including properties for the kernels and results for the
associated measure-valued processes. In Sections~\ref{S:web} and \ref{S:net}
we make the heuristics from Section~\ref{S:heur} rigorous. In
Section~\ref{S:web}, in particular in Theorem~\ref{T:HWconst}, we present our
construction of Howitt-Warren flows based on a `reference' Brownian web with a
Poisson marking, which is the main result of this paper. Along the way, we
will recall the necessary background on the Brownian web. In
Section~\ref{S:net}, we show that a special subclass of the Howitt-Warren
flows can be constructed alternatively as flows of mass in the Brownian
net. Along the way, we will recall the necessary background on the Brownian
net and establish some new results on a coupling between a Brownian web and a
Brownian net. Sections \ref{S:outline}--\ref{S:inf} are devoted to proofs. In
particular, we refer to Section~\ref{S:outline} for an outline of the
proofs. The paper concludes with a number of appendices and a list of
notation.

Our work leaves several open problems. One question, for example, is how to
characterize the measure-valued processes associated with a Howitt-Warren flow
(see (\ref{rho}) below) by means of a well-posed martingale problem. Other
questions (martingale problem formulation, path properties) refer to the duals
(in the sense of linear systems duality) of these measure-valued processes,
introduced in (\ref{zetat}) below, which we have not investigated in much
detail.

Moving away from the Brownian case, we note that it is an open problem whether
our methods can be generalized to other stochastic flows of kernels than those
introduced by Howitt and Warren. In particular, this applies to the stochastic
flows of kernels with $\al$-stable L\'evy $n$-point motions introduced in
\cite{LR04PTRF} for $1<\al<2$. A first step on this road would be the
construction of an $\al$-stable L\'evy web which should generalize the
presently known Brownian web. Some first steps in this direction have recently
been taken in \cite{EMS13}.

\section{Results for Howitt-Warren flows}\label{S:HWclass}

In this section, we recall the notion of a stochastic flow of kernels, define
the Howitt-Warren flows, and state our results on these Howitt-Warren flows,
which include almost sure path properties and ergodic theorems for the
associated measure-valued processes. The proofs of these results are based on
our graphical construction of the Howitt-Warren flows, which we postpone to
Sections~\ref{S:web}--\ref{S:net} due to the extensive background we need to
recall.

\subsection{Stochastic flows of kernels}\label{S:kernel}

In \cite{LR04AOP}, Le Jan and Raimond developed a theory of {\em stochastic
  flows of kernels}, which may admit versions that can be interpreted as the
random transition probability kernels of a Markov process in a stationary
random space-time environment. The notion of a stochastic flow of kernels
generalizes the usual notion of a stochastic flow, which is a family of random
mappings $(\phi^\om_{s,t})_{s\leq t}$ from a space $E$ to itself. In the
special case that all kernels are delta-measures, a stochastic flow of kernels
reduces to a stochastic flow in the usual sense of the word.

Since stochastic flows of kernels play a central role in our work, we take
some time to recall their defintion. For any Polish space $E$, we let $\Bi(E)$
denote the Borel \si-field on $E$ and write $\Mi(E)$ and $\Mi_1(E)$ for the
spaces of finite measures and probability measures on $E$, respectively,
equipped with the topology of weak convergence and the associated Borel
\si-field. By definition, a probability kernel on $E$ is a function
$K:E\times\Bi(E)\to\R$ such that the map $x\mapsto K(x,\,\cdot\,)$ from $E$ to
$\Mi_1(E)$ is measurable. By a {\em random probability kernel}, defined on
some probability space $(\Om,\Fi,\P)$, we will mean a function $K:\Om\times
E\times\Bi(E)\to\R$ such that the map $(\om,x)\mapsto K^\om(x,\,\cdot\,)$ from
$\Om\times E$ to $\Mi_1(E)$ is measurable. We say that two random probability
kernels $K,K'$ are equal in finite dimensional distributions if for each
$x_1,\ldots,x_n\in E$, the $n$-tuple of random probability measures
$\big(K(x_1,\,\cdot\,),\ldots,K(x_n,\,\cdot\,)\big)$ is equally distributed
with $\big(K'(x_1,\,\cdot\,),\ldots,K'(x_n,\,\cdot\,)\big)$. We say that two
or more random probability kernels are independent if their finite-dimensional
distributions are independent.

\begin{defi}{\bf(Stochastic flow of kernels)}\label{D:stochflow}
\index{stochastic flow of kernels}\index{flow!of kernels}
A {\em stochastic flow of kernels} on $E$ is a collection $(K_{s,t})_{s\leq
  t}$ of random probability kernels on $E$ such that
\footnote{For simplicity, we have omitted two regularity conditions on
  $(K_{s,t})_{s\leq t}$ from the original definition in
  \cite[Def.~2.3]{LR04AOP}, which are some form of weak continuity of
  $K_{s,t}(x,\cdot)$ in $x,s$ and $t$. It is shown in that paper that a
  stochastic flow of kernels on a compact metric space $E$ satisfies these
  regularity conditions if and only if it arises from a consistent family of
  Feller processes.}
\begin{itemize}
\item[{\rm(i)}] For all $s\leq t\leq u$ and $x\in E$,
a.s.\ $K_{s,s}(x,A)=\de_x(A)$ and
$\dis\int_EK_{s,t}(x,\di y)K_{t,u}(y,A)=K_{s,u}(x,A)$ for all $A\in\Bi(E)$.
\item[{\rm(ii)}] For each $t_0<\cdots<t_n$, the random probability kernels
  $(K_{t_{i-1},t_i})_{i=1,\ldots,n}$ are independent.
\item[{\rm(iii)}] $K_{s,t}$ and $K_{s+u,t+u}$ are equal in finite-dimensional
  distributions for each real $s\leq t$ and~$u$.
\end{itemize}
The finite-dimensional distributions of a stochastic flow of
kernels are the laws of $n$-tuples of random probability measures of the
form $\big(K_{s_1,t_1}(x_1,\,\cdot\,),\ldots,K_{s_n,t_n}(x_n,\,\cdot\,)\big)$,
where $x_i\in E$ and $s_i\leq t_i$, $i=1,\ldots,n$.
\end{defi}
{\bf Remark.} If the random set of probability $1$ on which Definition
\ref{D:stochflow}~(i) holds can be chosen uniformly for all $s\leq t\leq u$
and $x\in E$, then we can interpret $(K_{s,t})_{s\leq t}$ as bona fide
transition kernels of a random motion in random environment. For
the stochastic flows of kernels we are interested in, we will prove the
existence of a version of $K$ which satisfies this property (see
Proposition~\ref{P:conpath} below). To the best of our knowledge, it is not
known whether such a version always exists for general stochastic flows of
kernels, even if we restrict ourselves to those defined by a consistent family
of Feller processes.\med

If $(K_{s,t})_{s\leq t}$ is a stochastic flow of kernels and $\rho_0$ is a
finite measure on $E$, then setting
\be\label{rho}
\rho_t(\di y):=\int\rho_0(\di x)K_{0,t}(x,\di y)
\qquad(t\geq 0)
\ee
defines an $\Mi(E)$-valued Markov process $(\rho_t)_{t\geq 0}$. Moreover,
setting
\be\label{npoint}
P^{(n)}_{t-s}(\vec x,\di\vec y)
:=\E\big[K_{s,t}(x_1,\di y_1)\cdots K_{s,t}(x_n,\di y_n)\big]
\qquad(\vec x\in E^n,\ s\leq t)
\ee
defines a Markov transition function on $E^n$. We call the Markov process with
these transition probabilities the {\em $n$-point motion} \index{n-point motion} associated with the
stochastic flow of kernels $(K_{s,t})_{s\leq t}$. We observe that the
$n$-point motions of a stochastic flow of kernels satisfy a natural
consistency condition: namely, the marginal distribution of any $k$ components
of an $n$-point motion is necessarily a $k$-point motion for the flow. A
fundamental result of Le Jan and Raimond \cite[Thm~2.1]{LR04AOP} states that
conversely, any consistent family of Feller processes on a locally compact
space $E$ gives rise to a stochastic flow of kernels on $E$ which is unique in
finite-dimensional distributions.\footnote{In fact, \cite[Thm~2.1]{LR04AOP} is
stated only for compact metrizable spaces, but the extension to locally
compact $E$ is straightforward using the one-point compactification of $E$.}

\subsection{Howitt-Warren flows}\label{S:flow}

As will be proved in Proposition~\ref{P:nconv} below, under the condition
(\ref{mucon}), if space and time are rescaled respectively by $\eps_k$ and
$\eps_k^2$, then the $n$-point motions associated with the discrete
Howitt-Warren flow introduced in Section~\ref{S:dis} with characteristic
measure $\mu_k$ converge to a collection of Brownian motions with drift $\bet$
and some form of sticky interaction characterized by the measure $\nu$. These
Brownian motions solve a well-posed martingale problem, which we formulate
now.

Let $\bet\in\R$, $\nu$ a finite measure on $[0,1]$, and define constants
$(\bet_+(m))_{m\geq 1}$ by
\bc\label{betplusm}
\dis\bet_+(1)&:=&\dis\bet\quad\mbox{and}\\[5pt]
\dis\bet_+(m)&:=&\dis\bet+2\int\nu(\di q)\sum_{k=0}^{m-2}(1-q)^k
\qquad(m\geq 2).
\ec
We note that in terms of these constants, (\ref{mucon}) is equivalent to
\be\label{mucon2a}
\eps_k^{-1}\int\big(1-2(1-q)^m\big)\mu_k(\di q)\asto{k}\bet_+(m)\qquad(m\geq 1).
\ee
For $\emptyset\neq\De\sub\{1,\ldots,n\}$, we define
\be\label{fgDe}
f_\De(\vec x):=\max_{i\in\De}x_i\qquad\mbox{and}\qquad
g_\De(\vec x):=\big|\{i\in\De:x_i=f_\De(\vec x)\}\big|
\qquad(\vec x\in\R^n),
\ee
where $|\,\cdot\,|$ denotes the cardinality of a set.

The martingale problem we are about to formulate was invented by Howitt and
Warren \cite{HW09a}. We have reformulated their definition in terms of the
functions $f_\De$ in (\ref{fgDe}), which form a basis of the vector space of
test functions used in \cite[Def~2.1]{HW09a} (see Appendix~\ref{A:HWMP} for a
proof). This greatly simplifies the statement of the martingale problem and
also facilitates our proof of the convergence of the $n$-point motions of
discrete Howitt-Warren flows.

\index{Howitt-Warren!martingale problem}
\begin{defi}{\bf(Howitt-Warren martingale problem)}\label{D:HWMP2}
We say that an $\R^n$-valued process $\vec X=(\vec X(t))_{t\geq 0}$ solves the
Howitt-Warren martingale problem with drift $\bet$ and characteristic
measure $\nu$ \index{characteristic measure} if $\vec X$ is a continuous, square-integrable semimartingale,
the covariance process between $X_i$ and $X_j$ is given by
\be\label{MP1b}
\langle X_i,X_j\rangle(t)=\int_0^t1_{\{X_i(s)=X_j(s)\}}\di s
\qquad(t\geq 0,\ i,j=1,\ldots,n),
\ee
and, for each nonempty $\De\sub\{1,\ldots,n\}$,
\be\label{MP2b}
f_\De\big(\vec X(t)\big)-\int_0^t\bet_+\big(g_\De(\vec X(s))\big)\di s
\ee
is a martingale with respect to the filtration generated by $\vec X$.
\end{defi}
{\bf Remark.} We could have stated a similar martingale problem where instead
of the functions $f_\De$ from (\ref{fgDe}) we use the functions $\ti
f_\De(x):=\min_{i\in\De}x_i$ and we replace the $\bet_+(m)$ defined in
(\ref{betplusm}) by
\be\label{betminm}
\bet_-(1):=\bet\quad\mbox{and}\quad
\bet_-(m):=\bet-2\int\nu(\di q)\sum_{k=0}^{m-2}q^k\qquad(m\geq 2).
\ee
It is not hard to prove that both martingale problems are equivalent.\med

\noi
{\bf Remark.} When $n=2$, condition (\ref{MP2b}) is
equivalent to the condition that
\be\label{thetacouple}
X_1(t)-\beta t, \quad X_2(t)-\beta t, \quad
|X_1(t)-X_2(t)| - 4\nu([0,1])\int_0^t 1_{\{X_1(s)=X_2(s)\}}\di s
\ee
are martingales. In~\cite{HW09a}, such $(X_1,X_2)$ are called $\theta$-coupled
Brownian motions, with $\theta=2\nu([0,1])$. In this case, $X_1(t)-X_2(t)$ is
a Brownian motion with stickiness at the
origin.\index{Brownian motion!$\theta$-coupled}
\index{Brownian motion!sticky} Such a process can be constructed by
time-changing a standard Brownian motion
in such a way that it spends positive Lebesgue time at the origin. More
generally, for solutions to the Howitt-Warren martingale problem started in
$X_1(0)=\cdots=X_n(0)$, the set of times such that
$X_1(t)=X_2(t)=\cdots=X_n(t)$ is a nowhere dense set with positive Lebesgue
measure. The measure $\nu$ then determines a two-parameter family of constants
$(\theta(k,l))_{k,l\geq 1}$ (see formula (\ref{tetnu}) in the Appendix), which
can be interpreted as the rate, in a certain excursion theoretic sense, at
which $(X_1,\cdots, X_n)$ split into two groups, $(X_1,\cdots, X_k)$ and
$(X_{k+1},\cdots, X_{k+l})$, with $k+l=n$.\med

Howitt and Warren \cite[Prop.~8.1]{HW09a} proved that their martingale problem
is well-posed and its solutions form a consistent family of Feller
processes. Therefore, by the already mentioned result of Le Jan and Raimond
\cite[Thm~2.1]{LR04AOP}, there exists a stochastic flow of kernels
$(K_{s,t})_{s\leq t}$ on $\R$, unique in finite-dimensional distributions,
such that the $n$-point motions of $(K_{s,t})_{s\leq t}$ (in the sense of
(\ref{npoint})) are given by the unique solutions of the Howitt-Warren
martingale problem. \index{Howitt-Warren!flow} \index{flow!Howitt-Warren} We call this stochastic flow of kernels the {\em
  Howitt-Warren flow} with {\em drift} $\bet$ and {\em characteristic measure}
$\nu$. It can be shown that Howitt-Warren flows are the diffusive
scaling limits, in the sense of weak convergence of finite dimensional
distributions, of the discrete Howitt-Warren flows with characteristic
measures $\mu_k$ satisfying (\ref{mucon}). (Indeed, this is a direct
consequence of Proposition~\ref{P:nconv} below on the convergence of $n$-point
motions.)

We will show that it is possible to construct versions of Howitt-Warren flows
which are bona fide transition probability kernels of a random motion in a
random space-time environment, and the kernels have `regular' parameter
dependence.
\bp{\bf(Regular parameter dependence)}\label{P:conpath}
For each $\bet\in\R$ and finite measure $\nu$ on $[0,1]$,
there exists a version of the Howitt-Warren flow $(K_{s,t})_{s\leq t}$ with
drift $\bet$ and characteristic measure $\nu$ such that in addition to the
properties (i)--(iii) from Definition~\ref{D:stochflow}:
\begin{itemize}
\item[{\rm(i)'}] A.s., $\dis\int_EK_{s,t}(x,\di y)K_{t,u}(y,A)=K_{s,u}(x,A)$
  for all $s\leq t\leq u$, $x\in E$ and $A\in\Bi(E)$.
\item[{\rm(iv)}] A.s., the map $t\mapsto K_{s,t}(x,\,\cdot\,)$ from
  $[s,\infty)$ to $\Mi_1(\R)$ is continuous for each $(s,x)\in\R^2$.
\end{itemize}
\ep

When the characteristic measure $\nu=0$, solutions to the Howitt-Warren
martingale problem are coalescing Brownian motions. In this case, the
associated stochastic flow of kernels is a stochastic flow (in the usual
sense), which is known as the {\em Arratia flow}. \index{flow!Arratia} In the
special case that $\bet=0$ and $\nu$ is Lebesgue measure, the Howitt-Warren
flow and its $n$-point motions are reversible. This stochastic flow of kernels
has been constructed before (on the unit circle instead of $\R$) by Le Jan and
Raimond in \cite{LR04PTRF} using Dirichlet forms. We will call any stochastic
flow of kernels with $\nu({\rm d}x)=c\,{\rm d}x$ for some $c>0$ a {\it Le
  Jan-Raimond flow}. \index{flow!Le Jan-Raimond} In~\cite{HW09b}, Howitt and
Warren constructed a stochastic flow of kernels with $\bet=0$ and
$\nu=\ffrac{1}{2}(\de_0+\de_1)$, which they called the {\em erosion flow}. In
this paper, we will call this flow the {\em symmetric erosion flow} and more
generally, we will say that a Howitt-Warren flow is an {\em erosion flow} if
$\nu=c_0\delta_0+c_1\delta_1$ with $c_0+c_1>0$. \index{flow!erosion} The paper
\cite{HW09b} gives an explicit construction of the symmetric erosion flow
based on coupled Brownian webs. Their construction can actually be extended to
any erosion flow and can be seen as a precursor and special case of our
construction of general Howitt-Warren flows in this paper.

\subsection{Path properties}\label{S:path}

\begin{figure}[htb]
\begin{center}
\includegraphics[width=7cm]{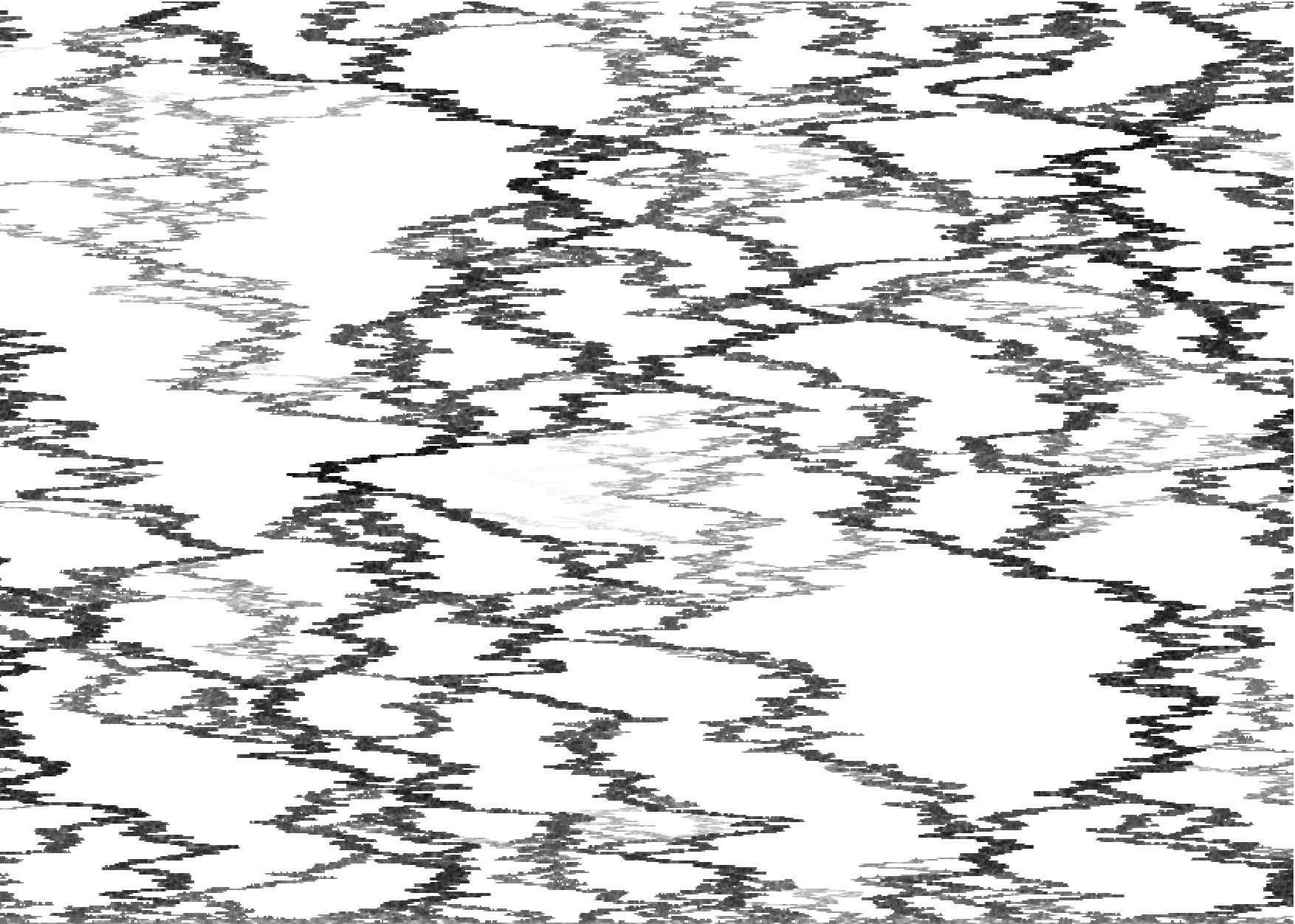} \quad
\includegraphics[width=7cm]{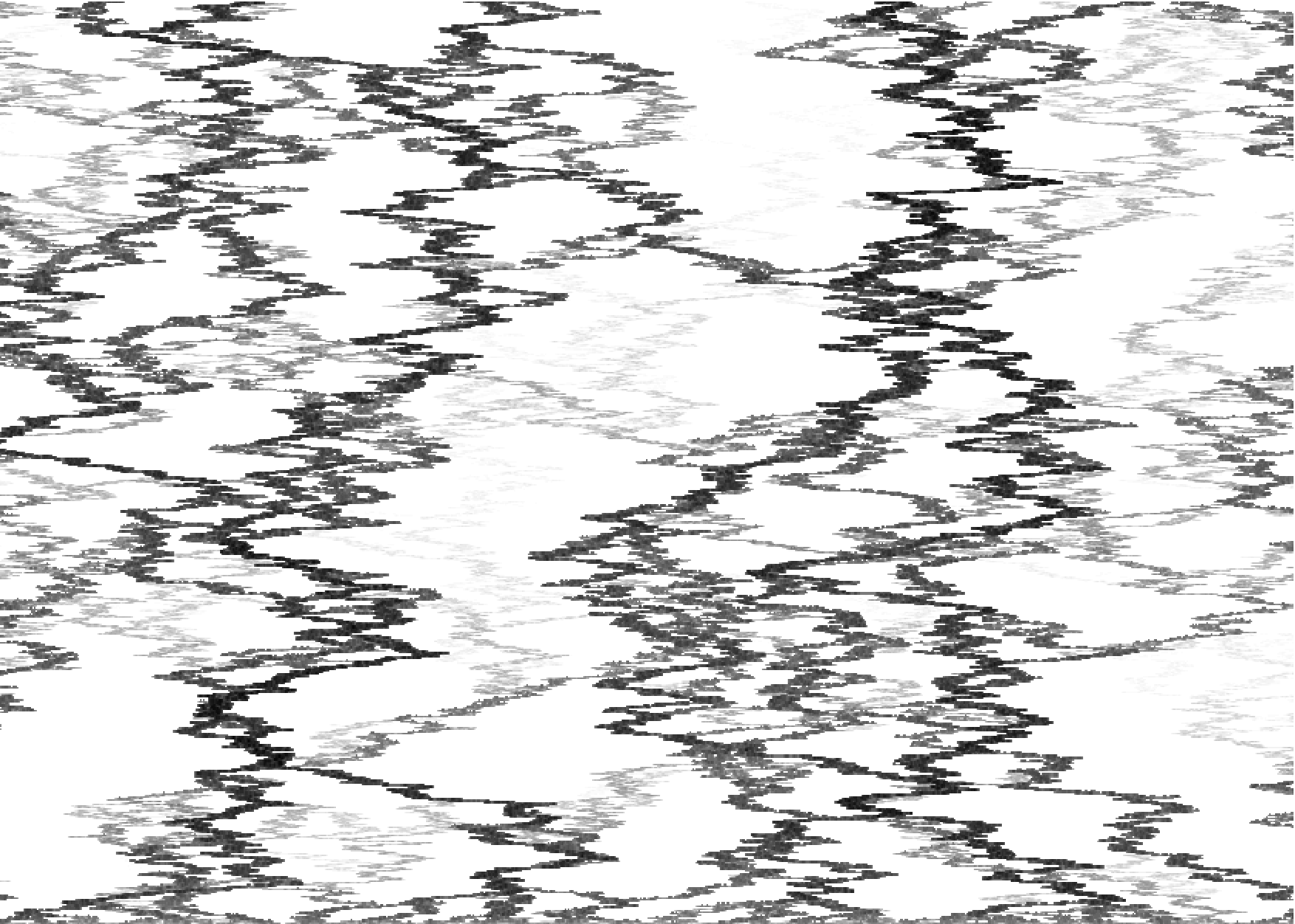}

\vspace{15pt}

\includegraphics[width=7cm]{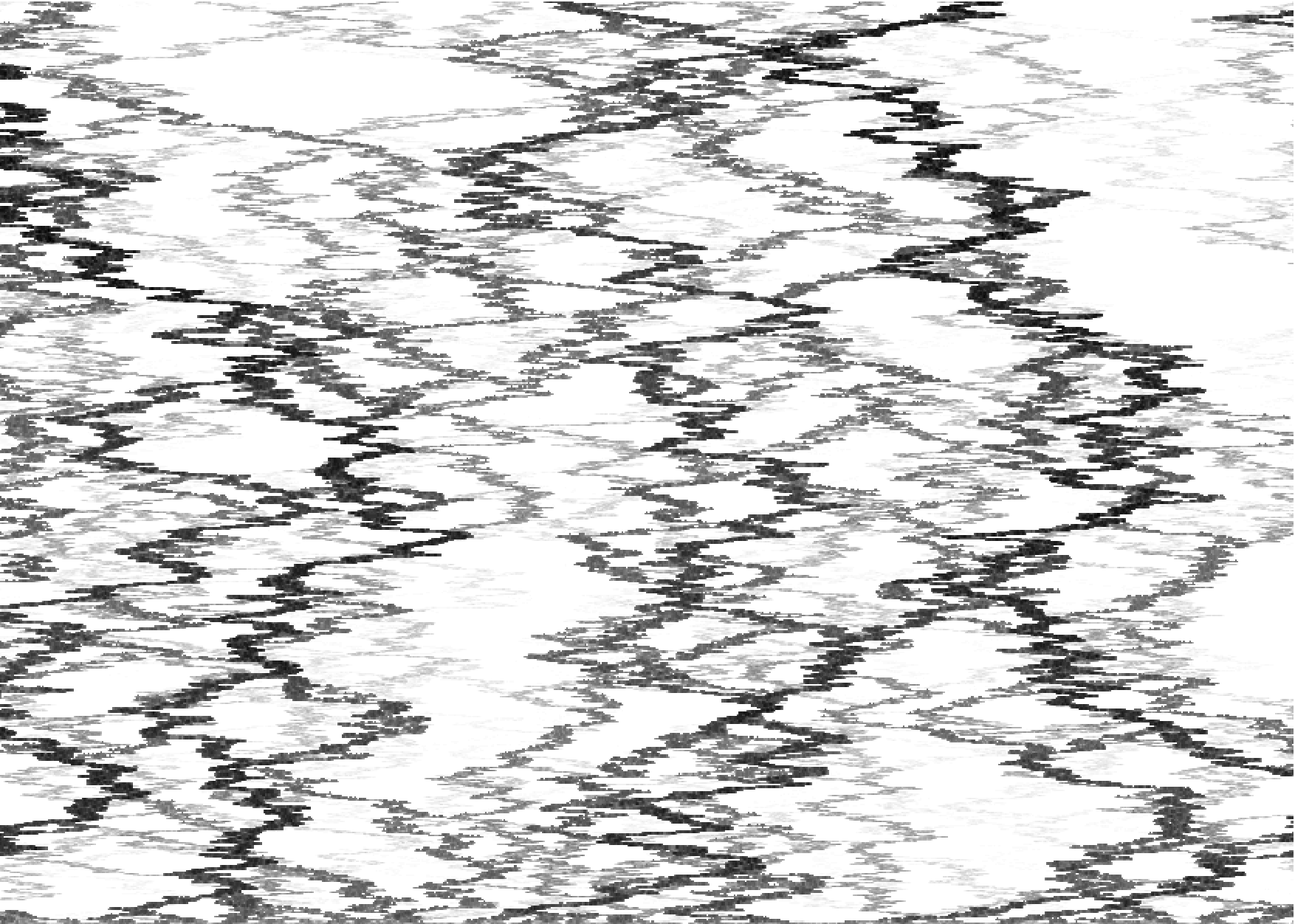}  \quad
\includegraphics[width=7cm]{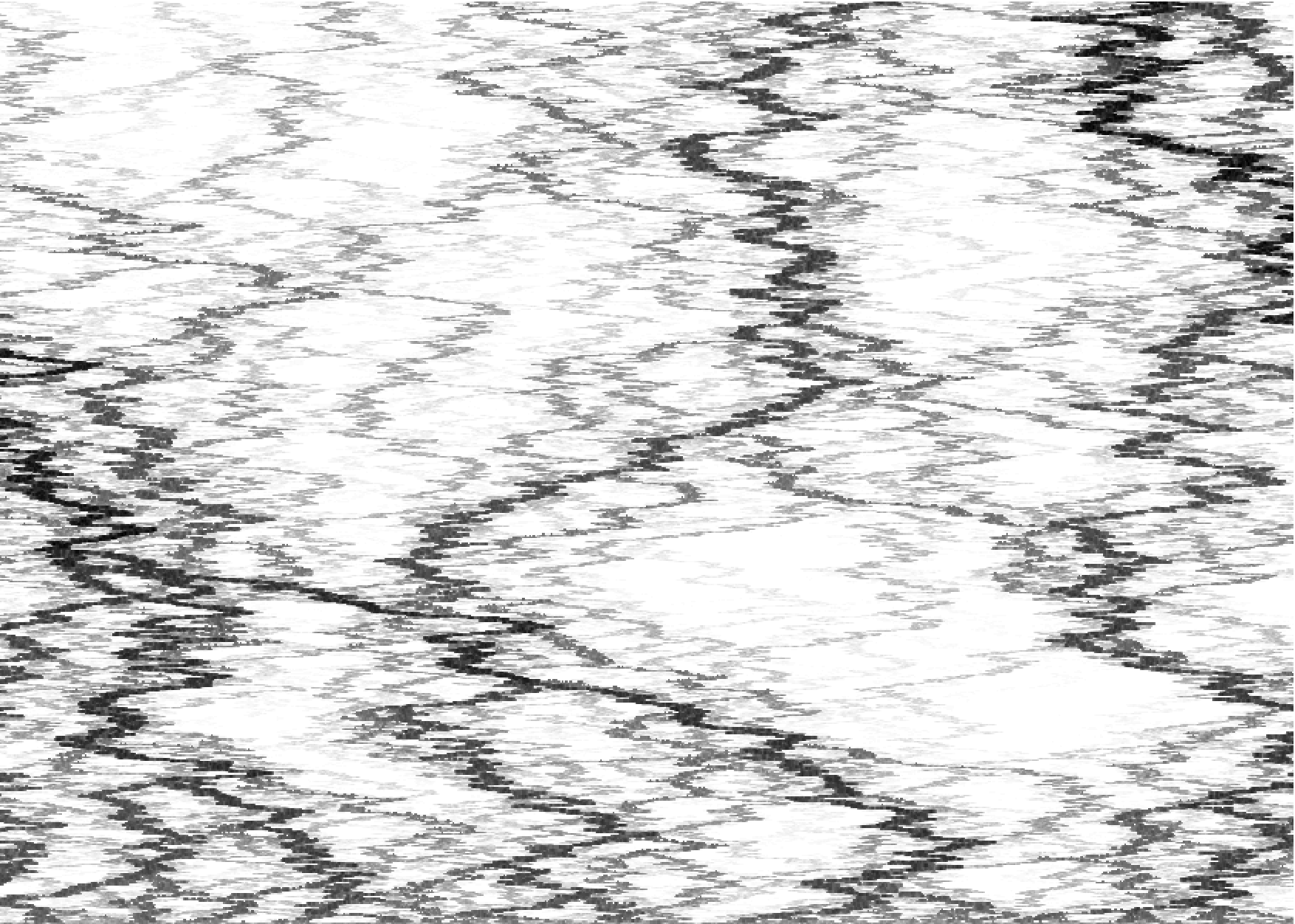}
\caption{Four examples of Howitt-Warren flows. All examples have drift
  $\bet=0$ and stickiness parameter $\tet=2$. From left to right and from top
  to bottom: 1.\ the equal splitting flow $\nu=\de_{1/2}$, 2.\ the `parabolic'
  flow $\nu(\di q)=6q(1-q)\di q$, 3.\ the Le Jan-Raimond flow $\nu(\di q)=\di
  q$, 4.\ the symmetric erosion flow $\nu=\ffrac{1}{2}\big(\de_0+\de_1)$. The
  first two flows have left and right speeds $\bet_-,\bet_+=\pm 4$ and
  $\bet_-,\bet_+=\pm 6$, respectively, while the last two flows have
  $\bet_-,\bet_+=\pm\infty$. Each picture shows a rectangle of $1.4$ units of
  space (horizontal) by $0.2$ units of time (vertical). The initial state
  is Lebesgue measure.}
\label{fig:flows}
\end{center}
\end{figure}

\begin{figure}[htb]
\begin{center}
\includegraphics[width=10cm]{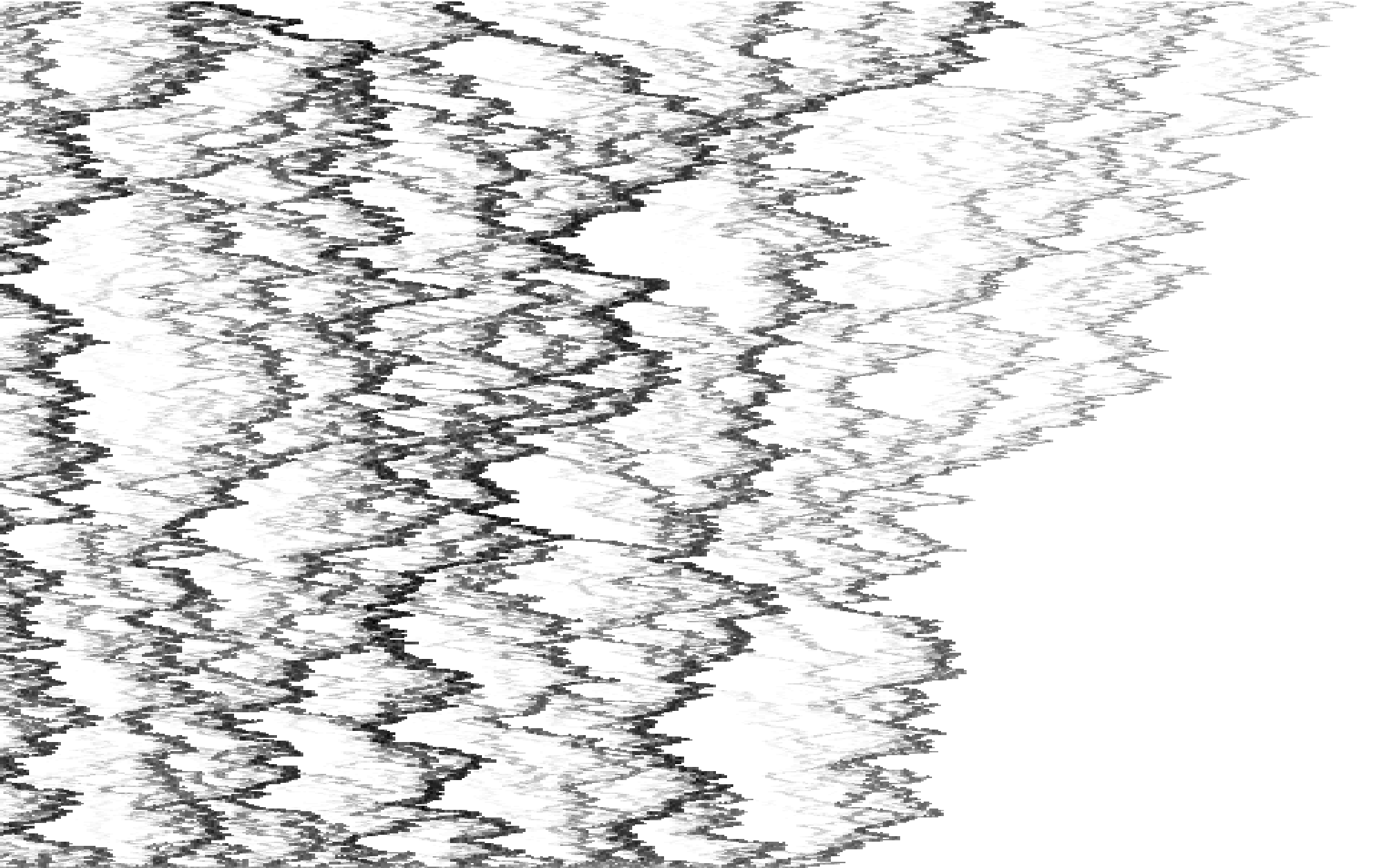}
\caption{Example of an asymmetric flow: the one-sided erosion flow with
  $\bet=0$ and $\nu=\de_1$. This flow has left and right speeds
  $\bet_-=-\infty$ and $\bet_+=2$, respectively. The picture shows a rectangle
  of $3.2$ units of space (horizontal) by $0.4$ units of time (vertical). The
  initial state is Lebesgue measure up to the point $1.9$ and zero from there
  onwards.}
\label{fig:onesid}
\end{center}
\end{figure}

In this subsection, we state a number of results on the almost sure path
properties of the measure-valued Markov process $(\rho_t)_{t\geq 0}$ defined
in terms of a Howitt-Warren flow by (\ref{rho}). Throughout this subsection,
we will assume that $\rho_0$ is a finite measure, and $\rho_t$ is defined
using a version of the Howitt-Warren flow $(K_{s,t})_{s\leq t}$, which
satisfies property (iv) in Proposition~\ref{P:conpath}, but not necessarily
property (i)'. Then it is not hard to see that for any $\rho_0\in\Mi(\R)$, the
Markov process $(\rho_t)_{t\geq 0}$ defined in (\ref{rho}) has continuous
sample paths in $\Mi(\R)$.  We call this process the {\em Howitt-Warren
  process} with drift $\bet$ and characteristic
measure~$\nu$. \index{Howitt-Warren!process}

See Figures~\ref{fig:flows} and \ref{fig:onesid} for some simulations of
Howitt-Warren processes for various choices of the characteristic measure
$\nu$. There are a number of parameters that are important for the behavior of
these processes. First of all, following \cite{HW09a}, we define
\be\label{tetdef}
\tet(k,l)=\int\nu(\di q)\,q^{k-l}(1-q)^{l-1}\qquad(k,l\geq 1).
\ee
In a certain excursion theoretic sense, $\tet(k,l)$ describes the rate
at which a group of $k+l$ coordinates of the $n$-point motion that are
at the same position splits into two groups consisting of $k$ and $l$
specified coordinates, respectively. In particular, following again
notation in \cite{HW09a}, we set
\be\label{stickyparam}
\tet:=2\tet(1,1)=2\int_{[0,1]}\nu(\di q),
\ee
and we call $\tet$ the {\em stickiness parameter} of the Howitt-Warren
flow. \index{stickiness parameter}
Note that when $\tet$ is increased, particles separate with a
higher rate, hence the flow is less sticky. The next proposition shows
that with the exception of the Arratia flow, by a simple
transformation of space-time, we can always scale our flow such that
$\bet=0$ and $\tet=2$. Below, for any $A\sub\R$ and $a\in\R$ we write
$aA:=\{ax:x\in A\}$ and $A+a:=\{x+a:x\in A\}$.

\bp{\bf(Scaling and removal of the drift)}\label{P:scale}
Let $(K_{s,t})_{s\leq t}$ be a Howitt-Warren flow with drift $\bet$ and
characteristic measure $\nu$. Then:
\begin{itemize}
\item[{\bf(a)}] For each $a>0$, the stochastic flow of kernels
  $(K'_{s,t})_{s\leq t}$ defined by $K'_{a^2s,a^2t}(ax,aA):=K_{s,t}(x,A)$ is a
  Howitt-Warren flow with drift $a^{-1}\bet$ and characteristic measure
  $a^{-1}\nu$.
\item[{\bf(b)}] For each $a\in\R$, the stochastic flow of kernels
  $(K'_{s,t})_{s\leq t}$ defined by $K'_{s,t}(x+as,A+at):=K_{s,t}(x,A)$ is a
  Howitt-Warren flow with drift $\bet+a$ and characteristic measure $\nu$.
\end{itemize}
\ep

There are two more parameters that are important for the behavior of a
Howitt-Warren flow. We define
\bc\label{speeds}
\dis\bet_-&:=&\dis\bet-2\int\nu(\di q)(1-q)^{-1},\\[5pt]
\dis\bet_+&:=&\dis\bet+2\int\nu(\di q)q^{-1}
\ec
Note that $\bet_+=\lim_{m\to\infty}\bet_+(m)$, where $(\bet_+(m))_{m\geq 1}$ are
the constants defined in (\ref{betplusm}). We call $\bet_-$ and $\bet_+$ the
{\em left speed} and {\em right speed} \index{left-right!speeds!of Howitt-Warren flow} \index{speeds!of Howitt-Warren flow} of a Howitt-Warren flow,
respectively. The next theorem shows that these names are justified. Below,
$\supp(\mu)$ denotes the support of a measure $\mu$, i.e., the smallest
closed set that contains all mass.

\bt{\bf(Left and right speeds)}\label{T:speed}
Let $(\rho_t)_{t\geq 0}$ be a Howitt-Warren process with drift $\bet$ and
characteristic measure $\nu$, and let $\bet_-,\bet_+$ be defined as in
(\ref{speeds}). Set $r_t:=\sup(\supp(\rho_t))$ $(t\geq 0)$. Then:
\begin{itemize}
\item[{\bf(a)}] If $\bet_+<\infty$ and $r_0<\infty$, then $(r_t)_{t\geq 0}$ is
  a Brownian motion with drift $\bet_+$. If $\bet_+<\infty$ and $r_0=\infty$,
  then $r_t=\infty$ for all $t\geq 0$.
\item[{\bf(b)}] If $\bet_+=\infty$, then $r_t=\infty$ for all $t>0$.
\end{itemize}
\noi
Analogue statements hold for $l_t:=\inf(\supp(\rho_t))$, with $\bet_+$
replaced by $\bet_-$.
\et

It turns out that the support of a Howitt-Warren process is itself a Markov
process. Let $\Kl(\R)$ be the space of closed subsets of $\R$. We equip
$\Kl(\R)$ with a topology such that $A_n\to A$ if and only if $\ov
A_n\stackrel{{\rm Haus}}{\longrightarrow}\ov A$, where $\ov A$ denotes the
closure of a set $A$ in $[-\infty,\infty]$ and $\stackrel{{\rm
    Haus}}{\longrightarrow}$ means convergence of compact subsets of
$[-\infty,\infty]$ in the Hausdorff topology. The {\em branching-coalescing
  point set} is a $\Kl(\R)$-valued Markov process that has been introduced in
\cite[Thm~1.11]{SS08}. Its definition involves the Brownian net; see formula
(\ref{braco}) below. The following proposition, which we cite from
\cite[Thm~1.11 and Prop.~1.15]{SS08} and \cite[Prop.~3.14]{SSS09}, lists some
of its elementary properties.

\index{branching-coalescing point set!standard}
\bp{\bf(Properties of the branching-coalescing point set)}\label{P:braco}
Let $\xi=(\xi_t)_{t\geq 0}$ be the branching-coalescing point set defined in
(\ref{braco}), started in any initial state $\xi_0\in\Kl(\R)$. Then:
\begin{itemize}
\item[{\bf(a)}] The process $\xi$ is a $\Kl(\R)$-valued Markov process with
  continuous sample paths.
\item[{\bf(b)}] If $\sup(\xi_0)<\infty$, then $(\sup(\xi_t))_{t\geq 0}$ is a
  Brownian motion with drift $+1$. Likewise, if $-\infty<\inf(\xi_0)$, then
  $(\inf(\xi_t))_{t\geq 0}$ is a Brownian motion with drift $-1$.
\item[{\bf(c)}] The law of a Poisson point set with intensity $2$ is a
  reversible invariant law for $\xi$ and the limit law of $\xi_t$ as $t\to\infty$
  for any initial state $\xi_0\neq\emptyset$.
\item[{\bf(d)}] For each deterministic time $t>0$, a.s.\ $\xi_t$ is a locally
  finite subset of $\R$.
\item[{\bf(e)}] Almost surely, there exists a dense set $\Ti\sub(0,\infty)$
  such that for each $t\in\Ti$, the set $\xi_t$ contains no isolated points.
\end{itemize}
\ep
Our next result shows how Howitt-Warren processes and the
branching-coalescing point set are related. Note that this result
covers all possible values of $\bet_-,\bet_+$, except the case
$\bet_-=\bet_+$ which corresponds to the Arratia flow. In
(\ref{supbraco}) below, we continue to use the notation
$aA+b:=\{ax+b:x\in A\}$.
\bt{\bf(Support process)}\label{T:supp}
Let $(\rho_t)_{t\geq 0}$ be a Howitt-Warren process with drift $\bet$ and
characteristic measure $\nu$ and let $\bet_-,\bet_+$ be defined as in
(\ref{speeds}). Then:
\begin{itemize}
\item[{\bf(a)}] If $-\infty<\bet_-<\bet_+<\infty$, then a.s.\ for all $t>0$,
\be\label{supbraco}
\supp(\rho_t)=\ffrac{1}{2}(\bet_+-\bet_-)\xi_t
+\ffrac{1}{2}(\bet_-+\bet_+)t,
\ee
where $(\xi_t)_{t\geq 0}$ is a branching-coalescing point set.
\item[{\bf(b)}] If $\bet_-=-\infty$ and $\bet_+<\infty$, then a.s.\ ${\rm
  supp}(\rho_t)=(-\infty,r_t]\cap\R$ for all $t>0$, where $r_t:=\sup({\rm
  supp}(\rho_t))$. An analogue statement holds when $\bet_->-\infty$ and
$\bet_+=\infty$.
\item[{\bf(c)}] If $\bet_-=-\infty$ and $\bet_+=\infty$, then a.s.\ ${\rm
  supp}(\rho_t)=\R$ for all $t>0$.
\end{itemize}
\et

Proposition~\ref{P:braco}~(d) and Theorem~\ref{T:supp}~(a) imply that if the
left and right speeds of a Howitt-Warren process are finite, then at
deterministic times the process is purely atomic. The next theorem generalizes
this statement to any Howitt-Warren process, but shows that if the
characteristic measure puts mass on the open interval $(0,1)$, then there are
random times when the statement fails to hold.

\bt{\bf(Atomicness)}\label{T:atom} Let
  $(\rho_t)_{t\geq 0}$ be a Howitt-Warren process with drift $\bet$ and
  characteristic measure $\nu$. Then:
\begin{itemize}
\item[{\bf(a)}] For each $t>0$, the measure $\rho_t$ is a.s.\ purely atomic.
\item[{\bf(b)}] If $\int_{(0,1)}\nu(\di q)>0$, then a.s.\ there exists a dense
  set of random times $t>0$ when $\rho_t$ is purely non-atomic.
\item[{\bf(c)}] If $\int_{(0,1)}\nu(\di q)=0$, then a.s.\ $\rho_t$ is purely
  atomic at all $t>0$.
\end{itemize}
\et
In the special case that $\nu$ is (a multiple of) Lebesgue measure, a weaker
version of part~(a) has been proved in \cite[Prop.~9~(c)]{LR04PTRF}. Part~(b)
is similar to Proposition~\ref{P:braco}~(e) and in fact, by
Theorem~\ref{T:supp}~(a), implies the latter. Note that parts~(b) and (c) of
the theorem reveal an interesting dichotomy between erosion flows (where $\nu$
is nonzero and concentrated on $\{0,1\}$) and all other Howitt-Warren flows
(except the Arratia flow, for which atomicness is trivial). The reason is that
atoms in erosion flows lose mass continuously (see the footnote in
Section~\ref{S:heur} and the construction in Section~\ref{S:flowcons} below),
while in all other flows atoms can be split into smaller atoms. This latter
mechanism turns out to be more effective at destroying atoms. For erosion
flows, we have an exact description of the set of space-time points where
$(\rho_t)_{t\geq 0}$ has an atom in terms of an underlying Brownian web, see
Theorem~\ref{T:erosion} below.

\subsection{Infinite starting measures and discrete
approximation}\label{S:infmass}

The ergodic behavior of the branching-coalescing point set is well-understood
(see Proposition \ref{P:braco}~(c)). As a consequence, by
Theorem~\ref{T:supp}~(a), it is known that if we start a Howitt-Warren process
with left and right speeds $\bet_-=-1$, $\bet_+=1$ in any nonzero initial
state, then its support will converge in law to a Poisson point process with
intensity $2$. This does not mean, however, that the Howitt-Warren process
itself converges in law. Indeed, since its 1-point motion is Brownian motion,
it is easy to see that any Howitt-Warren process started in a finite initial
measure satisfies $\lim_{t\to\infty}\E[\rho_t(K)]=0$ for any compact
$K\sub\R$. To find nontrivial invariant laws, we must start the process in
infinite initial measures.

To this aim, let $\Mi_{\rm loc}(\R)$ denote the space of locally finite
measures on $\R$, endowed with the vague topology. Let $(K_{s,t})_{s\leq t}$ be
a version of the Howitt-Warren flow with $-\infty<\bet_-$ and $\bet_+<\infty$,
which satisfies Proposition~\ref{P:conpath}~(iv). We will prove that for any
$\rho_0\in\Mi_{\rm loc}(\R)$,
\be\label{rhodef4}
\rho_t:=\int\rho_0(\di x)K_{0,t}(x,\,\cdot\,)\qquad(t\geq 0)
\ee
defines an $\Mi_{\rm loc}(\R)$-valued Markov process. If
$\bet_+-\bet_-=\infty$, then mass can spread infinitely fast, hence we cannot
define the Howitt-Warren process $(\rho_t)_{t\geq 0}$ for arbitrary
$\rho_0\in\Mi_{\rm loc}(\R)$. In this case, we will use the class
\be\label{MiG}
\Mi_{\rm g}(\R):=\Big\{\rho\in\Mi_{\rm loc}(\R):\int_\R e^{-cx^2}\rho(\di x)<\infty
\mbox{ for all }c>0\Big\},
\ee
endowed with the topology that $\mu_n\to\mu$ if and only if
$e^{-cx^2}\mu_n(\di x)$ converges weakly to $e^{-cx^2}\mu(\di x)$ for all
$c>0$, which can be seen to be equivalent to $\mu_n\to\mu$ in the vague
topology plus $\int e^{-cx^2}\mu_n(\di x)\to\int e^{-cx^2}\mu(\di x)$ for all
$c>0$. Note that $\Mi_{\rm loc}(\R)$ and $\Mi_{\rm g}(\R)$ are Polish spaces.

Observe that by Definition~\ref{D:stochflow}~(i), the Howitt-Warren process
$(\rho_t)_{t\geq 0}$ defined in (\ref{rhodef4}) satisfies
\be\label{rhodef5}
\rho_t = \int \rho_s(\di x) K_{s,t}(x,\cdot)
\quad{\rm a.s.}
\ee
for each deterministic $s<t$. We will also use (\ref{rhodef5}) to define Howitt-Warren processes starting at
any deterministic time $s\in\R$.

\bt{\bf(Infinite starting mass and continuous dependence)}\label{T:infmass}
Let $\bet\in\R$, let $\nu$ be a finite measure on $[0,1]$, and let
$(K_{s,t})_{s\leq t}$ be a version of the Howitt-Warren flow with drift $\bet$
and characteristic measure $\nu$ satisfying property~(iv) from
Proposition~\ref{P:conpath}. Then:\med

\noi
{\bf(a)} For any $\rho_0\in\Mi_{\rm g}(\R)$, formula (\ref{rhodef4}) defines
an $\Mi_{\rm g}(\R)$-valued Markov process with continuous sample paths,
satisfying
\be\label{EK}
\E[\rho_t(K)]<\infty\qquad(t\geq 0,\ K\sub\R\mbox{ compact}).
\ee
Moreover, if $(\rho^{\li n\re}_t)_{t\geq s_n}$ are processes started at times $s_n$ with deterministic initial data $\rho^{\li n\re}_{s_n}$, and $s_n\to 0$,
then for any $t>0$ and $t_n\to t$,
\be\label{rhonconv}
\rho^{\li n\re}_{s_n}\Asto{n}\rho_0\quad\mbox{implies}\quad
\rho^{\li n\re}_{t_n}\Asto{n}\rho_t\qquad a.s.,
\ee
where $\Rightarrow$ denotes convergence in $\Mi_{\rm g}(\R)$.\med

\noi
{\bf(b)} Assume moreover that $\bet_+-\bet_-<\infty$. Then, for any
$\rho_0\in\Mi_{\rm loc}(\R)$, formula (\ref{rhodef4}) defines an $\Mi_{\rm
  loc}(\R)$-valued Markov process with continuous sample paths.
Moreover, formula (\ref{rhonconv}) holds with convergence in
$\Mi_{\rm g}(\R)$ replaced by vague convergence in $\Mi_{\rm loc}(\R)$.
\et
{\bf Remark.} The convergence in (\ref{rhonconv}) implies the continuous dependence of the law of $\rho_t$ on the starting time and the initial law,
which is known as the {\it Feller property}. Note that when $\rho_0$ is a finite measure, the continuity in $t$ of $\rho_t$ in the space $\Mi(\R)$
already follows from Proposition~\ref{P:conpath}~(iv). However, for our purposes, we will only consider the spaces $\Mi_{\rm g}$ and $\Mi_{\rm loc}$ .
\medskip

\noi
{\bf Remark.} When $\beta_+-\beta_-=\infty$, $(\rho_t)_{t>0}$ may not be
well-defined if $\rho_0\notin\Mi_{\rm g}(\R)$. Indeed, by
Theorem~\ref{T:supp}, if $\beta_+-\beta_-=\infty$, then for any fixed $t>0$,
we can find $x_n\in\Z$ with $|x_n|\to\infty$ such that $\P(K_{0,t}(x_n,
[0,1])<\eps_n)<2^{-n}$ for some $\eps_n>0$. Therefore
$\rho_0:=\sum_n\eps_n^{-1} \delta_{x_n}$ has $\rho_0\in\Mi_{\rm loc}(\R)$, and
almost surely, $\rho_t([0,1])=\infty$.  \med

\noi
{\bf Remark.}
Theorems~\ref{T:speed}, \ref{T:supp}, and \ref{T:atom} carry over without change to the case of
infinite starting measures. To see this, note that it is easy to check from
(\ref{rhodef4}) that
\be
\rho_0\ll\tilde\rho_0\quad\mbox{implies}\quad\rho_t\ll\tilde\rho_t\qquad(t\geq 0),
\ee
where $\ll$ denotes absolute continuity. Since for each $\rho_0\in\Mi_{\rm
  loc}(\R)$, we can find a finite measure $\rho'_0$ that is equivalent to
$\rho$, statements about the support of $\rho_t$ and atomicness immediately
generalize to the case of locally finite starting measures.
\medskip

We also collect here a discrete approximation result for Howitt-Warren processes.

\bt{\bf {(Convergence of discrete Howitt-Warren processes)}}\label{T:disflow}
Let $\eps_k$ be positive constants converging to zero, and let $\mu_k$ be
probability measures on $[0,1]$ satisfying (\ref{mucon}) for some real $\bet$
and finite measure $\nu$ on $[0,1]$. Let $(\rho^{\langle k\rangle}_t)_{t\geq 0}$ be a
discrete Howitt-Warren process with characteristic measure $\mu_k$ defined as
in (\ref{disrhoK}), where $K^\omega_{0,t}(x,\cdot)$ therein is defined for
all $t>0$ by letting the random walk $(X_t)_{t\geq 0}$ in (\ref{disK})
be linearly interpolated between integer times.
Let $\bar\rho^{\langle k\rangle}_t(\di x) := \rho^{\langle k\rangle}_{\eps_k^{-2}t}(\eps_k^{-1}\di x)$.
If $\rho^{\langle k\rangle}_0$ is deterministic and $\bar\rho^{\langle k\rangle}_0\Rightarrow \rho_0$
in $\Mi_g(\R)$, then for any $T>0$,
\be\label{disflow}
(\bar\rho^{\langle k\rangle}_t)_{0\leq t\leq T} \Asto{k} (\rho_t)_{0\leq t\leq T},
\ee
where $\rho_t$ is a Howitt-Warren process with drift $\beta$, characteristic measure $\nu$, and
initial condition $\rho_0$, and $\Rightarrow$ denotes weak convergence in law of random variables taking
values in $\Ci([0,T], \Mi_g(\R))$, the space of continuous functions from $[0,T]$ to $\Mi_g(\R)$ equipped
with the uniform topology.
\et

\subsection{Ergodic properties}\label{S:ergo}

We are now ready to discuss the ergodic behavior of Howitt-Warren processes.
Note that for a given Howitt-Warren flow $(K_{s,t})_{s\leq t}$, the right-hand
side of (\ref{rhodef4}) is a.s.\ a linear function of the starting measure
$\rho_0$. In view of this, Howitt-Warren processes belong to the class of
so-called {\it linear systems}. The theory of linear systems on $\Z^d$ has
been developed by Liggett and Spitzer, see e.g.\ \cite{LS81} and
\cite[Chap.~IX]{Lig05}. We will adapt this theory to the continuum setting
here. First we define the necessary notion. \index{linear system}

We let $\Ii$ denote the set of invariant laws of a given Howitt-Warren
processes, i.e., $\Ii$ is the set of probability laws $\La$ on $\Mi_{\rm loc}(\R)$
(resp.\ $\Mi_{\rm g}(\R)$ if $\bet_+-\bet_-=\infty$) such that
$\P[\rho_0\in\cdot\,]=\La$ implies $\P[\rho_t\in\cdot\,]=\La$ for all $t\geq
0$. We let $\Ti$ denote the set of homogeneous (i.e., translation invariant) laws on $\Mi_{\rm loc}(\R)$
(resp.\ $\Mi_{\rm g}(\R)$), i.e., laws $\La$ such that
$\P[\rho\in\cdot\,]=\La$ implies $\P[T_a\rho\in\cdot\,]=\La$ for all $a\in\R$,
where $T_a\rho(A):=\rho(A+a)$ denotes the spatial shift map.  Note that $\Ii$
and $\Ti$ are both convex sets. We write $\Ii_{\rm e}$, $\Ti_{\rm e}$, and
$(\Ii\cap \Ti)_{\rm e}$ to denote respectively the set of extremal elements in
$\Ii$, $\Ti$, and $\Ii\cap\Ti$. Below, $\Ci_{\rm c}(\R)$ denotes the space of
continuous real function on $\R$ with compact support.

\index{homogeneous invariant laws}

\bt{\bf(Homogeneous invariant laws for Howitt-Warren processes)}
\label{T:HIL}
Let $\bet\in\R$, let $\nu$ be a finite measure on $[0,1]$ with $\nu\neq
0$. Then for the corresponding Howitt-Warren process $(\rho_t)_{t\geq 0}$,
we have:
\begin{itemize}
\item[{\bf(a)}] $(\Ii\cap\Ti)_{\rm e}$  is a
one-parameter familly $\{\La_c:c\geq 0\}$ of measures satisfying
$\Lambda_c(\di(c\rho))=\Lambda_1(\di \rho)$ for all $c\geq 0$, and
\begin{eqnarray}
\!\!\!\!\!\!\!\!\!\!\!\!\label{1stmom}
\int\!\Lambda_1(\di \rho)\int\!\rho(\di x)\,\phi(x)
&=&\int\phi(x)\,\di x,\\
\!\!\!\!\!\!\!\!\!\!\!\!\label{2ndmom}
\int\!\Lambda_1(\di \rho)\int\!\rho(\di x)\,\phi(x)\int\!\rho(\di y)\,\psi(y)
&=&\int\phi(x)\,\di x\int\psi(y)\,\di y
+\frac{\int\phi(x)\psi(x)\,\di x}{2\nu([0,1])}
\end{eqnarray}
for any $\phi,\psi\in\Ci_{\rm c}(\R)$.

\item[{\bf(b)}] If $\P[\rho_0\in\cdot\,]\in\Ti_{\rm e}$ and
  $\E[\rho_0([0,1])]=c\geq 0$, then $\P[\rho_t\in\cdot\,]$ converges weakly to
  $\Lambda_c$. Furthermore, if $\E[\rho_0([0,1])^2]<\infty$, then for any
  $\phi,\psi\in\Ci_{\rm c}(\R)$,
\be\label{2momconv}
\lim_{t\to\infty}\E\Big[\int\!\rho_t(\di x)\,\phi(x)
\int\!\rho_t(\di y)\,\psi(y)\Big]
=\int\!\Lambda_c(\di \rho)\int\!\rho(\di x)\,\phi(x)\int\!\rho(\di y)\,\psi(y).
\ee

\item[{\bf(c)}] If $\P[\rho_0\in\cdot\,]\in\Ti_{\rm e}$ and
  $\E[\rho_0([0,1])]=\infty$, then the laws $\P[\rho_t\in\cdot\,]$
have no weak cluster point as $t\to\infty$ which is supported on $\Mi_{\rm loc}(\R)$.

\item[{\bf(d)}] If $\Lambda\in\Ii\cap\Ti$, then there exists a probability
  measure $\gamma$ on $[0,\infty)$ such that $\Lambda=\int_0^\infty\gamma(\di c)
  \,\Lambda_c$.
\end{itemize}
\et
{\bf Remark.} When $\nu$ is Lebesgue measure, it is known that (see
\cite[Prop.~9~(b)]{LR04PTRF}) $\La_{\rm c}$ is the law of $c\rho^*$, where
$\rho^*=\sum_{(x,u)\in\cal P} u\delta_x$ for a Poisson point process $\cal P$
on $\R\times [0,\infty)$ with intensity measure $\di x\times u^{-1}e^{-u}\di
  u$.\med

\noi
Theorem \ref{T:HIL} shows that each Howitt-Warren process has a unique
(modulo a constant multiple) homogeneous invariant law, which by
(\ref{2ndmom}) has zero off-diagonal correlations. Moreover, any ergodic law
at time 0 with finite density converges under the dynamics to the unique
homogeneous invariant law with the same density.

In line with Theorems~\ref{T:supp} and \ref{T:atom} we have the following
support properties for $\Lambda_c$.
\bt{\bf(Support of stationary process)}\label{T:invsupp}
Let $c>0$ and let $\rho$ be an $\Mi_{\rm loc}(\R)$-valued random variable with
law $\Lambda_c$, the extremal homogeneous invariant law defined in
Theorem~\ref{T:HIL}. Then:
\begin{itemize}
\item[{\bf(a)}] If $\bet_+-\bet_-<\infty$, then $\supp(\rho)$ is a
Poisson point process with intensity $\bet_+-\bet_-$.
\item[{\bf(b)}] If $\bet_+-\bet_-=\infty$, then $\rho$ is a.s.\ atomic with
$\supp(\rho)=\R$.
\end{itemize}
\et

\section{Construction of Howitt-Warren flows in the Brownian web}
\label{S:web}

In this section, we make the heuristics in Section \ref{S:heur}
rigorous and give a graphical construction of the Howitt-Warren flows
using a procedure of Poisson marking of the Brownian web invented by
Newman, Ravishankar and Schertzer~\cite{NRS10}. The random environment
for the Howitt-Warren flow will turn out to be a Brownian web, which
we call the {\em reference web}, plus a marked Poisson point process on the
reference web. Given such an environment, we will then construct a
second coupled Brownian web, which we call the {\em sample web}, which is
constructed by modifying the reference web by switching the
orientation of marked points of type $(1,2)$. The kernels of the
Howitt-Warren flow are then constructed from the quenched law of the
sample web, conditional on the reference web and the associated marked
Poisson point process.

This construction generalizes the construction of the erosion flow
based on coupled Brownian webs given in \cite{HW09b}. For erosion
flows, the random environment consists only of a reference web
(without marked points) and the construction of the sample web can be done
by specifying the joint law of the reference web and the sample web by
means of a martingale problem. This is the approach taken in
\cite{HW09b}. In the general case, when the random environment also
contains marked points, this approach does not work. Therefore, in our
approach, even for erosion flows, we will give a graphical
construction of the sample web by marking and switching paths in the
reference web.

Discrete approximation will be an important tool in many of our proofs and is
helpful for understanding the continuum models. Therefore, in Section
\ref{S:disref}, we will first formulate the notion of a quenched law of sample
webs conditional on the random environment for discrete Howitt-Warren
flows. In Section~\ref{S:BMweb}, we then recall the necessary background on
the Brownian web and Poisson marking for the Brownian web. In
Section~\ref{S:modweb}, we show how coupled Brownian webs can be constructed
by Poisson marking and switching paths in a reference web. In
Section~\ref{S:flowcons}, we state our main result, Theorem~\ref{T:HWconst},
which is the construction of Howitt-Warren flows using the Poisson marking of
a reference Brownian web, and we also state some regularity properties for the
Howitt-Warren flows. Lastly, in Section~\ref{S:approx}, we state a convergence
result on the quenched law of discrete webs, which will be used to identify
the flows we construct in Theorem \ref{T:HWconst} as being, indeed, the
Howitt-Warren flows defined in Section~\ref{S:flow} earlier through their
$n$-point motions. The statements of this section are proved in
Sections~\ref{S:prelim} and \ref{S:main}.

\subsection{A quenched law on the space of discrete webs}\label{S:disref}

\begin{figure}[htb]
\begin{center}
\includegraphics[width=9cm]{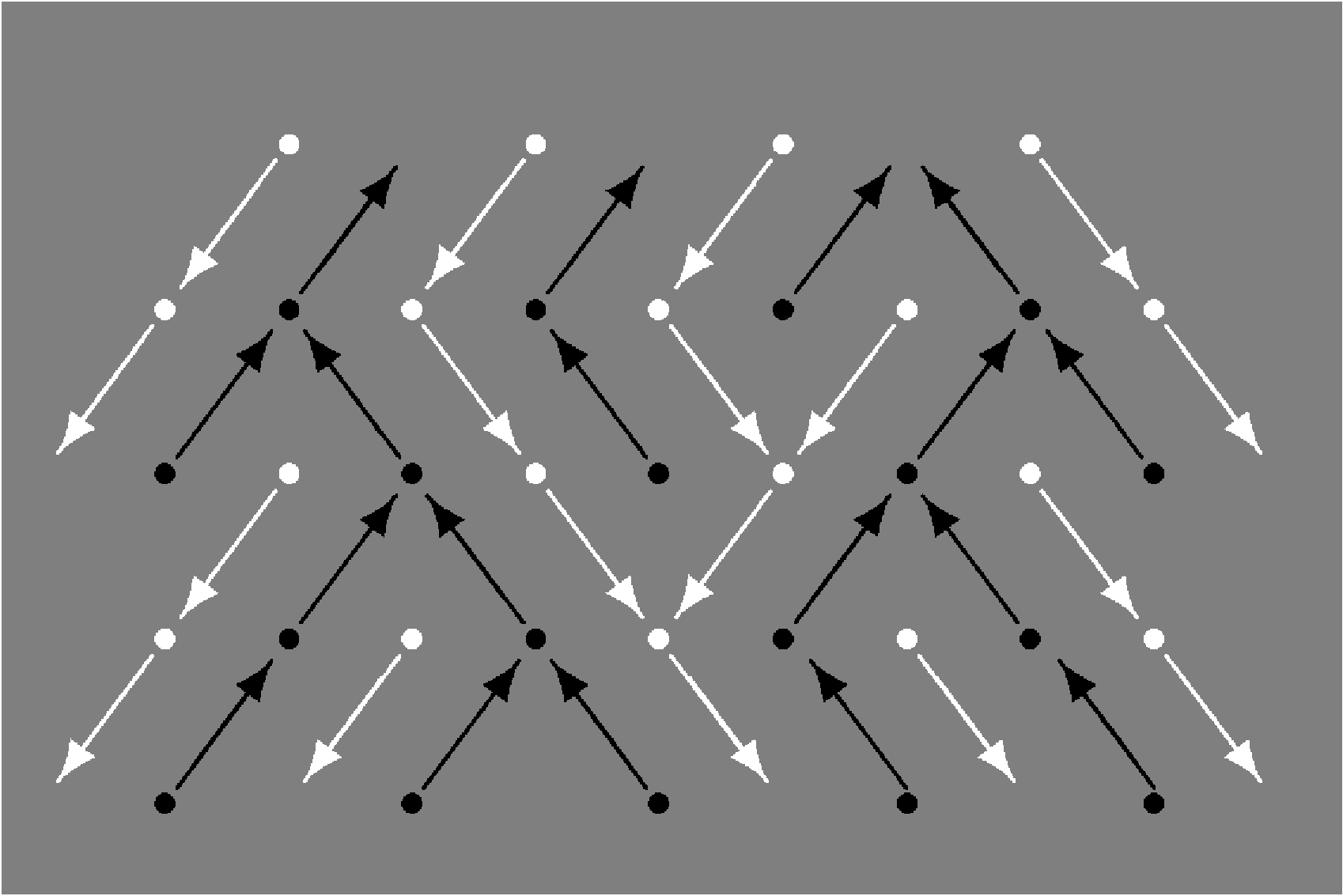}
\caption{A discrete web and its dual.}
\label{fig:disdual}
\end{center}
\end{figure}

As in Section~\ref{S:dis}, let $\om:=(\om_z)_{z\in \Zev}$ be
i.i.d.\ $[0,1]$-valued random variables with common distribution $\mu$.
Instead of using $\om$ as a random environment for a single random walk
started from one fixed time and position, as we did in Section~\ref{S:dis}, we
will now use $\om$ as a random environment for a collection of coalescing
random walks starting from each point in $\Zev$. To this aim, conditional on
$\om$, let $\al=(\al_z)_{z\in\Zev}$ be a collection of independent
$\{-1,+1\}$-valued random variables such that $\al_z=+1$ with probability
$\om_z$ and $\al_z=-1$ with probability $1-\om_z$. This $\al$ will play a
somewhat different role from the $\varal^{\li k\re}$ in Section~\ref{S:heur};
see the discussion below Theorem~\ref{T:HWconst}. For each $(x,s)\in\Zev$, we
let $\pad^\al_{(x,s)}:\{s,s+1,\ldots\}\to\Z$ be the function
$\pad^\al_{(x,s)}=p$ defined by
\be\label{pal}
p(s):=x\quad\mbox{and}\quad
p(t+1):=p(t)+\al_{(p(t),t)}\qquad(t\geq s).
\ee
Then $\pad^\al_{(x,s)}$ is the path of a random walk in the random environment
$\om$, started at time $s$ at position $x$. It is easy to see that paths
$\pad^\al_z,\pad^\al_{z'}$ starting at different points $z,z'$ coalesce when
they meet. We call the collection of paths
\be\label{Uial}
\Ui^\al:=\{\pad^\al_z:z\in\Zev\}
\ee
the {\em discrete web} \index{discrete!web} associated with $\al$ (see
Figure~\ref{fig:disdual}), where starting from this
section, for the rest of the paper, we will use different notation for
discrete webs and nets compared to Section~\ref{S:dis}, to avoid confusion
with certain other symbols that we will need. Let $\P$ denote the law of $\om$
and let
\be\label{quench}
\Qdis^\om:=\P\big[\Ui^\al\in\cdot\,\big|\,\om\big]
\ee
denote the conditional law of $\Ui^\al$ given $\om$. Then under the averaged law
$\int\P(\di\om)\Qdis^\om$, paths in $\Ui^\al$ are coalescing random walks that
in each time step jump to the right with probability $\int\mu(\di q)q$ and to
the left with the remaining probability $\int\mu(\di q)(1-q)$.

We will be more interested in the {\em quenched} law $\Qdis^\om$ defined in
(\ref{quench}). One has
\be
\Qdis^\om_z=\Qdis^\om[\pad^\al_z\in\cdot\,]
\ee
where $\Qdis^\om_z$ is the conditional law of the random walk in random
environment started from $z\in\Zev$ defined in Section~\ref{S:dis}. In
particular, by (\ref{disK}),
\be\label{KQ}
K^\om_{s,t}(x,y)=\Qdis^\om\big[p^\al_{(x,s)}(t)=y\big]
=\P\big[p^\al_{(x,s)}(t)=y\,\big|\,\om\big]
\qquad\big((x,s),(y,t)\in\Zev\big),
\ee
where $(K^\om_{s,t})_{s\leq t}$ is the discrete Howitt-Warren flow with
characteristic measure $\mu$. In view of this, the random law $\Qdis^\om$
contains all information that we are interested in. We call $\Qdis^\om$ the {\em
  discrete quenched law} with {\em characteristic measure} $\mu$. In the next
sections, we will construct a continuous analogue of this quenched law and use
it to define Howitt-Warren flows. \index{quenched law!discrete} \index{Howitt-Warren!quenched law!discrete}

\subsection{The Brownian web}\label{S:BMweb}

As pointed out in the previous subsection, under the averaged law
$\int\P(\di\om)\Qdis^\om$, the discrete web $\Ui^\al$ is a collection of
coalescing random walks, started from every point in $\Zev$. It turns out that
such discrete webs have a well-defined diffusive scaling limit, which is
basically a collection of coalescing Brownian motions, starting from each
point in space and time, and which is called a {\em Brownian web}. The
Brownian web arose from the work of Arratia \cite{Arr79,Arr81} and has since
been studied by T\'oth and Werner \cite{TW98}.  More recently, Fontes, Isopi,
Newman and Ravishankar \cite{FINR04} have introduced a by now standard
framework in which the Brownian web is regarded as a random compact set of
paths, and is an element of a suitable Polish space. \index{Brownian web}

It turns out that associated to each Brownian web, there is a dual Brownian
web, which is a collection of coalescing Brownian motions running backwards in
time. To understand this on a heuristic level, let $(\al_z)_{z\in\Zev}$ be an
i.i.d.\ collection of $\{-1,+1\}$-valued random variables. If for each
$z=(x,t)\in\Zev$, we draw an arrow from $(x,t)$ to $(x+\al_z,t+1)$, then paths
along these arrows form a discrete web as introduced in the previous section.
Now, if for each $z=(x,t)\in\Zev$, we draw in addition a {\em dual arrow} from
$(x,t+1)$ to $(x-\al_z,t)$, then paths along these dual arrows form a {\em
  dual discrete web} \index{dual!discrete web} of coalescing random walks running backwards in time,
which do not cross paths in the forward web (see Figure~\ref{fig:disdual}).
The dual Brownian web arises as the diffusive scaling limit of such a dual
discrete web.

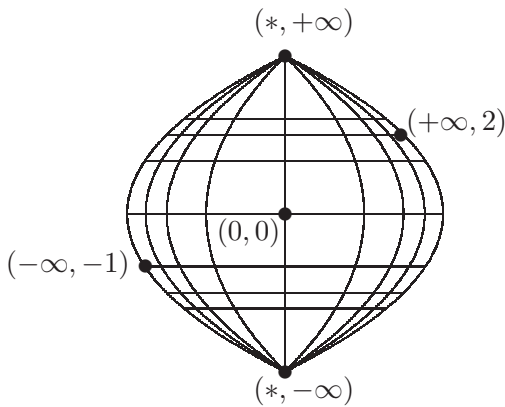
\begin{figure}[htb]
\begin{center}
\setlength{\unitlength}{.7cm}
\begin{picture}(10,8)(-5,-4)
\linethickness{.4pt}
\qbezier(0,-3)(0,0)(0,3)
\qbezier(-3,0)(0,0)(3,0)
\linethickness{.4pt}
\qbezier(0,-3)(-6,0)(0,3)
\qbezier(0,-3)(-5.3,0)(0,3)
\qbezier(0,-3)(-4.5,0)(0,3)
\qbezier(0,-3)(-3,0)(0,3)
\qbezier(0,-3)(3,0)(0,3)
\qbezier(0,-3)(4.5,0)(0,3)
\qbezier(0,-3)(5.3,0)(0,3)
\qbezier(0,-3)(6,0)(0,3)
\qbezier(-1.9,-1.8)(0,-1.8)(1.9,-1.8)
\qbezier(-2.2,-1.5)(0,-1.5)(2.2,-1.5)
\qbezier(-2.65,-1)(0,-1)(2.65,-1)
\qbezier(-2.65,1)(0,1)(2.65,1)
\qbezier(-2.2,1.5)(0,1.5)(2.2,1.5)
\qbezier(-1.9,1.8)(0,1.8)(1.9,1.8)

\put(0,-3){\circle*{.25}}
\put(0,3){\circle*{.25}}
\put(0,0){\circle*{.25}}
\put(2.2,1.5){\circle*{.25}}
\put(-2.65,-1){\circle*{.25}}

\put(-.6,-3.5){$(\ast,-\infty)$}
\put(-.6,3.5){$(\ast,+\infty)$}
\put(-1.3,-.5){$(0,0)$}
\put(2.3,1.6){$(+\infty,2)$}
\put(-5.3,-1.1){$(-\infty,-1)$}

\end{picture}
\caption{The compactification $\Rc$ of $\R^2$.}\label{fig:compac}
\end{center}
\end{figure}

We now introduce these objects formally. Let $\Rc$ be the compactification of
$\R^2$ obtained by equipping the set
$\Rc:=\R^2\cup\{(\pm\infty,t):t\in\R\}\cup\{(\ast,\pm\infty)\}$ with a
topology such that $(x_n,t_n)\to(\pm\infty,t)$ if $x_n\to\pm\infty$ and
$t_n\to t\in\R$, and $(x_n,t_n)\to(\ast,\pm\infty)$ if $t_n\to\pm\infty$
(regardless of the behavior of $x_n$). An explicit way to construct such a
compactification is as follows. Let $\Tet:\R^2\to\R^2$ be defined by
\be\label{Tet}
\Tet(x,t)=\big(\Tet_1(x,t),\Tet_2(t)\big)
:=\Big(\frac{\tanh(x)}{1+|t|},\tanh(t)\Big),
\ee
and let $\Tet(\R^2)$ denote the image of $\R^2$ under $\Tet$. Then the closure
of $\Tet(\R^2)$ in $\R^2$ is in a natural way isomorphic to $\Rc$ (see
Figure~\ref{fig:compac}).

By definition, a {\em path}\index{path} $\pi$ in $\Rc$ with {\em starting time} $\sig_\pi$
is a function $\pi:[\sig_\pi,\infty]\to[-\infty,\infty]\cup\{\ast\}$ such
that $t\mapsto(\pi(t),t)$ is a continuous map from $[\sig_\pi,\infty]$ to
$\Rc$. We will often view paths as subsets of $\Rc$, i.e., we identify a path
$\pi$ with its graph $\{(\pi(t),t):t\in[\sig_\pi,\infty]\}$. We let $\Pi$
denote the space of all paths in $\Rc$ with all possible starting times in
$[-\infty,\infty]$, equipped with the metric
\be\label{dPi}
d(\pi_1,\pi_2)
:=|\Tet_2(\sig_{\pi_1})-\Tet_2(\sig_{\pi_2})|
\vee\sup_{t\geq\sig_{\pi_1}\wedge\sig_{\pi_2}}
\big|\Tet_1\big(\pi_1(t\vee\sig_{\pi_1}),t)
-\Tet_1(\pi_2(t\vee\sig_{\pi_2}),t\big)\big|,
\ee
and we let $\Ki(\Pi)$ denote the space of all compact subsets $K\sub\Pi$,
equipped with the Hausdorff metric
\be\label{Hmet}
d_{\rm H}(K_1,K_2)=\sup_{x_1\in K_1}\inf_{x_2\in K_2}d(x_1,x_2)
\vee\sup_{x_2\in K_2}\inf_{x_1\in K_1}d(x_1,x_2).
\ee
Both $\Pi$ and $\Ki(\Pi)$ are complete separable metric spaces. The set
$\hat\Pi$ of all {\em dual paths} \index{path!dual} \index{dual!path}
$\hat\pi:[-\infty,\hat\sig_{\hat\pi}]\to[-\infty,\infty]\cup\{\ast\}$ with
{\em starting time} $\hat\sig_{\hat\pi}\in[-\infty,\infty]$ is defined
analoguously to $\Pi$.

We adopt the convention that if $f:\Rc\to\Rc$ and $A\sub\Rc$, then
$f(A):=\{f(z):z\in A\}$ denotes the image of $A$ under $f$. Likewise, if $\Ai$
is a set of subsets of $\Rc$ (e.g.\ a set of paths), then
$f(\Ai):=\{f(A):A\in\Ai\}$. This also applies to notation such as
$-A:=\{-z:z\in A\}$. If $\Ai\sub\Pi$ is a set of paths and $A\sub\Rc$, then we
let $\Ai(A):=\{\pi\in\Ai:(\pi(\sig_\pi),\sig_\pi)\in A\}$ denote the subspace
of all paths in $\Ai$ with starting points in $A$, and for $z\in\Rc$ we write
$\Ai(z):=\Ai(\{z\})$.

The next proposition, which follows from \cite[Theorem 2.1]{FINR04},
\cite[Theorem 3.7]{FINR06}, and \cite[Theorem 1.9]{SS08}, gives a
characterization of the Brownian web $\Wi$ and its dual $\hat\Wi$. Below, we
say that a path $\pi\in\Pi$ {\em crosses}\index{crossing!of dual paths} a dual path $\hat\pi\in\hat\Pi$ from
left to right if there exist $\sig_\pi\leq s<t\leq\hat\sig_{\hat\pi}$ such
that $\pi(s)<\hat\pi(s)$ and $\hat\pi(t)<\pi(t)$. Crossing from right to
left is defined analogously.

\index{Brownian web!dual}\index{dual!Brownian web}
\bp{\bf(Characterization of the Brownian web and its dual)}\label{P:char}
For each $\bet\in\R$, there exists a $\Ki(\Pi)\times\Ki(\hat\Pi)$-valued random
variable $(\Wi,\hat\Wi)$, called the double Brownian web with drift $\bet$,
whose distribution is uniquely determined by the following properties:
\begin{itemize}
\item[{\rm (a)}] For each deterministic $z\in\R^2$, almost surely there is a
  unique path $\pi_z\in\Wi(z)$ and a unique dual path $\hat \pi_z\in\hat\Wi(z)$.

\item[{\rm (b)}] For any deterministic countable dense subset $\Di
  \subset\R^2$, almost surely, $\Wi$ is the closure in $\Pi$ of
  $\{\pi_z:z\in\Di\}$ and $\hat\Wi$ is the closure in $\hat\Pi$ of $\{\hat
  \pi_z:z\in\Di\}$.

\item[{\rm (c)}] For any finite deterministic set of points
  $z_1,\ldots,z_k\in\R^2$, the paths $(\pi_{z_1},\ldots,\pi_{z_k})$ are
  distributed as a collection of coalescing Brownian motions, each with drift
  $\bet$.

\item[{\rm (d)}] For any deterministic $z\in\R^2$, the dual path $\hat \pi_z$ is
  the a.s.\ unique path in $\hat\Pi(z)$ that does not cross any path in $\Wi$.
\end{itemize}
\ep
If $(\Wi,\hat\Wi)$ is a double Brownian web as defined in
Proposition~\ref{P:char}, then we call $\Wi$ a Brownian web and $\hat\Wi$ the
associated {\it dual Brownian web}. Note that $\hat\Wi$ is a.s.\ uniquely determined
by $\Wi$. Although this is not obvious from the definition, the dual Brownian
web is indeed a Brownian web rotated by 180 degrees. Indeed,
$(\Wi,\hat\Wi)$ is equally distributed with $(-\hat\Wi,-\Wi)$.

\bd{\bf(Incoming and outgoing paths)}\label{D:inout}
\index{incoming and outgoing paths} \index{incoming and outgoing paths!strong equivalence} \index{equivalence!of paths, strong}
We say that a path $\pi\in\Pi$ is an {\em incoming path} at a point
$z=(x,t)\in\R^2$ if $\sig_\pi<t$ and $\pi(t)=x$. We say that $\pi$ is an {\em
  outgoing path} at $z$ if $\sig_\pi= t$ and $\pi(t)=x$. We say that two
incoming paths $\pi_1,\pi_2$ at $z$ are {\em strongly equivalent}, denoted as
$\pi_1=^z_{\rm in}\pi_2$, if $\pi_1=\pi_2$ on $[t-\eps,t]$ for some
$\eps>0$. For $z\in\R^2$, let $m_{\rm in}(z)$ denote the number of equivalence
classes of incoming paths in $\Wi$ at $z$ and let $m_{\rm out}(z)$ denote the
cardinality of $\Wi(z)$. Then $(m_{\rm in}(z),m_{\rm out}(z))$ is called the
{\em type} \index{special points!type} of the point $z$ in $\Wi$. The type $(\hat m_{\rm in}(z),\hat
m_{\rm out}(z))$ of a point $z$ in the dual Brownian web $\hat\Wi$ is defined
analogously.
\ed

\begin{figure}[htb]
\begin{center}
\includegraphics[width=9cm]{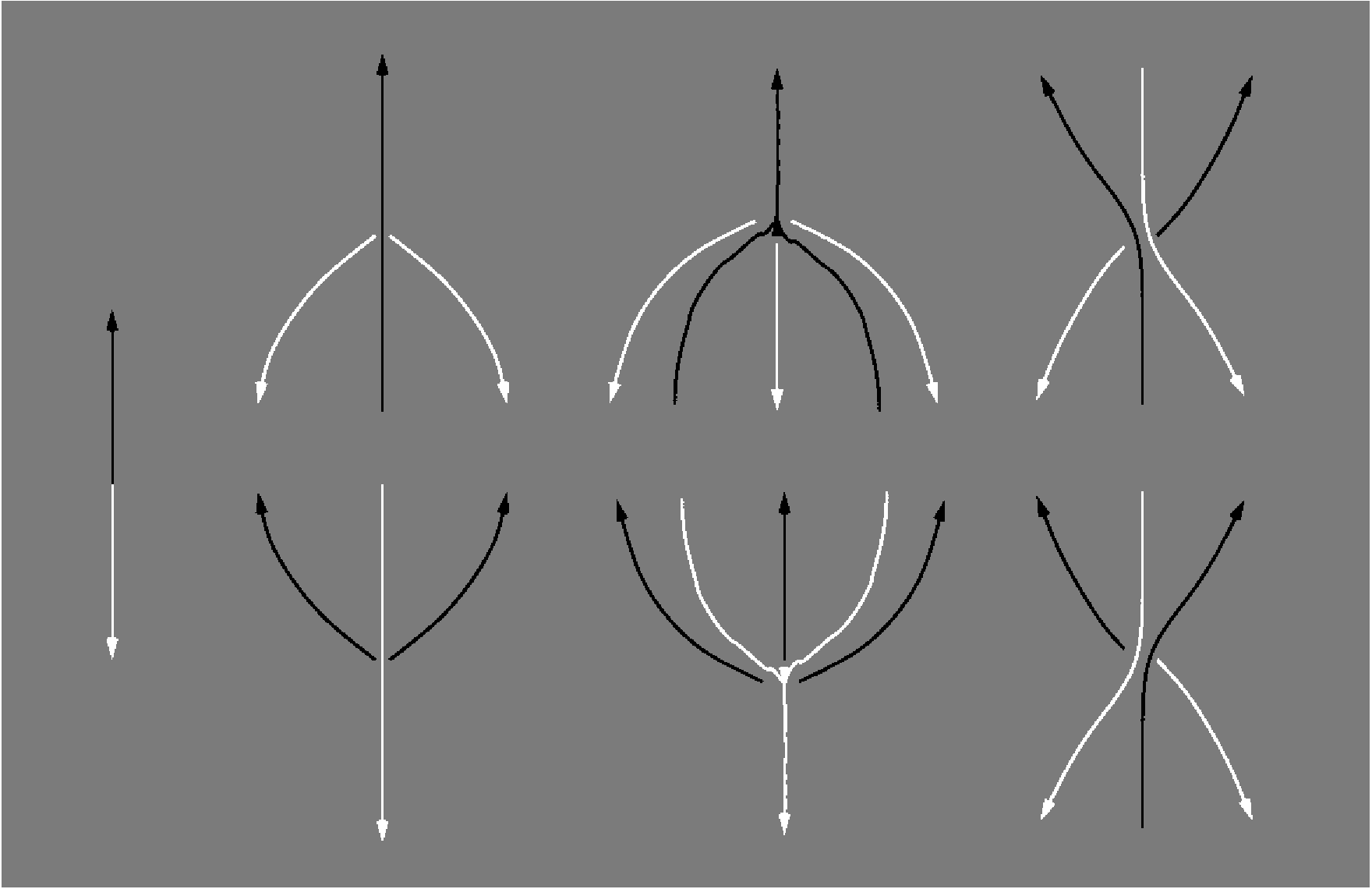}
\caption{Special points of the Brownian web. On the left: $(0,1)/(0,1)$. Top
  row: $(1,1)/(0,2)$, $(2,1)/(0,3)$, $(1,2)_{\rm l}/(1,2)_{\rm l}$. Bottom
  row: $(0,2)/(1,1)$, $(0,3)/(2,1)$, $(1,2)_{\rm r}/(1,2)_{\rm r}$.}
\label{fig:webpoints}
\end{center}
\end{figure}

We cite the following result from \cite[Proposition~2.4]{TW98} or
\cite[Theorems~3.11--3.14]{FINR06}. See Figure~\ref{fig:webpoints} for an
illustration.

\index{special points!Brownian web}
\bp\label{P:classweb}{\bf(Special points of the Brownian web)}
Almost surely, all points $z\in\R^2$ are of one of the following types in
$\Wi/\hat\Wi$: $(0,1)/(0,1)$, $(0,2)/(1,1)$, $(0,3)/(2,1)$, $(1,1)/(0,2)$,
$(1,2)/(1,2)$, and $(2,1)/(0,3)$. For each deterministic $t\in\R$, almost
surely, each point in $\R\times\{t\}$ is of type $(0,1)/(0,1)$, $(0,2)/(1,1)$,
or $(1,1)/(0,2)$. Deterministic points $z\in\R^2$ are a.s.\ of type
$(0,1)/(0,1)$.
\ep

For us, points of type $(1,2)/(1,2)$ are of special importance. Note that
these are the only points at which there are incoming paths both in $\Wi$ and
in $\hat\Wi$. Points of type $(1,2)$ in $\Wi$ are further distinguished into
points of type $(1,2)_{\rm l}$ and $(1,2)_{\rm r}$, according to whether the
left or the right outgoing path in $\Wi$ is the continuation of the (up to
equivalence unique) incoming path.

Proposition~\ref{P:classweb} shows that although for each deterministic
$z\in\R^2$, a.s.\ $\Wi(z)$ contains a single path, there exist random points
$z$ where $\Wi(z)$ contains up to three paths. Sometimes, it will be necessary
to choose a unique element of $\Wi(z)$ for each $z\in\R^2$. To that aim, for
each $z\in\R^2$, we let $\pi^+_z$ denote the right-most element of
$\Wi(z)$. We define $\pi^\up_z$ in the same way, except that at points of
type $(1,2)_{\rm l}$, we let $\pi^\up_z$ be the left-most element of
$\Wi(z)$. Note that as a consequence of this choice, whenever there are
incoming paths at $z$, the path $\pi^\up_z$ is the continuation of any
incoming path at $z$.

The next proposition, which follows from \cite[Prop.~3.1]{NRS10}, shows that
it is possible to define something like the intersection local time of $\Wi$ and
$\hat\Wi$. Below, $|I|$ denotes the Lebesgue measure of a set $I\sub\R$.

\index{intersection local time}
\bp{\bf(Intersection local time)}\label{P:refloc}
Let $(\Wi,\hat\Wi)$ be the double Brownian web. Then a.s.\ there exists a
unique measure $\ell$, concentrated on the set of points of type $(1,2)$ in
$\Wi$, such that for each $\pi\in\Wi$ and $\hat\pi\in\hat\Wi$,
\be\ba{l}\label{elldef}
\dis\ell\big(\big\{z=(x,t)\in\R^2:
\sig_\pi<t<\hat\sig_{\hat\pi},\ \pi(t)=x=\hat\pi(t)\big\}\big)\\[5pt]
\dis\qquad=\lim_{\eps\down 0}\eps^{-1}\big|\big\{t\in\R:
\sig_\pi<t<\hat\sig_{\hat\pi},\ |\pi(t)-\hat\pi(t)|\leq\eps\big\}\big|,
\ec
where the limit on the right-hand side exists and is finite. The measure
$\ell$ is a.s.\ non-atomic and $\si$-finite. We let $\ell_{\rm
  l}$ and $\ell_{\rm r}$ denote the restrictions of $\ell$ to the sets of
points of type $(1,2)_{\rm l}$ and $(1,2)_{\rm r}$, respectively.
\ep
We remark that $\ell(O)=\infty$ for every open nonempty subset $O\sub\R^2$,
but $\ell$ is \si-finite. To see the latter, for any path $\pi\in\Pi$, let
$\pi^\circ:=\{(\pi(t),t):t\in(\sig_\pi,\infty)\}$ denote its {\em interior},
and define the interior $\hat\pi^\circ$ of a dual path $\hat\pi\in\hat\Pi$
analogously. Let $\Di\sub\R^2$ be a deterministic countable dense set and for
$z\in\Di$, let $\pi_z$, resp.\ $\hat \pi_z$, denote the a.s.\ unique path in
$\Wi$, resp.\ $\hat\Wi$, starting from $z$. Then by
Proposition~\ref{P:refloc}, $\ell(\pi^\circ_z\cap\hat \pi^\circ_{\hat
  z})<\infty$ for each $z,\hat z\in\Di$, while by \cite[Lemma~3.4~(b)]{SS08},
$\ell$ is concentrated on
$\bigcup_{z,\hat z\in\Di}(\pi^\circ_z\cap\hat \pi^\circ_{\hat z})$.

\subsection{Sticky Brownian webs}\label{S:modweb}

We collect here some facts about a natural way to couple two Brownian
webs. Such coupled Brownian webs will then be used in the next subsection to
give a graphical construction of Howitt-Warren flows. We first start with a
`reference' Brownian web $\Wi$, which is then used to construct a second,
`modified' or `sample' Brownian web $\Wi'$ by `switching' a suitable Poisson
subset of points of type $(1,2)_{\rm l}$ of $\Wi$ into points of type
$(1,2)_{\rm r}$, and vice versa, using a marking procedure developed in
\cite{NRS10}.

To formulate this rigorously, let $z=(x,t)$ be a point of type
$(1,2)_{\rm l}$ in $\Wi$, and let
\be
\Wi_{\rm in}(z):=\{\pi\in\Wi:\sig_\pi<t,\ \pi(t)=x\}
\ee
denote the set of incoming paths in $\Wi$ at $z$. For any $\pi\in\Wi_{\rm
  in}(z)$, let $\pi^t:=\{(\pi(s),s):\sig_\pi\leq s\leq t\}$ denote the piece
of $\pi$ leading up to $z$, and let $\Wi(z)=\{l,r\}$ be the outgoing paths in
$\Wi$ at $z$, where $l<r$ on $(t,t+\eps)$ for some $\eps>0$. Since $z$ is of
type $(1,2)_{\rm l}$, identifying a path with its graph, we have
$\pi=\pi^t\cup l$ for each $\pi\in\Wi_{\rm in}(z)$. We define
\be\label{switch}\index{switching}
\switch_z(\Wi):=\big(\Wi\beh\Wi_{\rm in}(z)\big)
\cup\{\pi^t\cup r:\pi\in\Wi_{\rm in}(z)\}.
\ee
Then $\switch_z(\Wi)$ differs from $\Wi$ only in that $z$ is now of type
$(1,2)_{\rm r}$ instead of $(1,2)_{\rm l}$. In a similar way, if $z$ is of
type $(1,2)_{\rm r}$ in $\Wi$, then we let $\switch_z(\Wi)$ denote the
web obtained from $\Wi$ by switching $z$ into a point of type $(1,2)_{\rm
  l}$. If $z_1,\ldots,z_n$ are points of type $(1,2)$ in $\Wi$, then we let
$\switch_{\{z_1,\ldots,z_n\}}(\Wi)
:=\switch_{z_1}\circ\cdots\circ\switch_{z_n}(\Wi)$ denote the web
obtained from $\Wi$ by switching the orientation of the points
$z_1,\ldots,z_n$. Note that it does not matter in which order we perform the
switching. Recall that a point $z$ is of type $(1,2)_{\rm l}$
(resp.\ $(1,2)_{\rm r}$) in the dual Brownian web $\hat\Wi$ if and only if it
is of type $(1,2)_{\rm l}$ (resp.\ $(1,2)_{\rm r}$) in $\Wi$. We define
switching in $\hat\Wi$ analogously to switching in $\Wi$.

The next theorem, which is similar to \cite[Prop.~6.1]{NRS10}, shows how by
switching the orientation of a countable Poisson set of points of type
$(1,2)$, we can obtain a well-defined modified Brownian web. Recall the
definition of the intersection local time measure $\ell$ from
Proposition~\ref{P:refloc} and note that since $\ell$ is \si-finite, the set
$S$ below is a.s.\ a countable subset of the set of all points of type
$(1,2)$.

\bt{\bf(Modified Brownian web)}\label{T:webmod}
Let $\Wi$ be a Brownian web wih drift $\bet$, let $\ell$ be the intersection
local time measure between $\Wi$ and its dual and let $\ell_{\rm l},\ell_{\rm
  r}$ denote the restrictions of $\ell$ to the sets of points of type
$(1,2)_{\rm l}$ and $(1,2)_{\rm r}$ in $\Wi$, respectively. Let $c_{\rm
  l},c_{\rm r}\geq 0$ be constants and conditional on $\Wi$, let $S$ be a
Poisson point set with intensity $c_{\rm l}\ell_{\rm l}+c_{\rm r}\ell_{\rm
  r}$. Then, a.s., for any sequence of finite sets $\De_n\up S$, the limit
\be\label{webmod}
(\Wi',\hat\Wi'):=\lim_{\De_n\up S}
\big(\switch_{\De_n}(\Wi),\switch_{\De_n}(\hat\Wi)\big)
\ee
exists in $\Ki(\Pi)\times\Ki(\hat\Pi)$ and does not depend on the choice
of the sequence $\De_n\up S$. Moreover, $\Wi'$ is a Brownian web with drift
$\bet'=\bet+c_{\rm l}-c_{\rm r}$ and $\hat\Wi'$ is its dual.
\et

\noi
{\bf Remark.} We recall that a countable set $S\sub\R^2$ is a Poisson point
set with \si-finite intensity $\mu$ if $S\cap A_n$ a Poisson point set with
intensity $\mu(A_n\cap\,\cdot\,)$ for some, and hence for every sequence of
measurable sets $A_n$ such that $\mu(A_n)<\infty$ for all $n$ and
$\mu(\R^2\beh\bigcup_nA_n)=0$. We apply this to the case that $\mu$ is the
random measure $c_{\rm l}\ell_{\rm l}+c_{\rm r}\ell_{\rm r}$ and the $A_n$ are
finite unions of intersections of forward and dual paths, started from
deterministic points, as mentioned below Proposition~\ref{P:refloc}. In
particular, when we say that $S$ is Poisson with intensity $c_{\rm l}\ell_{\rm
  l}+c_{\rm r}\ell_{\rm r}$, this should be interpreted in this particular
sense. Some care is needed when talking about the conditional law of $S$ given
$\Wi$, since it is not clear whether $S$ (being a countable dense subset of
$\R^2$), on its own, can be viewed as a random variable with values in a
decent (at least measurable) space. Nevertheless, it is not hard to see that
the triple $(\Wi,\hat\Wi,S)$ (being a marked double Brownian web) can be
constructed as a legitimate random variable on a suitable probability space
and that $(\Wi',\hat\Wi')$ is a.s.\ a measurable function of
$(\Wi,\hat\Wi,S)$.\med

If $(\Wi,\Wi')$ are coupled as in Theorem~\ref{T:webmod}, then we say that
$(\Wi,\Wi')$ is a {\em pair of sticky Brownian webs} with drifts $\bet,\bet'$
and {\em coupling parameter} $\kappa:=\min\{c_{\rm l},c_{\rm
  r}\}$. \index{Brownian web!sticky} \index{Brownian web!modified}
\index{coupling parameter} In the special case that $\kappa=0$ and
$\bet\leq\bet'$, we call $(\Wi,\Wi')$ a {\em left-right Brownian web} with
drifts $\bet,\bet'$. \index{Brownian web!left-right}
\index{left-right!Brownian web}Left-right Brownian webs have been introduced
with the help of a `left-right stochastic differential equation' (instead of
the marking construction above) in \cite{SS08}. Pairs of sticky Brownian webs
with general coupling parameters $\kappa\geq 0$ have been introduced by means
of a martingale problem in \cite[Section~7]{HW09b}. They are, indeed, sticky in
the sense that a pair of paths, one from each web, are Brownian motions with
sticky interaction. We will prove in Lemma~\ref{L:lreq} below that the
constructions of left-right Brownian webs given above and in \cite{SS08} are
equivalent. We will not make use of the martingale formulation of sticky
Brownian webs developed in \cite{HW09b}.

\detail{For $\bet=\bet'$ the equivalence of our construction of sticky
Brownian webs and the construction given in \cite[Section~7]{HW09b}
follows from our results on convergence of 2-point motions. For
$\kappa=0$ is is probably straightforward to check that left-right
Brownian webs as defined in \cite{SS08} solve the martingale problem
in \cite[Section~7]{HW09b}. In general, the easiest way to prove that
two paths sampled from a sticky Brownian web as defined by us solve
the martingale problem of \cite[Section~7]{HW09b} is probably to use
discrete approximation. Howitt and Warren say that this martingale
problem `together with standard localozation arguments' can be used to
characterize sticky Brownian webs. Except in the case $\kappa=0$, I am
actually not sure these `standard localozation arguments' are so
straightforward as they make it sound.}

For any point $z$ of type $(1,2)$ in some Brownian web $\Wi$, we call
\be\label{signz}
\sign_\Wi(z):=\left\{\ba{ll}
-1\quad&\mbox{if $z$ is of type $(1,2)_{\rm l}$ in $\Wi$},\\
+1\quad&\mbox{if $z$ is of type $(1,2)_{\rm r}$ in $\Wi$}
\ea\right.
\ee
the {\em sign}\index{sign of point of type $(1,2)$} of $z$ in $\Wi$. If $(\Wi,\Wi')$ is a pair of sticky Brownian
webs, then it is known \cite[Thm.~1.7]{SSS09} that the set of points of type
$(1,2)$ in $\Wi$ in general does not coincide with the set of points of type
$(1,2)$ in $\Wi'$. The next proposition says that nevertheless, in the sense
of intersection local time measure, almost all points of type $(1,2)$ in $\Wi$
are also of type $(1,2)$ in $\Wi'$, and these point have the orientation one
expects.

\bp{\bf(Change of reference web)}\label{P:refchange}
In the setup of Theorem~\ref{T:webmod}, let
$\ell'$ be the intersection local time measure between $\Wi'$ and its dual and
let $\ell'_{\rm l},\ell'_{\rm r}$ denote the restrictions of $\ell'$ to the
sets of points of type $(1,2)_{\rm l}$ and $(1,2)_{\rm r}$ in $\Wi'$,
respectively. Then:
\begin{itemize}
\item[{\rm(i)}] Almost surely, $\ell'_{\rm l}=\ell_{\rm l}$ and $\ell'_{\rm
  r}=\ell_{\rm r}$.
\item[{\rm(ii)}] $S=\big\{z\in\R^2:z\mbox{ is of type $(1,2)$ in both $\Wi$ and $\Wi'$, and } \sign_\Wi(z)\neq\sign_{\Wi'}(z)\big\}$ a.s.
\item[{\rm(iii)}] Conditional on $\Wi'$, the set $S$ is a Poisson point set
  with intensity $c_{\rm r}\ell'_{\rm l}+c_{\rm l}\ell'_{\rm r}$ and
  $\Wi=\lim_{\De_n\up S}\switch_{\De_n}(\Wi')$.
\end{itemize}
\ep

Let $(\Wi_0,\Wi)$ be a pair of sticky Brownian webs with drifts $\bet_0,\bet$
and coupling parameter $\kappa$, and let $\pi^\up_z$ and $\pi^+_z$ denote the
special paths in $\Wi(z)$ defined below Proposition~\ref{P:classweb}. Let
\be\label{erosion}
K^\up_{s,t}(x,A):=\P\big[\pi^\up_{(x,s)}(t)\in A\,\big|\,\Wi_0\big]
\qquad\big(s\leq t,\ x\in\R,\ A\in\Bi(\R)\big),
\ee
and let $K^+_{s,t}(x,A)$ be defined similarly, with $\pi^\up_z$ replaced
by $\pi^+_z$. Then, as we will see in Theorem~\ref{T:HWconst} below,
$(K^\up_{s,t})_{s\leq t}$ and $(K^+_{s,t})_{s\leq t}$ are versions of the
Howitt-Warren flow with drift $\bet$ and characteristic measure $\nu=c_{\rm
  l}\de_1+c_{\rm r}\de_0$, where
\be
c_{\rm l}:=\kappa+\max\{0,\bet'-\bet\}
\quad\mbox{and}\quad
c_{\rm r}:=\kappa+\max\{0,\bet-\bet'\}.
\ee
In the special case that $c_{\rm l}=c_{\rm r}$, this was proved in
\cite[formula~(5)]{HW09b}. In the next subsection, we set out to give a similar
construction for {\em any} Howitt-Warren flow.

\detail{In \cite{NRS10}, the authors define
\be\ba{l}\label{Lidef}
\dis\Li\big(\big\{z=(x,t)\in\R^2:
\sig_\pi<t<\hat\sig_{\hat\pi},\ \pi(t)=\hat\pi(t)\big\}\big)\\[5pt]
\dis\qquad=\lim_{\eps\down 0}(2\eps)^{-1}\big|\big\{t\in\R:
\sig_\pi<t<\hat\sig_{\hat\pi},\ |\pi(t)-\hat\pi(t)|
\leq(\sqrt{2}\eps)\big\}\big|\\[5pt]
\dis\qquad=\lim_{\eps'\down 0}(\sqrt{2}\eps')^{-1}\big|\big\{t\in\R:
\sig_\pi<t<\hat\sig_{\hat\pi},\ |\pi(t)-\hat\pi(t)|\leq\eps'\big\}\big|\\[5pt]
\dis\qquad=\ffrac{1}{\sqrt{2}}\ell\big(\big\{z=(x,t)\in\R^2:
\sig_\pi<t<\hat\sig_{\hat\pi},\ \pi(t)=\hat\pi(t)\big\}\big).
\ec
Then they define a marking $\Mi(\tau)$ which is a Poisson point process on
$\R^2\times\half$ with intensity $\sqrt{2}\Li(\di z)\di\tau=\ell(\di
z)\di\tau$. Now their Prop.~5.3 shows that switching only the
$\Mi(\tau)$-marked points of type $(1,2)_{\rm l}$ yields a drift of size
$\tau$.}

\subsection{Marking construction of Howitt-Warren flows}\label{S:flowcons}

We now give a construction of a general Howitt-Warren flow based on two
coupled webs, which is the central result of this paper. More precisely, the
random environment of the Howitt-Warren flow will be represented by a
reference Brownian web $\Wi_0$ plus a set $\Mi$ of marked points of type
$(1,2)$.  Conditional on $(\Wi_0,\Mi)$, we will modify $\Wi_0$ in a similar
way as in Theorem~\ref{T:webmod} to construct a sample Brownian web $\Wi$,
whose law conditional on $(\Wi_0,\Mi)$ then defines the Howitt-Warren flow via
the continuous analogue of (\ref{KQ}). For erosion flows, the set $\Mi$ of
marked points is empty, hence our representation reduces to (\ref{erosion}).

\index{Brownian web!reference} \index{Brownian web!sample}

Let $\Wi_0$ be a Brownian web with drift $\bet_0$ and let $\nu_{\rm l}$ and
$\nu_{\rm r}$ be finite measures on $[0,1]$.  Let $\ell$, $\ell_{\rm l}$ and
$\ell_{\rm r}$ be defined for $\Wi_0$ as in Proposition~\ref{P:refloc}, and
conditional on $\Wi_0$, let $\Mi$ be a Poisson point set on $\R^2\times[0,1]$
with intensity
\be\label{Poisint}
\ell_{\rm l}(\di z)\otimes2\,1_{\{0<q\}}q^{-1}\nu_{\rm l}(\di q)
+\ell_{\rm r}(\di z)\otimes2\,1_{\{q<1\}}(1-q)^{-1}\nu_{\rm r}(\di q).
\ee
Elements of $\Mi$ are pairs $(z,q)$ where $z$ is a point of type $(1,2)$ in
$\Wi_0$ and $q\in[0,1]$. Since $\ell$ is non-atomic, for each point $z$ of type
$(1,2)$ there is at most one $q$ such that $(z,q)\in\Mi$, and we may write
$\Mi=\{(z,\om_z):z\in M\}$. We call points $z\in M$ {\em marked points} and we
call $\om_z$ the {\em mark} of $z$. \index{marked points, set of}

Conditional on the reference web $\Wi_0$ and the set of marked points $\Mi$, we
construct independent $\{-1,+1\}$-valued random variables
$(\al_z)_{z\in M}$ with $\P[\al_z=+1|(\Wi_0,\Mi)]=\om_z$, and we set
\be\label{Adef}
A:=\{z\in M:\al_z\neq\sign_{\Wi_0}(z)\big\}.
\ee
In addition, conditional on $(\Wi_0,\Mi)$, we let $B$ be a Poisson point set
with intensity $2\nu_{\rm l}(\{0\})\ell_{\rm l}+2\nu_{\rm r}(\{1\})\ell_{\rm
  r}$, independent of $A$. We observe that conditional on $\Wi_0$, but
integrating out the randomness of $\Mi$, the set $A\cup B$ is a Poisson point
set with intensity
\be\ba{l}\label{ABint}
\dis2\Big(\int_{(0,1]}qq^{-1}\nu_{\rm l}(\di q)
+\nu_{\rm l}(\{0\})\Big)\ell_{\rm l}
+2\Big(\int_{[0,1)}(1-q)(1-q)^{-1}\nu_{\rm r}(\di q)
+\nu_{\rm r}(\{1\})\Big)\ell_{\rm r}\\[15pt]
\qquad\dis=2\nu_{\rm l}\big([0,1]\big)\ell_{\rm l}+
2\nu_{\rm r}\big([0,1]\big)\ell_{\rm r}.
\ec
Therefore, by Theorem~\ref{T:webmod}, a.s.\ the limit
\be\label{sample}
\Wi:=\lim_{\De_n\up A\cup B}\switch_{\De_n}(\Wi_0)
\ee
exists in $\Ki(\Pi)$ and is a Brownian web with drift
\be\label{tib}
\bet:=\bet_0+2\nu_{\rm l}\big([0,1]\big)-2\nu_{\rm r}\big([0,1]\big).
\ee

\index{Howitt-Warren!flow!construction}
\index{flow!Howitt-Warren!construction}
\bt{\bf(Construction of Howitt-Warren flows)}\label{T:HWconst}
Let $\beta\in\R$ and let $\nu$ be a finite measure on $[0,1]$. For any finite
measures $\nu_{\rm l}$ and $\nu_{\rm r}$ on $[0,1]$ satisfying
\be\label{nuform}
\nu(\di q):=(1-q)\nu_{\rm l}(\di q)+q\nu_{\rm r}(\di q),
\ee
let $\beta_0$ be determined from $\beta$, $\nu_{\rm l}$ and $\nu_{\rm r}$ as
in (\ref{tib}). Let $\Wi_0$ be a reference Brownian web with drift $\bet_0$
and define a set of marked points $\Mi$ and sample Brownian web $\Wi$ as in
(\ref{Poisint}) and (\ref{sample}). Let $\pi^\up_z$ and $\pi^+_z$ denote the
special paths in $\Wi(z)$ defined below Proposition~\ref{P:classweb}. Set
\bc\label{HWconst}
K^\up_{s,t}(x,A):=\P\big[\pi^\up_{(x,s)}(t)\in A\,\big|\,(\Wi_0,\Mi)\big]
\qquad\big(s\leq t,\ x\in\R,\ A\in\Bi(\R)\big)
\ec
and define $K^+_{s,t}(x,A)$ similarly with $\pi^\up_{(x,s)}$ replaced by
$\pi^+_{(x,s)}$. Then $(K^\up_{s,t})_{s\leq t}$ and $(K^+_{s,t})_{s\leq t}$
are versions of the Howitt-Warren flow with drift $\bet$ and characteristic
measure $\nu$.  In the special case that $\nu_{\rm l}=\nu_{\rm r}$, the triple
$(\Wi_0,\Mi,\Wi)$ is equally distributed with $(\Wi,\Mi,\Wi_0)$.
\et

\noi
{\bf Remark.} Since $\Wi$ take values in the Polish space $\Ki(\Pi)$, we may
construct a regular version of the conditional probability
$\P[\Wi\in\cdot\,|(\Wi_0,\Mi)]$, which is a random probability measure on
$\Ki(\Pi)$. A.s., under this random law, the random set of paths $\Wi$ has the
same a.s.\ properties as a Brownian web. In particular, for a.e.\ $\om$ in our
underlying probability space, the paths $\pi^\up_z$ and $\pi^+_z$ in $\Wi$ are
well-defined for all $z\in\R^2$ and we obtain a version of the
conditional probabilities in formula (\ref{HWconst}) for all
$(x,s)\in\R^2$ and $t\geq s$ simultaneously. Interpreting formula
(\ref{HWconst}) in this way (as we will always do), we obtain versions
$(K^\up_{s,t})_{s\leq t}$ and $(K^+_{s,t})_{s\leq t}$ of our Howitt-Warren
flows with different properties, see Proposition~\ref{P:regul} below.\med

\noi
{\bf Remark.} Note that in (\ref{nuform}), it is always possible to choose
$\nu_{\rm l}=\nu_{\rm r}=\nu$. In this case, the construction of the reference
Brownian web and set of marked points arises as the scaling limit of the
discrete construction outlined in Section~\ref{S:heur}. If $\om^{\li
  k\re}=(\om^{\li k\re}_z)_{z\in\Zev}$ is a collection of independent
$[0,1]$-valued random variables with laws $\mu_k$ satisfying (\ref{mucon}) and
conditional on $\om^{\li k\re}$ we construct two independent collections
$\al^{\li k\re}=(\al^{\li k\re}_z)_{z\in\Zev}$ and
$\varal^{\li k\re}=(\varal^{\li k\re}_z)_{z\in\Zev}$ of
$\{-1,+1\}$-valued random variables with
$\P[\al^{\li k\re}_z=+1\,|\,\om^{\li k\re}]=\om^{\li k\re}_z$, and similar for
$\varal^{\li k\re}$, then the discrete webs corresponding to $\al^{\li k\re}$
and $\varal^{\li k\re}$ converge after diffusive rescaling to $\Wi$ and
$\Wi_0$, and $\{(z,\om_z):0<\om_z<1\}$ converges to the set of marked points
$\Mi$. The more general construction in Theorem~\ref{T:HWconst} where possibly
$\nu_{\rm l}\neq\nu_{\rm r}$ and $\bet_0$ is possibly different from $\bet$
arises as the diffusive scaling limit of discrete constructions where we first
choose a reference collection of random variables
$\varal^{\li k\re}$ with $\eps_k^{-1}\E[\varal^{\li k\re}_z]\to\bet_0$
and then conditional on $\varal^{\li k\re}$, we choose
independent $(\om^{\li k\re}_z)_{z\in\Zev}$ with
$\P[\om_z\in\di q\,|\,\varal^{\li k\re}=-1]=\mu^{\rm l}_k(\di q)$ and
$\P[\om_z\in\di q\,|\,\varal^{\li k\re}=+1]=\mu^{\rm r}_k(\di q)$, where
generalizing (\ref{mulr}), $\mu^{\rm l}_k$ and $\mu^{\rm r}_k$ are any laws
such that
\be
\P[\varal^{\li k\re}=-1]\mu^{\rm l}_k(\di q)+
\P[\varal^{\li k\re}=+1]\mu^{\rm r}_k(\di q)=\mu_k(\di q).
\ee
This more general construction will sometimes be handy. For example, for erosion flows where $\nu=c_0\de_0+c_1\de_1$ for some $c_0,c_1\geq 0$, it is most natural to choose $\nu_{\rm l}=c_0\de_0$ and $\nu_{\rm r}=c_1\de_1$. Also, for Howitt-Warren flows where one or both of the speeds $\bet_-,\bet_+$ defined in (\ref{speeds}) are finite, it is sometimes handy to choose either $\nu_{\rm l}=0$ or $\nu_{\rm r}=0$.\med

If $(\Wi_0,\Mi,\Wi)$ are a reference Brownian web, the set of marked points, and
the sample Brownian web as defined above Theorem~\ref{T:HWconst}, then we call the random probability measure $\Q$ on $\Ki(\Pi)$ defined by
\be\label{HWquen}
\Q:=\P\big[\Wi\in\cdot\,\big|\,(\Wi_0,\Mi)\big]
\ee
the {\em Howitt-Warren quenched law} with drift $\bet$ and characteristic
measure $\nu$. \index{Howitt-Warren!quenched law}\index{quenched law!Howitt-Warren} In Section~\ref{S:approx} below, we will show that in some
precisely defined way, these Howitt-Warren quenched laws are the diffusive
scaling limits of the discrete quenched laws defined in
Section~\ref{S:disref}.

Since at deterministic points in the Brownian web $\Wi$ there is a.s.\ only
one outgoing path, the stochastic flows of kernels $(K^\up_{s,t})_{s\leq t}$
and $(K^+_{s,t})_{s\leq t}$ from Theorem~\ref{T:HWconst} obviously have the
same finite-dimensional distributions. They are, however, not the same. Each
version has its own pleasant properties.
\bp{\bf(Regular parameter dependence)}\label{P:regul}
Let $(K^\up_{s,t})_{s\leq t}$ and $(K^+_{s,t})_{s\leq t}$ be defined as in
Theorem~\ref{T:HWconst}. Then, of the following properties,
$(K^+_{s,t})_{s\leq t}$ satisfies (a)--(c) and $(K^\up_{s,t})_{s\leq t}$
satisfies (a), (b) and (d).
\begin{itemize}
\item[{\bf(a)}] Setting $\R^2_\leq:=\{(s,t)\in\R^2:s\leq t\}$, the map
  $(s,t,x,\om)\mapsto K_{s,t}(x,\,\cdot\,)(\om)$ is a measurable map from
  $\R^2_\leq\times\R\times\Om$ to $\Mi_1(\R)$.

\item[{\bf(b)}] A.s., the map
  $t\mapsto K_{s,t}(x,\,\cdot\,)$ from $[s,\infty)$ to $\Mi_1(\R)$ is
  continuous for all $s\in\R$ and $x\in\R$.

\item[{\bf(c)}] A.s., $x\mapsto K_{s,t}(x,A)$ is a c\`adl\`ag function from
  $\R$ to $\R$ for each $s<t$ and $A\in \Bi(\R)$.

\item[{\bf(d)}] A.s., $\dis\int_{\R}K_{s,t}(x,\di
  y)K_{t,u}(y,A)=K_{s,u}(x,A)$ for all $s\leq t\leq u$, $x\in\R$, and
  $A\in\Bi(\R)$.
\end{itemize}
\ep
Proposition~\ref{P:regul}~(b) and (d) show that $(K^\up_{s,t})_{s\leq
t}$ yields a version of a Howitt-Warren flow with the properties
listed in Proposition~\ref{P:conpath}. In particular,
Proposition~\ref{P:regul}~(d) makes it a family of bona fide
transition probability kernels of a Markov process in a random
space-time environment.

\subsection{Discrete approximation}\label{S:approx}

Recall the definitions of the discrete quenched laws $\Qdis$ in (\ref{quench})
and the Howitt-Warren quenched laws $\Q$ in (\ref{HWquen}). In this
subsection, we formulate a convergence result which says that if $\mu_k$ is a
sequence of probability laws on $[0,1]$ satisfying (\ref{mucon}), then the
associated discrete quenched laws $\Qdis_{\li k\re}$, diffusively rescaled,
converge to the Howitt-Warren quenched law $\Q$ with drift $\bet$ and
characteristic measure $\nu$. This abstract result then implies other
convergence results such as the convergence of Howitt-Warren flows,
Howitt-Warren processes, and $n$-point motions. Since the $n$-point motions
of discrete Howitt-Warren flows will be shown in Proposition~\ref{P:nconv} to converge
to solutions of the Howitt-Warren martingale problems, this will also verify that
the flows we constructed in Theorem \ref{T:HWconst} are indeed versions of
the Howitt-Warren flow.

To formulate our convergence statement properly, we need to identify a
discrete quenched law $\Qdis$ with a random probability law on the space
$\Ki(\Pi)$ of compact subsets of the space $\Pi$ of paths defined in
Section~\ref{S:BMweb}. Recall the definition of the paths $p^\al_z$ in
(\ref{pal}) and the discrete webs $\Ui^\al$ in (\ref{Uial}). We wish to view
$\Ui^\al$ as a random variable with values in $\Ki(\Pi)$. To this aim, we {\em
  modify our definition of $\Ui^\al$ as follows}. First, for each
$z=(x,s)\in\Zev$ we make $p^\al_z(t)$ into a path in $\Pi$ by linear
interpolation between integer times and by setting
$p^\al_z(\infty):=\ast$. Next, we add to $\Ui^\al$ all trivial paths $\pi$,
with starting times $\sig_\pi\in\Z\cup\{-\infty,\infty\}$, such that $\pi$ is
identically $-\infty$ or $+\infty$ on $[\sig_\pi,\infty)\cap\R$.
With this modified definition, it can be checked that $\Ui^\al$ is indeed a
random compact subset of $\Pi$, as desired.

For $\eps>0$, we let $S_\eps:\Rc\to\Rc$ denote the scaling map
\be\label{Seps2}
S_\eps(x,t):=(\eps x,\eps^2 t)\qquad\big((x,t)\in\Rc\big).
\ee
As usual, we identify paths with their graphs; then $S_\eps(\pi)$ is the path
obtained by diffusively rescaling a path $\pi$ with $\eps$, and
$S_\eps(\Ui^\al)$ is the random collection of paths obtained by diffusively
rescaling paths in $\Ui^\al$. If $\Qdis^\om$ is a discrete quenched law as
defined in (\ref{quench}) and $\eps>0$, then we write
\be
S_\eps(\Qdis^\om):=\Qdis^\om\big[S_\eps(\Ui^\al)\in\cdot\,\big],
\ee
i.e., $S_\eps(\Qdis^\om)$ is the image under the scaling map $S_\eps$ of the
quenched law of $\Ui^\al$. Note that $S_\eps(\Qdis^\om)$, so defined, is a
random probability law on the space of compact subsets of $\Pi$, i.e., a
random variable with values in $\Mi_1(\Ki(\Pi))$.

\bt{\bf(Convergence of quenched laws)}\label{T:quench}
Let $\eps_k$ be positive constants, converging to zero and $\mu_k$ be
probability measures on $[0,1]$ satisfying (\ref{mucon}) for some real $\bet$
and finite measure $\nu$ on $[0,1]$. Let $\om^{\li k\re}=(\om^{\li
  k\re}_z)_{z\in\Zev}$ be i.i.d.\ $[0,1]$-valued random variables with
distibution $\mu_k$, let $\Qdis_{\li k\re}:=\Qdis^{\om^{\li k\re}}$ be the
discrete quenched law defined in (\ref{quench}), and let $\Q$ be the
Howitt-Warren quenched law with drift $\bet$ and characteristic measure $\nu$
defined in (\ref{HWquen}). Then
\be\label{quencon}
\P\big[S_{\eps_k}(\Qdis_{\li k\re})\in\cdot\,\big]\Asto{k}\P[\Q\in\cdot\,],
\ee
where $\Rightarrow$ denotes weak convergence of probability laws on
$\Mi_1(\Ki(\Pi))$.
\et

\section{Construction of Howitt-Warren flows in the Brownian net}\label{S:net}

In this section, we show that when a Howitt-Warren flow with drift $\beta$ and
characteristic measure $\nu$ has finite left and right speeds, or
equivalently, $b:=\int q^{-1}(1-q)^{-1}\nu(\di q)<\infty$, then we can
alternatively construct the flow as a random flow of mass in the Brownian
net. Analogous to Theorem \ref{T:HWconst}, the random environment will now be
represented as a Brownian net $\Ni$ plus a set of i.i.d.\ marks
$\bar\om:=(\bar\om_z)_{z\in S}$ attached to the separation points $S$ of
$\Ni$, each with law $\bar\nu(\di q):=b^{-1}q^{-1}(1-q)^{-1}\nu(\di
q)$. Conditional on $(\Ni, \bar\om)$, we can construct the sample web $\Wi$ by
choosing trajectories in $\Ni$ that turn in the `right' way at separation
points. The Howitt-Warren flow is then defined from the law of $\Wi$
conditional on $(\Ni, \bar\om)$ as in (\ref{HWconst}).

In Sections~\ref{S:netdef}--\ref{S:separ}, we recall the necessary background
on the Brownian net and properties of its separation points. In
Section~\ref{S:netmark}, we first state some coupling results between the
Brownian web and Brownian net, which will help shed more light on the marking
constructions of sticky Brownian webs in Theorem~\ref{T:webmod}. In
Section~\ref{S:marknet} we then give our main result on the alternative
construction of Howitt-Warren flows with finite left and right speeds using
the Brownian net. Lastly in Section~\ref{S:quensup}, we formulate what we call
Brownian half-nets, and state some properties for the support of the
Howitt-Warren quenched law defined in (\ref{HWquen}), which will imply
Theorems~\ref{T:speed} and \ref{T:supp}. We note that, apart from being used
to construct Howitt-Warren flows with finite left and right speeds, the theory
of the Brownian net will also play an important role for Howitt-Warren flows
with infinite left or right speed, such as in the proof of results in
Section~\ref{S:quensup}, as well as in the proof of Theorems \ref{T:webmod}
and \ref{T:quench}.

\subsection{The Brownian net}\label{S:netdef}

The Brownian net arises as the diffusive scaling limit of branching-coalescing
random walks in the limit of small branching probability. It was
first introduced by Sun and Swart in \cite{SS08} and independently by
Newman, Ravishankar and Schertzer in \cite{NRS10}. A further study of its
properties was carried out in \cite{SSS09}. We now recall the definition of the Brownian net
given in \cite{SS08}.

Recall that in Section~\ref{S:modweb}, we defined a left-right Brownian web to
be a pair of sticky Brownian webs $(\Wl,\Wr)$ with drifts $\bet_-\leq\bet_+$
and coupling parameter $\kappa=0$.
At present, we will need the original
definition of a left-right Brownian web given in \cite{SS08}. In
Lemma~\ref{L:lreq} below, we will prove that both definitions are equivalent.


Following \cite{SS08}, we call $(l_1,\ldots,l_m;r_1,\ldots,r_n)$ a collection
of {\em left-right coalescing Brownian motions} with drifts
$\bet_-\leq\bet_+$, if $(l_1,\ldots,l_m)$ and $(r_1,\ldots,r_n)$ are
distributed as collections of coalescing Brownian motions with drift $\bet_-$
and $\bet_+$, respectively, if paths in $(l_1,\ldots,l_m;r_1,\ldots,r_n)$
evolve independently when they are apart, and the interaction between $l_i$
and $r_j$ when they meet is described by the two-dimensional stochastic
differential equation
\bc\label{lrsde}
\dis\di L_t\ &=&\dis 1_{\{L_t\neq R_t\}}\di B^{\rm l}_t
+1_{\{L_t=R_t\}}\di B^{\rm s}_t+\bet_-\di t,\\[5pt]
\dis\di R_t\ &=&\dis 1_{\{L_t\neq R_t\}}\di B^{\rm r}_t
+1_{\{L_t=R_t\}}\di B^{\rm s}_t+\bet_+\di t,
\ec
where $B^{\rm l}_t,B^{\rm r}_t,B^{\rm s}_t$ are independent standard
Brownian motions, and $(L,R)$ are subject to the constraint that
\be\label{lrorder}
L_t\leq R_t\mbox{ for all }t\geq\inf\{s:L_s=R_s\}.
\ee
It can be shown that subject to the condition (\ref{lrorder}), solutions to
the SDE (\ref{lrsde}) are unique in distribution \cite[Proposition 2.1]{SS08}.

Let $\Wl,\Wr$ be two Brownian webs with drifts $\bet_-\leq\bet_+$, and for
determinstic $z\in\R^2$, let $l_z$ resp.\ $r_z$ denote the a.s.\ unique path
in $\Wl$ resp.\ $\Wr$ starting from $z$. Following \cite{SS08}, we say that
$(\Wl,\Wr)$ is a {\em left-right Brownian web}\index{Brownian web!left-right}
\index{left-right!Brownian web} if for any finite deterministic set of points
$z_1,\ldots,z_m, z'_1,\ldots,z'_n\in\R^2$, the collection
$(l_{z_1},\ldots,l_{z_m};r_{z'_1},\ldots,r_{z'_n})$ is distributed as
left-right coalescing Brownian motions. Elements of $\Wl$ (resp.\ $\Wr$) are
called {\em left-most}\index{left-right!paths} (resp.\ {\em right-most})
paths. It is known \cite[formula~(1.22)]{SS08} that if $(\Wl,\Wr)$ is a
left-right Brownian web and $\hat\Wl,\hat\Wr$ are the dual Brownian webs
associated with $\Wl,\Wr$, then $(-\hat\Wl,-\hat\Wr)$ is equally distributed
with $(\Wl,\Wr)$.

It was shown in \cite{SS08} that each left-right Brownian web a.s.\ determines
an associated Brownian net and vice versa. There, three different ways were
given to construct a Brownian net from its associated left-right Brownian web,
which are known as the {\em hopping construction} and the constructions using
{\em wedges} and {\em meshes}, which we recall now.\med

\noi
{\bf Hopping:} We call $t\in\R$ an {\em intersection time}\index{intersection time} of two paths
$\pi,\pi'\in\Pi$ if $\sig_\pi,\sig_{\pi'}<t<\infty$ and $\pi(t)=\pi'(t)$. If
$t$ is an intersection time of $\pi$ and $\pi'$, then we can define a new path
$\pi''$ by concatenating the piece of $\pi$ before $t$ with the piece of
$\pi'$ after $t$, i.e., by setting
$\pi'':=\{(\pi(s),s):s\in[\sig_\pi,t]\}\cup\{(\pi'(s),s):s\in[t,\infty]\}$. For
any collection of paths $\Ai\sub\Pi$, we let $\Hi_{\rm int}(\Ai)$ denote the
smallest set of paths containing $\Ai$ that is closed under such `hopping'
from one path onto another at intersection times, i.e., $\Hi_{\rm int}(\Ai)$
is the set of all paths $\pi\in\Pi$ of the form
\be\label{hopmeet}
\pi=\bigcup_{k=1}^m\big\{(\pi_k(s),s):s\in[t_{k-1},t_k]\big\},
\ee
where $\pi_1,\ldots,\pi_m\in\Ai$, $\sig_{\pi_1}=t_0<\cdots<t_m=\infty$,
and $t_k$ is an intersection time of $\pi_k$ and $\pi_{k+1}$ for each
$k=1,\ldots,m-1$.\med

\index{hopping}

\begin{figure}[htb] 
\centering
\includegraphics[width=10cm]{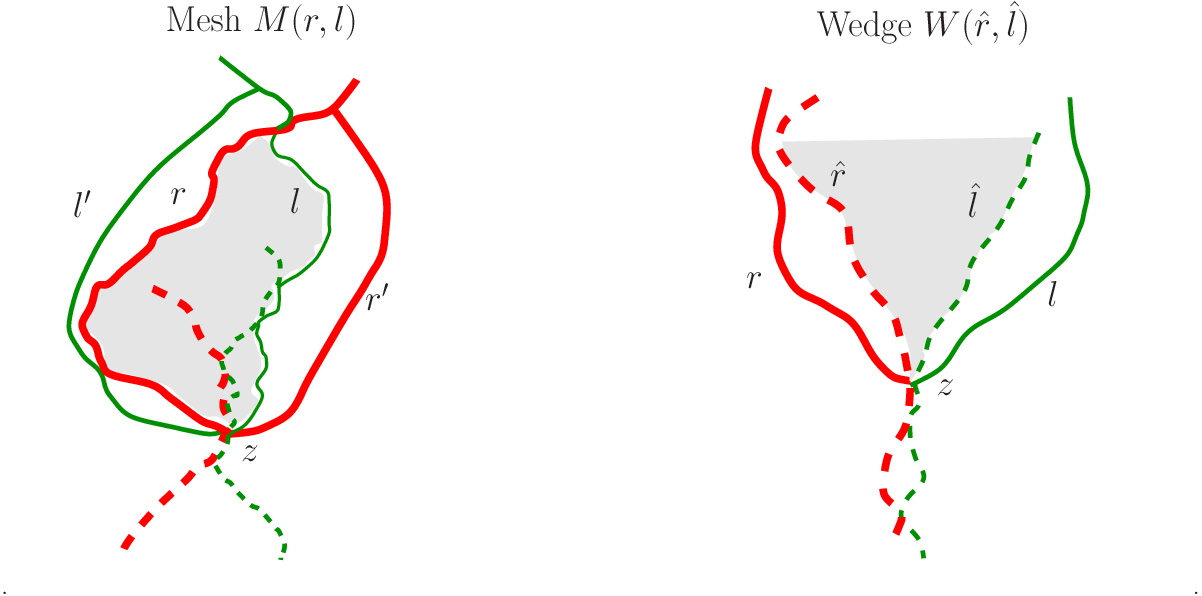}
\caption{A mesh $M(r,l)$ with bottom point $z$ and a wedge
$W(\hat r,\hat l)$ with bottom point $z$.}
\label{fig:meshwedge}
\end{figure}

\noi
{\bf Wedges:} Let $\hat\Wl,\hat\Wr$ be the dual Brownian webs associated with
a left-right Brownian web $(\Wl,\Wr)$. Any pair $\hat l\in\hat\Wl$, $\hat
r\in\hat\Wr$ with $\hat r(\hat\sig_{\hat l}\wedge \hat\sig_{\hat r}) <\hat
l(\hat\sig_{\hat l}\wedge\hat\sig_{\hat r})$ defines an open set (see
Figure~\ref{fig:meshwedge})
\be\label{wedge}
W(\hat r,\hat l):=\{(x,u)\in\R^2:
\hat\tau_{\hat r,\hat l}<u<\hat\sig_{\hat l}\wedge\hat\sig_{\hat r},
\ \hat r(u)<x<\hat l(u)\},
\ee
where $\tau_{\hat r,\hat l}:=\sup\{t<\hat\sig_{\hat l}\wedge\hat\sig_{\hat
  r}:\hat r(t)=\hat l(t)\}$ is the first (backward) hitting time of $\hat r$
and $\hat l$, which might be $-\infty$. Such an open set is called a {\em
  wedge} of $(\hat\Wl,\hat\Wr)$. If $\tau_{\hat r,\hat l}>-\infty$, then we call
$\tau_{\hat r,\hat l}$ the bottom time, and $(\hat l(\tau_{\hat r,\hat
  l}),\tau_{\hat r,\hat l})$ the bottom point of the wedge $W(\hat r,\hat l)$.
\med

\index{wedge}

\noi
{\bf Meshes:} By definition, a {\em mesh} \index{mesh} of $(\Wl,\Wr)$ (see Figure
\ref{fig:meshwedge}) is an open set of the form
\be\label{mesh}
M=M(r,l)=\{(x,t)\in\R^2:\sig_l<t<\tau_{l,r},\ r(t)<x<l(t)\},
\ee
where $l\in\Wl$, $r\in\Wr$ are paths such that $\sig_l=\sig_r$,
$l(\sig_l)=r(\sig_r)$ and $r(s)<l(s)$ on $(\sig_l,\sig_l+\eps)$ for some
$\eps>0$, and $\tau_{l,r}:=\inf\{t>\sigma_l : l(t)=r(t)\}$.
We call $(l(\sig_l),\sig_l)$ the bottom point, $\sig_l$ the
bottom time, $(l(\tau_{l,r}),\tau_{l,r})$ the top point, $\tau_{l,r}$ the top
time, $r$ the left boundary, and $l$ the right boundary of $M$.\med

Given an open set $A\subset \R^2$ and a path $\pi\in\Pi$, we say $\pi$ {\em
  enters} $A$ if there exist $\sig_\pi<s<t$ such that $\pi(s)\notin A$ and
$\pi(t)\in A$. We say $\pi$ {\em enters $A$ from outside} if there exists
$\sig_\pi <s<t$ such that $\pi(s)\notin\ov A$ and $\pi(t)\in A$. We now recall
the following characterization of the Brownian net from \cite[Theorems~1.3,
  1.7 and 1.10]{SS08}. Below, $\ov\Ai$ denotes the closure of a set of paths
$\Ai\sub\Pi$ in the topology on $\Pi$.

\bt\hspace{-2pt}{\bf(Brownian net associated with a left-right Brownian
web)}\hspace{-1pt}\label{T:net}
Let $(\Wl,\Wr)$ be a left-right Brownian web with drifts
$\bet_-\leq\bet_+$ and let $\hat\Wl,\hat\Wr$ be the dual
Brownian webs associated with $\Wl,\Wr$. Then there exists a random compact
set of paths $\Ni\in\Ki(\Pi)$, called the Brownian net, that is a.s.\ uniquely determined by any of the
following equivalent conditions:
\begin{itemize}
\item[{\rm(i)}] $\Ni=\ov{\Hi_{\rm int}(\Wl\cup\Wr)}$ a.s.
\item[{\rm(ii)}] $\Ni=\{\pi\in\Pi:\pi\mbox{ does not enter any wedge of
  $(\hat\Wl,\hat\Wr)$ from outside}\}$ a.s.
\item[{\rm(iii)}] $\Ni=\{\pi\in\Pi:\pi\mbox{ does not enter any mesh of
  $(\Wl,\Wr)$}\}$ a.s.
\end{itemize}
\index{Brownian net}
The set $\Ni$ is closed under hopping, i.e., $\Ni=\Hi_{\rm int}(\Ni)$ a.s.
Moreover, if $\Di\sub\R^2$ is a deterministic countable dense set,
then a.s., for each $z\in\Di$, the set $\Ni(z)$ contains a minimal element
$l_z$ and a maximal element $r_z$, and one has $\Wl=\ov{\{l_z:z\in\Di\}}$
  and $\Wr=\ov{\{r_z:z\in\Di\}}$~a.s.
\et

If $(\Wl,\Wr)$ and $\Ni$ are coupled as in Theorem~\ref{T:net}, then we call
$\Ni$ the {\em Brownian net} associated with $(\Wl,\Wr)$. We also call
$\bet_-,\bet_+$ the {\em left} and {\em right speed} of
$\Ni$. \index{left-right!speeds!of Brownian net} \index{speeds!of Brownian
  net} The Brownian net with left and right speeds $\bet_-=-1$ and $\bet_+=+1$
is called the {\em standard Brownian net}. We note that if $\bet_-=\bet_+$,
then $\Wl=\Wr=\Ni$, i.e., the Brownian net reduces to a Brownian web.  Since
$(\hat\Wl,\hat\Wr)$ is equal in law to a left-right Brownian web rotated over
180 degrees, such a dual left-right Brownian web defines an a.s.\ unique {\em
  dual Brownian net} $\hat\Ni$ in the same way as $(\Wl,\Wr)$ defines
$\Ni$. \index{Brownian net!dual}

If $A$ is any closed subset of $\R$ and $\Ni$ is a standard Brownian net, then
setting
\be\label{braco}
\xi_t:=\big\{\pi(t):\pi\in\Ni(A\times\{0\})\big\}\qquad(t\geq 0)
\ee
defines a Markov process taking values in the space of closed subsets of $\R$,
called the {\em branching-coalescing point set}. \index{branching-coalescing
  point set!standard} We refer to Proposition~\ref{P:braco} for some of its
basic properties.

\subsection{Separation points}\label{S:separ}

Loosely speaking, the {\em separation points} of a Brownian net are the limits of separation points of the
approximating branching-coalescing random walks, i.e., they are points
where paths in the Brownian net have a choice whether to `turn left' or `turn
right'. These points play an important role in our proofs. In this subsection, we recall some basic facts about them.

Recall from Definition~\ref{D:inout} the definition of strong equivalence of incoming paths. Following \cite{SSS09}, we adopt the following definition.

\index{incoming and outgoing paths!equivalence}\index{equivalence!of paths}
\bd{\bf(Equivalence of incoming and outgoing paths)}\label{D:equiv}
We call two incoming paths $\pi_1,\pi_2\in\Pi$ at a point $z=(x,t)\in\R^2$
{\em equivalent paths entering} $z$, denoted by $\pi_1\sim^{z}_{\rm in}\pi_2$,
if $\pi_1(t-\eps_n)=\pi_2(t-\eps_n)$ for a sequence $\eps_n\down 0$. We call
two outgoing paths $\pi_1,\pi_2$ at a point $z$ {\em equivalent paths leaving}
$z$, denoted by $\pi_1\sim^{z}_{\rm out}\pi_2$, if
$\pi_1(t+\eps_n)=\pi_2(t+\eps_n)$ for a sequence $\eps_n\down 0$.
\ed

In spite of the suggestive notation, these are not equivalence relations on
the spaces of all paths in $\Pi$ entering resp.\ leaving a point. However, it
is known that:
\be\ba{rl}\label{order}
{\rm(i)}&\mbox{If $\Wi$ is a Brownian web and $\pi_1,\pi_2\in\Wi$ satisfy
$\pi_1(t)=\pi_2(t)$}\\
&\mbox{for some $\sig_{\pi_1},\sig_{\pi_2}<t$, then $\pi_1=\pi_2$ on
$[t,\infty]$.}\\[5pt]
{\rm(ii)}&\mbox{If $(\Wl,\Wr)$ is a left-right Brownian web and
$l\in\Wl,r\in\Wr$ satisfy}\\
&\mbox{$l(t)=r(t)$ for some $\sig_l,\sig_r<t$, then $l\leq r$ on $[t,\infty]$.}
\ec
Using (\ref{order}), it is easy to see that if $(\Wl,\Wr)$ is a left-right
Brownian web, then a.s.\ for all $z\in\R^2$, the relations $\sim^{z}_{\rm in}$
and $\sim^{z}_{\rm out}$ are equivalence relations on the set of paths in
$\Wl\cup\Wr$ entering resp.\ leaving $z$, and the equivalence classes of paths
in $\Wl\cup\Wr$ entering resp.\ leaving $z$ are naturally ordered from left to
right.

In previous work \cite[Theorem 1.7]{SSS09}, we have given a complete
classification of points $z\in\R^2$ according to the structure of the
equivalence classes in $\Wl\cup\Wr$ entering resp.\ leaving $z$, in the spirit
of the classification of special points of the Brownian web in
Proposition~\ref{P:classweb}. It turns out there are 20 types of special
points in a left-right Brownian web. Here, we will only be interested in
separation points.

\begin{figure}[htb]
\begin{center}
\includegraphics[height=7cm]{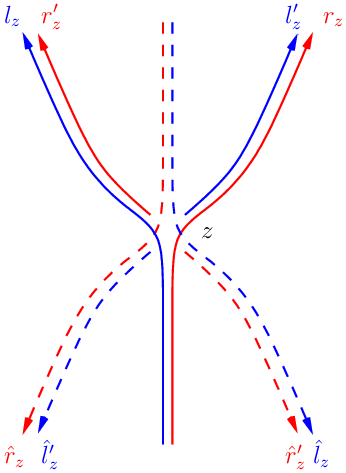}
\hspace{25pt}
\includegraphics[height=7cm]{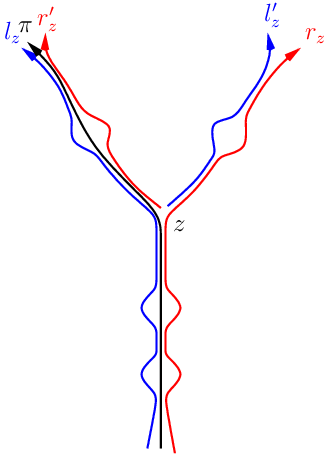}
\caption{Structure of separation points and a path $\pi$ in the Brownian net
turning left at a separation point $z$.}\label{fig:separ}
\end{center}
\end{figure}

By definition, we say that a point $z=(x,t)\in\R^2$ is a {\em separation
  point} of two paths $\pi_1,\pi_2\in\Pi$ if $\sig_{\pi_1},\sig_{\pi_2}<t$,
$\pi_1(t)=x=\pi_2(t)$, and $\pi_1\not\sim^z_{\rm out}\pi_2$. We say that $z$
is a separation point of some collection of paths $\Ai$ if there exist
$\pi_1,\pi_2\in\Ai$ such that $z$ is a separation point of $\pi_1$ and
$\pi_2$. Recall the definition of the dual Brownian net below
Theorem~\ref{T:net}. We cite the following proposition from \cite[Prop.~2.6
  and Thm.~1.12(a)]{SSS09}. See Figure~\ref{fig:separ}.

\index{separation points}
\bp{\bf(Separation points)}\label{P:separ}
Let $\Ni$ be a Brownian net with left and right speeds $\bet_-<\bet_+$ and let
$(\Wl,\Wr)$ be its associated left-right Brownian web. Then:
\begin{itemize}
\item[{\bf(a)}] A.s., $S:=\{z\in\R^2:z\mbox{ is a separation point
  of }\Ni\}=\{z\in\R^2:z\mbox{ is a separation}$ $\mbox{point of }\hat\Ni\}$,
  and $S$ is countable.
\item[{\bf(b)}] A.s., $S=\{z\in\R^2:z\mbox{ is of type $(1,2)_{\rm l}$
  in $\Wl$ and of type $(1,2)_{\rm r}$ in $\Wr$}\}$.
\item[{\bf(c)}] For given $z\in S$, let $l_z$ and $r_z$ denote the, up to
  strong equivalence unique, incoming paths in $\Wl$ resp.\ $\Wr$ at $z$ and
  let $l'_z$ and $r'_z$ be the elements of $\Wl(z)$ resp.\ $\Wr(z)$ that are
  not continuations of $l_z$ resp.\ $r_z$. Then, a.s.\ for each $z\in S$, one
  has $l_z\sim^z_{\rm in}r_z$, $l_z\sim^z_{\rm out}r'_z$ and $l'_z\sim^z_{\rm
    out}r_z$.
\item[{\bf(d)}] With the same notation as in (c), a.s.\ for each $z=(x,t)\in
  S$ and for each incoming path $\pi\in\Ni$ at $z$, there exists some $\eps>0$
  such that $l_z\leq\pi\leq r_z$ on $[t-\eps,\infty)$. Moreover, each path
  $\pi\in\Ni$ leaving $z$ satisfies either $l_z\leq\pi\leq r'_z$ on
  $[t,\infty)$ or $l'_z\leq\pi\leq r_z$ on $[t,\infty)$.
\end{itemize}
\ep
Note that part~(c) of this proposition says that at each separation point $z$
there is one pair of equivalent (in the sense of $\sim^z_{\rm in}$) incoming
paths $\{l_z,r_z\}$ and there are two pairs of equivalent (in the sense of
$\sim^z_{\rm out}$) outgoing paths: $\{l_z,r'_z\}$ and $\{l'_z,r_z\}$. By
part~(d), whenever a path $\pi\in\Ni$ enters $z$, it must do so squeezed
between $\{l_z,r_z\}$ and it must leave $z$ squeezed either between the pair
$\{l_z,r'_z\}$ or between the pair $\{l'_z,r_z\}$ (see
Figure~\ref{fig:separ}). By part~(b), these are the only points in $\R^2$
where paths in $\Ni$ can separate from each other, and by part~(a), there are
only countable many of these points.

\subsection{Switching and hopping inside a Brownian net}\label{S:netmark}

In this subsection, we show how it is possible to construct a Brownian web
inside a Brownian net by turning separation points into points of type $(1,2)$
with i.i.d.\ orientations. In the next subsection, this will be used to state
the main result of this section, which is an analogue of
Theorem~\ref{T:HWconst} and gives an alternative construction of Howitt-Warren
flows with finite left and right speeds based on a reference Brownian net.

Recall from Proposition~\ref{P:separ}~(d) (see also Figure~\ref{fig:separ})
that if $\Ni$ is a Brownian net and $\pi\in\Ni$ is some path entering a
separation point $z=(x,t)$ of $\Ni$, then $\pi$ must leave $z$ squeezed
between one of the two outgoing pairs $\{l_z,r'_z\}$ or $\{l'_z,r_z\}$. We
write
\be\label{signznet}
\sign_\pi(z):=\left\{\ba{ll}
-1\quad&\mbox{if }l_z\leq\pi\leq r'_z\mbox{ on }[t,\infty),\\[5pt]
+1\quad&\mbox{if }l'_z\leq\pi\leq r_z\mbox{ on }[t,\infty).
\ea\right.
\ee

Recall the definition of the dual Brownian net below Theorem~\ref{T:net} and
recall from Proposition~\ref{P:separ}~(a) that the set $S$ of separation
points of $\Ni$ coincides with the set of separation points of $\hat\Ni$. For
$\hat\pi\in\hat\Ni$, we define $\sign_{\hat\pi}(z)$ to be the sign of $-z$ in
$-\hat\Ni$. The next theorem, which will be proved by discrete approximation,
shows that it is possible to define a Brownian web `inside' a Brownian net.

\bt{\bf(Brownian web inside a Brownian net)}\label{T:webinnet}
Let $\Ni$ be a Brownian net with left and right speeds $\bet_-\leq\bet_+$, let
$\hat\Ni$ be its associated dual Brownian net, and let $r\in[0,1]$. Let $S$ be
the set of separation points of $\Ni$ and conditional on $\Ni$, let
$\al=(\al_z)_{z\in S}$ be a collection of i.i.d.\ $\{-1,+1\}$-valued random
variables such that $\P[\al_z=+1\,|\,\Ni]=r$. Then
\bc\label{webinnet}
\dis\Wi&:=&\dis\{\pi\in\Ni:\sign_\pi(z)=\al_z\ \forall z\in S
\mbox{ s.t.\ $\pi$ enters }z\},\\[5pt]
\dis\hat\Wi&:=&\dis\{\hat\pi\in\hat\Ni:\sign_{\hat\pi}(z)=\al_z\ \forall z\in S
\mbox{ s.t.\ $\hat\pi$ enters }z\}
\ec
defines a Brownian web $\Wi$ with drift $\bet:=(1-r)\bet_-+r\bet_+$ and its
associated dual Brownian web $\hat\Wi$. In particular, if $r=0$ resp.\ $r=1$,
then $\Wi$ is the left (resp.\ right) Brownian web associated with $\Ni$. In
general, if $\ell$ denotes the intersection local time measure between $\Wi$
and its dual and $\ell_{\rm l},\ell_{\rm r}$ are the restrictions of $\ell$ to
the sets of points of type $(1,2)_{\rm l}$ resp.\ $(1,2)_{\rm r}$, then,
conditional on $\Wi$, the sets $S_{\rm l}:=\{z\in S:\al_z=-1\}$ and $S_{\rm
  r}:=\{z\in S:\al_z=+1\}$ are independent Poisson point sets with intensities
$(\bet_+-\bet)\ell_{\rm l}$ and $(\bet-\bet_-)\ell_{\rm r}$, respectively.
\et

For any Brownian web $\Wi$, we may without loss of generality assume that
$\Wi$ is constructed `inside' some Brownian net $\Ni$ as in
Theorem~\ref{T:webinnet}. This will be very helpful in understanding marking
constructions based on $\Wi$ such as the `switching' construction of sticky
Brownian webs in Theorem~\ref{T:webmod} or the marking construction of a
Brownian net. To reap the full profit of Theorem~\ref{T:webinnet}, we need one
more result, which we formulate next.

Let $\Wi$ be a Brownian web with drift $\bet$. For each point $z$ of type
$(1,2)$, let $\switch_z(\Wi)$ denote the web obtained from $\Wi$ by
switching the orientation of $z$ as in (\ref{switch}), and let
\be
\hop_z(\Wi):=\Wi\cup\switch_z(\Wi)
\ee
be the compact set of paths obtained from $\Wi$ by allowing hopping at $z$,
i.e., by allowing incoming paths at $z$ to continue along any of the
outgoing paths. \index{hopping} More generally, if $\De$ is a finite set of points of type
$(1,2)$ in $\Wi$, then we set
\be
\hop_\De(\Wi):=\bigcup_{\De'\sub\De}\switch_{\De'}(\Wi),
\ee
where the union ranges over all subsets $\De'\sub\De$, with ${\rm
  switch}_\emptyset(\Wi):=\Wi$.

\bp{\bf(Switching and hopping inside a Brownian net)}\label{P:mark}
Let $\Ni$ be a Brownian net with left and right speeds $\bet_-\leq\bet_+$ and
set of separation points $S$. Conditional on $\Ni$, let $\al=(\al_z)_{z\in S}$
be a collection of i.i.d.\ $\{-1,+1\}$-valued random variables such that
$\P[\al_z=+1\,|\,\Ni]=r$ and let $\Wi$ be a Brownian web with drift
$\bet:=(1-r)\bet_-+r\bet_+$ defined inside $\Ni$ as in (\ref{webinnet}).
Then, a.s.\ for each subset $S'\sub S$ and for each sequence of finite sets
$\De_n\up S'$, the limits
\bc\label{WNac}
\dis\Wi'&:=&\dis\lim_{\De_n\up S'}\switch_{\De_n}(\Wi),\\[5pt]
\dis\Ni'&:=&\dis\lim_{\De_n\up S'}\hop_{\De_n}(\Wi)
\ec
exist in $\Ki(\Pi)$ and are given by
\bc\label{limid}
\dis\Ni'&=&\big\{\pi\in\Ni:\sign_\pi(z)=\al_z\ \forall z\in S\beh S'
\mbox{ s.t.\ $\pi$ enters }z\big\},\\[5pt]
\dis\Wi'&=&\Ni'\cap\big\{\pi\in\Ni:\sign_\pi(z)=-\al_z\ \forall z\in S'
\mbox{ s.t.\ $\pi$ enters }z\big\}.
\ec
\ep

\subsection{Construction of Howitt-Warren flows inside a Brownian net}\label{S:marknet}

By combining Theorem~\ref{T:webinnet} and Proposition~\ref{P:mark}, one can
give short proofs of some of the results we have seen before, such as the
marking construction of sticky Brownian webs (Theorem~\ref{T:webmod}),
Proposition~\ref{P:refchange} on changing the reference web, and the
equivalence of the definitions of a left-right Brownian web given in
Sections~\ref{S:modweb} and \ref{S:netdef}. From Theorem~\ref{T:webinnet} and
Proposition~\ref{P:mark}, one moreover easily deduces the following result,
which is similar to the marking construction of the Brownian net given in
\cite[Sec.~3.3.1 and Thm.~5.5]{NRS10}. For the proofs of all these results, we
refer to Section~\ref{S:markproof}.

\bt{\bf(Marking construction of the Brownian net)}\label{T:marknet}
Let $\Wi$ be a Brownian web with drift $\bet$ and let $c_{\rm l},c_{\rm r}\geq
0$. Let $\ell$ denote the intersection local time measure between $\Wi$ and
its dual, and let $\ell_{\rm l}$ and $\ell_{\rm r}$ denote the restrictions of
$\ell$ to the sets of points of type $(1,2)_{\rm l}$ and $(1,2)_{\rm r}$,
respectively. Conditional on $\Wi$, let $S_{\rm l}$ and $S_{\rm r}$ be
\index{Brownian net!marking construction}
independent Poisson point sets with intensities $c_{\rm l}\ell_{\rm l}$ and
$c_{\rm r}\ell_{\rm r}$, respectively. Then, for any sequence of finite sets
$\De^{\rm l}_n\up S_{\rm l}$ and $\De^{\rm r}_n\up S_{\rm r}$, the limits
\be\ba{rr@{\,}c@{\,}l}\label{marknet}
{\rm(i)}&\dis\Ni&:=&\dis\lim_{n\to\infty}
\hop_{\De^{\rm l}_n\cup\De^{\rm r}_n}(\Wi),\\[5pt]
{\rm(ii)}&\dis\Wl&:=&\dis\lim_{n\to\infty}\switch_{\De^{\rm r}_n}(\Wi),\\[5pt]
{\rm(iii)}&\dis\Wr&:=&\dis\lim_{n\to\infty}\switch_{\De^{\rm l}_n}(\Wi)
\ec
exist in $\Ki(\Pi)$ a.s.\ and do not depend on the choice of the sequences
$\De^{\rm l}_n\up S_{\rm l}$ and $\De^{\rm r}_n\up S_{\rm r}$. Moreover, $\Ni$
is a Brownian net with left and right speeds $\bet_-:=\bet-c_{\rm r}$ and
$\bet_+:=\bet+c_{\rm l}$, $(\Wl,\Wr)$ is its associated left-right
Brownian web, and $S:=S_{\rm l}\cup S_{\rm r}$ is its set of separation points.
If $c_{\rm l}+c_{\rm r}>0$, then conditional on $\Ni$, the random variables
$(\sign_\Wi(z))_{z\in S}$ are i.i.d.\ with $\P[\sign_\Wi(z)=+1\,|\,\Ni]=c_{\rm
  r}/(c_{\rm l}+c_{\rm r})$.
\et

Our final result of this subsection shows how the construction of Brownian webs
inside a Brownian net given in Theorem~\ref{T:webinnet} can be used to
construct Howitt-Warren flows with finite left and right speeds, providing an
alternative to Theorem~\ref{T:HWconst}. Recall that a Howitt-Warren quenched
law with drift $\bet$ and characteristic measure $\nu$ is a random probability
measure $\Q$ on $\Ki(\Pi)$ with law as defined in (\ref{HWquen}).

\index{Howitt-Warren!flow!construction}
\index{flow!Howitt-Warren!construction}
\bt{\bf (Construction of Howitt-Warren flows with finite speeds)}
\label{T:HWconst2}
Let $\beta\in\R$ and let $\nu$ be a finite measure on $[0,1]$ such that the
speeds $\beta_-,\beta_+$ defined in (\ref{speeds}) are finite. Let $\Ni$ be a
Brownian net with left and right speeds $\beta_-,\beta_+$ and let $S$ be
its set of separation points. Conditional on $\Ni$, let
$\om:=(\om_z)_{z\in S}$ be a collection of i.i.d.\ $[0,1]$-valued
random variables with law $\bar\nu(\di q):=b^{-1}q^{-1}(1-q)^{-1}\nu(\di q)$,
where $b:=\int q^{-1}(1-q)^{-1}\nu(\di q)$, and conditional on $(\Ni,\om)$,
let $(\alpha_z)_{z\in S}$ be a collection of independent
$\{-1,+1\}$-valued random variables such that
$\P[\alpha_z=1\,|\,(\Ni,\om)]=\om_z$. Set
\be\label{webinnet4}
\Wi:=\{\pi\in\Ni:\sign_\pi(z)=\al_z\ \forall z\in S
\mbox{ s.t.\ $\pi$ enters }z\}.
\ee
Then setting
\be\label{HWquen2}
\Q:=\P\big[\Wi\in\cdot\,\big|\,(\Ni,\om)\big]
\ee
yields a Howitt-Warren quenched law with drift $\bet$ and characteristic measure
$\nu$. In particular, setting
\bc\label{HWconst2}
K^\up_{s,t}(x,A):=\P\big[\pi^\up_{(x,s)}(t)\in A\,\big|\,(\Ni,\om)\big]
\qquad\big(s\leq t,\ x\in\R,\ A\in\Bi(\R)\big)
\ec
and defining $K^+_{s,t}(x,A)$ similarly with $\pi^\up_{(x,s)}$ replaced by
$\pi^+_{(x,s)}$ yields versions of the Howitt-Warren flow with drift $\bet$
and characteristic measure $\nu$ with properties as described in
Proposition~\ref{P:regul}.
\et

\subsection{Support of the quenched law}\label{S:quensup}

In this subsection, we formulate a theorem on the support of Howitt-Warren
quenched laws, which will imply Theorems~\ref{T:speed} and
\ref{T:supp}. Before we can do this, we need to introduce Brownian half-nets,
which are basically Brownian nets with either infinite left speed and finite
right speed, or vice versa. Recall that a path $\pi\in\Pi$ {\em crosses} a
dual path $\hat\pi\in\hat\Pi$ from left to right if there exist $\sig_\pi\leq
s<t\leq\hat\sig_{\hat\pi}$ such that $\pi(s)<\hat\pi(s)$ and
$\pi(t)>\hat\pi(t)$. Crossing from right to left and crossing of forward paths
are defined analogously.\index{crossing!of forward paths} We will prove the following analogue of
Theorem~\ref{T:net}.

\index{Brownian net!half-net}\index{half-net}
\bt{\bf(Brownian half-net associated with a Brownian web)}\label{T:halfnet}
Let $\Wi$ be a Brownian web with drift $\bet$ and let $\hat\Wi$ be its dual.
Then there exists a random closed set of paths $\Hi_-\sub\Pi$ that is a.s.\
uniquely determined by any of the following equivalent conditions:
\begin{itemize}
\item[{\rm(i)}] $\Hi_-=\{\pi\in\Pi:\pi\mbox{ does not cross any path of $\Wi$
  from left to right}\}$ a.s.
\item[{\rm(ii)}] $\Hi_-=\{\pi\in\Pi:\pi\mbox{ does not cross any path of
  $\hat\Wi$ from left to right}\}$ a.s.
\end{itemize}
Moreover, if $\Di\sub\R^2$ is a deterministic countable dense set, then a.s.,
for each $z\in\Di$, the set $\Hi_-(z)$ contains a maximal element $\pi_z$, and
one has $\Wi=\ov{\{\pi_z:z\in\Di\}}$. Analogue statements hold with $\Hi_-$
replaced by $\Hi_+$, `from left to right' replaced by `from right to left' and
`maximal element' replaced by `minimal element'.
\et

If $\Hi_-$ (resp.\ $\Hi_+$) and $\Wi$ are coupled as in
Theorem~\ref{T:halfnet}, then we call $\Hi_-$ (resp.\ $\Hi_+$) a {\em Brownian
  half-net} with {\em left} and {\em right speeds} $-\infty,\bet$
(resp.\ $\bet,+\infty$), and we call $\Wi$ the {\em right} (resp.\ {\em left})
{\em Brownian web}\index{left-right!Brownian web!for Brownian halfnet}\index{Brownian web!left-right!of Brownian halfnet} associated with $\Hi_-$ (resp.\ $\Hi_+$).

Let $\Q$ be a Howitt-Warren quenched law with drift $\bet$ and characteristic
measure $\nu$ as defined as in (\ref{HWquen}), or alternatively, in the case
of finite speeds, as in (\ref{HWquen2}). Then $\Q$ is a random probability law
on the space of webs. In particular, if $\Wi$ is a $\Ki(\Pi)$-valued random
variable with (random) law $\Q$, then for each $z\in\R^2$ we can define
special paths $\pi^\up_z$ and $\pi^+_z$ in $\Wi(z)$ as below
Proposition~\ref{P:classweb}. In analogy with the conditional law
$\Qdis^\om_{(x,s)}$ of the random walk in random environment defined in
Section~\ref{S:dis}, in the continuum setting, we define
\be\label{HWquenz}
\Q^+_z:=\Q\big[\pi^+_z\in\cdot\,\big]\qquad(z\in\R^2),
\ee
and we define $\Q^\up_z$ similarly, with $\pi^+_z$ replaced by $\pi^\up_z$.
In particular, if $\Q$ is defined as in (\ref{HWquen}) or as in
(\ref{HWquen2}), this says that $\Q^+_z:=\P[\pi^+_z\in\cdot\,|\,(\Wi_0,\Mi)]$
resp.\ $\Q^+_z:=\P[\pi^+_z\in\cdot\,|\,(\Ni,\om)]$. We note that since
typical points in $\R^2$ are of type $(0,1)$ in $\Wi$, for deterministic
$z\in\R^2$, the random variables $\Q^+_z$ and $\Q^\up_z$ are equal a.s.
It follows that for any deterministic finite measure $\mu$ on $\R^2$,
one has $\int\!\mu(\di z)\,\Q^+_z=\int\!\mu(\di z)\,\Q^\up_z$.

\bt{\bf(Support property)}\label{T:support}
Let $\Q$ be a Howitt-Warren quenched law with drift $\bet$ and characteristic
measure $\nu$, and let $\bet_-,\bet_+$ be the left and right speeds defined in
(\ref{speeds}). Then there exists a random, closed subset $\Ni\sub\Pi$
such that for any deterministic finite measure $\mu$ on $\R^2$,
\be
\supp\Big(\int\!\mu(\di z)\,\Q^+_z\Big)=\ov{\Ni\big(\supp(\mu)\big)}
\qquad{\rm a.s.}
\ee
If $-\infty<\bet_-\leq\bet_+<+\infty$, then $\Ni$ is Brownian net with left
and right speeds $\bet_-,\bet_+$. If either $-\infty=\bet_-<\bet_+<+\infty$ or
$-\infty<\bet_-<\bet_+=+\infty$, then $\Ni$ is a Brownian half-net with left
and right speeds $\bet_-,\bet_+$. If $-\infty=\bet_-<\bet_+=+\infty$, then
$\Ni=\Pi$.
\et

Note that above, $\supp(\mu)$ is a closed subset of $\R^2$, but not
necessarily of $\Rc$, which is why in general we need to take the closure of
$\Ni(\supp(\mu))$ in the space of paths $\Pi$. If
$(\ast,-\infty)\not\in\ov{\supp(\mu)}$ or if $\supp(\mu)=\R^2$, then it is
moreover true that $\ov{\Ni(\supp(\mu))}=\Ni\big(\,\ov{\supp(\mu)}\,\big)$,
where $\ov{\supp(\mu)}$ denotes the closure of $\supp(\mu)$ in $\Rc$; see
Lemma~\ref{L:NgenA} below.

We note, without proof, that in the setup of Theorem~\ref{T:support}, it can be
shown that $\Ni=\cup\,\supp(\Q)$, where $\cup\,{\rm
  supp}(\Q):=\{\pi:\pi\in\Ai\mbox{ for some }\Ai\in\supp(\Q)\}$ denotes the
union of all elements of $\supp(\Q)\subset\Ki(\Pi)$. We state as an open
problem to characterize $\supp(\Q)$ itself (rather than just
$\cup\,\supp(\Q)$).

\section{Outline of the proofs}\label{S:outline}

Our results are proved in Sections~\ref{S:prelim}--\ref{S:ergproof} below. In
Section~\ref{S:prelim} we collect some well-known and less well-known facts
about the Brownian web and net and prove some new results that we will need
further on. In particular, in Section~\ref{S:fingraph} we prove a `finite
graph representation' that gives a precise description of how paths in the
Brownian net move between deterministic
times. Sections~\ref{S:webap}--\ref{S:webnetap} then culminate in
Theorem~\ref{T:webnetap}, the central result of the section, which is about
discrete approximation of a Brownian web embedded in a Brownian net and
implies Theorem~\ref{T:webinnet}. In Section~\ref{S:markproof}, this is then
used, together with the finite graph representation, to prove
Theorem~\ref{T:webmod} and Proposition~\ref{P:refchange} on the construction
of sticky Brownian webs and related results such as Proposition~\ref{P:mark}
and Theorem~\ref{T:marknet}.

In Section~\ref{S:main} we prove our main results: Theorem~\ref{T:quench}
on the convergence of the quenched laws on the space of webs, and
Theorems~\ref{T:HWconst} and \ref{T:HWconst2} on the construction of
Howitt-Warren flows using a marked reference Brownian web or net. Here we also
harvest some immediate consequences of our construction, such as the existence
of regular versions of Howitt-Warren flows (Proposition~\ref{P:conpath} and
\ref{P:regul}) and scaling (Proposition~\ref{P:scale}).

In Section~\ref{S:support} we prove our results on the support of
Howitt-Warren flows. In Section~\ref{S:halfnet}, we prove a number of
preparatory results about generalized Brownian nets with possibly infinite
left and right speeds. In particular, we prove Theorem~\ref{T:halfnet} on
Brownian half-nets and prepare for the proof of Theorem~\ref{T:support} on the
support of the quenched law on the space of webs. In Section~\ref{S:HWsupp},
we prove Theorem~\ref{T:support} and use it to deduce Theorems~\ref{T:speed}
and \ref{T:supp} on the left and right speeds and the support of Howitt-Warren
processes.

In Section~\ref{S:atom}, we address questions of atomicness. In particular,
parts (a), (b) and (c) of Theorem~\ref{T:atom} are proved in Sections~\ref{S:detatom}, \ref{S:nonatom} and \ref{S:erosion}, respectively.

In Section~\ref{S:inf} we prove Theorems~\ref{T:infmass} and \ref{T:disflow}
on Howitt-Warren processes with infinite starting mass and the convergence of
rescaled discrete Howitt-Warren processes, while Section~\ref{S:ergproof}
contains the proofs of Theorems~\ref{T:HIL} and \ref{T:invsupp} on homogeneous
invariant laws.

The paper concludes with four appendices on the Howitt-Warren martingale problem and some other technical issues.

The table below gives a complete overview of where the proofs can be found of
all results stated so far. Further results stated in the following sections
will be proved on the spot. Below, {\em cited} means that the listed result is
cited from other sources.

\vspace{\baselineskip}

\newcommand{\cited}[1]{cited}

\begin{tabular}{|c|l||c|l||c|l|}
\hline
Result & Proved in & Result & Proved in & Result & Proved in \\ \hline
Prop.~\ref{P:conpath} & Sect.~\ref{S:immed} & Thm.~\ref{T:invsupp} & Sect.~\ref{S:HIL} & Thm.~\ref{T:net} & \cited{SS08}\\
Prop.~\ref{P:scale} & Sect.~\ref{S:immed} & Prop.~\ref{P:char} & \cited{FINR04,FINR06,SS08} & Prop.~\ref{P:separ} & \cited{SSS09}\\
Thm.~\ref{T:speed} & Sect.~\ref{S:HWsupp} & Prop.~\ref{P:classweb} & \cited{TW98,FINR06} & Thm.~\ref{T:webinnet} & Sect.~\ref{S:webnetap}\\
Prop.~\ref{P:braco} & \cited{SS08,SSS09} & Prop.~\ref{P:refloc} & \cited{NRS10} & Prop.~\ref{P:mark} & Sect.~\ref{S:markproof}\\
Thm.~\ref{T:supp} & Sect.~\ref{S:HWsupp} & Thm.~\ref{T:webmod} & Sect.~\ref{S:markproof} & Thm.~\ref{T:marknet} & Sect.~\ref{S:markproof}\\
Thm.~\ref{T:atom} & Sect.~\ref{S:atom} & Prop.~\ref{P:refchange} & Sect.~\ref{S:markproof} & Thm.~\ref{T:HWconst2} & Sect.~\ref{S:HWconst}\\
Thm.~\ref{T:infmass} & Sect.~\ref{S:Feller} & Thm.~\ref{T:HWconst} & Sect.~\ref{S:HWconst} & Thm.~\ref{T:halfnet} & Sect.~\ref{S:halfnet}\\
Thm.~\ref{T:disflow} & Sect.~\ref{S:disflow} & Prop.~\ref{P:regul} & Sect.~\ref{S:immed} & Thm.~\ref{T:support} & Sect.~\ref{S:HWsupp}\\
Thm.~\ref{T:HIL} & Sect.~\ref{S:HIL} & Thm.~\ref{T:quench} & Sect.~\ref{S:quench} & & \\ \hline
\end{tabular}

\vspace{\baselineskip}

\section{Coupling of the Brownian web and net}\label{S:prelim}

The main aim of this section is to prove Theorem~\ref{T:webinnet} and
Proposition~\ref{P:mark}, which will be our main tools for constructing
modified Brownian webs and nets by switching or hopping inside a reference
Brownian web or net. In particular, after proving these theorems, we will
apply them to prove Theorems~\ref{T:webmod} and \ref{T:marknet} on the
switching construction of sticky Brownian webs and the marking construction of
the Brownian net.

In order to prepare for the proofs of Theorem~\ref{T:webinnet} and
Proposition~\ref{P:mark}, we first need to take a closer look at the
separation points of a Brownian net, introduced in Section~\ref{S:separ}. It
has been proved in \cite{SSS09} that for deterministic times $S<U$, there are
only locally finitely many `$S,U$-relevant' separation points that decide
where paths in the Brownian net started at time $S$ end up at time $U$. After
recalling some basic facts about these relevant separation points in
Section~\ref{S:relev}, we use them in Section~\ref{S:fingraph} to give a
rather precise description, by means of a `finite graph representation', of
the way paths in the Brownian net move between time $S$ and $U$.

Since discrete approximation will play an important role in our proofs,
Sections~\ref{S:webap}--\ref{S:netap} are devoted to discrete approximation of
the Brownian web and net, and related objects such as intersection local times
and separation points. In Section~\ref{S:webnetap}, we then use these results
to prove a result about the convergence of a discrete web embedded in a
discrete net to analogue Brownian objects. This result then immediately yields
Theorem~\ref{T:webinnet} on the construction of a Brownian web inside a
Brownian net. In addition, it lays the basis for proofs of other convergence
results such as Theorem~\ref{T:quench} on the convergence of quenched laws.
In Section~\ref{S:markproof}, finally, we use the finite graph representation
developed in Section~\ref{S:fingraph} together with Theorem~\ref{T:webinnet}
to prove Proposition~\ref{P:mark} and we combine Theorem~\ref{T:webinnet} and
Proposition~\ref{P:mark} to prove Theorems~\ref{T:webmod} and \ref{T:marknet}
and some related results.

\subsection{Relevant separation points}\label{S:relev}

The set of separation points of a Brownian net $\Ni$ is dense in $\R^2$ and
also along any path $\pi\in\Ni$. It turns out, however, that for given
deterministic times $S<U$, the set of separation points that are relevant for
deciding where paths in the Brownian net started at time $S$ end up at time
$U$ is a locally finite subset of $\R\times[S,U]$.

\begin{figure}[htb]
\begin{center}
\includegraphics[height=6cm]{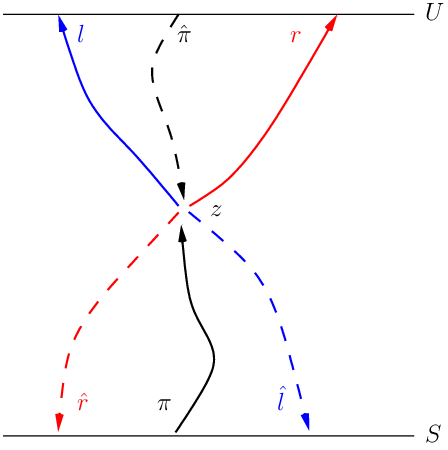}
\caption{An $S,U$-relevant separation point.}\label{fig:relev}
\end{center}
\end{figure}

Following \cite{SSS09}, we say that a separation point $z=(x,t)$ of a Brownian
net $\Ni$ is {\em $S,U$-relevant} for some $-\infty\leq S<t<U\leq\infty$, if
there exists $\pi\in\Ni$ such that $\sig_\pi=S$ and $\pi(t)=x$, and there exist
$l\in\Wl(z)$ and $r\in\Wr(z)$ such that $l<r$ on $(t,U)$. (Note that since we
are assuming that $z$ is a separation point, $l$ and $r$ have to be the paths
$l_z$ and $r_z$ from Proposition~\ref{P:separ}~(c). In particular, $l$ and $r$
are continuations of incoming paths at $z$.) The next proposition follows
easily, by Brownian scaling, from \cite[Lemma~2.8 and
  Prop.~2.9]{SSS09}. Part~(a) says that the definition of relevant separation
points is symmetric with respect to duality; see also Figure~\ref{fig:relev}.

\bp{\bf(Relevant separation points)}\label{P:relev}
Let $\Ni$ be a Brownian net with left and right speeds $\bet_-\leq\bet_+$.
Then:\med
\index{separation points!relevant}

\noi
{\bf(a)} A.s.\ for each $-\infty\leq S<U\leq\infty$, a separation point
$z=(x,t)$ with $S<t<U$ is $S,U$-relevant in $\Ni$ if and only if $-z$ is
$-U,-S$-relevant in the rotated dual Brownian net $-\hat\Ni$.\med

\noi
{\bf(b)} For each deterministic $-\infty\leq S<U\leq\infty$, if $R_{S,U}$
denotes the set of $S,U$-relevant separation points, then
\be\label{reldens}
\E\big[\big|R_{S,U}\cap A\big|\big]=2b\int_A\Psi_b(t-S)\Psi_b(U-t)\,\di x\,\di t
\qquad\big(A\in\Bi(\R\times(S,U))\big),
\ee
where $b:=(\bet_+-\bet_-)/2$,
\be\label{Psib}
\Psi_b(t):=\frac{e^{-b^2t}}{\sqrt{\pi t}}+2b\Phi(b\sqrt{2t})
\qquad(0<t\leq\infty),
\ee
and $\Phi(x):=\frac{1}{\sqrt{2\pi}}\int_{-\infty}^xe^{-y^2/2}\di y$. In
particular, if $-\infty<S$, $U<\infty$, then $R_{S,U}$ is a.s.\ a locally
finite subset of $\R\times[S,U]$.
\ep

We will need yet another characterization of relevant separation points. To
formulate this, we first need to recall the definition of crossing times from
\cite[Def.~2.4]{SSS09}.

\begin{defi}{\bf(Crossing and crossing points)}
We say that a forward path $\pi\in\Pi$ crosses a dual
path $\hat\pi\in\hat\Pi$ from left to right at time $t$ if there exist
$\sig_\pi\leq t_-<t<t_+\leq\hat\sig_{\hat\pi}$ such that
$\pi(t_-)<\hat\pi(t_-)$, $\hat\pi(t_+)<\pi(t_+)$, and
$t=\inf\{s\in(t_-,t_+):\hat\pi(s)<\pi(s)\}
=\sup\{s\in(t_-,t_+):\pi(s)<\hat\pi(s)\}$. Crossing from right to left is
defined analogously. We call $z=(x,t)\in\R^2$ a crossing point of $\pi\in\Pi$
and $\hat\pi\in\hat\Pi$ if $\pi(t)=x=\hat\pi(t)$ and $\pi$ crosses $\hat\pi$
either from left to right or from right to left at time $t$.
\end{defi}
\index{crossing!points} 

\bl{\bf(Relevant separation points and crossing points)}\label{L:relcros}
Almost surely for each $-\infty\leq S<U\leq\infty$ and
$z\in\R\times(S,U)$, the following statements are equivalent:
\begin{itemize}
\item[{\rm(i)}] $z$ is an $S,U$-relevant separation point.
\item[{\rm(ii)}] $z$ is a crossing point of some $\pi\in\Ni$ and
$\hat\pi\in\hat\Ni$ with $\sig_\pi=S$ and $U=\hat\sig_{\hat\pi}$.
\end{itemize}
\el
{\bf Proof.} If $z=(x,t)$ is an $S,U$-relevant separation point, then by
Proposition~\ref{P:relev} there exist $\pi'\in\Ni$ starting at time
$\sig_{\pi'}=S$ and $\hat\pi'\in\hat\Ni$ starting at time
$\hat\sig_{\hat\pi'}=U$ such that $\pi'$ and $\hat\pi'$ enter $z$. By
\cite[Prop.~2.6]{SSS09}, $z$ is a crossing point of some $r\in\Wr$ and $\hat
l\in\hat\Wl$. Let $\pi$ be the concatenation of $\pi'$ on $[S,t]$ and $r$
on $[t,U]$ and likewise, let $\hat\pi\in\Ni$ be the concatenation of
$\hat\pi'$ on $[t,U]$ and $\hat l$ on $[S,t]$. Since by Theorem~\ref{T:net},
$\Ni$ is closed under hopping, we see that $\pi\in\Ni$ and
$\hat\pi\in\hat\Ni$. By the structure of separation points
(Proposition~\ref{P:separ}~(d)), $z$ is a crossing point of $\pi$ and
$\hat\pi$, proving the implication (i)$\volgt$(ii).

Conversely, if $z$ is a crossing point of some $\pi\in\Ni$ and
$\hat\pi\in\hat\Ni$ with $\sig_\pi=S$ and $U=\hat\sig_{\hat\pi}$, then
by the classification of special points of the Brownian net
\cite[Thm.~1.7]{SSS09} and their structure \cite[Thm.~1.12~(d)]{SSS09}, $z$
must be a separation point of $\Ni$. By \cite[Lemma~2.7~(a)]{SSS09}, the
presence of the dual path $\hat\pi$ implies the existence of $l\in\Wl(z)$,
$r\in\Wr(z)$ such that $l<r$ on $(t,U)$, hence $z$ is $S,U$-relevant.\qed

\subsection{Finite graph representation}\label{S:fingraph}

In this section, we give a rather precise description of how paths in a
Brownian net move between deterministic times $S,U$. In particular, we will
construct an oriented graph whose internal vertices are relevant separation
points and whose directed edges are pairs consisting of a left-most and
right-most path, such that each path in the Brownian net starting at time
$S$ must between times $S$ and $U$ move through an oriented path in this
graph, and conversely, for each oriented path in the graph there exist paths
in the Brownian net following this path.

As a preparation, we need some results from \cite{SSS09} on the special points
of the Brownian net. Almost surely, there are 20 types of special points in the
Brownian net, but we will only need those that occur at deterministic times,
of which there are only three. Let $\Ni$ be a Brownian net with associated
left-right Brownian web $(\Wl,\Wr)$. Recall the notion of strong equivalence
of paths from Definition~\ref{D:inout} and the relations $\sim^z_{\rm in}$ and
$\sim^z_{\rm out}$ from Definition~\ref{D:equiv}. As remarked there, these are
equivalence relations on the set of paths in $\Wl\cup\Wr$ entering
resp.\ leaving a point $z$, and the corresponding equivalence classes are
naturally ordered from left to right. In general, such an equivalence class
may be of three types. If it contains only paths in $\Wl$ then we say it is of
type ${\rm l}$, if it contains only paths in $\Wr$ then we say it is of type
${\rm r}$, and if it contains both paths in $\Wl$ and $\Wr$ then we say it is
of type ${\rm p}$, standing for pair. To denote the type of a point $z\in\R^2$
in a Brownian net $\Ni$, we first list the incoming equivalence classes in
$\Wl\cup\Wr$ from left to right and then, separated by a comma, the outgoing
equivalence classes.

In our case, there are only three types of points of interest, namely the
types ${\rm (o,p)}$, ${\rm (p,p)}$ and ${\rm (o,pp)}$, where a $o$ means that
there are no incoming paths in $\Ni$ at $z$. We note that by property
(\ref{order})~(i), an outgoing equivalence class of type ${\rm p}$ at a point
$z$ contains exactly one path in $\Wl(z)$ and one path in $\Wr(z)$. By the
same property, at points of type ${\rm (p,p)}$, all incoming paths in $\Wl$
are strongly equivalent and likewise all incoming paths in $\Wr$ are strongly
equivalent. We cite the following result from \cite[Thms.~1.7 and 1.12]{SSS09}
and \cite[Prop.~1.8]{SS08}. Recall the definition of the dual Brownian net
below Theorem~\ref{T:net}.

\index{special points!Brownian net}
\bp{\bf(Special points at deterministic times)}\label{P:detspec}
Let $\Ni$ be a Brownian net, let $\hat\Ni$ be its dual, and let
$(\Wl,\Wr)$ and $(\hat\Wl,\hat\Wr)$ be the left-right Brownian web and the dual
left-right Brownian web associated with $\Ni$ and $\hat\Ni$. Then:
\begin{itemize}
\item[{\bf(a)}] For each deterministic $t\in\R$, almost surely, each point in
  $\R\times\{t\}$ is either of type ${\rm (o,p)}/{\rm (o,p)}$, ${\rm
  (p,p)}/{\rm (o,pp)}$ or ${\rm (o,pp)}/{\rm (p,p)}$ in $\Ni/\hat\Ni$, and all
  of these types occur.
\item[{\bf(b)}] Almost surely, for each point $z=(x,t)$ of type ${\rm (o,p)}$,
  ${\rm (p,p)}$ or ${\rm (o,pp)}$ in $\Ni$ and $\pi\in\Ni(z)$, there exist
  $l\in\Wl(z)$ and $r\in\Wr(z)$ such that $l\sim^z_{\rm out}r$ and
  $l\leq\pi\leq r$ on $[t,\infty)$.
\item[{\bf(c)}] Almost surely, for each point $z=(x,t)$ of type ${\rm
  (p,p)}$ in $\Ni$, for each $l\in\Wl$, $r\in\Wr$ and $\pi\in\Ni$ entering
  $z$, there exists an $\eps>0$ such that $l\leq\pi\leq r$ on $[t-\eps,\infty)$.
\end{itemize}
\ep

Let $-\infty<S<U<\infty$ be deterministic times,
let $R_{S,U}$ be the set of $S,U$-relevant separation points of $\Ni$ and set
\bc
\dis R_S&:=&\dis\R\times\{S\},\\[5pt]
\dis R_U&:=&\dis\big\{(x,U):x\in\R,\ \exists\pi\in\Ni\mbox{ with }
\sig_\pi=S\mbox{ s.t.\ }\pi(U)=x\big\}.
\ec
We make the set $R:=R_S\cup R_{S,U}\cup R_U$ into an oriented graph by
writing $z\to_{l,r}z'$ if $z,z'\in R$, $z\neq z'$, $l\in\Wl(z)$, $r\in\Wr(z)$,
$l\sim^z_{\rm out} r$, and $l\sim^{z'}_{\rm in} r$.

\begin{figure}[htb]
\begin{center}
\includegraphics[width=15.5cm]{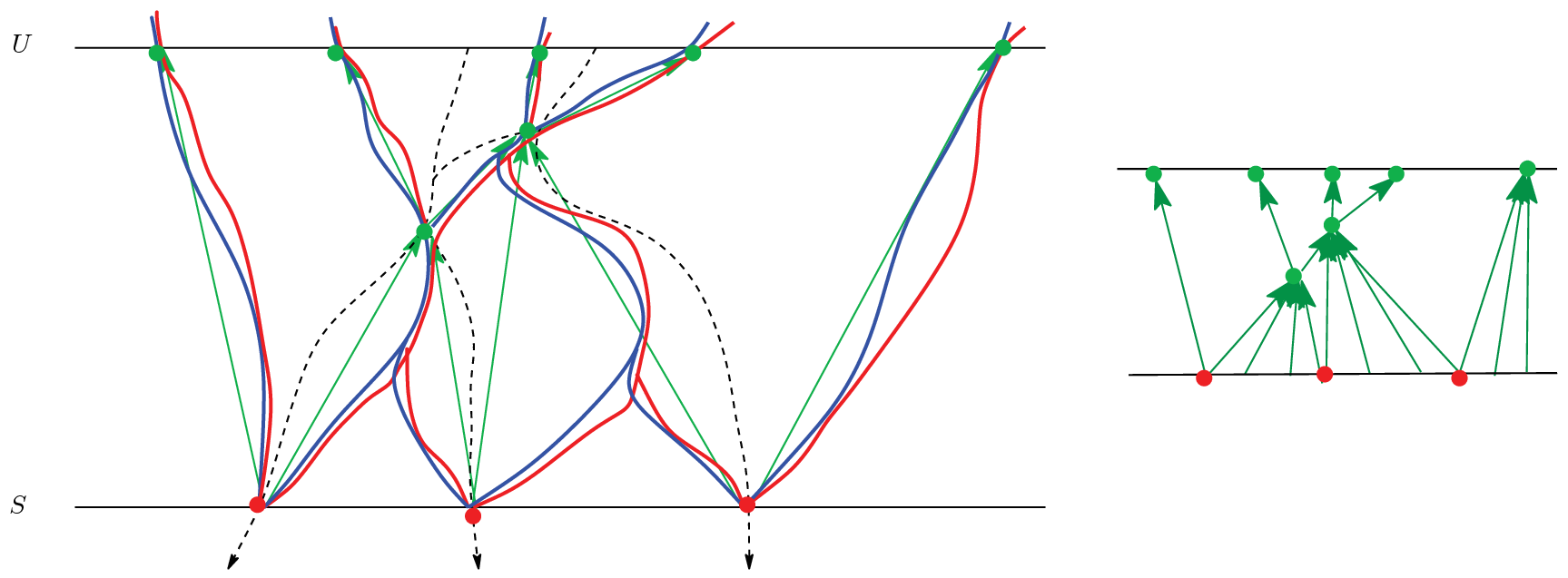}
\caption{Finite graph representation.}\label{fig:fingraph}
\end{center}
\end{figure}

\index{finite graph representation}
\bp{\bf(Finite graph representation)}\label{P:fingraph}
Let $\Ni$ be a Brownian net with associated left-right Brownian web
$(\Wl,\Wr)$ and let $-\infty<S<U<\infty$ be deterministic times.
Let $R:=R_S\cup R_{S,U}\cup R_U$ and the relation $\to_{l,r}$ be defined as
above. Then, a.s.\ (see Figure~\ref{fig:fingraph}):
\begin{itemize}
\item[{\bf(a)}] For each $z\in R_S$ that is not of type $\rm(o,pp)$, there
  exist unique $l\in\Wl(z)$, $r\in\Wr(z)$ and $z'\in R$ such that
  $z\to_{l,r}z'$.
\item[{\bf(b)}] For each $z=(x,t)$ such that either $z\in R_{S,U}$ or $z\in
  R_S$ is of type $\rm(o,pp)$, there exist unique $l,l'\in\Wl(z)$,
  $r,r'\in\Wr(z)$ and $z',z''\in R$ such that $l\leq r'<l'\leq r$ on
  $(t,t+\eps)$ for some $\eps>0$, $z\to_{l,r'}z'$ and $z\to_{l',r}z''$. For
  $z\in R_{S,U}$ one has $z'\neq z''$. For $z\in R_S$ of type $\rm(o,pp)$, one
  has $z'\neq z''$ if and only if there exists a dual path $\hat\pi\in\hat\Ni$
  with $\hat\sig_{\hat\pi}=U$ such that $\hat\pi$ enters $z$.
\item[{\bf(c)}] For each $\pi\in\Ni$ with $\sig_\pi=S$, there exist
  $z_i=(x_i,t_i)\in R$ $(i=0,\ldots,n)$ and $l_i\in\Wl(z_i)$, $r_i\in\Wl(z_i)$
  $(i=0,\ldots,n-1)$ such that $z_0\in R_S$, $z_n\in R_U$,
  $z_i\to_{l_i,r_i}z_{i+1}$ and $l_i\leq\pi\leq r_i$ on $[t_i,t_{i+1}]$
  $(i=0,\ldots,n-1)$.
\item[{\bf(d)}] If $z_i=(x_i,t_i)\in R$ $(i=0,\ldots,n)$ and
  $l_i\in\Wl(z_i)$, $r_i\in\Wl(z_i)$ $(i=0,\ldots,n-1)$ satisfy $z_0\in
  R_S$, $z_n\in R_U$, and $z_i\to_{l_i,r_i}z_{i+1}$ $(i=0,\ldots,n-1)$, then
  there exists a $\pi\in\Ni$ with $\sig_\pi=S$ such that $l_i\leq\pi\leq r_i$
  on $[t_i,t_{i+1}]$.
\end{itemize}
\ep
{\bf Proof.} By Proposition~\ref{P:detspec}~(a), each $z=(x,t)\in R_S$ that is
not of type $\rm(o,pp)$ must be of type $\rm(o,p)$ or $\rm(p,p)$, hence there
exists a unique pair $(l,r)$ consisting of one left-most path $l\in\Wl(z)$ and
one right-most path $r\in\Wr(z)$, such that $l\sim^z_{\rm out}r$. Likewise, by
Proposition~\ref{P:separ}, for each $z\in R_{S,U}$ there exist exactly two
pairs $(l,r')$ and $(l',r)$ such that $l,l'\in\Wl(z)$, $r,r'\in\Wr(z)$,
$l\sim^z_{\rm out}r'$ and $l'\sim^z_{\rm out}r$, and the same is true for
$z\in R_S$ that is of type $\rm(o,pp)$, by the properties of such
points. Therefore, in order to prove parts~(a) and (b), assume that $z\in
R_S\cup R_{S,U}$ and that $l\in\Wl(z)$ and $r\in\Wr(z)$ satisfy $l\sim^z_{\rm
  out}r$. We claim that there exists a unique $z'\in R_{S,U}\cup R_U$ such
that $z\to_{l,r}z'$.

To see this, let $\tau:=\sup\{u\in(t,U):l(u)=r(u)\}$ be the last time $l$ and
$r$ separate before time $U$. If $\tau=U$, then there exists some
$z'\in\R\times\{U\}$ such that $l$ and $r$ enter $z'$ and hence, by
Proposition~\ref{P:detspec}~(a), $l\sim^{z'}_{\rm in}r$. On the other hand, if
$\tau<U$, then we claim that $z'=(x',t'):=(l(\tau),\tau)$ is an $S,U$-relevant
separation point. To prove this, we must show that there exists some
$\pi\in\Ni$ with $\sig_\pi=S$ and $\pi(t')=x'$, the other parts of the
definition being obviously satisfied. If $z\in\R_S$ we may take $\pi=l$. If
$z\in R_{S,U}$, then there exists some $\pi\in\Ni$ with $\sig_\pi=S$ such that
$\pi(t)=x$. By Theorem~\ref{T:net}, the Brownian net is closed under hopping,
therefore we may concatenate $\pi$ with $l$ to find a path in $\Ni$ starting
at time $S$ and entering $z'$. This proves that $z'$ is an $S,U$-relevant
separation point. By Proposition~\ref{P:separ}~(b), the left-most and
right-most paths entering a separation point are up to strong equivalence
unique, and $l\sim^{z'}_{\rm in}r$. This proves the existence of a $z'\in
R_{S,U}\cup R_U$ such that $z\to_{l,r}z'$. The uniqueness of $z'$ follows from
the fact that only the last separation point of $l$ and $r$ before time $U$
can be $S,U$-relevant. This completes the proof of part~(a).

To complete the proof of part~(b), it suffices to show that $z'\neq z''$ if
and only if there exists a dual path $\hat\pi\in\hat\Ni$ with
$\hat\sig_{\hat\pi}=U$ such that $\hat\pi$ enters $z$. In particular, since by
Proposition~\ref{P:relev}~(a) such a dual path exists for each $z\in R_{S,U}$,
this then shows that $z'\neq z''$ for such $z$. We observe that in general,
if $z'=z''$, then the paths $l$ and $r$ starting at $z=(x,t)$ meet before time
$U$. Conversely, if $\tau_{l,r}:=\inf\{u>t:l(u)=r(u)\}<U$, then $l$ and $r$
cannot enter a $S,U$-relevant separation point before time $\tau_{l,r}$, while
after time $\tau_{l,r}$, by the arguments above, $l$ and $r$ must lead to the
same point in $R$. The statement now follows from the fact that by
\cite[Lemma~2.7]{SSS09}, there exists a dual path $\hat\pi\in\hat\Ni$ with
$\hat\sig_{\hat\pi}=U$ entering $z$ if and only if there exist $l\in\Wl(z)$
and $r\in\Wr(z)$ such that $l<r$ on $(t,U)$.

To prove part~(c), set $z_0:=(\pi(S),S)$. By Proposition~\ref{P:detspec}~(b)
there exist unique $l_0\in\Wl(z_0)$ and $r_0\in\Wl(z_0)$ such that
$l_0\sim^{z_0}_{\rm out}r_0$ and $l_0\leq\pi\leq r_0$, hence by what we have
just proved one has $z_0\to_{l_0,r_0}z_1$ for some unique $z_1=(x_1,t_1)\in
R$. If $t_1=U$ we are done. Otherwise, by what we have just proved, there
exist $l,l'\in\Wl(z_1)$, $r,r'\in\Wl(z_1)$, and $z',z''\in R$ such that
$z_0\to_{l,r'}z'$ and $z_0\to_{l',r}z''$. By Proposition~\ref{P:separ}, the
path $\pi$ must either turn left or right at $z_1$, so setting either
$(l_1,r_1)=(l,r')$ or $(l_1,r_1)=(l',r)$ and $z_2=(x_2,t_2):=z'$ or $z''$ we
have that $l_1\leq\pi\leq r_1$ on $[t_1,t_2]$. Continuing this process, which
terminates after a finite number of steps by Proposition~\ref{P:relev}, we
find a sequence of points $z_0,\ldots,z_n$ and paths
$l_0,r_0,\ldots,l_{n-1},r_{n-1}$ with the desired properties.

Finally, part~(d) follows from \cite[Thm.~1.12~(d)]{SSS09} which implies that
the concatenation of the paths $l_0,\ldots,l_{n-1}$ defines a path $\pi\in\Ni$
with all the desired properties.\qed

We will sometimes need the following extension of Proposition~\ref{P:fingraph}.

\bcor{\bf(Steering paths between deterministic times)}\label{C:steer}
Let $\Ni$ be a Brownian net with associated left-right Brownian web
$(\Wl,\Wr)$ and let $-\infty<T_1<\cdots<T_m<\infty$ be deterministic times.
Set 
\bc\label{RTk}
\dis R_{T_1}&:=&\dis\R\times\{T_1\},\\[5pt]
\dis R_{T_m}&:=&\dis\big\{(x,T_m):x\in\R,\ \exists\pi\in\Ni\mbox{ with }
\sig_\pi=T_1\mbox{ s.t.\ }\pi(T_m)=x\big\},\\[5pt]
\dis R_{T_k,T_{k+1}}&:=&\dis\big\{z\in\R^2:z\mbox{ is a
  $T_k,T_{k+1}$-relevant separation point,}\\
&&\dis\phantom{\big\{z\in\R^2:}
\exists\pi\in\Ni\mbox{ with }\sig_\pi=T_1\mbox{ s.t.\ }\pi
\mbox{ enters }z\big\}.
\ec
Then all conclusions of Proposition~\ref{P:fingraph} remain valid with $R_S$
replaced by $R_{T_1}$, $R_{S,U}$ replaced by
$\bigcup_{k=1}^{m-1}R_{T_k,T_{k+1}}$, and $R_U$ replaced by $R_{T_m}$, except
that in part~(b), it may happen that $z'=z''$ for some
$z\in\bigcup_{k=1}^{m-1}R_{T_k,T_{k+1}}$ or $z'\neq z''$ for some $z\in R_S$
even though there is no path $\hat\pi\in\hat\Ni$ starting at time $T_m$
entering $z$.

Moreover, a.s.\ for each $\pi,\pi'\in\Ni$ satisfying $\sig_\pi,\sig_{\pi'}\leq
T_1$, $\sig_\pi\wedge\sig_{\pi'}<T_1$, $\pi(T_1)\leq\pi'(T_1)$ and
$\sign_\pi(z)\leq\sign_{\pi'}(z)$ for all
$z\in\bigcup_{k=1}^{m-1}R_{T_k,T_{k+1}}$ such that both $\pi$ and $\pi'$ enter
$z$, one has $\pi(T_k)\leq\pi'(T_k)$ for $k=1,\ldots,m$.
\ecor
{\bf Proof.} This generalization of Proposition~\ref{P:fingraph} follows by
`pasting together' the finite graph representations for the consecutive time
intervals $[T_k,T_{k+1}]$, where we use that by Proposition~\ref{P:detspec}~(a),
if $\pi\in\Ni$ satisfies $\sig_\pi=T_1$, then the points $(\pi(T_k),T_k)$
$(k=2,\ldots,m)$ must be of type ${\rm(p,p)}$ in $\Ni$.

To prove the statement about the paths $\pi,\pi'$, by symmetry, we may assume
without loss of generality that $\sig_\pi<T_1$. In this case, by
Proposition~\ref{P:detspec}~(a), the point $(\pi(T_1),T_1)$ must be of type
${\rm(p,p)}$ in $\Ni$. Let $z_i=(x_i,t_i)$ $(i=1,\ldots,n)$ and $l_i,r_i$
$(i=1,\ldots,n-1)$ be defined for $\pi$ as in
Proposition~\ref{P:fingraph}~(c), and define $\pi_-,\pi_+:[T_1,T_m]\to\R$ by
\be
\pi_-:=l_i\quad\mbox{and}\quad\pi_+:=r_i\mbox{ on }[t_i,t_{i+1}]
\qquad(i=1,\ldots,n-1).
\ee
Then $\pi_-(T_k)=\pi(T_k)=\pi_+(T_k)$ for $k=1,\ldots,m$ and $\pi_-\leq\pi'$
on $[T_1,T_m]$.\qed

\subsection{Discrete approximation of the Brownian web}\label{S:webap}

In this and the next section, we recall known results about convergence of
discrete webs and nets to Brownian webs and nets and prove some related, new
results about convergence of intersection local times and relevant separation
points. In Section~\ref{S:webnetap}, we then use these results to prove a new
convergence result about Brownian webs embedded in Brownian nets, which will
form the basis for the proof of Theorem~\ref{T:quench} on the convergence of
quenched laws, which will be proved in Section~\ref{S:quench}.

Before we turn our attention to the details of the Brownian web, we first
explain two simple, general principles that we will be using several time in
what follows.

\bl{\bf(Weak convergence of coupled random variables)}\label{L:weakcouple}
\begin{itemize}
\item[{\bf(a)}] Let $E$ be a Polish space, let $(F_i)_{i\in I}$ be a finite or countable collection of Polish spaces and for each $i\in I$, let $f_i:E\to F_i$ be a measurable function. Let $X,X_k,Y_{k,i}$ be random variables $(k\geq 1,\ i\in I)$ such that $X,X_k$ take values in $E$ and $Y_{k,i}$ takes values in $F_i$. Then
\be\ba{l}\label{weakcouple}
\dis\P\big[\big(X_k,Y_{k,i}\big)\in\cdot\,\big]
\Asto{k}\P\big[\big(X,f_i(X)\big)\in\cdot\,\big]\quad\forall i\in I\\[10pt]
\dis\quad\mbox{implies}\quad
\dis\P\big[\big(X_k,(Y_{k,i})_{i\in I}\big)\in\cdot\,\big]
\Asto{k}\P\big[\big(X,(f_i(X))_{i\in I}\big)\in\cdot\,\big],
\ec
where $\Rightarrow$ denotes weak convergence of probability laws on the Polish
spaces $E\times F_i$ and $E\times\prod_{i\in I}F_i$, respectively.

\item[{\bf(b)}] Let $E,F,G$ be Polish spaces, let $f:E\to F$ and $g:F\to G$ be measurable functions, let $X,X_k$, $Y,Y_k$ and $Z_k$ be random variables taking values in $E$, $F$ and $G$, respectively $(k\geq 1)$. Then
\be\ba{l}\label{weakcouple2}
\P\big[(X_k,Y_k)\in\cdot\,\big]\Asto{k}\P\big[(X,f(X))\in\cdot\,\big]
\quad\mbox{and}\quad
\P\big[(Y_k,Z_k)\in\cdot\,\big]\Asto{k}\P\big[(Y,g(Y))\in\cdot\,\big]\\[10pt]
\dis\quad\mbox{implies}\quad
\P\big[(X_k,Y_k,Z_k)\in\cdot\,\big]\Asto{k}
\P\big[(X,f(X),g(f(X))\big)\in\cdot\,\big],
\ec
where $\Rightarrow$ denotes weak convergence of probability laws on the Polish
spaces $E\times F$, $F\times G$, and $E\times F\times G$, respectively.
\end{itemize}
\el
{\bf Proof.} For part~(a), we observe that the assumed weak convergence of $(X_k,Y_{k,i})$ for each $i\in I$ implies tightness of the laws of $(X_k,(Y_{k,i})_{i\in I})$. Let $(X,(Y_i)_{i\in I},\ldots)$ be any weak subsequential limit. Then $(X,Y_i)$ is equally distributed with $(X,f_i(X))$, hence $Y_i=f_i(X)$ a.s.\ for each $i\in I$. Similarly, in the set-up of part~(b), the weak convergence of $(X_k,Y_k)$ and $(Y_k,Z_k)$ implies tightness of the laws of $(X_k,Y_k,Z_k)$, while for each weak subsequential limit $(X,Y,Z)$ one has $Y=f(X)$ and $Z=g(Y)$ a.s.\qed

We note that by Skorohod's representation theorem (see
e.g.\ \cite[Theorem~6.7]{Bil99}) the left-hand side of (\ref{weakcouple})
implies that for each $i\in I$, we can find a coupling of the $X_k,Y_{k,i}$
and $X$ such that $(X_k,Y_{k,i})\to(X,f_i(X))$ a.s. By the right-hand side of
(\ref{weakcouple}), we can find a coupling that works for all $i\in I$
simultaneously. We will apply this principle many times, e.g.\ when $X$ is a
Browian web, $f(X)$ is its associated dual Brownian web, $g(X)$ is the set of
paths starting at a given time etc. We will not always be explicit in our
choice of the measurable maps $f,g$ but it is clear from the context
that they can be constructed.

Recall from Section~\ref{S:disref} that each i.i.d.\ collection
$\al=(\al_z)_{z\in\Zev}$ of $\{-1,+1\}$-valued random variables defines a
discrete web $\Ui^\al=\{p^\al_z:z\in\Zev\}$ as in (\ref{Uial}). As in
Section~\ref{S:approx}, by linear interpolation and by adding trivial paths
that are constantly $-\infty$ or $+\infty$, we view $\Ui^\al$ as a random
compact subset of the space of paths $\Pi$ introduced in
Section~\ref{S:BMweb}.

Let $\Zod:=\{(x,t):x,t\in\Z,\ x+t\mbox{ is odd}\}$ be the odd sublattice of
$\Z^2$. For each $(x,s)\in\Zod$, we let $\hat p^\al_{(x,s)}=\hat p$, defined
by (compare (\ref{pal}))
\be\label{dupal}
\hat p(s):=x\quad\mbox{and}\quad
\hat p(t-1):=\hat p(t)-\al_{(\hat p(t),t-1)}\qquad(t\leq s)
\ee
denote the dual path started at $(x,s)$ and we let $\hat\Ui^\al=\{\hat
p^\al_z:z\in\Zev\}$ denote the dual discrete web associated with $\Ui^\al$.
We view $\hat\Ui^\al$ as a random compact subset of the space of dual paths
$\hat\Pi$. In line with earlier notation, for any $A\sub\Zev$
(resp.\ $A\sub\Zod$), we let $\Ui^\al(A)$ (resp.\ $\hat\Ui^\al(A)$) denote the
set of paths in $\Ui^\al$ (resp.\ $\hat\Ui^\al$) starting from $A$. We define
diffusive scaling maps $S_\eps$ as in (\ref{Seps2}) and use
$S_\eps(\Ai_1,\ldots,\Ai_n)$ as a shorthand for
$(S_\eps(\Ai_1),\ldots,S_\eps(\Ai_n))$.

The following result follows easily from \cite[Theorem~6.1]{FINR04} on the
convergence of discrete webs to the Brownian web and Proposition~\ref{P:char}
on the characterization of the dual Brownian web.

\bt{\bf(Convergence to the double Brownian web)}\label{T:webcon}
Let $\eps_k$ be positive constants, tending to zero. For each $k$, let
$\al^{\li k\re}=(\al^{\li k\re}_z)_{z\in\Zev}$ be an i.i.d.\ collection of
$\{-1,+1\}$-valued random variables, let $\Ui_{\li k\re}:=\Ui^{\al^{\li
    k\re}}$ and $\hat\Ui_{\li k\re}:=\hat\Ui^{\al^{\li k\re}}$ be the discrete
web and dual discrete web associated with $\al^{\li k\re}$, and assume that
$\lim_{k\to\infty}\eps_k^{-1}\E[\al^{\li k\re}_z]=\bet$ for some
$\bet\in\R$. Then
\be\label{webcon}
\P\big[S_{\eps_k}(\Ui_{\li k\re},\hat\Ui_{\li k\re})\in\cdot\,\big]
\Asto{k}\P\big[(\Wi,\hat\Wi)\in\cdot\,\big],
\ee
where $\Rightarrow$ denotes weak convergence of probability laws on
$\Ki(\Pi)\times\Ki(\hat\Pi)$, $\Wi$ is a Brownian web with drift $\bet$ and
$\hat\Wi$ is its dual.
\et

For notational convenience, let us write
\be\label{SigT}
\Sig_T:=\big\{(x,t)\in\Rc:t=T\big\},
\ee
so that, e.g., $\Ni(\Sig_T)=\{\pi\in\Ni:\sig_\pi=T\}$. We use similar notation
for sets of discrete paths. The following strengthening of
Theorem~\ref{T:webcon} is sometimes handy.

\bl{\bf(Convergence of paths started at given times)}\label{L:webTcon}
In the setup of Theorem~\ref{T:webcon}, let $T_k\in\Z\cup\{-\infty,\infty\}$
be times such that $\eps_k^2T_k\to T$ for some $T\in[-\infty,+\infty]$. Then
\be\label{webKcon}
\P\big[S_{\eps_k}\big(\Ui_{\li k\re},
\Ui_{\li k\re}(\Sig_{T_k})\big)\in\cdot\,\big]
\Asto{k}\P\big[\big(\Wi,\Wi(\Sig_T)\big)\in\cdot\,\big].
\ee
\el
{\bf Proof.} It follows from the tightness of the $S_{\eps_k}\big(\Ui_{\li
  k\re})$ and Lemma~\ref{L:Hautight} that also the laws of the
$\Ki(\Pi)^2$-valued random variables $S_{\eps_k}\big((\Ui_{\li k\re},\Ui_{\li
  k\re}(\Sig_{T_k}))$ are tight. By going to a subsequence if necessary and invoking
Skorohod's representation theorem, we may assume that they converge to an
a.s.\ limit $(\Wi,\Ai)$. It is easy to see that $\Ai\sub\Wi(\Sig_T)$. If
$T=\pm\infty$, then $\Wi(\Sig_T)$ contains only trivial paths and it is easy
to check that also $\Ai\supset\Wi(\Sig_T)$. To get this inclusion for
$-\infty<T<\infty$, let $\Di$ be a deterministic countable dense subset of
$\R\times\{T\}$. Since $\Ai(z)$ is nonempty for each $z\in\Di$ and since
$\Wi(z)$ contains a single path for each $z\in\Di$, we conclude that
$\Ai\supset\Wi(\Di)$. Since $\Ai$ is compact and $\Wi(\Sig_T)$ is the closure
of $\Wi(\Di)$, it follows that $\Ai=\Wi(\Sig_T)$.\qed

We next formulate a result which says that the intersection local time measure
$\ell$ between a forward and dual Brownian web as defined in
Proposition~\ref{P:refloc} is the limit of the intersection local time
measures between approximating forward and dual discrete webs. Since $\ell$ is
locally infinite, such a statement on its own cannot make sense. Rather, we
will show that the restriction of $\ell$ to the intersection of finitely many
forward and dual paths is a.s.\ the weak limit of the analogue discrete
object.

For any $K\sub\Ki(\Rc)$, we let
\be\label{Img}
{\rm Img}(K):=\{z\in\Rc:\exists A\in K\mbox{ s.t.\ }z\in A\}
\ee
denote the union of all sets in $K$. We call ${\rm Img}(K)$ the {\em image
  set} (or {\em trace}) of $K$. In particular, if $\Ai$ is a set of paths
(which, as usual, we identify with their graphs), then ${\rm
  Img}(\Ai)=\{(\pi(t),t):t\geq\sig_\pi,\ \pi\in\Ai\}$.  Similarly, if $\Ai$ is
a set of discrete paths, then
\be
{\rm Img}(\Ai):=\{(x,t)\in\Zev:t\geq\sig_\pi,\ \pi\in\Ai\},
\ee
and we use similar notation for a set $\hat\Ai$ of discrete dual paths, where
in this case ${\rm Img}(\hat\Ai)$ is a subset of $\Zod$.

\bp{\bf(Convergence of intersection local time)}\label{P:ellcon}
Let $\eps_k$ be positive constants, tending to zero. Let $\al^{\li k\re}$ be
collections of i.i.d.\ $\{-1,+1\}$-valued random variables satisfying
$\lim_{k\to\infty}\eps_k^{-1}\E[\al^{\li k\re}_z]=\bet$ for some $\bet\in\R$,
and let $\Ui_{\li k\re},\hat\Ui_{\li k\re}$ be the discrete web and its dual
associated with $\al^{\li k\re}$. Let $\Wi$ be a Brownian web
with drift $\bet$, $\hat\Wi$ be its dual, $\ell$ be the intersection local time measure
between $\Wi$ and $\hat\Wi$, and $\ell_{\rm r}$ be the restriction of $\ell$ to the
set of points of type $(1,2)_{\rm r}$. Let
\be\ba{r@{\,}c@{\,}lr@{\,}c@{\,}l}
\dis\De_k&=&\dis\{z^k_1,\ldots,z^k_m\}\sub\Zev,
\quad&\dis\hat\De_k
&=&\dis\{\hat z^k_1,\ldots,\hat z^k_n\}\sub\Zod,\\[5pt]
\dis\De&=&\dis\{z_1,\ldots,z_m\}\sub\R^2,
\quad&\dis\hat\De
&=&\dis\{\hat z_1,\ldots,\hat z_n\}\sub\R^2
\ec
be finite sets such that $S_{\eps_k}(z^k_i)\to z_i$ and
$S_{\eps_k}(\hat z^k_j)\to\hat z_j$ as $k\to\infty$ for each
$i=1,\ldots,m$ and $j=1,\ldots,n$. Set
$$
\dis\ell_{{\rm r}}^{\li k\re} :=\eps_k\sum_{z\in I_k\cap Z_{{\rm r}}^{\li k\re}}\de_{S_{\eps_k}(z)},
$$
where
\be\label{IkZk}
\begin{aligned}
Z_{{\rm r}}^{\li k\re}&:=\{z\in\Zev:\al^{\li k\re}_{(x,t)}=+1\}, \\
I_k: &=\big\{(x,t)\in\Zev:(x,t)\in{\rm Img}(\Ui_{\li k\re}(\De_k)),\
(x,t+1)\in{\rm Img}(\hat\Ui_{\li k\re}(\hat\De_k))\big\}.
\end{aligned}
\ee
Let $\ell_{\rm r}(\De,\hat\De)$ denote the restriction of $\ell_{\rm r}$
to the set ${\rm Img}(\Wi(\De))\cap{\rm Img}(\hat\Wi(\hat\De))$. Then
\be\label{localconv}
\P\big[\big(S_{\eps_k}(\Ui_{\li k\re}),\ell_{{\rm r}}^{\li k\re}\big) \in\cdot\big]
\Asto{k} \P\big[\big(\Wi, \ell_{\rm r}(\De,\hat\De)\big) \in \cdot\big],
\ee
where $\Rightarrow$ denotes weak convergence of probability laws on $\Ki(\Pi)\times\Mi(\R^2)$, and $\Mi(\R^2)$ is the space of
finite measures on $\R^2$ equipped with the topology of weak convergence.
\ep
{\bf Proof.} Since
$$
\P\big[S_{\eps_k}(\Ui_{\li k\re}, \hat \Ui_{\li k\re}, \Ui_{\li k\re}(\Delta_k), \hat\Ui_{\li k\re}(\hat\Delta_k))\in\cdot\big] \Asto{k}
\P\big[(\Wi, \hat\Wi, \Wi(\Delta), \hat\Wi(\hat\Delta)) \in \cdot\big],
$$
and a.s.\ $(\Wi(\Delta), \hat\Wi(\hat\Delta))$ determines $\ell_{\rm r}(\De,\hat\De)$, by Lemma~\ref{L:weakcouple}~(b), proving (\ref{localconv}) reduces to proving
\be\label{localconv2}
\P\big[\big(S_{\eps_k}(\Ui_{\li k\re}(\Delta_k),\hat \Ui_{\li k\re}(\hat \Delta_k)), \ell_{{\rm r}}^{\li k\re}\big)\in\cdot\,\big]
\Asto{k}\P\big[\big(\Wi(\Delta), \hat\Wi(\hat\Delta)),\ell_{\rm r}(\De,\hat\De)\big)\in\cdot\,\big].
\ee
We will make a further reduction.

For $k\in\N$, $1\leq i\leq m$ and $1\leq j\leq n$, let $\Ui_{\li k\re}(z^k_i)= \{p^k_i\}$ and $\hat\Ui_{\li k\re}(\hat z^k_j)=\{\hat p^k_j\}$. Let $t^k_i$ and
$\hat t^k_j$ denote respectively the starting time of $p^k_i$ and $\hat p^k_j$. Similarly, let $\Wi(z_i)=\{\pi_i\}$ and $\hat\Wi(\hat z_j)=\{\hat\pi_j\}$,
with starting time $t_i$ for $\pi_i$ and $\hat t_j$ for $\hat \pi_j$. Let $(\tau^k_{uv})_{1\leq u<v\leq m}$ be the time of coalescence between $p^k_u$ and $p^k_v$, and let $(\hat\tau^k_{uv})_{1\leq u<v\leq n}$ be the time of coalescence between $\hat p^k_u$ and $\hat p^k_v$. Define $(\tau_{uv})_{1\leq u<v\leq m}$ and $(\hat\tau_{uv})_{1\leq u<v\leq n}$ similarly for $\Wi(\Delta)$ and $\hat\Wi(\hat\Delta)$. For $1\leq u\leq m$ and $1\leq v\leq n$, let
$$
\ell^{\li k\re}_{{\rm r}, uv}:=\eps_k\sum_{z:=(x,t)\in\Z^2_{\rm even}}\delta_{S_{\eps_k}(z)} 1_{\{p^k_u(t)=\hat p^k_v(t+1)=\hat p^k_v(t)+1=x\}}
$$
be the rescaled intersection local time measure of $p^k_u$ and $\hat p^k_v$ on points with $\alpha^{\li k\re}_{(x,t)}=1$, and similarly let $\ell_{{\rm r},uv}$ be the intersection local time measure of $\pi_u$ and $\hat \pi_v$ on points of type $(1,2)_{\rm r}$. We note that $\ell_{\rm r}(\Delta, \hat\Delta)$ can be uniquely constructed from
$(\tau_{uv})_{1\leq u<v\leq m}$, $(\hat \tau_{uv})_{1\leq u<v\leq n}$, and $(\ell_{{\rm r}, uv})_{1\leq u\leq m,1\leq v\leq n}$. For example, we can go through the
indices $(uv)_{1\leq u\leq m, 1\leq v\leq n}$ in numeric order, and at each step, we add the proper restriction of $\ell_{{\rm r}, uv}$ to the construction of
$\ell_{\rm r}(\Delta, \hat\Delta)$ so as to exclude overlaps among $(\ell_{{\rm r}, uv})_{1\leq u\leq m,1\leq v\leq n}$ due to the coalescence of paths. By the same procedure, $\ell^{\li k\re}_{\rm r}$ can be constructed from  $(\tau^k_{uv})_{1\leq u<v\leq m}$, $(\hat\tau^k_{uv})_{1\leq u<v\leq n}$, and $(\ell^{\li k\re}_{{\rm r}, uv})_{1\leq u\leq m,1\leq n\leq v}$. To prove (\ref{localconv2}), it then suffices to prove
\begin{eqnarray}
\!\!\!\!\!\!&& \!\!\!\!\!\!\! \P\Big[\Big(S_{\eps_k}\big((p^k_i)_{1\leq i\leq m}, (\hat p^k_j)_{1\leq j\leq n}\big),(\eps_k^2\tau^k_{uv})_{1\leq u<v\leq m}, (\eps_k^2\hat\tau^k_{uv})_{1\leq u<v\leq n}, (\ell^{\li k\re}_{{\rm r}, uv})_{1\leq u\leq m,1\leq n\leq v} \Big)\in\cdot \Big]  \nonumber\\
\!\!\!\!\!\!&\Asto{k} & \!\!\! \P\Big[\Big( (\pi_i)_{1\leq i\leq m}, (\hat \pi_j)_{1\leq j\leq n}, (\tau_{uv})_{1\leq u<v\leq m}, (\hat \tau_{uv})_{1\leq u<v\leq n}, (\ell_{{\rm r}, uv})_{1\leq u\leq m,1\leq v\leq n}\Big)\in\cdot \Big]. \label{localconv2.5}
\end{eqnarray}

It has been shown in the proof of \cite[Thm.~8]{STW00} that
\begin{eqnarray} \nonumber
&& \P\big[S_{\eps_k}\big( (p^k_i)_{1\leq i\leq m}, (\hat p^k_j)_{1\leq j\leq n},(\tau^k_{uv})_{1\leq u<v\leq m}, (\hat\tau^k_{uv})_{1\leq u<v\leq n}\big) \in\cdot \big]\\
&\Asto{k}& \P\big[ \big( (\pi_i)_{1\leq i\leq m}, (\hat \pi_j)_{1\leq j\leq n}, (\tau_{uv})_{1\leq u<v\leq m}, (\hat \tau_{uv})_{1\leq u<v\leq n} \big)\in\cdot \big], \label{STWconv}
\end{eqnarray}
where $(p^k_i)_{1\leq i\leq m}$ and $(\hat p^k_j)_{1\leq j\leq n}$ were constructed as a deterministic transformation (via Skorohod reflection and coalescence) of a collection of independent random walks $(W^k_i)_{1\leq i\leq m}$ and $(\hat W^k_j)_{1\leq j\leq n}$. The same transformation was used to construct $(\pi_i)_{1\leq i\leq m}$ and $(\hat\pi_j)_{1\leq j\leq n}$ from a collection of independent Brownian motions $(B_i)_{1\leq i\leq m}$ and $(\hat B_j)_{1\leq j\leq n}$. Furthermore, this transformation together with the times of coalescence $(\tau_{uv})_{1\leq u<v\leq m}$ and
$(\hat\tau_{uv})_{1\leq u<v\leq n}$ are a.s.\ continuous in $(B_i)_{1\leq i\leq m}$ and $(\hat B_j)_{1\leq j\leq n}$. The convergence in (\ref{STWconv}) then follows
from Donsker's invariance principle. Since $\ell_{{\rm r},uv}$ is uniquely determined by $\pi_u$ and $\hat \pi_v$, by Lemma~\ref{L:weakcouple}, to prove (\ref{localconv2.5}), it then suffices to show that for each $1\leq u\leq m$ and $1\leq v\leq m$,
\be\label{pairconv}
\P\big[\big(S_{\eps_k}(p^k_u, \hat p^k_v), \ell^{\li k\re}_{{\rm r}, uv}\big)\in\cdot \big] \Asto{k} \P\big[(\pi_u, \hat\pi_v, \ell_{{\rm r},uv})\in\cdot \big].
\ee

Without loss of generality, we may assume $u=v=1$ in (\ref{pairconv}). We may also assume that $z_1=(x_1,t_1)$ and
$\hat z_1=(\hat x_1, \hat t_1)$ satisfy $t_1<\hat t_1$, so that $\ell_{\rm r, 11}$ is not a.s.\ the zero measure, in which case (\ref{pairconv}) is trivial.
We recall from \cite[Thm.~3]{STW00} (see also \cite[(3.6) and Thm.~3.7]{FINR06}) that, conditional on $\hat \pi_1$, $\pi_1$ is distributed as an independent Brownian motion $B_1$ with drift $\beta$, starting from $z_1$, and Skorohod reflected away from $\hat\pi_1$. More precisely, $\pi_1$ admits the representation
\be\label{skoropair}
\pi_1(t) = \left\{
\begin{aligned}
& B_1(t) + L_{\rm r}(t) \qquad & \quad \mbox{if} \ \hat\pi_1(t_1)<x_1, \\
& B_1(t) - L_{\rm l}(t) \qquad & \quad \mbox{if} \ x_1<\hat\pi_1(t_1),
\end{aligned}
\right.
\ee
where
\be\label{refloc}
\begin{aligned}
L_{\rm r}(t) & =\sup_{t_1\leq s\leq t}\max\{0, \hat \pi_1(s)-B_1(s)\} \ \ \mbox{ for  }t\in [t_1, \hat t_1], \qquad & L_{\rm r}(t)=L_{\rm r}(\hat t_1) \mbox{  for  }t\geq \hat t_1,\\
L_{\rm l}(t) & =\sup_{t_1\leq s\leq t}\max\{0, B_1(s)- \hat\pi_1(s)\} \ \ \mbox{ for  }t\in [t_1, \hat t_1], \qquad & L_{\rm l}(t)=L_{\rm l}(\hat t_1) \mbox{  for  }t\geq \hat t_1.
\end{aligned}
\ee
It was shown in \cite[Prop.~3.1]{NRS10} and its proof\footnote{Our definition of $\ell$ in (\ref{elldef}) differs from the definition in \cite[(3.2)]{NRS10} by a factor of $\sqrt{2}$, which is compensated by the fact that we consider Poisson point process in $\R^2$ with intensity measure $\ell$ instead of $\sqrt{2}\ell$, as done in \cite[(3.8)]{NRS10}.}
that, with the construction of $\pi_1$ as in (\ref{skoropair}), almost surely $\ell_{\rm r, 11}(\R\times\cdot)={\rm d}L_{\rm r}(\cdot)$, or equivalently,
\be\label{ellrep}
\ell_{\rm r,11}(\R\times[t_1, t])=L_{\rm r}(t) \qquad \mbox{for all } t\in [t_1, \hat t_1].
\ee
Since $\ell_{\rm r, 11}$ is concentrated on the graph of $\pi_1$, it follows that
\be\label{imagemap1}
\ell_{\rm r, 11}=\di L_r\circ\pi_1^{-1},
\ee
i.e., $\ell_{\rm r, 11}$ is the image of the measure $\di L_r$ under the map $\pi_1$.

There is a similar representation for $p^k_1$, $\hat p^k_1$, and $\ell^{\li k\re}_{\rm r, 11}$. Indeed, if $W^k_1$ is an independent simple random walk on $\Z$ with drift $\E[\alpha^{\li k\re}_z]$ and starting from $z^k_1$, then conditional on $\hat p^k_1$, we can construct $p^k_1$ as (see e.g.\ \cite[Sec.~2.2.2]{STW00} or the proof of \cite[Lemma 2.1]{SSS09})
\be\label{skoropair2}
p^k_1(t) = \left\{
\begin{aligned}
W^k_1(t) + L^k_{\rm r}(t) \qquad & \quad \mbox{if} \ \hat p^k_1(t^k_1) <x^k_1, \\
W^k_1(t) - L^k_{\rm l}(t) \qquad & \quad \mbox{if} \ x^k_1< \hat p^k_1(t^k_1), \\
\end{aligned}
\right.
\ee
where
\be\label{refloc2}
\begin{aligned}
\!\!\!\! L^k_{\rm r}(t) & =\sup_{t^k_1\leq s\leq t}\max\{0, 1+\hat p^k_1(s)-W^k_1(s)\} \  \mbox{ for  }t\in [t^k_1, \hat t^k_1], \ \  & L^k_{\rm r}(t)=L^k_{\rm r}(\hat t^k_1) \mbox{  for  }t\geq \hat t^k_1,\\
\!\!\!\! L^k_{\rm l}(t) & =\sup_{t^k_1\leq s\leq t}\max\{0, 1+W^k_1(s)- \hat p^k_1(s)\} \  \mbox{ for  }t\in [t^k_1, \hat t^k_1], \ \ & L^k_{\rm l}(t)=L^k_{\rm l}(\hat t^k_1) \mbox{  for  }t\geq \hat t^k_1.
\end{aligned}
\ee
The constant $1$ arises because $\hat p^k_1$ is a walk on $\Z^2_{\rm odd}$, and $W^k_1$ is a walk on $\Z^2_{\rm even}$. Let
$$
\bar L^k_{\rm r}(t) = \sum_{i=t^k_1}^t 1_{\{p^k_1(i)=\hat p^k_1(i+1)=\hat p^k(i)+1\}} \qquad \mbox{for} \ t^k_1 \leq t <\hat t^k_1.
$$
Then, in analogy with (\ref{imagemap1}),
\be\label{imagemap2}
\ell^{\li k\re}_{\rm r, 11}={\rm d}(S_{\eps_k}\bar L^k_{\rm r})\circ(p^k_1)^{-1},
\ee
where ${\rm d}(S_{\eps_k}\bar L^k_{\rm r})(t)=\eps_k{\rm d}\bar L^k_{\rm r}(\eps_k^{-2} t)$. To relate $L^k_{\rm r}$ and $\bar L^k_{\rm r}$, note that conditional on $\hat p^k_1$,
\be\label{reflocerr1}
L^k_{\rm r}(t+1)-(1-\E[\alpha^{\li k\re}_z])\bar L^k_{\rm r}(t)= \sum_{i=t^k_1}^t \big(L^k_{\rm r}(i+1)-L^k_{\rm r}(i)-(1-\E[\alpha^{\li k\re}_z])1_{\{\bar L^k_{\rm r}(i)-\bar L^k_{\rm r}(i-1)=1\}}\big)
\ee
is a martingale, because $L^k_{\rm r}(i+1)-L^k_{\rm r}(i)\neq 0$ only when $p^k_1(i)=\hat p^k_1(i+1)=\hat p^k(i)+1$, and conditional on the later event,
$L^k_{\rm r}(i+1)-L^k_{\rm r}(i)=0$ with probability $\frac{1+\E[\alpha^{\li k\re}_z]}{2}$ and $L^k_{\rm r}(i+1)-L^k_{\rm r}(i)=2$ with probability $\frac{1-\E[\alpha^{\li k\re}_z]}{2}$. By Doob's maximal inequality, conditional on $\hat p^k_1$,
\begin{eqnarray}
\E\Big[\sup_{t^k_1\leq t <\hat t^k_1} \big|L^k_{\rm r}(t+1)-(1-\E[\alpha^{\li k\re}_z])\bar L^k_{\rm r}(t)\big|^2\Big]
&\leq& 4 \E\Big[ \big|L^k_{\rm r}(\hat t^k_1)-(1-\E[\alpha^{\li k\re}_z])\bar L^k_{\rm r}(\hat t^k_1)\big|^2\Big] \nonumber\\
&\leq& 16 \E[\bar L^k_{\rm r}(\hat t^k_1)]. \label{reflocerr2}
\end{eqnarray}

We are now ready to prove (\ref{pairconv}). By (\ref{imagemap1}) and (\ref{imagemap2}), we can replace $\ell_{\rm r,11}$ by $L_r$, and $\ell^{\li k\re}_{\rm r,11}$ by
$S_{\eps_k} \bar L^k_{\rm r}$. First let us extend the definition of all processes in discrete time to continuous time by linear interpolation. Note that (\ref{skoropair2}) and (\ref{refloc2}) remain valid. By Donsker's invariance principle, the pair of independent processes $S_{\eps_k}(\hat p^k_1, W^k_1)$ converge weakly to $(\hat \pi_1, B_1)$. By Skorohod's representation, we may assume from now on this convergence is almost sure by suitable coupling. If $x_1<\hat \pi_1(t_1)$, then trivially
\be\label{pairconv2}
S_{\eps_k}(p^k_1, \hat p^k_1, L^k_{\rm r}) \asto{k} (\pi_1, \hat \pi_1, L_r)
\ee
since $L_{\rm r}=0$, and so is $L^k_{\rm r}$ for all $k$ large. If $\hat \pi_1(t_1)<x_1$, then $S_{\eps_k}(L^k_{\rm r})\to L_{\rm r}$ uniformly on compacts, because the Skorohod reflection map which defines both $L^k_{\rm r}$
and $L_{\rm r}$ in (\ref{refloc2}) and respectively (\ref{refloc}) is continuous in its arguments. Therefore,
\be\label{pairconv3}
\P\big[S_{\eps_k}(p^k_1, \hat p^k_1, L^k_{\rm r})\in \cdot\big] \Asto{k} \P[(\pi_1, \hat \pi_1, L_{\rm r})\in\cdot].
\ee
On the other hand, by (\ref{reflocerr2}),
\be\label{pairconv4}
\E\Big[\sup_{t^k_1\leq t <\hat t^k_1} \eps_k^2\big|L^k_{\rm r}(t+1)-(1-\E[\alpha^{\li k\re}_z])\bar L^k_{\rm r}(t)\big|^2\Big] \leq 16 \eps_k^2\E[\bar L^k_{\rm r}(\hat t^k_1)]\asto{k} 0,
\ee
because the above inequality implies by triangle inequality that
$$
(1-\E[\alpha^{\li k\re}_z])\E[ (\eps_k \bar L^k_{\rm r}(\hat t^k_1))^2]^{\frac{1}{2}} \leq \E[(\eps_k L^k_{\rm r}(\hat t^k_1))^2]^{\frac{1}{2}} + 4\sqrt{\eps_k} \E[\eps_k \bar L^k_{\rm r}(\hat t^k_1)]^{\frac{1}{2}},
$$
and since $\E[(\eps_k L^k_{\rm r}(\hat t^k_1))^2]$ is uniformly bounded in $k$ as easily seen from the definition of $L^k_{\rm r}$, so are
$\E[(\eps_k \bar L^k_{\rm r}(\hat t^k_1))^2]$ and $\E[\eps_k \bar L^k_{\rm r}(\hat t^k_1)]$. Since $\E[\alpha^{\li k\re}_z]\to 0$, we conclude from (\ref{pairconv4}) that
\be\label{pairconv5}
\P[S_{\eps_k}(L^k_r, \bar L^k_r)\in\cdot ] \Asto{k} \P[(L_{\rm r}, L_{\rm r})\in \cdot].
\ee
By Lemma~\ref{L:weakcouple}, (\ref{pairconv3}) and (\ref{pairconv5}) imply that
\be
\P\big[S_{\eps_k}(p^k_1, \hat p^k_1, \bar L^k_{\rm r})\in\cdot \big] \Asto{k} \P\big[(\pi_1, \hat\pi_1, L_{\rm r})\in\cdot\big],
\ee
which in turn implies (\ref{pairconv}) and concludes our proof.
\qed

\subsection{Discrete approximation of the Brownian net}\label{S:netap}

It has been shown in \cite{SS08} that the Brownian net arises as the limit of
systems of branching-coalescing random walks, in the limit of small
branching probability and after diffusive rescaling. In this section, we
review this result and add some additional results on the approximation of
(relevant) separation points by discrete separation points.

Let $\bet_-\leq\bet_+$ be real constants. Let $\eps_k$ be
positive constants, converging to zero, and for each $k$, let
\be\label{alr}
(\al^{\li k\re\,{\rm l}},\al^{\li k\re\,{\rm r}})
=(\al^{\li k\re\,{\rm l}}_z,\al^{\li k\re\,{\rm r}}_z)_{z\in\Zev}
\ee
be an i.i.d.\ collection of $\{-1,+1\}^2$-valued random variables such that
$\al^{\li k\re\,{\rm l}}_z\leq\al^{\li k\re\,{\rm r}}_z$ and
\be\label{alr2}
\eps_k^{-1}\E[\al^{\li k\re\,{\rm l}}_z]\asto{k}\bet_-
\quad\mbox{and}\quad\eps_k^{-1}\E[\al^{\li k\re\,{\rm r}}_z]\asto{k}\bet_+.
\ee
We let $\Ui^{\rm l}_{\li k\re}$ and $\Ui^{\rm r}_{\li k\re}$ denote the
discrete webs associated with $\al^{\li k\re\,{\rm l}}$ and $\al^{\li
  k\re\,{\rm r}}$, respectively. Then $(\Ui^{\rm l}_{\li k\re},\Ui^{\rm
  r}_{\li k\re})$ is a discrete analogue of a left-right Brownian web as
introduced in Section~\ref{S:modweb}. We call the collection of discrete paths
\be\label{disnet1}
\Vi_{\li k\re}
:=\big\{p:p(t+1)-p(t)\in
\{\al^{\li k\re\,{\rm l}}_{(p(t),t)},\al^{\li k\re\,{\rm r}}_{(p(t),t)}\}
\ \forall t\geq\sig_p\big\}
\ee
the {\em discrete net}\index{discrete!net} associated with $(\Ui^{\rm l}_{\li k\re},\Ui^{\rm
  r}_{\li k\re})$. We observe that except for a rotation by 180 degrees and
a shift from $\Zev$ to $\Zod$, the discrete dual left-right web $(\hat\Ui^{\rm
  l}_{\li k\re},\hat\Ui^{\rm r}_{\li k\re})$ is equally distributed with
$(\Ui^{\rm l}_{\li k\re},\Ui^{\rm r}_{\li k\re})$. In view of this, we define
a {\em dual discrete net}\label{dual!discrete net} $\hat\Vi_{\li k\re}$ analogously to $\Vi_{\li
  k\re}$. As in Section~\ref{S:approx}, we view the sets of discrete
paths $\Ui^{\rm l}_{\li k\re},\Ui^{\rm r}_{\li k\re},\Vi_{\li k\re}$ as random
compact subsets of the space of continuous paths $\Pi$.

We cite the following result from \cite[Thm.~5.4]{SS08}.

\bt{\bf(Convergence to the Brownian net)}\label{T:netcon}\ \\
Let $\eps_k$ and $\Vi_{\li k\re},\Ui^{\rm l}_{\li k\re},\Ui^{\rm
  r}_{\li k\re},\hat\Vi_{\li k\re},\hat\Ui^{\rm l}_{\li k\re},\hat\Ui^{\rm
  r}_{\li k\re}$ be as above and let $\Ni,\Wl,\Wr,\hat\Ni,\hat\Wl,\hat\Wr$ be
a Brownian net with left and right speeds $\bet_-\leq\bet_+$, its associated
left-right Brownian web, and their duals. Then
\be\label{netcon}
\P\big[S_{\eps_k}(\Vi_{\li k\re},\Ui^{\rm l}_{\li k\re},\Ui^{\rm r}_{\li
  k\re},\hat\Vi_{\li k\re},\hat\Ui^{\rm l}_{\li k\re},\hat\Ui^{\rm r}_{\li
  k\re})\in\cdot\,\big]
\Asto{k}\P\big[(\Ni,\Wl,\Wr,\hat\Ni,\hat\Wl,\hat\Wr)\in\cdot\,\big],
\ee
where $\Rightarrow$ denotes weak convergence of probability laws on
$\Ki(\Pi)^3\times\Ki(\hat\Pi)^3$.
\et

The following analogue of Lemma~\ref{L:webTcon} is sometimes handy.

\bl{\bf(Convergence of paths started at a given time)}\label{L:netTcon}
In the setup of Theorem~\ref{T:netcon}, let
$T_k\in\Z\cup\{-\infty,+\infty\}$ satisfy $\eps_k^2T_k\to T$ for some
$T\in[-\infty,\infty]$. Then
\be\label{netTcon}
\P\big[S_{\eps_k}\big(\Vi_{\li k\re},\Vi_{\li k\re}(\Sig_{T_k})\big)
\in\cdot\,\big]
\Asto{k}\P\big[\big(\Ni,\Ni(\Sig_T)\big)\in\cdot\,\big].
\ee
\el
{\bf Proof.} By Lemma~\ref{L:Hautight} in the appendix, the tightness of the
$S_{\eps_k}\big(\Vi_{\li k\re})$ implies that also the laws of the
$\Ki(\Pi)^2$-valued random variables $S_{\eps_k}\big((\Vi_{\li k\re},\Vi_{\li
  k\re}(\Sig_{T_k}))$ are tight. By going to a subsequence if necessary and
invoking Skorohod's representation theorem, we may assume that they converge
to an a.s.\ limit $(\Ni,\Ai)$. It is easy to see that $\Ai\sub\Ni(\Sig_T)$. To
get the other inclusion, we distinguish three cases. The case $T=+\infty$ is
trivial. If $-\infty<T<+\infty$, let $\Di^{\rm l},\Di^{\rm r}$ be
deterministic countable dense subsets of $\R^2$ such that $\Di^{\rm l}$ is
also dense in $\R\times\{T\}$. Let $\Hi_{\rm cros}(\Wl(\Di^{\rm
  l}),\Wr(\Di^{\rm r}))$ be the set of paths that can be obtained by
concatenating finitely many paths in $\Wl(\Di^{\rm l})$ and $\Wr(\Di^{\rm r})$
at crossing times between left and right paths. Arguing as in the proof of
\cite[Thm.~5.4]{SS08}, we obtain that
\be
\Ai\supset
\Pi(\Sig_T)\cap\Hi_{\rm cros}(\Wl(\Di^{\rm l}),\Wr(\Di^{\rm r})),
\ee
hence by \cite[Lemma~8.1]{SS08} we conclude that
$\Ai\supset\Ni(\Sig_T)$. Finally, if $T=-\infty$, then let $\Vi_{\li
  k\re}(\Sig_{-\infty})\big|_{T_k}^\infty$ denote the set of all restrictions
of paths in $\Vi_{\li k\re}(\Sig_{-\infty})$ to the time interval
$[T_k,\infty]$. Since $\Vi_{\li
  k\re}(\Sig_{-\infty})\big|_{T_k}^\infty\sub\Vi_{\li k\re}(\Sig_{T_k})$, it
then suffices to prove the claim if $T_k=-\infty$ for all $k$. But this is
just \cite[Lemma~9.2]{SS08}.\qed

We will need one more result that is very close in spirit to
Lemma~\ref{L:netTcon} and can in fact be seen as a strengthening of the
latter. If $\Ni$ is a Brownian net with left and right speeds
$\bet_-\leq\bet_+$, then, generalizing (\ref{braco}), for any closed
$A\sub\R$, we may define a Markov process taking values in the closed subsets
of the real line by
\be\label{xiA}
\xi^A_t:=\big\{\pi(t):\pi\in\Ni(A\times\{0\})\big\}\qquad(t\geq 0).
\ee
We call $\xi^A$ the {\em branching-coalescing point set with left and right
  speeds $\bet_-,\bet_+$}. \index{branching-coalescing point set!with finite speeds} By combining \cite[Prop.~1.12]{SS08}, Brownian
scaling, and the well-known density of the Arratia flow (see
\cite[equation~(1.6)]{FINR02}), it is easy to check that the density of
$\xi^\R_t$ is given by
\be
\E\big[\big|\xi^\R\cap[x,y]\big|\big]=(y-x)\Psi_b(t)\qquad(t>0,\ x<y),
\ee
where $b:=(\bet_+-\bet_-)/2$ and $\Psi_b$ is the function in (\ref{Psib}).

If $\Vi$ is a discrete net defined from an i.i.d.\ collection of random
variables $(\al^{\rm l}_z,\al^{\rm r}_z)_{z\in\Zev}$ as in (\ref{disnet1}),
then we can define a discrete branching-coalescing point set in analogy with
(\ref{xiA}). In particular, we let
\be\label{disdens}
\Psi_{b_-,b_+}(t)
:=\P\big[\exists\pi\in\Vi\mbox{ s.t.\ }\sig_\pi=0,\ \pi(t)=x\big]
\qquad\big(t\geq 0,\ (x,t)\in\Zev\big)
\ee
denote its density, which is a function of $t$ and the speeds
$b_-:=\E[\al^{\rm l}_z]$ and $b_+:=\E[\al^{\rm r}_z]$ of the discrete net
$\Vi$. In what follows, we will need the following fact.
\bl{\bf(Convergence of the density)}\label{L:denscon}
Let $\eps_k$ be positive constants, converging to zero and assume that $-1\leq
b_{k,\,-}\leq b_{k,\,+}\leq 1$ satisfy
\be\label{bk+-}
\eps_k^{-1}b_{k,\,-}\asto{k}\bet_-
\quad\mbox{and}\quad
\eps_k^{-1}b_{k,\,+}\asto{k}\bet_+
\ee
for some $\bet_-\leq\bet_+$. Let $\Psi_b(t)$ be the function in (\ref{Psib})
with $b:=(\bet_+-\bet_-)/2$. Then
\be\label{Psib1}
\lim_{k\to\infty}\sup_{\de\leq t\leq\de^{-1}}\big|
(2\eps_k)^{-1}\Psi_{b_{k,\,-},b_{k,\,+}}(\lfloor\eps_k^{-2}t\rfloor)
-\Psi_b(t)\big|=0\qquad(\de>0)
\ee
and
\be\label{Psib2}
\lim_{\de\down 0}\limsup_{k\to\infty}\int_0^\de\di t\,
(2\eps_k)^{-1}\Psi_{b_{k,\,-},b_{k,\,+}}(\lfloor\eps_k^{-2}t\rfloor)=0.
\ee
\el
{\bf Proof.} Fix $\delta>0$. First we derive a formula for $\Psi_{b_-, b_+}(t)$, defined as in (\ref{disdens}).
For $(x,t)\in\Z^2_{\rm even}$, let $\hat p^{\rm r}$ (resp.\ $\hat p^{\rm l}$) be the path
starting from $(x-1,t)$ (resp.\ $(x+1,t)$) in the dual discrete rightmost (resp.\ leftmost) web
$\hat\Ui^{\rm r}$ (resp.\ $\hat\Ui^{\rm l}$) associated with the discrete dual net $\hat\Vi$.
Then by the discrete analogue of the wedge characterization of the Brownian net in Theorem~\ref{T:net},
the event in the RHS of (\ref{disdens}) occurs if and only if $\hat p^{\rm r}$ and
$\hat p^{\rm l}$ do not intersect on the time interval $[0,t]$. Before $\hat p^{\rm r}$ and $\hat p^{\rm l}$
intersect, the two paths evolve independently, with $\frac{\hat p^{\rm l}(t-\cdot)-\hat p^{\rm r}(t-\cdot)}{2}$ distributed
as a random walk $(D_i)_{i\geq 0}$ with $D_0=1$ and increment distribution
$\P(\Delta D=1)= \gamma_+:=\frac{(1-b_-)(1+b_+)}{4}$, $\P(\Delta D=-1)= \gamma_-:= \frac{(1+b_-)(1-b_+)}{4}$,
and $\P(\Delta D=0) = \gamma_0:=1-\gamma_--\gamma_+$. Therefore
\be\label{Psib-+}
\Psi_{b_-,b_+}(t)=\P^D_1(\tau_0 >t),
\ee
where $\P^D_1(\cdot)$ denotes probability w.r.t.\ $D$ with $D_0=1$, and $\tau_0:=\inf\{i\geq 0: D_i=0\}$.

Let $N_t$ be the number of non-zero increments of $D$ up to time $t$, and let $\bar D$ be a random walk on $\Z$ with $\bar D_0=1$ and
increment distribution $\P(\Delta \bar D=\pm 1)=\bar\gamma_{\pm}:=\frac{\gamma_{\pm}}{\gamma_++\gamma_-}$. Then
\be\label{Psibrep}
\Psi_{b_-,b_+}(t)=\E^D_1[\P^D_1(\tau_0 >t|N_t)] = \E^D_1[\P^{\bar D}_1(\tau_0>N_t)].
\ee

Note that for any $n\in\N$, the law of $(\bar D_i-1)_{0\leq i\leq n}$ is absolutely continuous w.r.t.\ the law of a simple
symmetric random walk $(X_i)_{0\leq i\leq n}$ with $X_0=0$, and the Radon-Nikodym derivative is given by
$(2\bar\gamma_+)^{\frac{n+X_n}{2}}(2\bar\gamma_-)^{\frac{n-X_n}{2}}$. Therefore
\begin{eqnarray}
&& \P^{\bar D}_1(\tau_0>n)=1-\P^{\bar D}_1(\tau_0\leq n) = 1-\P^{\bar D}_1(\bar D_n\leq 0) - \P^{\bar D}_1(\bar D_n\geq 1, \tau_0\leq n) \nonumber \\
&=&\!\!\!\!\!\! \sum_{m=0}^\infty \big(\P^{\bar D}_1(\bar D_n=m+1)-\P^{\bar D}_1(\bar D_n=m+1, \tau_0\leq n)\big) \nonumber \\
&=&\!\!\!\!\!\! \sum_{m=0}^\infty \E^X_0\Big[(2\bar\gamma_+)^{\frac{n+m}{2}}(2\bar\gamma_-)^{\frac{n-m}{2}}\big(1_{\{X_n=m\}}-1_{\{X_n=m,\tau_{-1}\leq n\}}\big)\Big] \nonumber\\
&=&\!\!\!\!\!\! \sum_{m=0}^\infty \E^X_0\Big[(2\bar\gamma_+)^{\frac{n+m}{2}}(2\bar\gamma_-)^{\frac{n-m}{2}}\big(1_{\{X_n=m\}}-1_{\{X_n=-m-2\}}\big)\Big] \nonumber\\
&=&\!\!\!\!\!\! \sum_{m=0}^1 (2\bar\gamma_+)^{\frac{n+m}{2}}(2\bar\gamma_-)^{\frac{n-m}{2}} \P^X_0(X_n=m)\\
&&+ \!\!\! \sum_{m=2}^\infty (2\bar\gamma_+)^{\frac{n+m}{2}}(2\bar\gamma_-)^{\frac{n-m}{2}}\Big(1-\frac{\bar\gamma_-}{\bar\gamma_+}\Big)\P^X_0(X_n=m), \qquad\quad \label{hitprobrep}
\end{eqnarray}
where we applied the reflection principle to $X$ and used $\P^X_0(X_n=-m-2)=\P^X_0(X_n=m+2)$.

We now specialize to the calculation of $\Psi_{b_{k,-}, b_{k,+}}(t_k)$ for $t_k:=\lfloor \eps_k^{-2}t\rfloor$, where $b_{k,-}, b_{k,+}, \eps_k$ satisfy (\ref{bk+-}).
Note that to prove (\ref{Psib2}), it suffices to restrict the integral to $t\in [\eps_k^{3/2}, \delta^{-1}]$. Therefore we assume
$t\in [\eps_k^{3/2}, \delta^{-1}]$ from now on, which implies in particular that $t_k\to\infty$ uniformly in $t$ as $k\to\infty$. By (\ref{Psibrep}), where we replace $D,\bar D, N_{\cdot}$ by $D^k, \bar D^k, N^k_{\cdot}$, we have
\be\label{Psibkdecomp}
\Psi_{b_{k,-}, b_{k,+}}(t_k) = \!\!\!\!\!\!\sum_{|n-t_k/2|\leq t_k^{3/4}}\!\!\!\!\!\! \P^{D^k}_1(N^k_{t_k}=n) \P^{\bar D^k}_1(\tau_0>n) + \!\!\!\!\!\!\sum_{|n-t_k/2|> t_k^{3/4}}\!\!\!\!\!\! \P^{D^k}_1(N^k_{t_k}=n) \P^{\bar D^k}_1(\tau_0>n),
\ee
where for $k$ large, the second term is bounded by
\be\label{Hoeffbd1}
\P^{D^k}_1\big(|N^k_{t_k}-t_k/2|>t_k^{3/4}\big) \leq \P^{D^k}_1\big(|N^k_{t_k}-(\gamma_{k,+}+\gamma_{k,-})t_k|>t_k^{2/3}\big)\leq 2 e^{-\frac{t_k^{4/3}}{2t_k}}=2e^{-\frac{1}{2}t_k^{1/3}},
\ee
where we applied Hoeffding's concentration of measure inequality~\cite{Hoe63} to $N^k_{t_k}$, which is a sum of $t_k$ i.i.d.\ $\{0,1\}$-valued random variables with
mean $\gamma_{k,+}+\gamma_{k,-}=\frac{1}{2}+\frac{\beta_+\beta_-}{2}\eps_k^2(1+o(1))$. Since as $k\to\infty$, $\eps_k^{-1}e^{-\frac{1}{2}t_k^{1/3}}\to 0$ uniformly in $t\in[\eps_k^{3/2}, \delta^{-1}]$, we can safely neglect the second term in (\ref{Psibkdecomp}) when proving (\ref{Psib1})--(\ref{Psib2}).

Note that in the first sum in (\ref{Psibkdecomp}), $n=\eps_k^{-2}t(\frac{1}{2}+o(1))$ uniformly in $n$ and $t\in [\eps_k^{3/2}, \delta^{-1}]$ as $k\to\infty$.
For $n_k:=\eps_k^{-2}t(\frac{1}{2}+o(1))$, we have a representation for $\P^{\bar D^k}_1(\tau_0>n_k)$ as in (\ref{hitprobrep}), where the first sum
in (\ref{hitprobrep}) gives
\begin{eqnarray}
\sum_{m=0}^1 (2\bar\gamma_{k,+})^{\frac{n_k+m}{2}}(2\bar\gamma_{k,-})^{\frac{n_k-m}{2}} \P^X_0(X_{n_k}=m) &=& \Big(\frac{4\gamma_{k,+}\gamma_{k,-}}{(\gamma_{k,+}+\gamma_{k,-})^2}\Big)^{\frac{n_k}{2}}\frac{2}{\sqrt{2\pi n_k}}(1+o(1)) \nonumber\\
&=& \frac{2\eps_k e^{-b^2t(1+o(1))}}{\sqrt{\pi t}} (1+o(1)), \label{Psibloc}
\end{eqnarray}
where we used
$$
\begin{aligned}
4\gamma_{k,\pm} & = & (1\pm b_{k,+})(1\mp b_{k,-}) & =1\pm \eps_k(2b+o(1)), \\
16\gamma_{k,+}\gamma_{k,-}& = & (1-b_{k,+}^2)(1-b_{k,-}^2) & = 1-\eps_k^2(\beta_-^2+\beta_+^2+o(1)), \\
2\gamma_{k,+}+2\gamma_{k,-} & = & 1-b_{k,+}b_{k,-}\qquad & =1-\eps_k^2(\beta_+\beta_-+o(1)),
\end{aligned}
$$
and we applied the local central limit theorem, a strong version of which we need later is
\be\label{LCLT}
\P^X_0(X_s=x) = 1_{\{s+x \mbox{ is even}\}} \frac{2e^{-\frac{x^2}{2s}}}{\sqrt{2\pi s}}(1+o(1))
\ee
uniformly for all $|x|\leq s^{\frac{3}{4}}$ as $s\to\infty$. This can be deduced from~\cite[Theorem 3]{Sto67}.

Analogously, the second term in (\ref{hitprobrep}) gives a contribution to $\P^{\bar D^k}_1(\tau_0>n_k)$ of
\begin{eqnarray}
&& \Big(\frac{4\gamma_{k,+}\gamma_{k,-}}{(\gamma_{k,+}+\gamma_{k,-})^2}\Big)^{\frac{n_k}{2}}\Big(1-\frac{\gamma_{k,-}}{\gamma_{k,+}}\Big)\sum_{m=2}^\infty
\Big(\frac{\gamma_{k,+}}{\gamma_{k,-}}\Big)^{\frac{m}{2}} \P^X_0(X_{n_k}=m)  \nonumber \\
&=& \eps_k (4b+o(1))e^{-b^2t(1+o(1))} \sum_{m=2}^\infty e^{\eps_k m (2b+o(1))} \P^X_0(X_{n_k}=m), \label{Psibint1}
\end{eqnarray}
where we note that the sum is bounded by
$$
\sum_{m\in\Z} e^{3b\eps_k m} \P^X_0(X_{n_k}=m) = \E^X_0\big[e^{3b\eps_k X_{n_k}}\big] = \Big(\frac{e^{3b\eps_k}+e^{-3b\eps_k}}{2}\Big)^{\eps_k^{-2}t(\frac{1}{2}+o(1))}
=O(1)
$$
uniformly in $t\in [\eps_k^{\frac{3}{2}}, \delta^{-1}]$ as $k\to\infty$. Combined with (\ref{Psibkdecomp})--(\ref{Psibloc}), this implies
(\ref{Psib2}).

To prove (\ref{Psib1}), we now restrict to $t\in [\delta, \delta^{-1}]$ and estimate the sum in (\ref{Psibint1}) more preicsely. By Hoeffding's inequality~\cite{Hoe63},
\be\label{Hoeffding}
\P^X_0(X_s\geq m) \leq e^{-\frac{m^2}{2s}}.
\ee
Substituting this bound into (\ref{Psibint1}) then gives
$$
\sum_{m>n_k^{3/4}}^\infty \!\!\!\!\!\! e^{\eps_k m (2b+o(1))} \P^X_0(X_{n_k}=m) \leq \!\!\!\!\!\! \sum_{m>n_k^{3/4}}^\infty \!\!\!\!\! e^{\eps_k m (2b+o(1))-\frac{m^2}{2n_k}}
\leq \!\!\!\!\!\! \sum_{m>n_k^{3/4}}^\infty \!\!\!\!\!\! e^{\eps_k m (2b+o(1))-\delta\eps_k^2m^2(1+o(1))} =o(1)
$$
uniformly in $t\in [\delta, \delta^{-1}]$ as $k\to\infty$. On the other hand, by (\ref{LCLT}),
\begin{eqnarray}
&& \sum_{2\leq m\leq n_k^{3/4}}^\infty \!\!\!\!\!\! e^{\eps_k m (2b+o(1))} \P^X_0(X_{n_k}=m) = (1+o(1))\!\!\!\!\!\! \sum_{2\leq m\leq n_k^{3/4}\atop 2| (m+n_k)}^\infty \!\!\!\!\! \frac{2e^{\eps_k m (2b+o(1))-\frac{m^2}{2n_k}}}{\sqrt{2\pi n_k}} \nonumber \\
&=& (1+o(1))\!\!\!\!\!\! \sum_{2\leq m\leq n_k^{3/4}\atop 2| (m+n_k)}^\infty \!\!\!\!\! \frac{e^{(2b\sqrt{t}+o(1))\frac{\eps_k}{\sqrt t}m -(1+o(1))\big(\frac{\eps_k }{\sqrt t}m\big)^2}}{\sqrt{\pi}} \frac{2\eps_k}{\sqrt{t}} \nonumber \\
&=& \frac{1+o(1)}{\sqrt \pi} \int_0^\infty e^{2b\sqrt{t}x-x^2}{\rm d}x= (1+o(1))e^{b^2t}\int_{-b\sqrt{2t}}^\infty \frac{e^{-\frac{x^2}{2}}}{\sqrt{2\pi}} {\rm d}x
\end{eqnarray}
by Riemann sum approximation. Substituting the last two estimates into (\ref{Psibint1}) and combining with (\ref{Psibkdecomp})--(\ref{Psibloc}) then gives (\ref{Psib1}).
\qed

\bp{\bf(Convergence of relevant separation points)}\label{P:relcon}
Let $\Vi_{\li k\re}$ be a sequence of discrete nets as defined in
(\ref{alr})--(\ref{disnet1}) and let $\Ni$ be a Brownian net with left and
right speeds $\bet_-\leq\bet_+$. Let $-\infty\leq S<U\leq\infty$ and let
$S_k,U_k\in\Z\cup\{-\infty,+\infty\}$ be such that $\eps_k^2S_k\to S$ and
$\eps_k^2U_k\to U$. Let $R_{S,U}$ denote the set of $S,U$-relevant separation
points of $\Ni$ and let $R^{\li k\re}_{S_k,U_k}$ denote the set of
$S_k,U_k$-relevant separation points of $\Vi_{\li k\re}$. Then it is possible to
couple the $\Vi_{\li k\re}$ and $\Ni$ in such a way that
\be\label{ViNi}
S_{\eps_k}(\Vi_{\li k\re})\asto{k}\Ni\quad{\rm a.s.}
\ee
and moreover
\be\label{relcon}
\sum_{z\in R^{\li k\re}_{S_k,U_k}}\de_{S_{\eps_k}(z)}
\Asto{k}
\sum_{z\in R_{S,U}}\de_z\quad{\rm a.s.},
\ee
where $\Rightarrow$ denotes vague convergence of locally finite measures
on~$\R^2$.
\ep
{\bf Remark.} The convergence in (\ref{relcon}) is stronger than the statement
that for each $z\in R_{S,U}$ there exist $z_k\in R^{\li k\re}_{S_k,U_k}$ such
that $S_{\eps_k}(z_k)\to z$. Indeed, since the counting measure on the
right-hand side of (\ref{relcon}) has no double points, such an approximating
sequence is eventually unique, a fact that wil be important in the proof of
Theorem~\ref{T:webnetap} below.\med

\noi
{\bf Proof of Proposition~\ref{P:relcon}.} By Theorem~\ref{T:netcon},
Lemma~\ref{L:netTcon}, Lemma~\ref{L:weakcouple}~(a) and the remarks below it, we
can couple our random variables such that
\be\label{SVV}
S_{\eps_k}\big(\Vi_{\li k\re},\hat\Vi_{\li k\re},
\Vi_{\li k\re}(\Sig_{S_k}),\hat\Vi_{\li k\re}(\Sig_{U_k})\big)
\asto{k}\big(\Ni,\hat\Ni,\Ni(\Sig_S),\hat\Ni(\Sig_U)\big)\quad{\rm a.s.}
\ee
We claim that with this coupling,
for each $z\in R_{S,U}$ there exist $z_k\in R^{\li k\re}_{S_k,U_k}$ with
$S_{\eps_k}(z_k)\to z$. To see this, note that by Lemma~\ref{L:relcros}, each
$z\in R_{S,U}$ is a crossing point of some $\pi\in\Ni$ and $\hat\pi\in\hat\Ni$
with $\sig_\pi=S$ and $U=\hat\sig_{\hat\pi}$. By (\ref{SVV}), there exist
$p_k\in\Vi_{\li k\re}(\Sig_{S_k})$ and $\hat p_k\in\hat\Vi_{\li
  k\re}(\Sig_{U_k})$ such that $p_k\to\pi$ and $\hat p_k\to\hat\pi$. It
follows from the definition of crossing points that for $k$ sufficiently
large, there must exist points $z_k\in\Zev$ such that $S_{\eps_k}(z_k)\to z$
and $p_k$ crosses $\hat p_k$ in $z_k$. In particular, this implies that the
$z_k$ must be $S_k,U_k$-relevant in $\Vi_{\li k\re}$.

We next claim that for each $-\infty<T_-<T_+<\infty$ and
$-\infty<x_-<x_+<+\infty$,
\be\label{denscon}
\E\big[\big|S_{\eps_k}(R^{\li k\re}_{S_k,U_k})\cap A\big|\big]
\asto{k}\E\big[\big|R_{S,U}\cap A\big|\big],
\quad\mbox{where}\quad
A=(x_-,x_+)\times(T_-,T_+).
\ee
To see this, recall that the discrete nets $\Vi_{\li k\re}$ are defined from
i.i.d.\ collections of random variables $(\al^{\li k\re\,{\rm l}}_z,\al^{\li
  k\re\,{\rm r}}_z)_{z\in\Zev}$.  We observe that for all $z=(x,t)\in\Zev$
with $S_k\leq t<U_k$,
\be\ba{l}\label{denscalc}
\dis\P\big[z\mbox{ is $S_k,U_k$-relevant in $\Vi_{\li k\re}$}\big]\\[5pt]
\dis\quad=\P[\al^{\rm l}_z<\al^{\rm r}_z]
\P\big[\exists\pi\in\Vi_{\li k\re}\mbox{ s.t.\ }\sig_\pi=S_k,\ \pi(t)=x\big]
\P\big[\exists\hat\pi\in\hat\Vi_{\li k\re}\mbox{ s.t.\ }
\hat\sig_{\hat\pi}=U_k,\ \hat\pi(t+1)=x\big]\\[5pt]
\dis\quad=\ffrac{1}{2}(b_{k,\,+}-b_{k,\,-})
\Psi_{b_{k,\,-},b_{k,\,+}}(t-S_k)\Psi_{b_{k,\,-},b_{k,\,+}}(U_k-(t+1)),
\ec
where
\be
b_{k,\,-}:=\E[\al^{\li k\re\,{\rm l}}_z]
\quad\mbox{and}\quad
b_{k,\,+}:=\E[\al^{\li k\re\,{\rm r}}_z],
\ee
and $\Psi_{b_-,b_+}(t)$ is the function in (\ref{disdens}). We claim that
(\ref{denscon}) now follows from Proposition~\ref{P:relev}~(b),
Lemma~\ref{L:denscon}, and Riemann sum approximation. Without going through
the details, note that after diffusive rescaling, the per unit density of
points of $S_\eps(\Zev)$ in the plane is $\ffrac{1}{2}\eps_k^{-3}$, and
therefore, by Lemma~\ref{L:denscon}, formula (\ref{denscalc}) says that after
diffusive rescaling, the per unit density of relevant separation points at
time $S<t<U$ is approximately given by
\be\ba{l}
\dis\ffrac{1}{2}\eps_k^{-3}\cdot\ffrac{1}{2}(\eps_k\bet_+-\eps_k\bet_-)
\cdot(2\eps_k)\Psi_b(t-S)\cdot(2\eps_k)\Psi_b(U-t)\\[5pt]
\dis\quad=2b\Psi_b(t-S)\Psi_b(u-t)
\quad\mbox{where}\quad b:=(\bet_+-\bet_-)/2,
\ec
which agrees with (\ref{reldens}).

To prove the existence of a coupling such that (\ref{relcon}) holds,
let
\be
\nu:=\!\sum_{z\in R_{S,U}\cap A}\de_z
\quad\mbox{and}\quad
\nu_k:=\!\!\!\sum_{z\in S_{\eps_k}(R^{\li k\re}_{S_k,U_k})\cap A}\de_z
\ee
be random counting measures with atoms at the positions of the sets in
(\ref{denscon}). By (\ref{denscon}), the laws of the $\nu_k$'s are tight, so
by going to a subsequence if necessary and invoking Skorohod's representation
theorem, we can find a coupling such that in addition to (\ref{SVV}), also
$\nu_k\Rightarrow\nu^\ast$, where $\Rightarrow$ denotes weak convergence and
$\nu^\ast$ is some finite counting measure on the closure $\ov A$ of
$A$. Since for each $z\in R_{S,U}$ there exist $z_k\in R^{\li k\re}_{S_k,U_k}$
such that $S_{\eps_k}(z_k)\to z$, we know that $\nu\leq\nu^\ast$. By
(\ref{denscon}), we see that moreover $\E[\nu(A)]=\E[\nu^\ast(\ov A)]$, so we
conclude that $\nu=\nu^\ast$.

By Lemma~\ref{L:weakcouple}~(a) and the remarks below it, we can find a coupling
such that the measures in (\ref{relcon}) converge weakly on
$(x_{n,-},x_{n,+})\times(T_{n,-},T_{n,+})$ for each $n$, where
$x_{n,-},T_{n,-}\down-\infty$ and $x_{n,+},T_{n,+}\up+\infty$, proving the
vague convergence in (\ref{relcon}).\qed

\subsection{Discrete approximation of a coupled Brownian web and net}\label{S:webnetap}

In this section, we prove a convergence result for discrete webs that are
defined `inside' a discrete net. As a result, we will obtain
Theorem~\ref{T:webinnet}. Our convergence result also prepares for the proof
of Theorem~\ref{T:quench} which will be given in Section~\ref{S:quench}.

For $k\geq 1$, let $(\al^{\li k\re\,{\rm l}},\al^{\li k\re\,{\rm r}})$ be a
collection of $\{-1,+1\}^2$-valued random variables indexed by $\Zev$ as in
(\ref{alr})--(\ref{alr2}) and let $\Vi_{\li k\re}$ and $\Ui^{\rm l}_{\li
  k\re},\Ui^{\rm r}_{\li k\re}$ be the associated discrete net (as defined in
(\ref{disnet1})) and discrete left-right web.

Let $r\in[0,1]$ and, conditional on $(\al^{\li k\re\,{\rm l}},\al^{\li
  k\re\,{\rm r}})$, let $\al^{\li k\re}=(\al^{\li k\re}_z)_{z\in\Zev}$ be a
collection of independent $\{-1,+1\}$-valued random variables such that
$\al^{\li k\re\,{\rm l}}_z\leq\al^{\li k\re}_z\leq\al^{\li k\re\,{\rm r}}_z$
a.s.\ and
\be
\P\big[\al^{\li k\re}_z=\al^{\li k\re\,{\rm r}}_z
\,\big|\,(\al^{\li k\re\,{\rm l}},\al^{\li k\re\,{\rm r}})\big]=r
\qquad(z\in\Zev).
\ee
Then, obviously, under the unconditioned law the collection $\al^{\li
  k\re}=(\al^{\li k\re}_z)_{z\in\Zev}$ is i.i.d.\ with
\be
\eps_k^{-1}\E[\al^{\li k\re\,{\rm l}}_z]\asto{k}\bet:=(1-r)\bet_-+r\bet_+.
\ee
We let $\Ui_{\li k\re}$ denote the discrete web associated with $\al^{\li
  k\re}$. The following theorem implies Theorem~\ref{T:webinnet}.

\bt{\bf(Convergence to a coupled Brownian web and net)}\label{T:webnetap}
Let $\Ui_{\li k\re}$ and $\Vi_{\li k\re}$ be coupled discrete webs and nets as
above. Then
\be\label{webnetcon}
\P\big[S_{\eps_k}(\Ui_{\li k\re},\Vi_{\li k\re})\in\cdot\,\big]
\Asto{k}\P\big[(\Wi,\Ni)\in\cdot\,\big],
\ee
where $\Ni$ is a Brownian net with left and right speeds $\bet_-\leq\bet_+$
and $\Wi$ is a Brownian web with drift $\bet$. Letting $S$ denote the set of
separation points of $\Ni$, one has a.s.:
\begin{itemize}
\item[{\rm(i)}] $\Wi\sub\Ni$ and each point $z\in S$ is of type $(1,2)$ in
  $\Wi$.
\item[{\rm(ii)}] Conditional on $\Ni$, the random variables
  $(\sign_\Wi(z))_{z\in S}$ are i.i.d.\ with $\P[\sign_\Wi(z)=+1\,|\,\Ni]=r$.
\item[{\rm(iii)}] Conditional on $\Wi$, the sets $S_{\rm l}:=\{z\in
  S:\sign_\Wi(z)=-1\}$ and $S_{\rm r}:=\{z\in S:\sign_\Wi(z)=+1\}$ are
  independent Poisson point sets with intensities $(\bet_+-\bet)\ell_{\rm l}$
  and $(\bet-\bet_-)\ell_{\rm r}$, respectively.
\end{itemize}
Moreover,
\be\label{webinnet2}
\Wi=\{\pi\in\Ni:\sign_\pi(z)=\sign_\Wi(z)\ \forall z\in S
\mbox{ s.t.\ $\pi$ enters }z\}.
\ee
In the special case that $r=0$ (resp.\ $r=1$), the Brownian web $\Wi$ is the
left (resp.\ right) Brownian web associated with $\Ni$.
\et
{\bf Proof.} By Theorems~\ref{T:webcon} and
\ref{T:netcon}, the random variables $\Ui_{\li k\re}$ and $\Vi_{\li k\re}$,
diffusively rescaled with $\eps_k$, converge weakly in law to a Brownian web
with drift $\bet$ and Brownian net with left and right speeds $\bet_-,\bet_+$,
respectively. It follows that the laws on the left-hand side of
(\ref{webnetcon}) are tight, so by going to a subsequence if necessary we can
assume that they converge weakly in law to a limit $(\Wi,\Ni)$. We will show
that each such limit point has the properties (i)--(iii) and satisfies
moreover (\ref{webinnet2}). Since property~(ii) and formula (\ref{webinnet2})
determine the joint law of $(\Wi,\Ni)$ uniquely, this then proves the
convergence in (\ref{webnetcon}). Note that if $r=0$ (resp.\ $r=1$), then it
follows from Theorem~\ref{T:netcon} that $\Wi$ is the left (resp.\ right)
Brownian web associated with $\Ni$.

{\em Proof of property~(i)}. The fact that $\Wi\sub\Ni$ is
immediate from the fact that $\Ui_{\li k\re}\sub\Vi_{\li k\re}$ for each
$k$. To prove that each point $z\in S$ is of type $(1,2)$ in $\Wi$, we first
claim that if $l\in\Wl(z')$ and $r\in\Wr(z')$ satisfy $l\sim^{z'}_{\rm out}r$
for some $z'=(x',t')$, then there exists a path $\pi\in\Wi(z')$ such that
$l\leq\pi\leq r$ on $[t',\infty)$. This follows from the fact that by
\cite[Prop.~3.6~(b)]{SS08}, there exist $t_n\down t'$ such that
$l(t_n)<r(t_n)$, while by \cite[Prop.~1.8]{SS08}, any path in $\Ni$ started at
time $t_n$ at a position in $(l(t_n),r(t_n))$ is contained between $l$ and
$r$. Using the compactness of $\Wi$, we find a path $\pi\in\Wi(z')$ with the
desired property. Applying our claim to the one incoming and two outgoing
left-right pairs at a separation point $z$ of $\Ni$, we find that there must
be at least one incoming path and at least two outgoing paths in $\Wi$ at each
such point. By the classification of special points in $\Wi$
(Proposition~\ref{P:classweb}), it follows that $z$ is of type $(1,2)$ in
$\Wi$.

{\em Intermezzo}. Before we turn to the proofs of properties (ii) and (iii),
we first prove some prepatarory results. Since we are assuming that
$S_{\eps_k}(\Ui_{\li k\re},\Vi_{\li k\re})$ converges weakly in law to
$(\Wi,\Ni)$, by Skorohod's representation theorem, we can find a coupling such
that the convergence is a.s. Let $\Ti\sub\R$ be a deterministic countable
dense set of times and for $T\in\Ti$, set
$T_{[k]}:=\lfloor\eps_k^{-2}T\rfloor$. By Proposition~\ref{P:relcon}, Lemma~\ref{L:weakcouple}~(a) and the remarks below it, we can
improve our coupling such that
\be\label{relcon3}
\sum_{z\in R^{\li k\re}_{S_{[k]},U_{[k]}}}\de_{S_{\eps_k}(z)}
\Asto{k}
\sum_{z\in R_{S,U}}\de_z\quad\forall S,U\in\Ti,\ S<U,
\ee
where $\Rightarrow$ denotes vague convergence of locally finite measures on
$\R^2$, and $R^{\li k\re}_{S_{[k]},U_{[k]}}$ and $R_{S,U}$ denote the sets of
$S_{[k]},U_{[k]}$-relevant and $S,U$-relevant separation points of
$\Vi_{\li k\re}$ and $\Ni$, respectively.

It follows from (\ref{relcon3}) that for each $z\in R_{S,U}$, there exist
$z_k\in R^{\li k\re}_{S_{[k]},U_{[k]}}$ such that $S_{\eps_k}(z_k)\to z$.  We
claim that such an approximating sequence is eventually unique. To see this,
assume that $z'_k\in R^{\li k\re}_{S_{[k]},U_{[k]}}$ satisfy
$S_{\eps_k}(z'_k)\to z$. We can choose $\de>0$ such that the ball of radius
$\de$ around $z$ does not contain any other $S,U$-relevant separation points
of $\Ni$ except for $z$. Then (\ref{relcon3}) shows that for $k$ sufficiently
large, there is exactly one $S_{[k]},U_{[k]}$-relevant separation point in the
ball of radius $\de$ around $z$, hence $z_k=z'_k$ for $k$ sufficiently large.

Now let $S,U\in\Ti$, $S<U$, $z\in R_{S,U}$, and let $z_k$ be the eventually
unique sequence of points in $R^{\li k\re}_{S_{[k]},U_{[k]}}$ such that
$S_{\eps_k}(z_k)\to z$. We claim that
\be\label{apsign}
\sign_\Wi(z)=+1\quad\mbox{if and only if}\quad\al^{\li k\re}_{z_k}=+1
\quad\mbox{eventually}.
\ee
(Here, ``$\al^{\li k\re}_{z_k}=+1$ eventually'' means that there exists a $K$
such that $\al^{\li k\re}_{z_k}=+1$ for all $k\geq K$.) It suffices to prove
that $\sign_\Wi(z)=+1$ implies that $\al^{\li k\re}_{z_k}=+1$
eventually. By symmetry, this then also shows that $\sign_\Wi(z)=-1$
implies that $\al^{\li k\re}_{z_k}=-1$ eventually, proving (\ref{apsign}).

If $\sign_\Wi(z)=+1$, then there exist $\pi\in\Wi$ and $\hat l\in\hat\Wl$
such that $\pi$ crosses $\hat l$ in $z=(x,t)$. It follows that there exist
$p_k\in\Ui_{\li k\re}$ and $\hat l_k\in\hat\Ui^{\rm l}_{\li k\re}$ such that
$S_{\eps_k}(p_k)\to\pi$ and $S_{\eps_k}(\hat l_k)\to\hat l$, and points
$z'_k\in\Zev$ with $S_{\eps_k}(z'_k)\to z$ such that $p_k$ crosses $\hat l_k$
in $z'_k$. Let $\sig_\pi<S'<t<U'<\hat\sig_{\hat l}$. Then the $z'_k$ are
$S'_{[k]},U'_{[k]}$-relevant for $k$ large enough, so by the principle of
eventual uniqueness applied to the times $S',U'$ we see that $z'_k=z_k$
eventually. Since $\al^{\li k\re}_{z'_k}=+1$ for each $k$, this proves
(\ref{apsign}).

Next, let $\ti z\in\R^2$ and $\ti z_k\in\Zev$ be deterministic points such
that $S_{\eps_k}(\ti z_k)\to\ti z$, let $p_k$ be the unique element of
$\Ui_{\li k\re}(\ti z_k)$ and let $\pi$ be the a.s.\ unique element of
$\Wi(\ti z)$. Since $S_{\eps_k}(\Ui_{\li k\re},p_k)$ converges weakly in law
to $(\Wi,\pi)$, by Lemma~\ref{L:weakcouple}~(a) and the remarks below it, we can
improve our coupling such that $S_{\eps_k}(p_k)\to\pi$ a.s. Now as before let
$S,U\in\Ti$, $S<U$, $z\in R_{S,U}$, and let $z_k$ be the eventually unique
sequence of points in $R^{\li k\re}_{S_{[k]},U_{[k]}}$ such that
$S_{\eps_k}(z_k)\to z$. We claim that
\be\label{apenter}
\mbox{$\pi$ enters $z$ if and only if $p_k$ enters $z_k$
eventually as $k\to\infty$.}
\ee
Indeed, if $\pi$ does not enter $z$, then $\pi$ dos not enter some open ball
around $z$, so it is clear that for $k$ sufficiently large, $p_k$ does not
enter $z_k$. On the other hand, if $\pi$ enters $z$, then since by
property~(i), $z$ is either of type $(1,2)_{\rm l}$ or of type $(1,2)_{\rm r}$
in $\Wi$, there must be either some $\hat r\in\hat\Wr$ such that $\pi$ crosses
$\hat r$ from right to left or some $\hat l\in\hat\Wl$ such that $\pi$ crosses
$\hat k$ from left to right. By symmetry, it suffices to consider only the
first case. In this case, there must exist $\hat r_k\in\hat\Ui^{\rm r}_{\li
  k\re}$ and $z'_k\in\Zev$ with $S_{\eps_k}(z'_k)\to z$ such that for $k$
sufficiently large, $p_k$ crosses $\hat r_k$ in $z'_k$. By the same argument
as in the proof of (\ref{apsign}) we see that $z'_k=z_k$ eventually, hence
$p_k$ enters $z_k$ eventually.

{\em Proof of Property~(ii)}. Let $\Theta\sub\R$ be a deterministic finite
set, say $\Theta=\{T_1,\ldots,T_m\}$ with $T_1<\cdots<T_m$, and set
\bc
\dis R_\Theta&:=&\dis\big\{z\in\R^2:z\mbox{ is a $T_i,T_{i+1}$-relevant
separation point for some $1\leq i\leq m-1$}\big\},\\[5pt]
\dis R^+_\Theta&:=&\dis\big\{z\in\R_\Theta:\sign_\Wi(z)=+1\big\}.
\ec
We let
\be
\nu_\Theta:=\sum_{z\in R_\Theta}\de_z
\quad\mbox{and}\quad
\nu^+_\Theta:=\sum_{z\in R^+_\Theta}\de_z
\ee
be counting measures with atoms at each point of $R_\Theta$ and $R^+_\Theta$,
respectively. We wish to show that $R^+_\Theta$ is an $r$-thinning of
$R_\Theta$. By formula (\ref{thindef}) and Lemma~\ref{L:pclosunit} of
Appendix~\ref{A:thin}, it suffices to show that (in notation introduced there)
\be
\E\big[(1-f)^{\txt\nu^+_\Theta}\,\big|\,\Ni\big]=(1-rf)^{\txt\nu_\Theta}
\qquad{\rm a.s.}
\ee
for each deterministic $f:\R^2\to[0,1]$ that is continuous and has
compact support. Equivalently, we may show that
\be\label{rtin}
\E\big[(1-f)^{\txt\nu^+_\Theta}g(\Ni)\big]
=\E\big[(1-rf)^{\txt\nu_\Theta}g(\Ni)\big]
\ee
for each $f$ as before and bounded continuous $g:\Ki(\Pi)\to\R$. We now choose
deterministic $T_{k,i}\in\Z$ with $\eps_k^2T_{k,i}\to T_i$ for each $1\leq
i\leq m$, and we set
\bc
\dis R^{\li k\re}
&:=&\dis\big\{z\in\Zev:z\mbox{ is a $T_{k,i},T_{k,i+1}$-relevant separation
point}\\
&&\dis\phantom{\big\{z\in\Zev:}
\:\mbox{of $\Vi_{\li k\re}$ for some $1\leq i\leq m-1$}\big\},\\[5pt]
\dis R^{+\,\li k\re}&:=&\dis\big\{z\in\R^{\li k\re}:\al_z=+1\big\},
\ec
and
\be
\nu^{\li k\re}:=\sum_{z\in R^{\li k\re}}\de_{S_{\eps_k}(z)}
\quad\mbox{and}\quad
\nu^{+\,\li k\re}:=\sum_{z\in R^{+\,\li k\re}}\de_{S_{\eps_k}(z)}.
\ee
By construction, $\nu^{+\,\li k\re}$ is an $r$-thinning of $\nu^{\li k\re}$,
so
\be\label{disrtin}
\E\big[(1-f)^{\txt\nu^{+\,\li k\re}}g(\Vi_{\li k\re})\big]
=\E\big[(1-rf)^{\txt\nu^{\li k\re}}g(\Vi_{\li k\re})\big].
\ee
Now (\ref{rtin}) will follow by taking the limit in (\ref{disrtin}), provided
we show that there exists a coupling such that
\be\label{nuscon}
{\rm(i)}\ \nu^{\li k\re}\Asto{k}\nu_\Theta,
\qquad
{\rm(ii)}\ \nu^{+\,\li k\re}\Asto{k}\nu^+_\Theta,
\ee
where $\Rightarrow$ denotes vague convergence on $\R^2$. The existence of a
coupling such that (\ref{nuscon})~(i) holds follows from
Proposition~\ref{P:relcon}. By (\ref{apsign}), we can improve this coupling
such that also (\ref{nuscon})~(ii) holds. This completes our proof that
$R^+_\Theta$ is an $r$-thinning of $R_\Theta$. Since $\Theta$ is arbitrary and
since each separation point $(x,t)$ is $S,U$-relevant for some $S<t<U$,
Property~(ii) follows.

{\em Proof of Property~(iii)}. Let $\De,\hat\De\sub\R^2$ be deterministic,
finite sets and let $\ell_{\rm l}(\De,\hat\De)$ and $\ell_{\rm
  r}(\De,\hat\De)$ be the restrictions of $\ell_{\rm l}$ and $\ell_{\rm r}$,
respectively, to the set $I:={\rm Img}(\Wi(\De))\cap{\rm
  Img}(\hat\Wi(\hat\De))$.  We note that if $\De_n,\hat\De_n$ are finite sets
increasing to countable limits $\De_\infty$ and $\hat\De_\infty$ that are
dense in $\R^2$, then ${\rm Img}(\Wi(\De_n))\cap{\rm Img}(\hat\Wi(\hat\De_n))$
increases to the set of all points of type $(1,2)$ in $\Wi$, so Property~(iii)
will follow provided we show that for any deterministic finite
$\De,\hat\De\sub\R^2$, conditional on $\Wi$, the sets $S_{\rm l}\cap I$ and
$S_{\rm r}\cap I$ are independent Poisson point sets with intensities
$(\bet_+-\bet)\ell_{\rm l}(\De,\hat\De)$ and $(\bet-\bet_-)\ell_{\rm
  r}(\De,\hat\De)$, respectively.

Equivalently, this says that the set $\{(z,-1):z\in S_{\rm l}\cap
I\}\cup\{(z,+1):z\in S_{\rm r}\cap I\}$ is a Poisson point set on
$\R^2\times\{-1,+1\}$ with intensity $(\bet_+-\bet)\int\ell_{\rm l}(\di
z)\de_{(z,-1)}+(\bet-\bet_-)\int\ell_{\rm r}(\di z)\de_{(z,+1)}$. Thus, by formula
(\ref{Poisdef}) and Lemma~\ref{L:pclosunit} of Appendix~~\ref{A:thin}, it
suffices to show that (in notation introduced there)
\be
\E\big[(1-f)^{\txt\nu_{\rm l}}(1-g)^{\txt\nu_{\rm r}}\,\big|\,\Wi\big]
=\ex{-(\bet_+-\bet)\!\int\! f\,\di\ell_{\rm l}(\De,\hat\De)
-(\bet-\bet_-)\!\int\! g\,\di\ell_{\rm r}(\De,\hat\De)}
\ee
for any deterministic, continuous $f,g:\R^2\to[0,1]$, where
\be
\nu_{\rm l}:=\sum_{z\in S_{\rm l}\cap I}\de_z
\quad\mbox{and}\quad
\nu_{\rm r}:=\sum_{z\in S_{\rm r}\cap I}\de_z.
\ee
Equivalently, we may show that
\be\label{Pois}
\E\big[(1-f)^{\txt\nu_{\rm l}}(1-g)^{\txt\nu_{\rm r}}h(\Wi)\big]
=\E\big[\ex{-(\bet_+-\bet)\!\int\! f\,\di\ell_{\rm l}(\De,\hat\De)
-(\bet-\bet_-)\!\int\! g\,\di\ell_{\rm r}(\De,\hat\De)}h(\Wi)\big],
\ee
for each $f,g$ as before and bounded continuous $h:\Ki(\Pi)\to\R$. Let
$\De_k\sub\Zev$ and $\hat\De_k\sub\Zod$ approximate $\De$ and $\hat\De$ as in
Proposition~\ref{P:ellcon}, let $Z^{\li k\re}_{\rm r}$, $I_k$ and $\ell^{\li k\re}_{\rm r}$ be as defined in (\ref{IkZk}) and let $\ell^{\li k\re}_{\rm l}$ be defined similarly, with $Z^{\li k\re}_{\rm r}$ replaced by
$Z^{\li k\re}_{\rm l}:=\{z\in\Zev:\al^{\li k\re}_z=-1\}$. Let $S_{\li k\re}$
be the set of separation points of $\Vi_{\li k\re}$, and set
\be
\nu^{\li k\re}_{\rm l}:=\!\!\!
\sum_{z\in S_{\li k\re}\cap I_k\cap Z^{\li k\re}_{\rm l}}
\!\!\!\de_{S_{\eps_k}(z)}
\quad\mbox{and}\quad
\nu^{\li k\re}_{\rm r}:=\!\!\!
\sum_{z\in S_{\li k\re}\cap I_k\cap Z^{\li k\re}_{\rm r}}
\!\!\!\de_{S_{\eps_k}(z)}.
\ee
We know that conditional on $\Ui_{\li k\re}$, the sets $S_{\li k\re}\cap
I_k\cap Z^{\li k\re}_{\rm l}$ and $S_{\li k\re}\cap I_k\cap Z^{\li k\re}_{\rm
  r}$ are independent thinnings of the sets $I_k\cap Z^{\li k\re}_{\rm l}$ and
$I_k\cap Z^{\li k\re}_{\rm r}$, with thinning probabilities $b_{k,{\rm l}}$
and $b_{k,{\rm r}}$, respectively, which satisfy
\be\ba{r@{\,}c@{\,}c@{\,}c@{\,}c@{\,}c@{\,}l}
\dis b_{k,{\rm l}}
&:=&\dis\P[\al^{\li k\re\,{\rm r}}_z=+1\,|\,\al^{\li k\re}_z=-1]
&=&\dis\frac{\E[\al^{\li k\re\,{\rm r}}_z-\al^{\li k\re}_z]}
{2\P[\al^{\li k\re}_z=-1]}&\dis\:\symto{k}\:&(\bet_+-\bet)\eps_k\\[15pt]
\dis b_{k,{\rm r}}
&:=&\dis\P[\al^{\li k\re\,{\rm l}}_z=-1\,|\,\al^{\li k\re}_z=+1]
&=&\dis\frac{\E[\al^{\li k\re}_z-\al^{\li k\re\,{\rm l}}_z]}
{2\P[\al^{\li k\re}_z=+1]}&\dis\:\symto{k}\:&(\bet-\bet_-)\eps_k.
\ec
Therefore, by formula (\ref{thindef}) of Appendix~~\ref{A:thin}, we have that
\be\label{Poisap}
\E\big[(1-f)^{\txt\nu^{\li k\re}_{\rm l}}(1-g)^{\txt\nu^{\li k\re}_{\rm r}}
h(S_{\eps_k}(\Ui_{\li k\re})\big]
=\E[(1-b_{k,{\rm l}}f)^{\txt\eps_k^{-1}\ell^{\li k\re}_{\rm l}}
(1-b_{k,{\rm r}}g)^{\txt\eps_k^{-1}\ell^{\li k\re}_{\rm r}}
h(S_{\eps_k}(\Ui_{\li k\re})\big].
\ee
Recall that we are assuming throughout that our random variables are coupled
in such a way that $\Ui_{\li k\re}$ and $\Vi_{\li k\re}$, diffusively rescaled
with $\eps_k$, converge to a Brownian web $\Wi$ and Brownian net $\Ni$,
respectively. By Proposition~\ref{P:ellcon}, Lemma~\ref{L:weakcouple}~(a) and
the remarks below it we can improve our coupling such that moreover
\be
\ell^{\li k\re}_{\rm l}\Asto{k}\ell_{\rm l}(\De,\hat\De)
\quad\mbox{and}\quad
\ell^{\li k\re}_{\rm r}\Asto{k}\ell_{\rm r}(\De,\hat\De).
\ee
Thus, (\ref{Pois}) will follow by taking the limit $k\to\infty$ in
(\ref{Poisap}), provided we show that our coupling can be further improved such
that also
\be\label{nuco}
\nu^{\li k\re}_{\rm l}\Asto{k}\nu_{\rm l}
\quad\mbox{and}\quad
\nu^{\li k\re}_{\rm r}\Asto{k}\nu_{\rm r},
\ee
where $\Rightarrow$ denotes weak convergence of finite measures on
$\R^2$. Since $z\in S_{\rm l}\cap I$ if and only if $z$ is a separation point
of $\Ni$, $\sign_\Wi(z)=-1$, and $z$ is entered by a path
$\pi\in\Wi(\De)$ and a path $\hat\pi\in\hat\Wi(\hat\De)$, formula
(\ref{nuco}) follows from (\ref{apsign}) and (\ref{apenter}).

{\em Proof of formula~(\ref{webinnet2})} Let $\ti\Wi$ be defined by
\be\label{webinnet3}
\ti\Wi:=\{\pi\in\Ni:\sign_\pi(z)=\sign_\Wi(z)\ \forall z\in S
\mbox{ s.t.\ $\pi$ enters }z\}.
\ee
Then $\Wi\sub\ti\Wi$ by the fact that $\Wi\sub\Ni$. To prove the other
inclusion, let $\Ti$ be some deterministic countable dense subset of $\R$. Fix
$\ti\pi\in\ti\Wi$. Choose $\sig_{\ti\pi}<s_n\in\Ti$ with
$s_n\down\sig_{\ti\pi}$. For each $n$, we may choose some $\pi_n\in\Wi$ with
$\sig_{\pi_n}=s_n$ and $\pi_n(s_n)=\ti\pi(s_n)$. Since $\ti\pi\in\Ni$ is an
incoming path at $(\ti\pi(s_n),s_n)$, by Proposition~\ref{P:detspec}~(a), this
point is of type ${\rm(p,p)}$ in $\Ni$. Using this and the finite graph
representation (Proposition~\ref{P:fingraph}), we see that
$\pi_n(t)=\ti\pi(t)$ for each $s_n\leq t\in\Ti$, hence $\pi_n=\ti\pi$ on
$[s_n,\infty)$. It follows that $\pi_n\to\tilde\pi$, and therefore, by the
compactness of $\Wi$, that $\ti\pi\in\Wi$.\qed
\med

\noi
{\bf Proof of Theorem~\ref{T:webinnet}.} Let $\Ui_{\li k\re}$ and $\Vi_{\li
  k\re}$ be as in Theorem~\ref{T:webnetap} and let $\hat\Ui_{\li k\re}$ and
$\hat\Vi_{\li k\re}$ be their associated dual discrete web and net. Then, by
Theorems~\ref{T:webcon} and \ref{T:netcon}
\be
\P\big[S_{\eps_k}(\Ui_{\li k\re},\Vi_{\li k\re},
\hat\Ui_{\li k\re},\hat\Vi_{\li k\re})\in\cdot\,\big]
\Asto{k}\P\big[(\Wi,\Ni,\hat\Wi,\hat\Ni)\in\cdot\,\big],
\ee
where $(\Wi,\Ni)$ are a coupled Brownian web and net as in
Theorem~\ref{T:webnetap} and $\hat\Wi,\hat\Ni$ are the duals of
$\Wi,\Ni$. Since $(-\hat\Ui_{\li k\re},-\hat\Vi_{\li k\re})$ is equally
distributed with $(\Ui_{\li k\re},\Vi_{\li k\re})$, we see that
$(-\hat\Wi,-\hat\Ni)$ is equally distributed with $(\Wi,\Ni)$. Now all
statements in Theorem~\ref{T:webinnet} follow from
Theorem~\ref{T:webnetap}.\qed

\subsection{Switching and hopping in the Brownian web and net}
\label{S:markproof}

In this section, we apply Theorem~\ref{T:webinnet} together with the finite
graph representation developed in Section~\ref{S:fingraph} to prove
Proposition~\ref{P:mark} on switching and hopping inside a Brownian net. We
then apply Theorem~\ref{T:webinnet} and Proposition~\ref{P:mark} to give short
proofs of the marking construction of sticky Brownian webs
(Theorem~\ref{T:webmod}) and the Brownian net (Theorem~\ref{T:marknet}),
Proposition~\ref{P:refchange} on changing the reference web, and the
equivalence of the definitions of a left-right Brownian web given in
Sections~\ref{S:modweb} and \ref{S:netdef}. We also formulate and prove a
result on the construction of sticky Brownian webs inside a Brownian net,
analogous to Theorem~\ref{T:webinnet}.\med

\noi
{\bf Proof of Proposition~\ref{P:mark}} For each set $\De\sub S$, set
\bc
\dis\Ni_\De&=&\big\{\pi\in\Ni:\sign_\pi(z)=\al_z\ \forall z\in S\beh\De
\mbox{ s.t.\ $\pi$ enters }z\big\},\\[5pt]
\dis\Wi_\De&=&\big\{\pi\in\Ni:\sign_\pi(z)=-\al_z\ \forall z\in\De
\mbox{ s.t.\ $\pi$ enters }z\big\}\cap\Ni_\De.
\ec
Since $\Ni_{S'},\Ni_{\De_n},\Wi_{S'},\Wi_{\De_n}$ are contained in the compact
set $\Ni$, it suffices to prove the following statements:
\begin{itemize}
\item[{\rm 1.}] $\Wi_{\De_n}=\switch_{\De_n}(\Wi)$ and $\Ni_{\De_n}={\rm
  hop}_{\De_n}(\Wi)$.
\item[{\rm 2.}] $\Ni_{S'},\Ni_{\De_n},\Wi_{S'},\Wi_{\De_n}$ are closed sets.
\item[{\rm 3.}] $\Ni_{\De_n}\to\Ni_{S'}$ and $\Wi_{\De_n}\to\Wi_{S'}$.
\end{itemize}
1.\ Since $\hop_{\De_n}(\Wi)=\bigcup_{\De'\sub\De_n}{\rm
  switch}_{\De_n}(\Wi)$ and $\Ni_{\De_n}=\bigcup_{\De'\sub\De_n}\Wi_{\De_n}$
it suffices to prove that $\Wi_\De=\switch_\De(\Wi)$ for each finite
$\De\sub S$. By induction, it suffices to prove that $\Wi_{\De\cup\{z\}}={\rm
  switch}_z(\Wi_\De)$ for each finite $\De\sub S$ and $z\in S$. Here, by
induction, we have that $\Wi_\De$ is a subset of $\Ni$ such that for each
$z\in S$ the set $\Wi_\De(z)$ contains exactly two paths, say $\pi_1,\pi_2$,
of which exactly one, say $\pi_1$, is the continuation of a path in the set
$\Wi_{\De,\,{\rm in}}(z)$ of paths in $\Wi_\De$ entering $z$. By definition,
writing $z=(x,t)$, one has
\be
\switch_z(\Wi_\De)
=\big(\Wi_\De\beh\Wi_{\De,\,{\rm in}}(z)\big)
\cup\{\pi^t\cup\pi_2:\pi\in\Wi_{\De,\,{\rm in}}(z)\}.
\ee
By the structure of separation points (Proposition~\ref{P:separ}~(c) and (d))
and the fact that the net is closed under hopping between paths at
intersection times \cite[Prop.~1.4]{SS08}, it follows that each path of the
form $\pi':=\pi^t_1\cup\pi_2$ with $\pi\in\Wi_{\De,\,{\rm in}}(z)$ is an
element of $\Ni$ and satisfies $\sign_{\pi'}(z)=-\sign_{\pi_1}(z)$,
proving that $\switch_z(\Wi_\De)\sub\Wi_{\De\cup\{z\}}$. Conversely, each
$\pi'\in\Wi_{\De\cup\{z\}}$ that enters $z$ is of the form $\pi'=\pi^t\cup\pi_2$
where $\pi:={\pi'}^t\cup\pi_1\in\Wi_{\De,\,{\rm in}}(z)$, showing that
$\switch_z(\Wi_\De)\supset\Wi_{\De\cup\{z\}}$.

2.\ It suffices to prove that if $\pi\in\Ni$ enters some point $z=(x,t)\in S$
and $\pi_n\in\Ni$ satisfy $\pi_n\to\pi$ and $\sign_{\pi_n}(z)=\al$
whenever $\pi_n$ enters $z$, then $\sign_\pi(z)=\al$. By symmetry, it
suffices to treat the case $\al=-1$. We start by noting that there exists an
$N$ such that for each $n\geq N$, the path $\pi_n$ enters $z$. This follows
from the fact that, by Proposition~\ref{P:separ}, there exist dual paths $\hat
l'_z$ and $\hat r'_z$ forming a dual mesh $\hat M(\hat l'_z,\hat r'_z)$, and
each path in $\Ni$ starting in $\hat M(\hat l'_z,\hat r'_z)$ must enter $z$
\cite[Lemma~3.3]{SSS09}. Since $\sign_{\pi_n}(z)=-1$ for all $n\geq N$,
we have that $\pi_n\leq r'_z$ on $[t,\infty)$ for all $n\geq N$, hence the
same holds for $\pi$ and $\sign_\pi(z)=-1$.

3.\ Since $\Ni_{\De_n},\Wi_{\De_n}\sub\Ni$ and $\Ni$ is compact, by
Lemma~\ref{L:Haucomp} in the appendix, the sets $\{\Ni_{\De_n}\}$ and
$\{\Wi_{\De_n}\}$ are precompact, so by going to a subsequence if necessary,
we may assume that $\Ni_{\De_n}\to\Ni^\ast$ and $\Wi_{\De_n}\to\Wi^\ast$ for
some $\Ni^\ast,\Wi^\ast\in\Ki(\Pi)$. We need to show that $\Ni^\ast=\Ni_{S'}$
and $\Wi^\ast=\Wi_{S'}$. We observe that $\Wi_{\De_n},\Ni_{\De_n}\sub\Ni_{S'}$
so $\Wi^\ast,\Ni^\ast\sub\Ni_{S'}$. Set
\be
\ti\Wi_\De:=\big\{\pi\in\Ni:\sign_\pi(z)=-\al_z\ \forall z\in\De
\mbox{ s.t.\ $\pi$ enters }z\big\}\cap\Ni_{S'}.
\ee
Then $\Wi_{\De_n}\sub\ti\Wi_{\De_m}$ for all $n\geq m$, so letting
$n\to\infty$ we see that $\Wi^\ast\sub\ti\Wi_{\De_m}$ for each $m$, hence
$\Wi^\ast\sub\bigcap_m\ti\Wi_{\De_m}=\Wi_{S'}$.

To prove the opposite inclusions, we must show that for each $\pi\in\Ni_{S'}$
there exist $\pi_n\in\Ni_{\De_n}$ such that $\pi_n\to\pi$ and likewise, for
each $\pi\in\Wi_{S'}$ there exist $\pi_n\in\Wi_{\De_n}$ such that
$\pi_n\to\pi$. Let $\Ti$ be some deterministic countable dense subset of $\R$
and let $\pi\in\Ni_{S'}$. Let $T_1,\ldots,T_m\in\Ti$ be such that
$\sig_\pi<T_1<\cdots<T_m$ and let $R_{T_k,T_{k+1}}$ be as in (\ref{RTk}).  We
observe that at each point in $\R^2$ there starts at least one path in
$\Ni_{\De_n}=\bigcup_{\De'\sub\De_n}\switch_{\De'}(\Wi)$. Therefore, for
each $n$ we can find some $\pi_n\in\Ni_{\De_n}$ such that $\sig_{\pi_n}=T_1$
and $\pi_n(T_1)=\pi(T_1)$. Provided $n$ is sufficiently large, we may moreover
choose $\pi_n$ with the property that $\sign_{\pi_n}(z)={\rm
  sign}_\pi(z)$ for each point $z\in\bigcup_{k=1}^{m-1}R_{T_k,T_{k+1}}$ such
that both $\pi_n$ and $\pi$ enter $z$, hence by Corollary~\ref{C:steer}, we
conclude that $\pi_n(T_k)=\pi(T_k)$ for $k=1,\ldots,m$.

Thus, we have shown that for each finite set $T\sub(\sig_\pi,\infty)\cap\Ti$
there exists an $N$ such that for all $n\geq N$ there exists some
$\pi_n\in\Ni_{\De_n}$ with $\pi_n=\pi$ on $T$.  Choosing
$T_m\up(\sig_\pi,\infty)\cap\Ti$, using the compactness of $\Ni$, going to a
subsequence if necessary, we can find $\pi_{n_m}\in\Ni_{\De_{n_m}}$ such that
$\pi_{n_m}\to\pi$ locally uniformly on $(\sig_\pi,\infty)$. Cutting off a
piece of $\pi_{n_m}$ if necessary to make the starting times converge, we have
found $\Ni_{\De_{n_m}}\ni\pi_{n_m}\to\pi$, proving that
$\Ni^\ast\supset\Ni_{S'}$. The proof that $\Wi^\ast\supset\Wi_{S'}$ is
completely analogous.\qed
\med

\noi
{\bf Proof of Theorem~\ref{T:webmod}.} Let $\bet\in\R$ and $c_{\rm l},c_{\rm
  r}\geq 0$. In Theorem~\ref{T:webinnet}, set $\bet_-:=\bet-c_{\rm r}$,
$\bet_+:=\bet+c_{\rm l}$, let $r:=c_{\rm r}/(c_{\rm l}+c_{\rm r})$ if $c_{\rm
  l}+c_{\rm r}>0$, and choose some arbitrary $r\in[0,1]$ otherwise. Then
$\Wi$, defined in (\ref{webinnet}), is a Brownian web with drift $\bet$ and
conditional on $\Wi$, the set $S$ is a Poisson point set with intensity
$c_{\rm l}\ell_{\rm l}+c_{\rm r}\ell_{\rm r}$. In Theorem~\ref{T:webmod}, we
may without loss of generality assume that $\Wi$ and $S$ are constructed in
this way. Then Proposition~\ref{P:mark} tells us that the limit
$\Wi'=\lim_{\De_n\up S}\switch_{\De_n}(\Wi)$ exists, does not depend on the
choice of the $\De_n$, and is given by
\be
\Wi'=\{\pi\in\Ni:\sign_{\hat\pi}(z)=-\al_z\ \forall z\in S
\mbox{ s.t.\ $\hat\pi$ enters }z\}.
\ee
By Theorem~\ref{T:webinnet}, the dual webs $\hat\Wi,\hat\Wi'$ associated with
$\Wi,\Wi'$ are given by
\bc
\dis\hat\Wi&=&\dis\{\hat\pi\in\hat\Ni:\sign_{\hat\pi}(z)=\al_z\ \forall z\in S
\mbox{ s.t.\ $\hat\pi$ enters }z\},\\[5pt]
\dis\hat\Wi'&=&\dis\{\hat\pi\in\hat\Ni:\sign_{\hat\pi}(z)=\al'_z\ \forall z\in S
\mbox{ s.t.\ $\hat\pi$ enters }z\},
\ec
so Proposition~\ref{P:mark} tells us that $\hat\Wi'=\lim_{\De_n\up
  S}\switch_{\De_n}(\hat\Wi)$. Since conditional on $\Ni$, the $(-\al_z)_{z\in
  S}$ are i.i.d.\ with parameter $1-r$, by Theorem~\ref{T:webinnet}, the
Brownian web $\Wi'$ has drift $\bet'=r\bet_-+(1-r)\bet_+=\bet+c_{\rm l}-c_{\rm
  r}$.\qed
\med

\noi
{\bf Proof of Proposition~\ref{P:refchange}~(ii) and (iii).} We continue to
assume that $\Wi$ and $\Wi'$ are defined inside a Brownian net $\Ni$ as in the
proof of Theorem~\ref{T:webmod}. Set $S_{\rm l}:=\{z\in S:\al_z=-1\}$ and
$S_{\rm r}:=\{z\in S:\al_z=+1\}$. Then, by Theorem~\ref{T:webinnet},
conditional on $\Wi$, the set $S_{\rm l}$ is a Poisson point set with
intensity $c_{\rm l}\ell_{\rm l}$ and the set $S_{\rm r}$ is a Poisson point
set with intensity $c_{\rm r}\ell_{\rm r}$, and likewise, conditional on
$\Wi'$, the set $S_{\rm l}$ is Poisson with intensity $c_{\rm r}\ell'_{\rm l}$
and $S_{\rm r}$ is Poisson with intensity $c_{\rm l}\ell'_{\rm r}$. In
particular, this implies that a.s., each point $z\in S_{\rm l}$ is of type
$(1,2)_{\rm l}$ in $\Wi$ and of type $(1,2)_{\rm r}$ in $\Wi'$. Conversely, if
$z\in\R^2$ is of type $(1,2)_{\rm l}$ in $\Wi$ and of type $(1,2)_{\rm r}$ in
$\Wi'$, then $z$ is a separation point of some paths $\pi\in\Wi$ and
$\pi'\in\Wi'$ and therefore, by the definition of separation points of $\Ni$
given before Proposition~\ref{P:separ}, $z\in S$.\qed
\med

\noi
{\bf Proof of Theorem~\ref{T:marknet}.} As in the previous two proofs, without
loss of generality, we assume that $\Wi$ is embedded in a Brownian net $\Ni$
as in (\ref{webinnet}) and that $S_{\rm l}:=\{z\in S:\al_z=-1\}$ and $S_{\rm
  r}:=\{z\in S:\al_z=+1\}$. Then, by Proposition~\ref{P:mark},
(\ref{marknet})~(i) holds and the limits in (\ref{marknet})~(ii) and (iii)
exist and are given by
\bc
\dis\Wl&=&\dis\{\pi\in\Ni:\sign_{\hat\pi}(z)=-1\ \forall z\in S
\mbox{ s.t.\ $\hat\pi$ enters }z\},\\[5pt]
\dis\Wr&=&\dis\{\pi\in\Ni:\sign_{\hat\pi}(z)=+1\ \forall z\in S
\mbox{ s.t.\ $\hat\pi$ enters }z\}.
\ec
By Theorem~\ref{T:webinnet}, $(\Wl,\Wr)$ is the left-right Brownian web
associated with $\Ni$. Since $S_{\rm l}$ and $S_{\rm r}$ are Poisson point
sets with intensities $c_{\rm l}\ell_{\rm l}$ and $c_{\rm r}\ell_{\rm r}$,
respectively, each $z\in S$ is of type $(1,2)$ in $\Wi$ and
$\sign_\Wi(z)=\al_z$, so by construction, conditional on $\Ni$, the random
variables $(\sign_\Wi(z))_{z\in S}$ are i.i.d.\ with
$\P[\sign_\Wi(z)=+1\,|\,\Ni]=r=c_{\rm r}/(c_{\rm l}+c_{\rm r})$.\qed

To prepare for the proof of Proposition~\ref{P:refchange}~(i), we need a lemma.

\bl{\bf(Sticky Brownian webs inside a Brownian net)}\label{L:stickynet}
Let $\Ni^\ast$ be a Brownian net with left and right speeds
$\bet^\ast_-\leq\bet^\ast_+$ and set of separation points
$S^\ast$. Conditional on $\Ni^\ast$, let $(\al_z,\al'_z)_{z\in S}$ be an
i.i.d.\ collection of random variables with values in $\{-1,+1\}^2$. Set
\bc
\dis p_{--}&:=&\dis\P\big[(\al_z,\al'_z)=(-1,-1)\,\big|\,\Ni\big],\\[5pt]
\dis S_{--}&:=&\dis\big\{z\in S^\ast:(\al_z,\al'_z)=(-1,-1)\big\},
\ec
and let $p_{-+},p_{+-},p_{++}$ and $S_{-+},S_{+-},S_{++}$ be defined
analogously. Set
\be\ba{rr@{\,}c@{\,}l}\label{stickdef}
{\rm(i)}&\dis\Wi&:=&\dis
\big\{\pi\in\Ni^\ast:\sign_\pi(z)=\al_z\ \forall z\in S^\ast
\mbox{ s.t.\ $\pi$ enters }z\big\},\\[5pt]
{\rm(ii)}&\dis\Wi'&:=&\dis
\big\{\pi\in\Ni^\ast:\sign_\pi(z)=\al'_z\ \forall z\in S^\ast
\mbox{ s.t.\ $\pi$ enters }z\big\}.
\ec
Then $\Wi$ is a Brownian web with drift
$\bet:=(p_{--}+p_{-+})\bet^\ast_-+(p_{+-}+p_{++})\bet^\ast_+$ and $\Wi'$ is a
Brownian web with drift
$\bet':=(p_{--}+p_{+-})\bet^\ast_-+(p_{-+}+p_{++})\bet^\ast_+$.

Let $\ell$ denote the intersection local time measure between $\Wi$ and its
dual, let $\ell_{\rm l},\ell_{\rm r}$ denote the restrictions of $\ell$ to the
sets of points of type $(1,2)_{\rm l}$ and $(1,2)_{\rm r}$ in $\Wi$,
respectively, and let $\ell',\ell'_{\rm l},\ell'_{\rm r}$ be the same objects
defined for $\Wi'$. Then, conditional on $\Wi$, the sets
\be\label{SSSS}
S_{--},\quad S_{-+},\quad S_{+-},\quad S_{++}
\ee
are independent Poisson point sets with respective intensities
\be\label{pppp}
p_{--}(\bet^\ast_+-\bet^\ast_-)\ell_{\rm l},\quad
p_{-+}(\bet^\ast_+-\bet^\ast_-)\ell_{\rm l},\quad
p_{+-}(\bet^\ast_+-\bet^\ast_-)\ell_{\rm r},\quad
p_{++}(\bet^\ast_+-\bet^\ast_-)\ell_{\rm r},
\ee
while conditional on $\Wi'$, the sets in (\ref{SSSS}) are independent Poisson
point sets with respective intensities
\be\label{ppppa}
p_{--}(\bet^\ast_+-\bet^\ast_-)\ell'_{\rm l},\quad
p_{-+}(\bet^\ast_+-\bet^\ast_-)\ell'_{\rm r},\quad
p_{+-}(\bet^\ast_+-\bet^\ast_-)\ell'_{\rm l},\quad
p_{++}(\bet^\ast_+-\bet^\ast_-)\ell'_{\rm r}.
\ee
Moreover, one has
\be\ba{rr@{\,}c@{\,}l}\label{stickmod}
{\rm(i)}&\dis\Wi'
&=&\dis\lim_{\De_n\up S_{-+}\cup S_{+-}}\switch_{\De_n}(\Wi),\\[5pt]
{\rm(ii)}&\dis\Wi
&=&\dis\lim_{\De_n\up S_{-+}\cup S_{+-}}\switch_{\De_n}(\Wi').
\ec
\el
{\bf Proof.} By Theorem~\ref{T:webinnet}, formulas (\ref{stickdef})~(i) and
(ii) define Brownian webs with drifts as claimed. By Proposition~\ref{P:mark},
the limits in (\ref{stickmod}) exist and coincide with the objects defined in
(\ref{stickdef}). By Theorem~\ref{T:webinnet}, conditional on $\Wi$, the sets
$S_{\rm l}:=S_{--}\cup S_{-+}$ and $S_{\rm r}:=S_{+-}\cup S_{++}$ are
independent Poisson point sets with intensities
$(p_{--}+p_{-+})(\bet^\ast_+-\bet^\ast_-)\ell_{\rm l}$ and
$(p_{+-}+p_{++})(\bet^\ast_+-\bet^\ast_-)\ell_{\rm r}$, respectively. In
particular, this implies that each $z\in S_{\rm l}$ (resp.\ $z\in S_{\rm r}$)
is of type $(1,2)_{\rm l}$ (resp.\ $(1,2)_{\rm r}$) in $\Wi$.

We claim that the \si-fields generated by, on the one hand, $\Wi$ and
$S^\ast$, and, on the other hand, $\Ni^\ast$ and the collection of random
variables $\al=(\al_z)_{z\in S^\ast}$ are identical. To see this, we note that
by Proposition~\ref{P:mark}, $\Ni=\lim_{\De_n\up
  S^\ast}\hop_{\De_n}(\Wi)$. Since moreover $\al_z=\sign_\Wi(z)$ for all $z\in
S^\ast$, this shows that $\Ni^\ast$ and $\al$ are a.s.\ uniquely determined by
$\Wi$ and $S^\ast$. Conversely, since $\Wi$ is given by (\ref{stickdef})~(i)
and $S^\ast$ is the set of separation points of $\Ni^\ast$, we see that $\Wi$
and $S^\ast$ are a.s.\ uniquely determined by $\Ni^\ast$ and $\al$.

Conditional on $\Ni^\ast$ and $\al$, the random variables $(\al'_z)_{z\in
  S^\ast}$ are independent, where $\P[\al'_z=+1\,|\,(\Ni^\ast,\al)]$ equals
$p_{-+}/(p_{--}+p_{-+})$ if $\al_z=-1$ and $p_{++}/(p_{+-}+p_{++})$ if
$\al_z=+1$. It follows that conditional on $\Wi$ and $S^\ast$, the set
$S_{-+}$ is obtained from $S_{--}\cup S_{-+}$ by independent thinning with
probability $p_{-+}/(p_{--}+p_{-+})$ and likewise, the set $S_{++}$ is
obtained from $S_{+-}\cup S_{++}$ by independent thinning with probability
$p_{++}/(p_{+-}+p_{++})$. Since independent thinning splits a Poisson point
set in two independent Poisson point sets, we conclude that conditional on
$\Wi$, the sets in (\ref{SSSS}) are independent Poisson point sets with
intensities given in (\ref{pppp}). By symmetry, an analogue statement holds
for $\Wi'$, i.e., conditional on $\Wi'$, the sets in (\ref{SSSS}) are
independent Poisson point sets with intensities given in (\ref{ppppa}).\qed
\med

\noi
{\bf Proof of Proposition~\ref{P:refchange}~(i).} By symmetry, it suffices to show that $\ell_{\rm l}=\ell'_{\rm l}$. Let $\bet\in\R$ and $c_{\rm l},c_{\rm r}\geq 0$. In Lemma~\ref{L:stickynet}, set $\bet^\ast_-:=\bet-c_{\rm r}$, $\bet^\ast_+:=\bet+c_{\rm l}+1$, and let $p_{--}:=1/(1+c_{\rm l}+c_{\rm r})$, $p_{-+}:=c_{\rm l}/(1+c_{\rm l}+c_{\rm r})$, $p_{+-}:=c_{\rm r}/(1+c_{\rm l}+c_{\rm r})$, and $p_{++}:=0$. Let $\Wi,\Wi'$ be as in (\ref{stickdef}) and set $S:=S_{-+}\cup S_{+-}$. Then conditional on $\Wi$, the set $S$ is a Poisson point set with intensity $c_{\rm l}\ell_{\rm l}+c_{\rm r}\ell_{\rm r}$ and $\Wi'=\lim_{\De_n\up S}\switch_{\De_n}(\Wi)$. Without loss of generality, we may assume that the sticky Brownian webs in Proposition~\ref{P:refchange} are constructed in this way.

It follows from (\ref{stickmod}) that the \si-fields generated by, on the one hand $\Wi$ and $S$, and, on the other hand, $\Wi'$ and $S$ coincide. By (\ref{pppp}) and (\ref{ppppa}), conditional on this \si-field, the set $S_{--}$ is a Poisson point set with intensity $\ell_{\rm l}$ and also a Poisson point set with intensity $\ell'_{\rm l}$, i.e., the conditional law $\P[S_{--}\in\cdot\,|\,\Wi,S]$ is the law of a Poisson point set with intensity $\ell_{\rm l}$ and also the law of a Poisson point set with intensity $\ell'_{\rm l}$. This is possible only if $\ell_{\rm l}=\ell'_{\rm l}$.\qed

The following lemma sometimes comes in handy.

\bl{\bf(Commutativity of switching)}\label{L:switchcom}
Let $\Wi$ be a Brownian web with drift $\bet$, let $\ell$ be the intersection local time measure between $\Wi$ and its dual and let $\ell_{\rm l},\ell_{\rm r}$ denote the restrictions of $\ell$ to the sets of points of type $(1,2)_{\rm l}$ and $(1,2)_{\rm r}$ in $\Wi$, respectively. Let $c_{\rm l},c_{\rm r},c'_{\rm l},c'_{\rm r}\geq 0$ be constants and conditional on $\Wi$, let $S,S'$ be independent Poisson point sets with intensities $c_{\rm l}\ell_{\rm l}+c_{\rm r}\ell_{\rm r}$ and $c'_{\rm l}\ell_{\rm l}+c'_{\rm r}\ell_{\rm r}$, respectively. Then
\be\label{switchcom}
\lim_{\De_n\up S}\switch_{\De_n}
\big(\lim_{\De'_m\up S'}\switch_{\De'_m}(\Wi)\big)
=\lim_{\De''_k\up S\cup S'}\switch_{\De''_k}(\Wi).
\ee
\el
{\bf Proof.} Choose $\bet^\ast_-\leq\bet^\ast_+$ and $p_{--},\ldots,p_{++}$, summing up to one, such that $c_{\rm l}=p_{--}(\bet^\ast_+-\bet^\ast_-)$, $c'_{\rm l}=p_{-+}(\bet^\ast_+-\bet^\ast_-)$, $c_{\rm r}=p_{+-}(\bet^\ast_+-\bet^\ast_-)$, $c'_{\rm l}=p_{++}(\bet^\ast_+-\bet^\ast_-)$, and $\bet=(p_{--}+p_{-+})\bet^\ast_-+(p_{+-}+p_{++})\bet^\ast_+$. Then, without loss of generality, we may assume that $\Wi$ is constructed inside a Brownian net $\Ni^\ast$ as in Lemma~\ref{L:stickynet} and that $S=S_{--}\cup S_{+-}$ and $S'=S_{-+}\cup S_{++}$. Now Proposition~\ref{P:mark} tells us that both sides of (\ref{switchcom}) are well-defined and given by
\be
\Wi''=\big\{\pi\in\Ni^\ast:\sign_\pi(z)\neq\al_z\ \forall z\in S^\ast
\mbox{ s.t.\ $\pi$ enters }z\big\}.
\ee
\qed

The following Lemma has been announced in Section~\ref{S:modweb}.

\bl{\bf(Equivalent definitions of left-right Brownian web)}\label{L:lreq}
A pair of Brownian webs $(\Wl,\Wr)$ is a left-right Brownian web with drifts
$\bet_-,\bet_+$ as defined in Section~\ref{S:netdef} if and only if
$(\Wl,\Wr)$ is a pair of sticky Brownian webs with drifts
$\bet_-,\bet_+$ and coupling parameter $\kappa=0$, as defined in
Section~\ref{S:modweb}.
\el
{\bf Proof.} Let $\Ni$ be a Brownian net with left and right speeds $\bet_-,\bet_+$ and let $S$ be its set of separation points. Then, by Theorem~\ref{T:webinnet}, the left-right Brownian web $(\Wl,\Wr)$ associated with $\Ni$ is given by
\bc\label{WlWr}
\dis\Wl&=&\dis\big\{\pi\in\Ni:\sign_\pi(z)=-1\ \forall z\in S
\mbox{ s.t.\ $\pi$ enters }z\big\},\\[5pt]
\dis\Wr&=&\dis\big\{\pi\in\Ni:\sign_\pi(z)=+1\ \forall z\in S
\mbox{ s.t.\ $\pi$ enters }z\big\}.
\ec
Moreover, by the same theorem, if $\ell$ denotes the intersection local time
measure between $\Wl$ and its dual, and let $\ell_{\rm l}$ and $\ell_{\rm r}$
denote the restrictions of $\ell$ to the sets of points of type $(1,2)_{\rm
  l}$ and $(1,2)_{\rm r}$ in $\Wl$, respectively, then conditional on $\Wl$,
the set $S$ is a Poisson point set with intensity $(\bet_+-\bet_-)\ell_{\rm
  l}$. By Proposition~\ref{P:mark}, it follows that $\Wr=\lim_{\De_n\up
  S}\switch_{\De_n}(\Wl)$, hence $(\Wl,\Wr)$ is a pair of sticky Brownian webs
with drifts $\bet_-,\bet_+$ and coupling parameter $\kappa=0$ as defined in
Section~\ref{S:modweb}.

Conversely, if $\Wl$ is a Brownian web with drift $\bet_-$ and if conditional
on $\Wl$, the set $S$ is a Poisson point set with intensity
$(\bet_+-\bet_-)\ell_{\rm l}$, and $\Wr=\lim_{\De_n\up S}{\rm
  switch}_{\De_n}(\Wl)$, then by Theorem~\ref{T:webinnet}, we may assume
without loss of generality that $\Wl$ is defined inside a Brownian net $\Ni$
such that $S$ is the set of separation points of $\Ni$. Now
Proposition~\ref{P:mark} tells us that $\Wr$ has the representation in
(\ref{WlWr}), hence by Theorem~\ref{T:webinnet} $(\Wl,\Wr)$ is the left-right
Brownian web associated with $\Ni$.\qed

\section{Construction and convergence of Howitt-Warren flows}\label{S:main}

In this section, we prove our main results. We start in Section~\ref{S:quench}
with the proof of Theorem~\ref{T:quench} on the convergence of the quenched
laws on the space of webs. In Section~\ref{S:HWconst}, we then use this to
show that the $n$-point motions of the sample web constructed in
Theorem~\ref{T:HWconst} solve the Howitt-Warren martingale problem, thereby
identifying the stochastic flow of kernels there as a Howitt-Warren flow. Here
we also prove the construction of Howitt-Warren flows inside a Brownian net
(Theorem~\ref{T:HWconst2}) and a result on the exchangeability of the
reference and sample Brownian webs from Theorem~\ref{T:HWconst} if $\nu_{\rm
  l}=\nu_{\rm r}$. In Section~\ref{S:immed}, finally, we harvest some
immediate consequences of our construction, such as scaling
(Proposition~\ref{P:scale}) and the existence of regular versions of
Howitt-Warren flows (Proposition~\ref{P:conpath} and \ref{P:regul}).

\subsection{Convergence of quenched laws}\label{S:quench}

In this section, we prove Theorem~\ref{T:quench}. The measures
$S_{\eps_k}(\Qdis_{\li k\re})$ and $\Q$ from Theorem~\ref{T:quench} are random
probability measures on the Polish space $\Ki(\Pi)$. Therefore, by
\cite[Thm.~3.2.9]{Daw91}, the convergence in (\ref{quencon}) is equivalent
to the convergence of the moment measures of $S_{\eps_k}(\Qdis_{\li k\re})$ to
the moment measures of $\Q$.

We start by describing these moment measures. Let $(\Wi_0,\Mi)$ be a marked
reference web as in Section~\ref{S:flowcons} and conditional on $(\Wi_0,\Mi)$,
let $\Wi_1,\Wi_2,\ldots$ be an i.i.d.\ sequence of sample webs constructed as in
(\ref{sample}). Then the unconditional law
\be\label{momeas}
\P\big[(\Wi_1,\ldots,\Wi_n)\in\cdot\,\big]
\ee
is the $n$-th moment measure of $\Q$. Similarly, for each $k$, conditional on
an i.i.d.\ collection of $[0,1]$-valued random variables
$\om^{\li k\re}=(\om^{\li k\re}_z)_{z\in\Zev}$ with law $\mu_k$ satisfying (\ref{mucon}), let
$\al^{\li k\re\,1},\ldots,\al^{\li k\re\,n}$ be independent collections
$\al^{\li k\re\,i}=(\al^{\li k\re\,i}_z)_{z\in\Zev}$ of $\{-1,+1\}$-valued
random variables with $\P[\al^{\li k\re\,i}_z=+1\,|\,\om^{\li k\re}]=\om^{\li
  k\re}_z$, and let $\Ui^i_{\li k\re}:=\Ui^{\al^{\li k\re\,i}}$ be the
discrete web associated with $\al^{\li k\re\,i}$ as defined in
(\ref{Uial}). Then the averaged law
\be\label{dismom}
\P\big[S_{\eps_k}\big(\Ui^1_{\li k\re},\ldots,\Ui^n_{\li k\re}\big)
\in\cdot\,\big]
\ee
is the $n$-th moment measure of $S_{\eps_k}(\Qdis_{\li k\re})$. We need to prove
weak convergence of the laws in (\ref{dismom}) to those in (\ref{momeas}).

Our strategy will be to embed the Brownian webs $\Wi_1,\ldots,\Wi_n$ in a
Brownian net $\Ni$, and similarly for the rescaled discrete webs. We will then
prove weak convergence in law for the discrete net and webs to
$(\Ni,\Wi_1,\ldots,\Wi_n)$ much in the same way as we have proved
Theorem~\ref{T:webnetap}.

We start by recalling how the sample Brownian webs $\Wi_1,\Wi_2,\ldots$ are constructed in terms of the marked reference Brownian web $(\Wi_0,\Mi)$. The basic ingredients of the construction are the drift $\bet_0$ of the reference web $\Wi_0$ and finite measures $\nu_{\rm l},\nu_{\rm r}$ on $[0,1]$. Given $\Wi_0$, the set of marked points $\Mi=\{(z,\om_z):z\in M\}$ is then a Poisson point set with intensity as in (\ref{Poisint}). To construct $\Wi_1,\Wi_2,\ldots$, conditional on $(\Wi_0,\Mi)$, independently for $i=1,2,\ldots$, we let $(\al^i_z)_{z\in M}$ be a collection of independent $\{-1,+1\}$-valued random variables with $\P[\al^i_z=+1\,|\,(\Wi_0,\Mi)]=\om_z$, we set $A_i:=\{z\in M:\al^i_z\neq\sign_{\Wi_0}(z)\}$, we let $B_i$ be a Poisson point set with intensity $2\nu_{\rm l}(\{0\})\ell_{\rm l}+2\nu_{\rm r}(\{1\})\ell_{\rm r}$, independent of $A_i$, and as in (\ref{sample}), we set
\be\label{sample2}
\Wi_i:=\lim_{\De_n\up A_i\cup B_i}\switch_{\De_n}(\Wi_0)
\qquad(i=1,2,\ldots).
\ee
Then the $\Wi_1,\Wi_2,\ldots$ are conditionally i.i.d.\ given $(\Wi_0,\Mi)$ and $\P[\Wi_i\in\cdot\,|\,(\Wi_0,\Mi)]$ is the Howitt-Warren quenched law with drift $\bet$ and characteristic measure $\nu$ given by (\ref{tib}) and (\ref{nuform}).

We wish to show that for each $n\geq 1$, the Brownian webs $\Wi_1,\ldots,\Wi_n$ from (\ref{sample2}) are in a natural way embedded in a Brownian net. To that aim, for any set of paths $\Ai\sub\Pi$ and set of times $\Ti\sub\R$, we let $\Hi_\Ti(\Ai)$ denote the set of paths that can be obtained from $\Ai$ by hopping finitely often at times in $\Ti$, i.e., $\Hi_\Ti(\Ai)$ contains all paths of the form
\be\ba{l}
\dis\pi=\bigcup_{i=1}^n\big\{(\pi_i(t),t):t_{i-1}\leq t\leq t_i\big\}
\quad\mbox{where}\quad
\pi_1,\ldots,\pi_n\in\Ai,\ t_1,\ldots,t_{n-1}\in\Ti,\\[5pt]
\dis\hspace{2cm} t_0<\cdots<t_n=\infty,\ \sig_{\pi_1}=t_0,
\ \sig_{\pi_{i+1}}\leq t_i,\ \pi_{i+1}(t_i)=\pi_i(t_i)\quad(1\leq i\leq n-1),
\ec
where as usual we identify a path $\pi$ with its graph $\{(\pi(t),t):t\geq\sig_\pi\}$. Moreover, we set
\be\label{Gan}
\Ga_n:=\{-1,+1\}^n\beh\{(-1,\ldots,-1),(+1,\ldots,+1)\}
\ee
and for each $n\geq 2$, we define a probability measure $\La_n$ on $\Ga_n$ by
\be\label{Lan}
\La_n(\vec\ga):=\frac{1}{Z}\int q^{k-1}(1-q)^{n-k-1}\nu(\di q),
\quad\mbox{where}\quad k:=|\{i:\ga^i=+1\}|,
\ee
where $Z$ is the normalization constant given by
\be\label{Zdef}
Z:=\int\frac{1-q^n-(1-q)^n}{q(1-q)}\,\nu(\di q),
\ee
with the convention that the integrand in (\ref{Zdef}) takes on the value $n$ at the points $q=0,1$.

Below, if $\Ai$ is a set of paths, then $\ov\Ai$ denotes the closure of $\Ai$ in the topology on the path space $\Pi$. We note that in (\ref{Nin}), if $\Di$ is moreover dense in $\R^2$, then $\ov{\Ni_n(\Di)}=\Ni_n$. (This follows, for example, from \cite[Thm.~1.3]{SS08}.)

\index{moment measures}
\bl{\bf(Construction of moment measures)}\label{L:momcon}
Conditional on the marked reference Brownian web $(\Wi_0,\Mi)$, let
$\Wi_1,\Wi_2,\ldots$ be an i.i.d.\ sequence of sample Brownian webs defined as
in (\ref{sample2}). Then, for each $n\geq 1$, there exists an a.s.\ unique
Brownian net $\Ni_n$ with left and right speeds $\bet_-(n),\bet_+(n)$, defined
in (\ref{betminm}) and (\ref{betplusm}) with $\bet$ and $\nu$ given by
(\ref{tib}) and (\ref{nuform}), such that for any deterministic countable set $\Di\sub\R^2$ and countable dense set of times $\Ti\sub\R$,
\be\label{Nin}
\ov{\Ni_n(\Di)}=\ov{\Hi_\Ti(\Wi_1\cup\cdots\cup\Wi_n)(\Di)}\quad{\rm a.s.}
\ee
Let $S_n$ be the set of separation points of $\Ni_n$. Then each $z\in S_n$ is of type $(1,2)$ in $\Wi_1,\ldots,\Wi_n$ and conditional on $\Ni_n$, the random variables $(\vec\al_z)_{z\in S_n}$ defined by
\be\label{alzdef}
\vec\al_z:=\big(\sign_{\Wi_1}(z),\ldots,\sign_{\Wi_n}(z)\big)
\ee
are i.i.d.\ with law $\La_n$ defined in (\ref{Lan}). Moreover, one has
\be\label{WinN}
\Wi_i=\{\pi\in\Ni_n:\sign_\pi(z)=\al^i_z\ \forall z\in S_n
\mbox{ s.t.\ $\pi$ enters }z\}
\qquad(i=1,\ldots,n).
\ee
\el
{\bf Proof.} Set $C:=\bigcup_{i=1}^n(A_i\cup B_i)$ and let $C_{\rm l}$ and $C_{\rm r}$ denote the restrictions of $C$ to the sets of points of type $(1,2)_{\rm l}$ and $(1,2)_{\rm r}$ in $\Wi_0$, respectively. For each $\vec\ga\in\{-1,+1\}^n$, set
\be
C(\vec\ga)
:=\big\{z\in C:\sign_{\Wi_i}(z)=\ga_i\ \forall i=1,\ldots,n\big\}
\ee
and define $C_{\rm l}(\vec\ga)$ and $C_{\rm r}(\vec\ga)$ similarly, with $C$ replaced by $C_{\rm l}$ resp.\ $C_{\rm r}$. By our definition of the sample Brownian webs, conditional on $\Wi_0$, the sets $\{C(\vec\ga):\vec\ga\in\{-1,+1\}^n\}$ are independent Poisson point sets with intensity $c_{\rm l}(\vec\ga)\ell_{\rm l}+c_{\rm r}(\vec\ga)\ell_{\rm r}$, where (compare (\ref{ABint}))
\be\ba{l}\label{clcr}
\dis c_{\rm l}(\vec\ga)
=2\,1_{\{0<k\}}\!\int_{(0,1]}q^k(1-q)^{n-k}q^{-1}\nu_{\rm l}(\di q)
+2\,1_{\{k=1\}}\nu_{\rm l}(\{0\})\\[5pt]
\left.\ba{@{\!}l}
\dis\phantom{c_{\rm r}(\vec\ga)}
=2\,1_{\{0<k\}}\int q^{k-1}(1-q)^{n-k}\nu_{\rm l}(\di q),\\[5pt]
\dis c_{\rm r}(\vec\ga)
=2\,1_{\{k<n\}}\int q^k(1-q)^{n-k-1}\nu_{\rm r}(\di q)
\ea\quad\right\}\quad\mbox{where }k:=|\{i:\ga_i=+1\}|.
\ec
We modify our reference web by setting
\be\label{newref}
\Wi'_0:=\lim_{\De_m\up C^\ast}
\switch_{\De_m}(\Wi_0)
\quad\mbox{where}\quad
C^\ast:=C_{\rm l}(+1,\ldots,+1)\cup C_{\rm r}(-1,\ldots,-1).
\ee
By Proposition~\ref{P:refchange}~(i), $\ell_{\rm l}$ and $\ell_{\rm r}$ are also the intersection local time measures for the modified reference web $\Wi'_0$. Since $\Wi'_0$ is a.s.\ uniquely determined by $\Wi_0$ and the set $C^\ast$, and since conditional on $\Wi_0$ and $C^\ast$, the sets $C(\vec\ga)$ with $\ga\neq(-1,\ldots,-1),(+1,\ldots,+1)$ are independent Poisson point sets with intensity $c_{\rm l}(\vec\ga)\ell_{\rm l}+c_{\rm r}(\vec\ga)\ell_{\rm r}$, by Theorem~\ref{T:marknet}, we can define a Brownian net $\Ni_n$ with set of separation points $S_n$ by (recall (\ref{Gan}))
\be\label{Nndef}
\Ni_n:=\lim_{\De_m\up S_n}\hop_{\De_m}(\Wi'_0)
\quad\mbox{where}\quad
S_n:=\bigcup_{\vec\ga\in\Ga_n}C(\vec\ga).
\ee
By Theorem~\ref{T:marknet}, conditional on $\Ni_n$, the random variables $(\sign_{\Wi'_0}(z))_{z\in S_n}$ are i.i.d.\ with
\be
\P[\sign_{\Wi'_0}(z)=+1\,|\,\Ni_n]
=\frac{\sum_{\vec\ga\in\Ga_n}c_{\rm r}(\vec\ga)}
{\sum_{\vec\ga\in\Ga_n}(c_{\rm l}(\vec\ga)+c_{\rm r}(\vec\ga))}.
\ee
Using this and the independence of the Poisson point sets $\{C(\vec\ga):\vec\ga\in\Ga_n\}$, it is straightforward to check from (\ref{nuform}) and (\ref{clcr}) that conditional on $\Ni_n$, the random variables in (\ref{alzdef}) are i.i.d.\ with law $\La_n$ defined in (\ref{Lan}). By Lemma~\ref{L:switchcom},
\be
\Wi_i=\lim_{\De_m\up C^i}\switch_{\De_m}(\Wi'_0)
\quad\mbox{where}\quad
C^i:=\!\!\!\bigcup\subb{\vec\ga\in\Ga_n}{\ga_i=+1}\!\!\!\!C_{\rm l}(\vec\ga)
\;\;\cup\;
\!\!\!\bigcup\subb{\vec\ga\in\Ga_n}{\ga_i=-1}\!\!\!\!C_{\rm r}(\vec\ga).
\ee
Therefore, by Proposition~\ref{P:mark}, we see that (\ref{WinN}) holds.

The speed of $\Wi'_0$ is given by
\be
\bet_0+c_{\rm l}(+1,\ldots,+1)-c_{\rm r}(-1,\ldots,-1),
\ee
and therefore the left speed of $\Ni_n$ is given by
\be\ba{l}\label{betncalc}
\dis\bet_0+c_{\rm l}(+1,\ldots,+1)-c_{\rm r}(-1,\ldots,-1)
-\sum_{\vec\ga\in\Ga_n}c_{\rm r}(\vec\ga)\\[5pt]
\dis\quad=\bet-2\nu_{\rm l}([0,1])+2\nu_{\rm r}([0,1])
+c_{\rm l}(+1,\ldots,+1)
-\sum_{\vec\ga\neq(+1,\ldots,+1)}c_{\rm r}(\vec\ga)\\[5pt]
\dis\quad=\bet-2\int\nu_{\rm l}(\di q)(1-q^{n-1})
+2\int\nu_{\rm r}(\di q)
\Big(1-\sum_{k=0}^{n-1}\left({n\atop k}\right)q^k(1-q)^{n-k-1}\Big)\\[5pt]
\dis\quad=\bet-2\int(1-q)\nu_{\rm l}(\di q)\sum_{k=0}^{n-2}q^k
-2\int q\nu_{\rm r}(\di q)\sum_{k=0}^{n-2}q^k
=\bet-2\int\nu(\di q)\sum_{k=0}^{n-2}q^k=\bet_-(n),
\ec
where we have used (\ref{tib}), (\ref{nuform}), (\ref{betplusm}) and the fact that
\be
\sum_{k=0}^{n-1}\left({n\atop k}\right)q^k(1-q)^{n-k-1}
=(1-q)^{-1}(1-q^n)=\sum_{k=0}^{n-1}q^k,
\ee
which is true even for $q=1$, even though the intermediate step is not defined
in this case. The calculation for $\bet_+(n)$ is completely analogous.

We are left with the task to prove (\ref{Nin}). The inclusion $\Ni_n(\Di)\supset\Hi_\Ti(\Wi_1\cup\cdots\cup\Wi_n)(\Di)$ follows from the fact that $\Ni_n$ is closed under hopping at deterministic times, see~\cite[Lemma~8.3]{SS08}. To prove the converse inclusion, by the compactness of $\Ni_n$, it suffices to prove that for each $t_1<\cdots<t_m$ with $t_1,\ldots,t_m\in\Ti$, $z=(x,t_0)\in\Di$ with $t_0<t_1$, and $\pi\in\Ni_n(z)$, we can find $\pi'\in\Hi_\Ti(\Wi_1\cup\cdots\cup\Wi_n)$ starting from $z$ such that $\pi(t_i)=\pi'(t_i)$ for $i=1,\ldots,m$. By the finite graph representation (in particular, by Corollary~\ref{C:steer}) and the fact that for each separation point $z$ of $\Ni_n$, there exists $1\leq i,j\leq n$ such that $\sign_{\Wi_i}(z)=-1$ and $\sign_{\Wi_j}(z)=+1$, we can find a $\pi''$ starting at $z$ and satisfying $\pi(t_i)=\pi''(t_i)$ for $i=1,\ldots,m$ that is obtained by concatenating finitey many paths in $\Wi_1,\ldots,\Wi_n$ at separation points of $\Ni_n$. By the fact that $\Ti$ is dense and the structure of separation points (see Proposition~\ref{P:separ}), we can modify $\pi''$ a bit such that the concatenation takes place at times in $\Ti$.\qed

\noi
{\bf Proof of Theorem~\ref{T:quench}} Let $\Ui^i_{\li k\re}$ $(i=1,\ldots,n)$
be the discrete webs in (\ref{dismom}), and let $\Vi_{\li k\re}$ be the
discrete net defined by
\be\label{disnet2}
\Vi_{\li k\re}:=\big\{p: p(t+1)-p(t)
\in\{\al^{\li k\re\,1}_{(p(t),t)},\ldots,\al^{\li k\re\,n}_{(p(t),t)}\}
\ \forall\, t\geq\sig_p\big\}.
\ee
By Theorem~\ref{T:netcon}, $\Vi_{\li k\re}$, diffusively rescaled,
converges to a Brownian net with left and right speeds given by
\bc
\dis\bet_-(n)&=&\dis\lim_{k\to\infty}\eps_k^{-1}
\E\big[\al^{\li k\re\,1}_z\wedge\cdots\wedge\al^{\li k\re\,n}_z\big]
=\lim_{k\to\infty}\eps_k^{-1}\int\mu_k(\di q)\big(q^n-(1-q^n)\big)\\[5pt]
&=&\dis\lim_{k\to\infty}\eps_k^{-1}\int\mu_k(\di q)
\Big((2q-1)-2q(1-q)\sum_{k=0}^{n-2}q^k\Big)
=\bet-2\int\nu(\di q)\sum_{k=0}^{n-2}q^k,\\[5pt]
\bet_+(n)&=&\dis\lim_{k\to\infty}\eps_k^{-1}
\E\big[\al^{\li k\re\,1}_z\vee\cdots\vee\al^{\li k\re\,n}_z\big]
=\bet+2\int\nu(\di q)\sum_{k=0}^{m-2}(1-q)^k,
\ec
where we have used (\ref{mucon}). Let $S_{\li k\re}$ be the set of separation
points of $\Vi_{\li k\re}$ and set $\vec\al^{\li k\re}_z:=(\al^{\li
  k\re\,1}_z,\ldots,\al^{\li k\re\,1}_z)$. Then, conditional on $\Vi_{\li
  k\re}$, the random variables $(\vec\al^{\li k\re}_z)_{z\in S_{\li k\re}}$
are i.i.d.\ with
\be
\P\big[\vec\al^{\li k\re}_z=\vec\ga\,\big|\,\Vi_{\li k\re}\big]
=\frac{1}{Z_k}\int\mu_k(\di q)q^l(1-q)^{n-l},
\quad\mbox{where}\quad l:=|\{i:\ga^i=+1\}|,
\ee
and $Z_k$ is a normalization constant. Using (\ref{mucon}), it is easy to check that this conditional law converges as $k\to\infty$ to the law in (\ref{Lan}).

For each $i=1,\ldots,n$, the pairs $(\Vi_{\li k\re},\Ui^i_{\li k\re})$ are distributed as the discrete nets and webs in Theorem~\ref{T:webnetap}, so by that theorem,
and going to a subsequence if necessary,
we can couple our random variables in such a way that
\be
S_{\eps_k}\big(\Vi_{\li k\re},\Ui^1_{\li k\re},\ldots,\Ui^n_{\li k\re}\big)
\asto{n}(\Ni,\Wi_1,\ldots,\Wi_n),
\ee
where $\Ni$ is a Brownian net with left and right speeds $\bet_-(n),\bet_+(n)$, $\Wi_1,\ldots,\Wi_n$ are Brownian webs with drift $\bet$, such that each separation point of $\Ni$ is of type $(1,2)$ in each $\Wi_i$ and
\be
\Wi_i:=\{\pi\in\Ni:\sign_\pi(z)=\sign_{\Wi_i}(z)\ \forall z\in S
\mbox{ s.t.\ $\pi$ enters }z\}\qquad(i=1,\ldots,n),
\ee
where $S$ is the set of separation points of $\Ni$. Much in the same way as in the proof of property~(ii) of Theorem~\ref{T:webnetap}, we find that conditional on $\Ni$, the random variables
\be
\big(\sign_{\Wi_1}(z),\ldots,\sign_{\Wi_n}(z)\big)_{z\in S}
\ee
are i.i.d.\ with common law as in (\ref{Lan}). By Lemma~\ref{L:momcon}, this
proves the convergence of the moment measures in (\ref{dismom}) to those in
(\ref{momeas}) and hence, by \cite[Thm.~3.2.9]{Daw91}, the convergence in
(\ref{quencon}).\qed

\subsection{Proof of the marking constructions of Howitt-Warren
flows}\label{S:HWconst}

In this section, we prove our main results, Theorems~\ref{T:HWconst} and
\ref{T:HWconst2} on the construction of Howitt-Warren flows inside a Brownian
web and net. It turns out that we already have most ingredients of the
proofs. The main point that still needs to be settled is to verify that our
construction agrees with the original definition of Howitt-Warren flows
based on $n$-point motions and the Howitt-Warren martingale problem.

\bp{\bf(Identification of $n$-point motions)}\label{P:niden}
Let $\bet_0\in\R$ and let $\nu_{\rm l},\nu_{\rm r}$ be finite measures on
$[0,1]$. Let $(\Wi_0,\Mi)$ be a marked reference web as in
Theorem~\ref{T:HWconst} and conditional on $(\Wi_0,\Mi)$, let
$\Wi_1,\ldots,\Wi_n$ be $n$ independent sample webs constructed as in
(\ref{sample}). For each deterministic $z\in\R^2$, let $\pi^i_z$ denote the
a.s.\ unique element of $\Wi_i(z)$. Then, for each $\vec x\in\R^n$ and
$s\in\R$, the process
\be
\big(\pi^1_{(x_1,s)}(s+t),\ldots,\pi^n_{(x_n,s)}(s+t)\big)_{t\geq 0}
\ee
solves the Howitt-Warren martingale problem with drift $\bet$ and
characteristic measure $\nu$ given by (\ref{tib}) and (\ref{nuform}).
\ep
{\bf Proof.} Instead of attempting a direct proof we will use discrete
approximation. It is easy to verify that Theorem~\ref{T:quench} implies the
convergence of the $n$-point motions of diffusively rescaled discrete
Howitt-Warren flows to the $n$-point motions of the quenched law $\Q$,
while by Proposition~\ref{P:nconv}, the same discrete $n$-point motions
converge to a solution of the Howitt-Warren martingale problem. The
proposition then follows.\qed

\detail{To formulate this somewhat more precisely, let $\eps_k\to 0$, let
$\mu_k$ be a sequence of probability laws on $[0,1]$ satisfying (\ref{mucon}),
let $\Ui^1_{\li k\re},\ldots,\Ui^n_{\li k\re}$ be conditionally independent
discrete webs as in (\ref{dismom}). Theorem~\ref{T:quench} tells us that
\be
\P\big[S_{\eps_k}
\big(\Ui^1_{\li k\re},\ldots,\Ui^n_{\li k\re}\big)\in\cdot\,\big]
\Asto{k}\P\big[(\Wi_1,\ldots,\Wi_n)\in\cdot\,\big],
\ee
where $\Rightarrow$ denotes weak convergence of probability measures on
$\Ki(\Pi)^n$. By Skorohod's representation theorem, we can find a coupling such
that the convergence is a.s. Let $p^{\li k\re\,i}_z$ denote the unique paths
in $\Ui^i_{\li k\re}$ started from a point $z\in\Zev$. Choose deterministic
$\vec x^{\langle k\rangle}\in\Z_{\rm even}^2$ such that $\eps_k\vec x^{\li
  k\re}\to\vec x$ for some $\vec x\in\R^n$. By Lemmas~\ref{L:Haucomp} and
\ref{L:Hauconv}~(i) in the appendix, each subsequence of
\be\label{p1n}
S_{\eps_k}\big(p^{\li k\re\,1}_{(x^{\li k\re}_1,0)},\ldots,
p^{\li k\re\,n}_{(x^{\li k\re}_n,0)}\big)
\ee
contains a further subsequence that converges, and the limit must be a
collection of paths in $\Wi_1,\ldots,\Wi_n$ starting from
$(x_1,0),\ldots,(x_n,0)$. Since in a Brownian web, at each deterministic point
there starts only one path, we conclude that the paths in (\ref{p1n}) converge
a.s.\ to the corresponding paths in $\Wi_1,\ldots,\Wi_n$. We conclude that
\be\label{npath}
\P\big[S_{\eps_k}\big(p^{\li k\re\,1}_{(x^{\li k\re}_1,0)},\ldots,
p^{\li k\re\,n}_{(x^{\li k\re}_n,0)}\big)\in\cdot\,\big]
\Asto{k}\P\big[(\pi^1_{(x_1,0)},\ldots,\pi^n_{(x_1,0)})\in\cdot\,\big]
\ee
where $\Rightarrow$ denotes weak convergence of probability measures on
$\Pi^n$ and $\pi^i_z$ is the a.s.\ unique element of $\Wi_i(z)$.  By
Proposition~\ref{P:nconv}, the left-hand side of (\ref{npath}) converges, in
the same topology, to the unique solution of the Howitt-Warren martingale
problem with drift $\bet$ and characteristic measure $\nu$. This proves our
claim for $s=0$. The general case follows by translation invariance.}

\noi
{\bf Proof of Theorem~\ref{T:HWconst}} We start by checking that the random
kernels $K^+_{s,t}$ defined as in (\ref{HWconst}) form a stochastic flow of
kernels on $\R$ as in Definition~\ref{D:stochflow}. Indeed, Property~(i)
follows from the fact that
\be\ba{l}
\dis\int_\R K^+_{s,t}(x,\di y)K^+_{t,u}(y,\di z)
=\int_\R\P[\pi^+_{(x,s)}(t)\in\di y\,|\,(\Wi_0,\Mi)]
\,\P[\pi^+_{(y,t)}(u)\in\di z\,|\,(\Wi_0,\Mi)]\\[5pt]
\dis\quad=\int_\R\P[\pi^+_{(x,s)}(t)\in\di y\,|\,(\Wi_0,\Mi)]
\,\P[\pi^+_{(y,t)}(u)\in\di z\,|\,(\Wi_0,\Mi),\ \pi^+_{(x,s)}(t)=y]\\[5pt]
\dis\quad= \P[\pi^+_{(\pi^+_{(x,s)}(t),t)}(u)\in\di z\,|\,(\Wi_0,\Mi)]
= \P[\pi^+_{(x,s)}(u)\in\di z\,|\,(\Wi_0,\Mi)]\quad{\rm a.s.},
\ec
where we have used that $\pi^+_{(x,s)}(t)$ and $\pi^+_{(y,t)}(u)$ are
conditionally independent given $(\Wi_0,\Mi)$ and
$\pi^+_{(\pi^+_{(x,s)}(t),t)}(u)=\pi^+_{(x,s)}(u)$ a.s., which follows from
the fact that for deterministic $t$, a.s.\ every point in $\R\times\{t\}$ is
of type $(0,1),(0,2)$ or $(1,1)$ (see
Proposition~\ref{P:classweb}). Property~(ii) of Definition~\ref{D:stochflow}
follows from the fact that the restrictions of $(\Wi,\Mi,\Wi_0)$ to disjoint
time intervals are independent, which follows from the analogue property for a
single Brownian web which is proved by discrete approximation. Property~(iii),
finally, is obvious from the translation invariance of our definitions. Since
$K^\up_{s,t}(x,\,\cdot\,)=K^+_{s,t}(x,\,\cdot\,)$ a.s.\ for deterministic
$s\leq t$ and $x\in\R$, the same conclusions can be drawn for $K^\up_{s,t}$.

To identify $(K^+_{s,t})_{s\leq t}$ (and likewise $(K^\up_{s,t})_{s\leq t}$)
as a Howitt-Warren flow with drift $\bet$ and characteristic measure $\nu$,
therefore, it suffices to check that for each deterministic $\vec x\in\R^n$
and $s\leq t$, one has (compare (\ref{npoint}))
\be
\E\big[K^+_{s,t}(x_1,\cdot\,)\cdots K^+_{s,t}(x_n,\cdot\,)\big]
=\P[\vec X^{\vec x}_{t-s}\in\cdot\,],
\ee
where $\vec X^{\vec x}$ is a solution the Howitt-Warren martingale problem
with drift $\bet$ and characteristic measure $\nu$, started in $\vec X^{\vec
  x}_0=\vec x$. Since
\be
\E\big[K^+_{s,t}(x_1,A_1\,)\cdots K^+_{s,t}(x_n,A_n)\big]
=\P\big[\pi^1_{(x_1,s)}(t)\in A_1,\ldots,\pi^n_{(x_n,s)}(t)\in A_n\big],
\ee
where $\pi^1_{(x_1,s)},\ldots,\pi^n_{(x_1,s)}$ are as in
Proposition~\ref{P:niden}, our claim follows from that result. The fact that
$(\Wi_0,\Mi,\Wi)$ and $(\Wi,\Mi,\Wi_0)$ are equally distributed if $\nu_{\rm
  l}=\nu_{\rm r}$ follows from the somewhat stronger
Proposition~\ref{P:exchange} below.\qed

\bp{\bf(Exchangeability of reference web)}\label{P:exchange}
Let $(\Wi_0,\Mi)$ be a marked reference web as in Theorem~\ref{T:HWconst} and conditional on $(\Wi_0,\Mi)$, let $(\Wi_1,\Wi_2,\ldots)$ be independent sample webs constructed as in (\ref{sample}). Assume that $\nu_{\rm l}=\nu_{\rm r}$. Then the sequence of Brownian webs $(\Wi_0,\Wi_1,\Wi_2,\ldots)$ is exchangeable.
\ep
{\bf Proof.} In the set-up of Lemma~\ref{L:momcon}, we will show that if $\nu_{\rm l}=\nu_{\rm r}$, then the joint law of $(\Wi_0,\ldots,\Wi_n)$ is equal to the law of $(\Wi_1,\ldots,\Wi_{n+1})$, which is clearly exchangeable. To see this, let $C$ be defined as in the proof of Lemma~\ref{L:momcon} and in analogy with (\ref{Nndef}), set
\be
\Ni'_n:=\lim_{\De_m\up C}\hop_{\De_m}(\Wi_0).
\ee
Then $\Ni'_n$ is a Brownian net with set of separation points $C$. For $i=0,\ldots,n$ and $z\in C$, set $\al^i_z:=\sign_{\Wi_i}(z)$ and as in (\ref{alzdef}) let $\vec\al_z:=(\al^1_z,\ldots,\al^n_z)$. In a similar way as in the proof of Lemma~\ref{L:momcon}, we check that conditional on $\Ni'_n$, the random variables $(\al^0_z,\vec\al_z)_{z\in C}$ are i.i.d.\ with
\[\ba{r@{\,}c@{\,}ll}
\dis\P\big[(\al^0_z,\vec\al_z)=(\ga_0,\vec\ga)\,\big|\,\Ni'_n\big]
&=&\dis\frac{c_{\rm l}(\vec\ga)}
{\sum_{\vec\ga\neq(-1,\ldots,-1)}c_{\rm l}(\vec\ga)
+\sum_{\vec\ga\neq(+1,\ldots,+1)}c_{\rm r}(\vec\ga)}
\quad&\dis\mbox{if }\ga_0=-1,\\[15pt]
\dis\P\big[(\al^0_z,\vec\al_z)=(\ga_0,\vec\ga)\,\big|\,\Ni'_n\big]
&=&\dis\frac{c_{\rm r}(\vec\ga)}
{\sum_{\vec\ga\neq(-1,\ldots,-1)}c_{\rm l}(\vec\ga)
+\sum_{\vec\ga\neq(+1,\ldots,+1)}c_{\rm r}(\vec\ga)}
\quad&\dis\mbox{if }\ga_0=+1,
\ea\]
where $c_{\rm l}(\vec\ga),c_{\rm r}(\vec\ga)$ are defined in (\ref{clcr}). In particular, if $\nu_{\rm l}=\nu_{\rm r}$, then
\be
c_{\rm l}(\vec\ga)=c(-1,\vec\ga)\quad\mbox{and}\quad
c_{\rm r}(\vec\ga)=c(+1,\vec\ga),
\ee
where we define
\[
c(\ga_0,\vec\ga):=21_{\{0<k<n+1\}}\int q^{k-1}(1-q)^{(n+1)-k-1}\nu(\di q)
\quad\mbox{with}\quad k:=|\{i:0\leq i\leq n,\ \ga_i=+1\}|.
\]
From this, it is easy to check that $\Ni'_n$ has left and right speeds $\bet_-(n+1),\bet_+(n+1)$. By Lemma~\ref{L:momcon}, since $\Wi_0,\ldots,\Wi_n$ can be constructed inside $\Ni'_n$ as in (\ref{WinN}), it follows that $(\Wi_0,\ldots,\Wi_n)$ is equally distributed with $(\Wi_1,\ldots,\Wi_{n+1})$.\qed

\noi
The next lemma implies Theorem~\ref{T:HWconst2}.

\bl{\bf(Limit of moment measures)}\label{L:infmom}
Let $\bet_0\in\R$, let $\nu_{\rm l},\nu_{\rm r}$ be finite measures on $[0,1]$, let $\bet$ and $\nu$ be given by (\ref{tib}) and (\ref{nuform}), and assume the left and right speeds $\bet_-,\bet_+$ defined in (\ref{speeds}) satisfy $-\infty<\bet_-$, $\bet_+<\infty$. Let $(\Wi_0,\Mi)$, with $\Mi=\{(z,\om_z):z\in M\}$, be a marked reference web as in Theorem~\ref{T:HWconst} and conditional on $(\Wi_0,\Mi)$, let $(\Wi_i)_{i\geq 1}$ be independent sample webs constructed as in (\ref{sample}). Then there exists a Brownian net $\Ni_\infty$, which is determined  a.s.\ uniquely by $(\Wi_0,\Mi)$ and has left and right speeds $\bet_-,\bet_+$, such that for any deterministic countable set $\Di\sub\R^2$ and countable dense set of times $\Ti\sub\R$,
\be\label{Ninf}
\ov{\Ni_\infty(\Di)}=\ov{\Hi_\Ti\big(\bigcup_{i\geq 1}\Wi_i\big)(\Di)}\quad{\rm a.s.}
\ee
Let $S_\infty$ be the set of separation points of $\Ni_\infty$. Then $S_\infty\sub M$ and conditional on $\Ni_\infty$, the collection of random variables $\om=(\om_z)_{z\in S_\infty}$ is i.i.d.\ with law $\bar\nu$ defined in Theorem~\ref{T:HWconst2}. A.s., each $z\in S_\infty$ is of type $(1,2)$ in each $\Wi_i$ $(i\geq 1)$. Conditional on $(\Ni_\infty,\om)$, the random variables $(\al^i_z)^{i\geq 1}_{z\in S_\infty}$ defined by $\al^i_z:=\sign_{\Wi_i}(z)$ are independent with $\P[\al^i_z=+1\,|\,(\Ni_\infty,\om)]=\om_z$, and one has
\be\label{Winf}
\Wi_i=\{\pi\in\Ni_\infty:\sign_\pi(z)=\al^i_z\ \forall z\in S_\infty
\mbox{ s.t.\ $\pi$ enters }z\}
\qquad(i\geq 1).
\ee
\el
{\bf Proof.} This is very similar to the proof of Lemma~\ref{L:momcon} so we will only sketch the main line of proof. By the assumption that the speeds $\bet_-,\bet_+$ are finite, conditional on $\Wi_0$, the set $M$ is a Poisson point set with intensity $c_{\rm l}\ell_{\rm l}+c_{\rm r}\ell_{\rm r}$, where $c_{\rm l},c_{\rm r}<\infty$ are given by
\be
c_{\rm l}:=2\int q^{-1}\nu_{\rm l}(\di q)
\quad\mbox{and}\quad
c_{\rm r}:=2\int(1-q)^{-1}\nu_{\rm r}(\di q).
\ee
In analogy with (\ref{newref}), we set
\be\label{Cast}
\Wi'_0:=\lim_{\De_m\up C^\ast}{\rm switch}_{\De_m}(\Wi_0)
\quad\mbox{with}\quad
C^\ast:=\{z\in M_{\rm l}:\om_z=1\}\cup\{z\in M_{\rm r}:\om_z=0\},
\ee
where $M_{\rm l},M_{\rm r}$ denote the restrictions of $M$ to the sets of points of type $(1,2)_{\rm l}$ and $(1,2)_{\rm r}$ in $\Wi_0$, respectively. Then, conditional on $\Wi'_0$, the set $\{(z,\om_z):z\in M\beh C^\ast\}$ is a Poisson point set on $\R^2\times(0,1)$ with intensity
\be\label{acint}
\ell_{\rm l}(\di z)\otimes 2q^{-1}1_{\{q<1\}}\nu_{\rm l}(\di q)
+\ell_{\rm r}(\di z)\otimes 2(1-q)^{-1}1_{\{0<q\}}\nu_{\rm r}(\di q).
\ee
Next, in analogy with (\ref{Nndef}), we define a Brownian net $\Ni_\infty$ with set of separation points $S_\infty:=M\beh C^\ast$ by $\Ni_\infty:=\lim_{\De_m\up S_\infty}\hop_{\De_m}(\Wi'_0)$. Then Proposition~\ref{P:mark} implies (\ref{Winf}), while a calculation similar to (\ref{betncalc}) shows that the left and right speeds of $\Ni_\infty$ are the constants $\bet_-,\bet_+$ from (\ref{speeds}).

By Theorem~\ref{T:marknet} and (\ref{acint}), conditional on $\Ni_\infty$, the random variables $(\sign_{\Wi'_0}(z))_{z\in M\beh C^\ast}$ are i.i.d., where $\P[\sign_{\Wi'_0}(z)=+1\,|\,\Ni_\infty]=c'_{\rm r}/(c'_{\rm l}+c'_{\rm r})$ and $c'_{\rm l}:=2\int q^{-1}1_{\{q<1\}}\nu_{\rm l}(\di q)$, $c'_{\rm r}:=2\int(1-q)^{-1}1_{\{0<q\}}\nu_{\rm r}(\di q)$. Using this, (\ref{nuform}) and (\ref{acint}), we see that conditional on $\Ni_\infty$, the collection of random variables $\om=(\om_z)_{z\in M}$ is i.i.d.\ with common law $\bar\nu$ from Theorem~\ref{T:HWconst2}. Since $(\Ni_\infty,\om)$ is determined a.s.\ by $(\Wi_0,\Mi)$, and since conditional on $(\Wi_0,\Mi)$, the random variables $(\al^i_z)^{i\geq 1}_{z\in S_\infty}$ are independent with $\P[\al^i_z=+1\,|\,(\Wi_0,\Mi)]=\om_z$, the same statement holds for the conditional law given $(\Ni_\infty,\om)$. The proof of (\ref{Ninf}), finally, is completely analogous to the proof of formula (\ref{Nin}) from Lemma~\ref{L:momcon}.\qed
\med

\noi
{\bf Proof of Theorem~\ref{T:HWconst2}.} In Lemma~\ref{L:infmom}, $(\Ni_\infty,\om)$ is determined a.s.\ by $(\Wi_0,\Mi)$. Moreover, by (\ref{Winf}), the conditional law $\P[\Wi_i\in\cdot\,\,|\,(\Wi_0,\Mi)]$ is a function of $(\Ni_\infty,\om)$ only. It follows that $\P[\Wi_i\in\cdot\,\,|\,(\Wi_0,\Mi)]=\P[\Wi_i\in\cdot\,\,|\,(\Ni_\infty,\om)]$.\qed

\subsection{Some immediate consequences of our construction}\label{S:immed}

{\bf Proof of Proposition~\ref{P:regul}~(b)--(d).} Let $\Q$ be the Howitt-Warren quenched law defined in (\ref{HWquen}). Then $\Q$ is a random probability law on the space $\Ki(\Pi)$ of compact subsets of the space of paths. We will be interested in a.s.\ properties of $\Q$ that hold for almost every realization of $(\Wi_0,\Mi)$. Let $\Wi$ be a random variable with law $\Q$. Since averaged over the law of $(\Wi_0,\Mi)$, the set $\Wi$ is distributed as a Brownian web, it follows that for a.e.\ realization of $(\Wi_0,\Mi)$, the compact set $\Wi\sub\Pi$ will satify all the a.s.\ properties of a Brownian web, such as the classification of special points. In particular, we can define the special paths $\pi^+_z$ and $\pi^\up_z$ for each $z\in\R^2$. Then
\be
K^+_{s,t}(x,A):=\Q[\pi^+_{(x,s)}(t)\in A]
\quad\mbox{and}\quad
K^\up_{s,t}(x,A):=\Q[\pi^\up_{(x,s)}(t)\in A]
\ee
are well-defined for {\em every} $s\leq t$, $x\in\R$, $A\in\Bi(\R)$.

To prove part~(b), it then suffices to note that for every $s\leq t_n$, $t_n\to t$, $x\in\R$ and continuous $f:\R\to\R$,
\be
\int K^+_{s,t_n}(x,\di y)f(y)=\Q[f(\pi^+_{(x,s)}(t_n))]
\asto{n}\Q[f(\pi^+_{(x,s)}(t))]=\int K^+_{s,t}(x,\di y)f(y),
\ee
where we also use the symbol $\Q$ to denote expectation with respect to the probability law $\Q$, and we have used the continuity of $t\mapsto\pi^+_{(x,s)}$. The same proof works for $K^\up_{s,t}$.

The proof of part~(c) is similar, where this time, for any $s<t$, $\R\ni x_n\down x$ and $A\in\Bi(\R)$,
\be
K^+_{s,t}(x_n,A)=\Q[\pi^+_{(x_n,s)}(t)\in A]\to
\Q[\pi^+_{(x,s)}(t)\in A]=K^+_{s,t}(x_n,A),
\ee
where we have used that by \cite[Lemma~3.4~(a)]{SS08}, under $\Q$, there is a random $m$ such that $\pi^+_{(x_n,s)}(t)=\pi^+_{(x,s)}(t)$ for all $n\geq m$. The existence of left limits follows in the same way.

To prove part~(d), we observe that for all $s\leq t\leq u$, $x\in\R$ and $A\in\Bi(\R)$,
\be\ba{l}\label{upprop}
\dis\Q\big[\pi^\up_{(x,s)}(u)\in A\big]
=\int\Q\big[\pi^\up_{(x,s)}(u)\in A\,\big|\,\pi^\up_{(x,s)}(t)=y\big]
\Q\big[\pi^\up_{(x,s)}(t)\in\di y\big]\\[5pt]
\dis\quad=\int\Q\big[\pi^\up_{(y,t)}(u)\in A\,\big|\,\pi^\up_{(x,s)}(t)=y\big]
\Q\big[\pi^\up_{(x,s)}(t)\in\di y\big]\\[5pt]
\dis\quad=\int\Q\big[\pi^\up_{(y,t)}(u)\in A\big]
\Q\big[\pi^\up_{(x,s)}(t)\in\di y\big],
\ec
where we have conditioned on the value of $\pi^\up_{(x,s)}(t)$, used the fact that $\pi^\up_{(y,t)}$ is the continuation of any incoming path at $(y,t)$, and in the last step we have used that under the law $\Q$, for any $y\in\R$ that may depend on the marked reference web $(\Wi_0,\Mi)$ but not on the sample web $\Wi$, the path $\pi^\up_{(y,t)}$ is independent of $\pi^\up_{(x,s)}(t)$. To prove this independence, for any $t_1\leq t_2$, let $\Wi|_{t_1}^{t_2}$ denote the restriction of $\Wi$ to the time interval $[t_1,t_2]$, i.e., $\Wi|_{t_1}^{t_2}:=\{\pi|_{t_1}^{t_2}:\pi\in\Wi\}$ where $\pi|_{t_1}^{t_2}:=\{(\pi(u),u):u\in[t_1,t_2]\cap[\sig_\pi,\infty]\}$ is the restriction of a path $\pi$ to $[t_1,t_2]$. It follows from the marking construction in Theorem~\ref{T:HWconst} that for all $t_1\leq t_2\leq t_3$, $\Wi|_{t_1}^{t_2}$ and $\Wi|_{t_2}^{t_3}$ are independent under $\Q$. Since $\pi^\up_{(y,t)}$ is a function of $\Wi|_{t-\eps}^\infty$ for each $\eps>0$, we conclude that $\pi^\up_{(x,s)}(t-\eps)$ is independent of $\pi^\up_{(y,t)}$ under $\Q$ for each $\eps>0$. Since $\pi^\up_{(x,s)}(t)=\lim_{\eps\to 0}\pi^\up_{(x,s)}(t-\eps)$, it follows that $\pi^\up_{(x,s)}(t)$ is independent of $\pi^\up_{(y,t)}$ under $\Q$.\qed

The proof of Proposition~\ref{P:regul}~(a) needs a bit of preparation. We start by proving the statement for the Arratia flow.

\bp{\bf(Measurability of special paths)}\label{P:specmeas}
There exists a measurable function $\Ki(\Pi)\times\R^2\ni(\Ai,z)\mapsto\pi^+_z(\Ai)\in\Pi$ such that if $\Wi$ is a Brownian web, then almost surely, for all $z\in\R^2$ the path $\pi^+_z(\Wi)=\pi^+_z$ is the special path in $\Wi(z)$ defined below Proposition~\ref{P:classweb}. An analogue statement holds for $\pi^\up_z$.
\ep
{\bf Proof.} It suffices to define $\pi^+_z(\Ai)$ on the measurable set of all
$\Ai\in\Ki(\Pi)$ such that $\Ai(z)$ contains a single path $\pi_z$ for each
$z\in\Q^2$. For any $r\in\R$, define $\lfloor r\rfloor_n\up r$ by $\lfloor
r\rfloor_n:=\sup\{r'\in\Z/n:r'<r\}$ and similarly, set $\lceil
r\rceil_n:=\inf\{r'\in\Z/n:r'>r\}$. Then $(\Ai,(x,t))\mapsto\pi_{(\lceil
  x\rceil_n,\lfloor t\rfloor_m)}$ is a measurable function. If $\Wi$ is a
Brownian web, then applying Lemma~\ref{L:nosimu} to the dual web $\hat\Wi$ we
see that for each $t\in\R$, there exists at most one $x\in\Q$ such that
$\Wi(x,t)$ contains more than one path. It follows that for each
$(x,t)\in\R^2$, there is at most one $n$ for which the limit
$\lim_{m\to\infty} \pi_{(\lceil x\rceil_n,\lfloor t\rfloor_m)}$ does not
exist, hence the double limit
\be
\pi^+_{(x,t)}=\lim_{n\to\infty}\lim_{m\to\infty}
\pi_{(\lceil x\rceil_n,\lfloor t\rfloor_m)}
\ee
is well-defined and gives the right-most path in $\Wi(x,t)$. Since pointwise
limits of measurable functions are meaurable, restricting ourselves to a
suitable measurable subset of $\Ki(\Pi)$, we see that $\pi^+_z$ depends
measurably jointly on $z$ and the Brownian web $\Wi$.

To also prove the statement for $\pi^\up_z$, we note that the dual $\hat\Wi$ of a Brownian web $\Wi$ is a measurable function of $\Wi$ and that by what we have just proved, both the left-most and right-most dual paths $\hat\pi^-_z$ and $\hat\pi^+_z$ depend measurably jointly on $z$ and $\hat\Wi$. For any $z=(x,t)\in\R^2$, set $\tau_z:=\inf\{s:s>0,\ \hat\pi^-_z(t-s)=\hat\pi^+_z(t-s)\}$ and
\[
z':=\Big(\ffrac{1}{2}\big(\hat\pi^-_z(t-\ffrac{1}{2}\tau_z)
+\hat\pi^+_z(t-\ffrac{1}{2}\tau_z)\big),t-\ffrac{1}{2}\tau_z\Big).
\]
Then $z'$ depends measurably jointly on $z$ and $\Wi$ and $\pi^\up_z$ is the restriction of $\pi^+_{z'}$ to $[t,\infty)$.\qed

\bl{\bf(No simultaneous incoming paths)}\label{L:nosimu}
Let $\Wi$ be a Brownian web and let $x,y\in\R$, $x\neq y$ be deterministic positions. Then a.s., there exist no time $t\in\R$ such that there exists paths $\pi,\pi'\in\Wi$ with $\sig_\pi,\sig_{\pi'}<t$, $\pi(t)=x$, $\pi'(t)=y$.
\el
{\bf Proof.} By \cite[Lemma~3.4~(b)]{SS08} it suffices to prove the statement for paths $\pi,\pi'$ started at deterministic points. The statement then follows from the fact that two-dimensional Brownian motion a.s.\ does not hit deterministic points.\qed

\bp{\bf(Measurability of quenched laws on path space)}\label{P:Qzmeas}
Let $\Q^+_z$ and $\Q^\up_z$ be the quenched laws on path space defined in (\ref{HWquenz}). Then $\R^2\times\Om\ni(z,\om)\mapsto\Q^+_z(\om)\in\Mi_1(\Pi)$ is a measurable map. An analogue statement holds for $\Q^\up_z$.
\ep
{\bf Proof.} The quenched law $\Q$ is a random variable taking values in $\Mi_1(\Ki(\Pi))$, i.e., a measurable map $\Om\ni\om\mapsto\Q(\om)\ni\Mi_1(\Ki(\Pi))$, where $(\om,\Fi,\P)$ is our underlying probability space. Since $\Q^+_z=\Q\circ{\pi^+_z}^{-1}$, where $(\Ai,z)\mapsto\pi^+_z(\Ai)$ is the measurable map from Proposition~\ref{P:specmeas}, the statement follows from Lemma~\ref{L:imeas} in the appendix. The same argument applies to $\Q^\up_z$.\qed
\med

\noi
{\bf Proof of Proposition~\ref{P:regul}~(a).} Define a continuous map $\Pi\times\R\ni(\pi,t)\mapsto\psi_t(\pi)\in\R$ by $\psi_t(\pi):=\pi(\sig_\pi\vee t)$. Then, since
\be
K^+_{s,t}(x,\,\cdot\,)=\Q^+_{(x,s)}\circ\psi_t^{-1}
=\Q\circ(\psi_t\circ\pi^+_{(x,s)})^{-1}
\qquad(s,t,x\in\R,\ s\leq t),
\ee
the statement follows from Lemma~\ref{L:imeas} in the appendix. The same argument applies to $K^\up_{s,t}(x,\,\cdot\,)$.\qed
\med

\noi
{\bf Proof of Proposition~\ref{P:conpath}.} Immediate from Proposition~\ref{P:regul}~(b) and (d).\qed
\med

Proposition~\ref{P:scale} is a direct consequence of the following proposition, which is formulated on the level of quenched laws on the space of webs.

\bp{\bf(Scaling of quenched laws)}\label{P:quenscale}
Let $\Q$ be a Howitt-Warren quenched law with drift $\bet$ and characteristic measure $\nu$. Define scaling maps $S_a:\Rc\to\Rc$ $(a>0)$ as in (\ref{Seps2}) and let $T_a:\Rc\to\Rc$ $(a\in\R)$ be defined by $T_a(x,t):=(x+at,t)$. Then:
\begin{itemize}
\item[{\bf(a)}] For each $a>0$, $S_a(\Q)$ is a Howitt-Warren quenched law with drift $a^{-1}\bet$ and characteristic measure $a^{-1}\nu$.
\item[{\bf(b)}] For each $a\in\R$, $T_a(\Q)$ is a Howitt-Warren quenched law with drift $\bet+a$ and characteristic measure $\nu$.
\end{itemize}
\ep
{\bf Proof.} To prove part~(a), choose $\eps_k$ and $\mu_k$ such that (\ref{mucon}) holds and set $\eps'_k:=a\eps_k$. Then, by Theorem~\ref{T:quench}, $S_{\eps_k}(\Qdis_{\li k\re})$ converges weakly in law to $\Q$ while, since the $\mu_k$ satisfy (\ref{mucon}) with $\eps_k$ replaced by $\eps'_k$ and $\bet$ and $\nu$ replaced by $a^{-1}\bet$ and $a^{-1}\nu$, respectively, $S_{\eps'_k}(\Qdis_{\li k\re})$ converges weakly in law to a Howitt-Warren quenched law with this drift and characteristic measure. Obviously, $S_{\eps'_k}(\Qdis_{\li k\re})=S_a(S_{\eps_k}(\Qdis_{\li k\re}))$ also converges to $S_a(\Q)$, so the latter is a Howitt-Warren quenched law with drift $a^{-1}\bet$ and characteristic measure $a^{-1}\nu$.

To prove part~(b), let $(\Wi_0,\Mi,\Wi)$ be a marked reference Brownian web and sample Brownian web as in Theorem~\ref{T:HWconst}. Then $T_a(\Wi_0)$ is a Brownian web with drift $\bet_0+a$. It follows from Proposition~\ref{P:refloc} that $T_a(\ell)$ is the intersection local time measure of $T_a(\Wi_0)$ and its dual. It follows that conditional on $\Wi_0$, the set $T_a(\Mi):=\{(T_a(z),\om_z):(z,\om_z)\in\Mi\}$ is a Poisson point set with intensity as in (\ref{Poisint}), with $\ell_{\rm l}$ and $\ell_{\rm r}$ replaced by $T_a(\ell_{\rm l})$ and $T_a(\ell_{\rm r})$. Since $T_a(\Wi)$ is constructed from $T_a(\Wi_0)$ and $T_a(\Mi)$ in the same way as $\Wi$ is constructed from $(\Wi_0,\Mi)$, in particular, it follows that $\P[T_a(\Wi)\in\cdot\,|\,(\Wi_0,\Mi)]$ is a Howitt-Warren quenched law with drift $\bet+a$ and characteristic measure $\nu$.\qed
\med

\noi
{\bf Proof of Proposition~\ref{P:scale}.} Immediate from Proposition~\ref{P:quenscale}.\qed

\section{Support properties}\label{S:support}

In this section, we will first prove Theorem~\ref{T:halfnet} on the
characterization of Brownian half-nets, establish its connection to
Howitt-Warren flows and prove some of its basic properties. We will then prove
Theorem~\ref{T:support} on the image set of the support of the Howitt-Warren
quenched law $\Q$, from which Theorems~\ref{T:speed} and \ref{T:supp} on the
support properties of Howitt-Warren processes follow immediately.

\subsection{Generalized Brownian nets}\label{S:halfnet}

{\bf Proof of Theorem~\ref{T:halfnet}.} Let $\gamma_n\in \Hi_-$ and
$\gamma_n\to \gamma$ in $\Pi$. If $\gamma\notin \Hi_-$ so that it crosses some
$\pi\in\Wi$ from left to right, i.e., $\gamma(s)<\pi(s)$ and
$\pi(t)<\gamma(t)$ for some $s<t$, then $\gamma_n$ crosses $\pi$ from left to
right for all $n$ large, a contradiction. Therefore $\Hi_-$ is a.s.\ closed.

The two characterizations of $\Hi_-$ in Theorem~\ref{T:halfnet}~(i) and (ii)
are equivalent, because if $\gamma\in \Pi$ crosses some $\pi\in\Wi$ from left
to right, then by the non-crossing property of paths in $\Wi$ and $\hat\Wi$,
$\gamma$ must also cross some $\hat\pi\in\hat\Wi$ from left to right, and the
same is true if $\Wi$ and $\hat\Wi$ are interchanged.

It only remains to show that for each deterministic $z=(x,t)\in\R^2$, if
$\pi_z$ denotes the a.s.\ unique path in $\Wi$ starting from $z$, then $\pi_z$
is the maximal element in $\Hi_-(z)$. Certainly $\pi_z\in \Hi_-(z)$. Note that
$z$ is a.s.\ of type $(0,1)$ in $\Wi$ by Prop.~\ref{P:classweb}, and hence for
any positive sequence $\eps_n\downarrow 0$, $\pi_{(x+\eps_n, t)} \to \pi_z$ as
$n\to\infty$.  If $\gamma\in \Hi_-(z)$, then it cannot cross
$\pi_{(x+\eps_n,t)}\in\Wi((x+\eps_n,t))$ from left to right. Therefore
$\gamma\leq \pi_{(x+\eps_n,t)}$ on $[t,\infty)$ for all $n\in\N$, which
implies that $\gamma \leq \pi_z$ on $[t,\infty)$. Therefore $\pi_z$ is the
maximal element in $\Hi_-(z)$.\qed

For any $-\infty\leq\bet_-\leq\bet_+\leq\infty$ with $\bet_-<\infty$ and $-\infty<\bet_+$, we define a {\em generalized Brownian net} with speeds $\bet_-,\bet_+$ to be a Brownian net with these speeds if $\bet_-,\bet_+$ are both finite, a Brownian half-net with these speeds if one of $\bet_-,\bet_+$ is infinite, and the space of all paths $\Pi$ if both speeds $\bet_-,\bet_+$ are infinite.

Consider a reference Brownian web $\Wi_0$ and set of marked points $\Mi$ as in Theorem~\ref{T:HWconst} and conditional on $(\Wi_0,\Mi)$, construct an i.i.d.\ sequence $\Wi_1,\Wi_2,\ldots$ of sample Brownian webs as in (\ref{sample2}). For each $n\geq 1$, let $\Ni_n$ denote the Brownian net containing $\Wi_1,\ldots,\Wi_n$ introduced in Lemma~\ref{L:momcon}. Recall that $\Ni_n$ has left and right speeds $\bet_-(n),\bet_+(n)$ given by (\ref{betplusm}), which converge as $n\to\infty$ to the speeds $\bet_-,\bet_+$ given by (\ref{speeds}).

\index{Brownian net!generalized}
\bl{\bf(Generalized Brownian net associated with
Howitt-Warren flow)}\label{L:Ngen}
Let $\bet_0\in\R$ and let $\nu_{\rm l},\nu_{\rm r}$ be finite measures on $[0,1]$. Let $(\Wi_0,\Mi)$, with $\Mi=\{(z,\om_z):z\in M\}$, be a marked reference web as in Theorem~\ref{T:HWconst} and conditional on $(\Wi_0,\Mi)$, let $(\Wi_i)_{i\geq 1}$ be independent sample webs constructed as in (\ref{sample}). Then there exists a generalized Brownian net $\Ni_\infty$ with left and right speeds given by (\ref{speeds}), which is a.s.\ uniquely defined by $(\Wi_0,\Mi)$, such that for any deterministic countable set $\Di\sub\R^2$ and countable dense set of times $\Ti\sub\R$
\be\label{infhop}
\ov{\Ni_\infty(\Di)}=\ov{\Hi_\Ti\big(\bigcup_{i\geq 1}\Wi_i\big)(\Di)}\quad{\rm a.s.}
\ee
\el
{\bf Proof.} We treat the cases when $\Ni_\infty$ is a Brownian net, a Brownian halfnet or the space of all paths $\Pi$ separately. If the speeds $\bet_-,\bet_+$ from (\ref{speeds}) are both finite, then the statements follow from Lemma~\ref{L:infmom}.

If only one of the speeds $\bet_-,\bet_+$ from (\ref{speeds}) is finite, then by symmetry, we may without loss of generality assume that $-\infty<\bet_-$ and $\bet_+=\infty$. We claim that without loss of generality, we may further assume that $\nu_{\rm r}=0$. To see this, let $M_{\rm l}$ and $M_{\rm r}$ be the restrictions of $M$ to the sets of points of type $(1,2)_{\rm l}$ and $(1,2)_{\rm r}$ in $\Wi_0$, respectively. As a first step, we reduce our problem to the case that $\nu_{\rm l}(\{1\})=0=\nu_{\rm r}(\{0\})$. If this is not yet the case, then define $\Wi'_0$ and $C^\ast$ as in (\ref{Cast}) and replace $\Wi_0$ by $\Wi'_0$ and $M$ by $M\beh C^\ast$.

Next, we observe that conditional on $\Wi_0$, the set $M_{\rm r}$ is a Poisson point set with intensity $c_{\rm r}\ell_{\rm r}$, where $c_{\rm r}:=2\int(1-q)^{-1}\nu_{\rm r}(\di q)<\infty$ by our assumption that $-\infty<\bet_-$. Let
\be
\Wi'_0:=\lim_{\De_m\up M_{\rm r}}\switch_{\De_m}(\Wi_0).
\ee
Using Proposition~\ref{P:refchange}, it is not hard to see that conditional on $\Wi'_0$, the set $\{(z,\om_z):z\in M\}$ is a Poisson point set on $\R^2\times(0,1]$ with intensity
\be\label{rint}
\ell_{\rm l}(\di z)\otimes\big(21_{\{0<q\}}q^{-1}\nu_{\rm l}(\di q)+
2(1-q)^{-1}\nu_{\rm r}(\di q)\big)
=\ell_{\rm l}(\di z)\otimes 2q^{-1}(1-q)^{-1}1_{\{0<q\}}\nu(\di q),
\ee
where we have used (\ref{nuform}) and the fact that $\nu_{\rm l}(\{1\})=0=\nu_{\rm r}(\{0\})$. By Proposition~\ref{P:refchange} and Lemma~\ref{L:switchcom}, conditional on $(\Wi'_0,\om)$, the random variables $(\al^i_z)^{i\geq 1}_{z\in M}$ defined by $\al^i_z:=\sign_{\Wi_i}(z)$ are independent with $\P[\al^i_z=+1\,|\,(\Wi'_0,\om)]=\om_z$, and
\be\label{rswitch}
\Wi_i=\lim_{\De_m\up A_i\cup B_i}\switch_{\De_m}(\Wi'_0)\qquad(i\geq 1),
\ee
where $A_i=\{z\in M:\al^i_z\neq\sign_{\Wi'_0}(z)\}$ and $B_i$ is an independent Poisson point set with intensity $2\nu_{\rm l}(\{0\})$. Thus, replacing our reference web $\Wi_0$ by $\Wi'_0$, we have reduced our problem to the case that $\nu_{\rm r}=0$ and $\nu_{\rm l}(\di q)=(1-q)^{-1}\nu(\di q)$.

In light of this, assume from now on that $\nu_{\rm r}=0$. Then (\ref{tib}) tells us that $\bet=\bet_0+2\nu_{\rm l}([0,1])=\bet_0+2\int(1-q)^{-1}\nu(\di q)$, hence $\bet_0=\bet_-$, the left speed from (\ref{speeds}). Let $\Ni_\infty$ be the Brownian halfnet with left Brownian web $\Wi_0$. Let $\Ni_n$ $(n\geq 1)$ be the Brownian nets defined in Lemma~\ref{L:momcon}. We recall from formulas (\ref{newref}) and (\ref{Nndef}) in the proof of that lemma that $\Ni_n$ is constructed by switching and then allowing hopping at subsets of the set $C:=\bigcup_i(A_i\cup B_i)$. Since points in $A_i\cup B_i$ are of type $(1,2)_{\rm l}$ in $\Wi_0$, it follows that paths in $\Ni_n$ cannot cross paths in $\Wi_0$ from right to left, hence $\Ni_n\sub\Ni_\infty$.

Let $(\Wl_n,\Wr_n)$ be the left-right Brownian web associated with $\Ni_n$. Let $z=(x,s)\in\R^2$ be deterministic and let $l^n_z,r^n_z$ denote the a.s.\ unique elements of $\Wl_n(z),\Wr_n(z)$, respectively. Since $\Ni_n\sub\Ni_{n+1}$, one has $l^n_z\down l^\infty_z$ and $r^n_z\up r^\infty_z$ for some functions $l^\infty_z:[s,\infty)\to[-\infty,\infty)$ and $r^\infty_z:[s,\infty)\to(-\infty,\infty]$. Since $\Ni_n\sub\Ni_\infty$ for each $n\geq 1$, we have $\pi^0_z\leq l^\infty_z$, where $\pi^0_z$ denotes the a.s.\ unique element of $\Wi_0$ starting at $z$. Since $\pi^0_z(t)$ and $l^\infty_z(t)$ are normally distributed with the same mean for each $t>s$, we must have $\pi^0_z=l^\infty_z$. On the other hand, since the $r^n_z$ are Brownian motions with drifts tending to infinity we must have $r^\infty_z(t)=\infty$ for all $t>s$.

We are now ready to prove (\ref{infhop}). Since $\Wi_i\sub\Ni_n\sub\Ni_\infty$ for each $i\leq n$ and since by \cite[Lemma~8.3]{SS08}, a Brownian net is closed under hopping at deterministic times, we see that $\Hi_\Ti(\bigcup_{i\geq 1}\Wi_i)=\bigcup_{n\geq 1}\Hi_\Ti(\Wi_1\cup\cdots\cup\Wi_n)\sub\bigcup_{n\geq 1}\Ni_n\sub\Ni_\infty$ and therefore $\ov{\Hi_\Ti(\bigcup_{i\geq 1}\Wi_i)(\Di)}\sub\ov{\Ni_\infty(\Di)}$. To prove the other inclusion, we first observe that $\ov{\Ni_n(\Di)}=\ov{\Hi_\Ti(\Wi_1\cup\cdots\cup\Wi_n)(\Di)}\sub\ov{\Hi_\Ti\big(\bigcup_{i\geq 1}\Wi_i\big)(\Di)}$ for all $n\geq 1$, hence $\bigcup_{n\geq 1}\Ni_n(\Di)\sub\ov{\Hi_\Ti\big(\bigcup_{i\geq 1}\Wi_i\big)(\Di)}$. Since $\Ni_1\sub\Ni_2\sub\cdots$ and since by Theorem~\ref{T:net}, each $\Ni_n$ is closed under hopping at intersection times, it follows that also $\bigcup_{n\geq 1}\Ni_n$ is closed under hopping at intersection times, i.e., $\Hi_{\rm int}(\bigcup_{n\geq 1}\Ni_n)=\bigcup_{n\geq 1}\Ni_n$. In view of this, we have $\ov{\Hi_{\rm int}(\bigcup_{n\geq 1}\Ni_n)(\Di)}=\ov{\bigcup_{n\geq 1}\Ni_n(\Di)}\sub\ov{\Hi_\Ti\big(\bigcup_{i\geq 1}\Wi_i\big)(\Di)}$, so it suffices to show that $\Ni_\infty(\Di)\sub\ov{\Hi_{\rm int}(\bigcup_{n\geq 1}\Ni_n)(\Di)}$. By Lemma~\ref{L:finpath} below, it suffices to show that each path $\pi\in\Ni_\infty(\Di)$ with $-\infty<\pi<\infty$ on $[\sig_\pi,\infty)$ can be approximated by paths in $\ov{\Hi_{\rm int}(\bigcup_{n\geq 1}\Ni_n)(\Di)}$.

\begin{figure}[htb]
\begin{center}
\includegraphics[width=12cm]{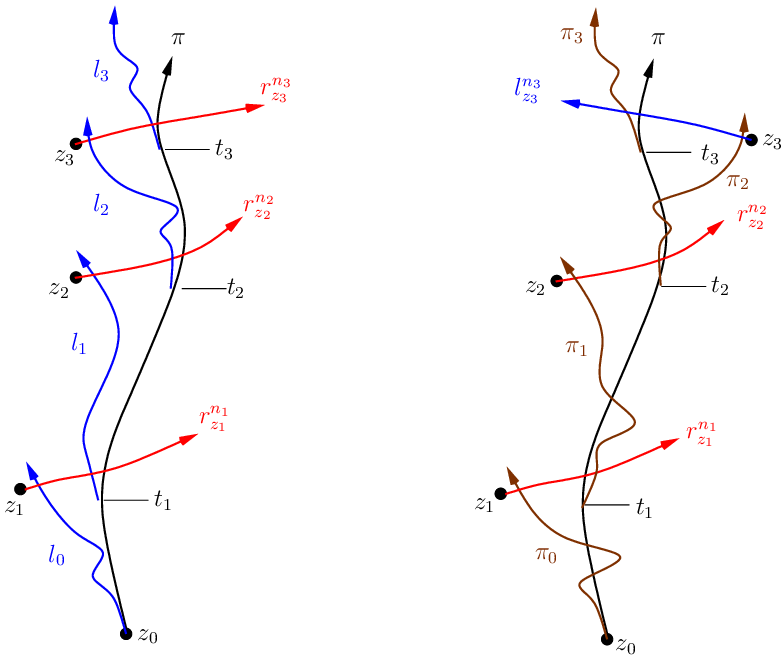}
\caption{Approximation of paths in generalized Brownian nets. On the left: a Brownian halfnet with finite left speed. On the right: the case when both the left and right speeds are infinite.}
\label{fig:Ninf}
\end{center}
\end{figure}

Let $\pi$ be such a path, $\eps>0$ and $T<\infty$. Let $z_0=(x_0,t_0)$ denote the starting point of $\pi$ and inductively choose times $t_k$ $(k\geq 1)$ and paths $l_k\in\Wi_0(\pi(t_k),t_k)$ $(k\geq 0)$ such that $t_k=\inf\{t\geq t_{k-1}:\pi(t)-l_{k-1}(t)\geq\eps\}$. (See Figure~\ref{fig:Ninf}.) Since paths in $\Ni_\infty$ do not cross paths in $\Wi_0$ from right to left, we can moreover choose the $l_k$ such that $l_k\leq\pi$ for each $k\geq 0$. The continuity of $\pi$ and the equicontinuity of $\Wi_0$ imply that $t_m\geq T$ for some $m\geq 1$. Let $\Di'\sub\R^2$ be a deterministic countable dense set. Choose $z_k=(x_k,t'_k)\in\Di'$ $(k\geq 1)$ such that $t_k<t'_k<t_{k+1}$, $x_k<l_{k-1}(t'_k)$, and $\pi-l_{k-1}<2\eps$ on $[t_{k-1},t'_k]$ and choose $n_k$ such that the right-most path $r^{n_k}_{z_k}$ in $\Ni_{n_k}(z_k)$ crosses $l_k$ before time $t_{k+1}$ and before $\pi-l_{k-1}$ exceeds $2\eps$. Then the concatenation of the paths $l_0,r^{n_1}_{z_1},l_1,r^{n_2}_{z_2},\ldots,r^{n_m}_{z_m},l_m$ approximates the path $\pi$ on $[t_0,T]$ within distance $2\eps$. We are not quite done yet, however, since the paths $l_k\in\Wi_0$ are not elements of $\bigcup_{n\geq 1}\Ni_n$.

To finish the argument, we will show that each path in $\Wi_0$ can be approximated by paths in $\bigcup_{n\geq 1}\Ni_n$ and in case $z_0$ is a deterministic point, that $l_0$ can be approximated by paths in $\bigcup_{n\geq 1}\Ni_n(z_0)$. Indeed, this follows from the fact, proved above, that for each deterministic $z\in\R^2$, the left-most path $l^{n}_z$ in $\Ni_n(z)$ converges to the a.s.\ unique path in $\Wi_0(z)$, and the fact that $\Wi_0=\ov{\Wi_0(\Di')}$. Replacing the paths $l_0,\ldots,l_m$ by sufficiently close approximating paths $l'_0,\ldots,l'_m\in\bigcup_{n\geq 1}\Ni_n$, with $l'_0$ also starting at $z_0$, we see that these approximating paths are crossed by the paths $r^{n_1}_{z_1},\ldots,r^{n_m}_{z_m}$ and hence $\pi$ can be approximated by a concatenation of $l'_0,r^{n_1}_{z_1},l'_1,r^{n_2}_{z_2},\ldots,r^{n_m}_{z_m},l'_m$. This completes the proof for case $-\infty<\bet_-$ and $\bet_+=\infty$, where we have identified $N_\infty$ as a Brownian half-net.

The proof for the case $\bet_-=-\infty$ and $\bet_+=\infty$ is similar, but easier. In this case, for an arbitrary path $\pi\in\Pi$ with starting point $z_0=(x_0,t_0)\in\Di$, we inductively choose times $t_k$ $(k\geq 1)$ and paths $\pi_k\in\Wi_1(\pi(t_k),t_k)$ $(k\geq 0)$ such that $t_k=\inf\{t\geq t_{k-1}:|\pi(t)-\pi_{k-1}(t)|\geq\eps\}$. If $\pi_k(t_k)<\pi(t_k)$, then we use a right-most path $r^{n_k}_{z_k}$ of the Brownian net $\Ni_n:=\ov{\Hi_\Ti(\Wi_1\cup\cdots\cup\Wi_n)}$ to connect $\pi_k$ to $\pi_{k+1}$ and if $\pi(t_k)<\pi_k(t_k)$, then we use a left-most path $l^{n_k}_{z_k}$ of $\Ni_n$.\qed

The proof of Lemma~\ref{L:Ngen} has a useful corollary.
\bcor{\bf(Paths going to infinity)}\label{C:infpath}
Let $\Ni$ be a generalized Brownian net with infinite right speed. Then, for each deterministic $z=(x,t)\in\R^2$, there a.s.\ exist paths $r^n_z\in\Ni(z)$ with $r^n_z(u)\up\infty$ for all $u>t$.
\ecor
{\bf Proof.} If $\Ni$ is a Brownian half-net, then we may take for $r^n_z$ the a.s.\ unique right-most paths in the Brownian nets $\Ni_n$ as in the proof of Lemma~\ref{L:Ngen}. If $\Ni=\Pi$ the statement is trivial.\qed

\bl{\bf(Hopping with left-most paths)}\label{L:lefthop}
Let $\Ni$ be a generalized Brownian net with left speed $-\infty<\bet_-$ and let $\Wl$ be its associated left Brownian web. Then, a.s.\ for each $\pi\in\Ni$ and $l\in\Wl$ with $\pi(t)<l(t)$ at $t:=\sig_\pi\vee\sig_l$, the following statements hold:
\begin{itemize}
\item[\rm(i)] The path $\pi'$ defined by $\sig_{\pi'}:=\sig_\pi$, $\pi':=\pi$ on $[\sig_\pi,t]$ and $\pi':=\pi\wedge l$ on $[t,\infty)$ satisfies $\pi'\in\Ni$.
\item[\rm(ii)] The path $\pi''$ defined by $\sig_{\pi''}:=\sig_l$, $\pi'':=l$ on $[\sig_l,t]$ and $\pi'':=\pi\vee l$ on $[t,\infty)$ satisfies $\pi''\in\Ni$.
\end{itemize}
\el
{\bf Proof.} Set $\tau:=\inf\{u\geq t:\pi(u)>l(u)\}$. Since paths in $\Ni$ cannot cross paths in $\Wl$ from right to left, one has $\pi'=\pi$ on $[\sig_\pi,\tau]$ and $\pi'=l$ on $[\tau,\infty]$, and also, $\pi''=l$ on $[\sig_l,\tau]$ and $\pi''=\pi$ on $[\tau,\infty]$. Therefore, if $\Ni$ is a Brownian net, then the statements follow from the fact that $\Ni$ is closed under hopping at intersection times (see Theorem~\ref{T:net}). We do not know if Brownian half-nets are closed under hopping at intersection times, so in this case we check directly that the paths $\pi'$ and $\pi''$ do not cross paths in $\Wl$ from right to left and hence (by Theorem~\ref{T:halfnet}) are elements of $\Ni$. If $l'\in\Wl$ satisfies $l'(s)<\pi'(s)$ for some $s\in[\sig_\pi,\infty)$, then $l'\leq\pi'$ on $[s,\tau]$ by the fact that $\pi\in\Ni$ while $l'\leq\pi'$ on $[\tau,\infty)$ by the fact that $l'$ cannot cross $l$. If $l'\in\Wl$ satisfies $l'(s)<\pi''(s)$ for some $s\in[\sig_l,\tau]$, then $l'\leq l\leq\pi''$ on $[s,\infty)$ by the fact that $l'$ cannot cross $l$, while if $l'\in\Wl$ satisfies $l'(s)<\pi''(s)$ for some $s\in[\tau,\infty]$, then $l'\leq\pi\leq\pi''$ on $[s,\infty)$ by the fact that $l'$ cannot cross $\pi$ from left to right.\qed

\bl{\bf(Finite paths)}\label{L:finpath}
Let $\Ni$ be a generalized Brownian net and let $\Ni_{\rm fin}:=\{\pi\in\Ni:-\infty<\pi<\infty\mbox{ on }[\sig_\pi,\infty)\}$. Then a.s.\ for each $z\in\R^2$, $\ov{\Ni_{\rm fin}(z)}=\Ni(z)$.
\el
{\bf Proof.} We need to show that for each $\pi\in\Ni(z)$ there exist $\pi_n\in\Ni_{\rm fin}(z)$ with $\pi_n\to\pi$. In case both the left and right speeds of $\Ni$ are infinite, and hence $\Ni=\Pi$, we may simply take $\pi'_n:=-n\vee(n\wedge\pi_n)$. Otherwise, by symmetry, we may without loss of generality assume that the left speed of $\Ni$ is finite. Let $\Wl$ be the left Brownian web. Write $z=(x,t)$ and for $n$ large enough such that $-n\leq x\leq n$ and $t\leq n$, choose $l^-_n,l^+_n\in\Wl$ with $\sig_{l^-_n}=\sig_{l^+_n}:=t$ such that $l^-_n\leq-n$ and $n\leq l^+_n$ on $[t,n]$ and set $\pi_n:=l^-_n\vee(l^+_n\wedge\pi)$. Then $\pi_n\in\Ni$ by Lemma~\ref{L:lefthop}, $-\infty<\pi_n<\infty$ on $[t,\infty)$, and $\pi_n\to\pi$.\qed

In Lemmas~\ref{L:momcon}, \ref{L:infmom} and \ref{L:Ngen}, we have seen objects of the form $\ov{\Ni(\Di)}$ where $\Ni$ is a generalized Brownian net and $\Di\sub\R^2$ is a deterministic countable set. Naively, one might guess that $\ov{\Ni(\Di)}=\Ni(\,\ov\Di\,)$, where $\ov\Di$ denotes the closure of $\Di$ in $\Rc$. It turns out, however, that this is not always true. In particular, it may happen that $(\ast,-\infty)\in\ov\Di$ but $\Ni(\ast,-\infty)\not\sub\ov{\Ni(\Di)}$. Our next result shows that this is indeed all that can go wrong. It is not very difficult, but rather tedious, to give a precise decription of $\ov{\Ni(\Di)}(\ast,-\infty)$ in terms of the shape of $\Di$ near $(\ast,-\infty)$. Since we will not need such a precise description, we will settle for the following lemma.

\bl{\bf(Closure of paths started from a countable set)}\label{L:NgenA}
Let $\Ni$ be a generalized Brownian net, let $\Di\sub\R^2$ be a deterministic countable set, let $\ov{\Ni(\Di)}$ denote the closure of $\Ni(\Di)$ in $\Pi$ and let $\ov\Di$ denote the closure of $\Di$ in $\Rc$. Then
\be\label{ovNi}
\Ni\big(\ov\Di\beh\{(\ast,-\infty)\}\big)\sub\ov{\Ni(\Di)}\sub\Ni(\,\ov\Di\,)
\quad{\rm a.s.}
\ee
If $\Di$ is dense in $\R^2$, then moreover $\Ni(\ast,-\infty)\sub\ov{\Ni(\Di)}$.
\el
{\bf Proof.} Since $\Ni$ is closed and since convergence of paths implies convergence of their starting points, the inclusion $\ov{\Ni(\Di)}\sub\Ni(\ov{\Di})$ is trivial. As a next step, we will show that $\Ni_{\rm fin}(\hat\Di)\sub\ov{\Ni(\Di)}$, where $\hat\Di$ denotes the closure of $\Di$ in $\R^2$ and $\Ni_{\rm fin}$ is defined in Lemma~\ref{L:finpath}. If $\Ni$ is a Brownian net, then the statement follows from \cite[Lemma~8.1]{SS08}. If $\Ni=\Pi$, the statement is trivial, so it suffices to treat the case when $\Ni$ is a Brownian half-net with, say, a finite left speed. We need to show that any path in $\Ni_{\rm fin}(\hat\Di)$ can be approximated by paths in $\Ni(\Di)$. We use the steering argument from the proof of Lemma~\ref{L:Ngen} (see Figure~\ref{fig:Ninf}), where in this case, for any $z_0=(x_0,t_0)\in\hat\Di$, we choose for $l_0$ the left-most path in $\Wi_0(z_0)$, which by the fact that paths in $\Ni$ cannot cross paths in $\Wi_0$ from right to left is guaranteed to stay on the left of any path $\pi\in\Ni$ starting from $z_0$. To apply the proof of Lemma~\ref{L:Ngen}, we only need to show that such an $l_0$ can be approximated by paths in $\bigcup_{n\geq 1}\Ni_n$ that moreover start in $\Di$. Since in the proof of Lemma~\ref{L:Ngen} it has already been shown that any path in $\Wi_0(\Di)$ can be approximated by paths in $\bigcup_{n\geq 1}\Ni_n(\Di)$, by a diagonal argument, it suffices to show that $l_0$ can be approximated by paths in $\Wi_0(\Di)$. By the structure of special points in a Brownian web, such an approximation is possible unless there exists a dual path $\hat l\in\hat\Wi_0$ entering $z_0$ and an $\eps>0$ such that $\{z=(x,t):|z-z_0|<\eps,\ x\leq\hat l(t)\}\cap\Di=\emptyset$. But by Lemma~\ref{L:noskim} below, such a dual path a.s.\ does not exist for any $z_0\in\hat\Di$.

This completes the proof that $\Ni_{\rm fin}(\hat\Di)\sub\ov{\Ni(\Di)}$. By Lemma~\ref{L:finpath}, it follows that also $\Ni(\hat\Di)\sub\ov{\Ni(\Di)}$. Obviously, the trivial path starting at time $(\ast,\infty)$ is an element of $\ov{\Ni(\Di)}$ if and only if there exists $\Di\ni z_n\to(\ast,\infty)$. To complete the proof that $\Ni(\ov\Di\beh\{(\ast,-\infty)\})\sub\ov{\Ni(\Di)}$, by symmetry, it therefore suffices to show that if $\pi\in\Ni(\ov\Di)$ starts at some $z\in\{-\infty\}\times\R$, then $\pi$ can be approximated by $\pi_n\in\Ni(\Di)$. Let $\Di\ni z_n=(x_n,t_n)\to z=(-\infty, t)$. If the left and right speeds of $\Ni$ are $-\infty$ and $+\infty$, then $\Ni=\Pi$ and the claim is trivial. If the right speed of $\Ni$ is finite, then $\Ni(z)$ contains a single path $\pi$ with $\pi(\cdot)=-\infty$ on $[t,\infty)$, which can be approximated by choosing $\pi_n$ to be the rightmost path in $\Ni(z_n)$. Now suppose that the left speed of $\Ni$ is finite while the right speed is $+\infty$. If either $t_n<t$ or $\pi(t_n)<x_n$, then we can simply take $\pi_n:=\pi_n'\vee \pi$, where $\pi_n'$ is the leftmost path in $\Ni(z_n)$ and $\pi_n\in\Ni(z_n)$ by Lemma~\ref{L:lefthop}. If $t_n>t$ and $\pi(t_n)\geq x_n$, then we construct a $\pi_n'\in\Ni(z_n)$ by first following a path in $\Ni(z_n)$ which crosses to the right of $\pi$ before time $t_n+\delta_n$ for some $\delta_n$ (with $\delta_n\downarrow 0$ as $n\to\infty$), and then hop onto a leftmost path $l_n$ with $\sigma_{l_n}<t_n+\delta_n$ and $l_n(t_n+\delta_n)>\pi(t_n+\delta_n)$. By Lemma~\ref{L:lefthop}, $\pi_n'\in\Ni(z_n)$ . We can then define $\pi_n\in\Ni(z_n)$ to equal $\pi'_n$ up to time $t_n+\delta_n$, and define $\pi_n:=\pi'_n\vee \pi=l_n\vee \pi$ from time $t_n+\delta_n$ onward. Then $\pi_n\in\Ni(z_n)$ and $\pi_n\to\pi$ as $n\to\infty$. This completes the proof of (\ref{ovNi}).


To prove that also $\Ni(\ast,-\infty)\sub\ov{\Ni(\Di)}$ if $\Di$ is dense in $\R^2$, by (\ref{ovNi}), it suffices to prove that each path $\pi\in\Ni(\ast,-\infty)$ can be approximated by paths $\pi_n\in\Ni(\Rc\beh\{(\ast,-\infty)\})$. In view of this, taking for $\pi_n$ the restriction of $\pi$ to $[-n,\infty]$ completes the proof of the lemma.\qed

A statement very similar to the lemma below has been demonstrated in the proof of \cite[Lemma~8.1]{SS08}.

\bl{\bf(The Brownian web does not skim closed sets)}\label{L:noskim}
Let $\Wi$ be a Brownian web and let $K\sub\R^2$ be a deterministic closed set. Then a.s.\ there exist no $z_0=(x_0,t_0)\in K$ and $\pi\in\Wi$ entering $z_0$, such that $\{z=(x,t)\in\R^2:|z-z_0|<\eps,\ x<\pi(t)\}\cap K=\emptyset$ for some $0<\eps<t-\sig_\pi$.
\el
{\bf Proof.} It suffice to prove the statement for paths $\pi$ started from a deterministic point and for deterministic $\eps>0$. By cutting $K$ into countably any pieces of diameter at most $\eps/4$ and using translation invariance, we can reduce the problem to the following statement: let $\pi$ be a Brownian motion started at time zero in the origin, let $K$ be a deterministic closed subset of $\R^2$, let $\eps>0$ and let $U:=\{(x,t):t\geq\eps,\ x<\pi(t)\}$. Then the event that $U\cap K=\emptyset$ but $\ov U\cap K\neq\emptyset$ has probability zero. To see that this is the case, set $\pi^\eps(t):=\pi(t)-\pi(\eps)$ $(t\geq\eps)$. Then, conditional on the path $(\pi^\eps(t))_{t\geq\eps}$, there is at most one value of $\pi(\eps)$ for which the event $\{U\cap K=\emptyset,\ \ov U\cap K\neq\emptyset\}$ occurs. Since $\pi(\eps)$ is independent of $(\pi^\eps(t))_{t\geq\eps}$ and normally distributed, we see that the conditional probability of $\{U\cap K=\emptyset,\ \ov U\cap K\neq\emptyset\}$ given $(\pi^\eps(t))_{t\geq\eps}$ is zero, hence integrating over the distribution of $(\pi^\eps(t))_{t\geq\eps}$ yields the desired result.\qed

Let $\Ni$ be a generalized Brownian net with left and right speeds $\bet_-,\bet_+$ (which may be infinite). Then, generalizing (\ref{braco}) and (\ref{xiA}), for any closed subset $A\sub\R$, setting
\be\label{braco2}
\xi^A_t:=\big\{\pi(t):\pi\in\Ni(A\times\{0\})\big\}\qquad(t\geq 0)
\ee
defines a Markov process taking values in the space of closed subsets of $\R$,
which we call the {\em branching-coalescing point set} with {\em left} and
{\em right speeds} $\bet_-,\bet_+$.

\index{branching-coalescing point set!with infinite speeds}
\bl{\bf(Edge of a branching-coalescing point set)}\label{L:edge}
Let $A$ be a deterministic non\-empty closed subset of the real line and let $(\xi^A_t)_{t\geq 0}$ be the branching-coalescing point set with left and right speeds $\bet_-,\bet_+$ defined in (\ref{braco2}). Set $r_t:=\sup(\xi^A_t)$ $(t\geq 0)$. Then:\med

\noi
{\bf(a)} If $\bet_+<\infty$ and $r_0<\infty$, then $(r_t)_{t\geq 0}$ is a Brownian motion with drift $\bet_+$. If $\bet_+<\infty$ and $r_0=\infty$, then $r_t=\infty$ for all $t\geq 0$.\med

\noi
{\bf(b)} If $\bet_+=\infty$, then $r_t=\infty$ for all $t>0$.\med

\noi
{\bf(c)} If $\bet_-=-\infty$ and $\bet_+<\infty$, then $\xi^A_t=(-\infty,r_t]\cap\R$ for all $t>0$.\med

\noi
{\bf(d)} If $\bet_-=-\infty$ and $\bet_+=\infty$, then $\xi^A_t=\R$ for all $t>0$.
\el
{\bf Proof.} To prove part~(a), let $\Wr$ be the right Brownian web associated with $\Ni$. If $\sup(A)<\infty$, then let $r$ be the a.s.\ unique path in $\Wr$ started from $r_0=\sup(A)$. Now $r$ is a Brownian motion with drift $\bet_+$, $r\in\Ni$, and $\pi\leq r$ for any path in $\Ni$ started from $A\times\{0\}$ by the fact that paths in $\Ni$ cannot cross paths in $\Wr$ (see Theorem~\ref{T:halfnet}~(i) and \cite[Prop.~1.8]{SS08}). If $\sup(A)=\infty$, then choose $x_n\in A$ with $x_n\up\infty$ and let $r_n$ be the a.s.\ unique paths in $\Wr$ started from $(x_n,0)$. Then $\P[\inf\{r^n(t):0\leq t\leq T\}\leq N]\to 0$ as $n\up\infty$ for each $N,T<\infty$, hence $\sup_nr^n_t=\infty$ for all $t\geq 0$ a.s.

Part~(b) is an immediate consequence of Corollary~\ref{C:infpath}.

To prove part~(c), we first consider the case that $\sup(A)<\infty$. Let $\Wr$ be the right Brownian web of $\Ni$ and let $r$ be the a.s.\ unique path in $\Wr$ started from $r_0=\sup(A)$. It has been shown in the proof of part~(a) that $\xi^A_t\sub(-\infty,r_t]$ for all $t\geq 0$. To prove the other inclusion, it suffices to show that for each $r'\in\Wr$ started from $\{(x,t):t\geq 0,\ x<r_t\}$ and $\eps>0$ there exists a path $\pi\in\Ni$ started at time zero from $\pi(0)=\sup(A)$ such that $\pi=r'$ on $[\sig_{r'}+\eps,\infty)$. Let $\Di\sub\R^2$ be a deterministic countable dense set. By Corollary~\ref{C:infpath}, we can find some path $\pi'\in\Ni(\Di)$ that starts on the right of $r$ and crosses both $r$ and $r'$ before time $\sig_{r'}+\eps$. Let $\pi$ be the concatenation of $r$, $\pi'$ and $r'$. Then $\pi\in\Ni$ by Lemma~\ref{L:lefthop}, $\pi(0)=\sup(A)$, and $\pi=r'$ on $[\sig_{r'}+\eps,\infty)$. The proof in the case $\sup(A)=\infty$ is similar, where now instead of $r$ we use a sequence of paths $r^n$ started from points $(x_n,0)$ with $A\ni x_n\up\infty$.

Part~(d), finally, is trivial since in this case $\Ni=\Pi$.\qed

\subsection{Support properties of Howitt-Warren flows and quenched laws}\label{S:HWsupp}

We are now ready to prove Theorems~\ref{T:speed} and \ref{T:supp} on the support properties of Howitt-Warren processes and Theorem~\ref{T:support} on the support of quenched laws. We start by preparing for the proof of the latter. In line with notation introduced in Section~\ref{S:quensup}, we set $\Q^+_z:=\P[\pi^+_z\in\cdot\,|(\Wi_0,\Mi)]$, where $(\Wi_0,\Mi)$ is the marked reference Brownian web as in Theorem~\ref{T:HWconst} and $\pi^+_z$ denotes the right-most path in the sample Brownian web $\Wi$ started from $z$.

\bl{\bf(Support of quenched laws on the space of paths)}\label{L:Qzsup}
Conditional on a marked reference Brownian web $(\Wi_0,\Mi)$, let $(\Wi_i)_{i\geq 1}$ be an i.i.d.\ sequence of sample Brownian webs as in (\ref{sample2}), and let $\Ni_\infty$ be the generalized Brownian net associated with $(\Wi_0,\Mi)$ defined in Lemma~\ref{L:Ngen}. Then, for any deterministic $z\in\R^2$,
\be
\supp(\Q^+_z)=\Ni_\infty(z)\qquad{\rm a.s.}
\ee
\el
{\bf Proof.} Given the marked reference Brownian web $(\Wi_0,\Mi)$, we note that $\Q^+_z$ is the conditional law of the a.s.\ unique path in $\Wi_n$ starting at $z$,
for each $n\in\N$. Since $(\Wi_n(z))_{n\in\N}$ are i.i.d., and $\Wi_n(z)\sub\Ni_\infty(z)$ a.s., the inclusion $\supp(\Q^+_z)\sub\Ni_\infty(z)$ is trivial.

To prove the other inclusion, by Lemma~\ref{L:Ngen}, it suffices to show that
\be\label{Hsup}
\Hi_\Ti(\Wi_1\cup\cdots\cup\Wi_n)(z)\sub\supp(\Q^+_z)\quad{\rm a.s.}
\ee
for each $n\geq 1$. Fix $1\leq i_0,\ldots,i_m\leq n$ and $t_1,\ldots,t_m\in\Ti$ with $t_1<\cdots<t_m$, and set $t_0:=-\infty$, $t_{m+1}:=+\infty$. Let $\Wi$ be the concatenation of $\Wi_{i_0},\ldots,\Wi_{i_m}$ on the time intervals $[t_0,t_1],\ldots,[t_m,t_{m+1}]$, i.e., $\Wi$ is the set of all paths $\pi\in\Pi$ such that for each $k=0,\ldots,m$ there is a $\pi'\in\Wi_{i_k}$ with $\pi=\pi'$ on $[\sig_\pi,\infty]\cap[t_k,t_{k+1}]$. Since conditional on the marked reference web $(\Wi_0,\Mi)$, restrictions of a sample Brownian web $\Wi_i$ to disjoint space-time regions are independent, we see that the conditional distribution of $\Wi$ equals that of the $\Wi_i$'s. In particular, conditional on $(\Wi_0,\Mi)$, the a.s.\ unique path in $\Wi(z)$ is distributed with law $\Q^+_z$. Since $1\leq i_0,\ldots,i_m\leq n$ and $t_1,\ldots,t_m\in\Ti$ are arbitrary, this proves (\ref{Hsup}).\qed

\noi
{\bf Proof of Theorem~\ref{T:support}.} Without loss of generality, we may assume that $\mu$ is a probability measure. In the set-up of Lemma~\ref{L:Qzsup}, let $(Z_j)_{j\geq 1}$ be an i.i.d.\ sequence of $\R^2$-valued random variables with law $\mu$, independent of the marked reference Brownian web $(\Wi_0,\Mi)$ and sequence of sample Brownian webs $(\Wi_i)_{i\geq 1}$. Then, conditional on $(\Wi_0,\Mi)$, for each $i,j\geq 1$, the random variable $\pi^{+\,i}_{Z_j}$ has law $\int\mu(\di z)\Q^+_z$, where $\pi^{+\,i}_z$ denotes the rightmost path in $\Wi_i$ started at $z$. Since $\Wi_i\sub\Ni_\infty$ and $Z_j\in\supp(\mu)$ a.s., we see that $\pi^{+\,i}_{Z_j}\in\Ni_\infty(\supp(\mu))$ a.s.\ and hence $\supp(\int\mu(\di z)\Q^+_z)\sub\ov{\Ni_\infty(\supp(\mu))}$ a.s.

To prove the other inclusion, set $\Di:=\{Z_j:j\geq 0\}$. Since $\pi^{+\,i}_{Z_j}\in\supp(\int\mu(\di z)\Q^+_z)$ a.s.\ for each $i,j\geq 1$, and conditional on $(\Wi_0,\Mi,(Z_j)_{j\geq 1})$, the random variable $\pi^{+\,i}_{Z_j}$ has law $\Q^+_{Z_j}$, we see that $\supp(\Q^+_{Z_j})\sub\supp(\int\mu(\di z)\Q^+_z)$ a.s.\ for each $j\geq 1$. On the other hand, by Lemma~\ref{L:Qzsup}, we have $\supp(\Q^+_{Z_j})=\Ni_\infty(Z_j)$. Therefore $\Ni_\infty(\Di)=\bigcup_j \Ni_\infty(Z_j)\sub\supp(\int\mu(\di z)\Q^+_z)$ a.s., which by Lemma~\ref{L:NgenA} implies that $\Ni_\infty(\supp(\mu))\sub\supp(\int\mu(\di z)\Q^+_z)$ a.s. (Recall that $\supp(\mu)$ is the support of $\mu$ in $\R^2$, not $\Rc$, which is why Lemma~\ref{L:NgenA} can be applied here.)\qed

\bp{\bf(Support of Howitt-Warren process)}\label{P:rhosupp}
Let $(\Wi_0,\Mi)$ be a marked reference Brownian web, let $\Ni_\infty$ be the generalized Brownian net associated with $(\Wi_0,\Mi)$ defined in Lemma~\ref{L:Ngen}, and for each closed $A\sub\R$, let $(\xi^A_t)_{t\geq 0}$ be the branching-coalescing point set associated with $\Ni_\infty$ defined in (\ref{braco2}). Then, for each deterministic finite measure $\rho_0$ on $\R$, almost surely
\be
\supp(\rho_t)=\xi^{\supp(\rho_0)}_t\qquad(t\geq 0),
\ee
where $(\rho_t)_{t\geq 0}$ is the Howitt-Warren process defined as in (\ref{rho}) with $K_{0,t}=K^+_{0,t}$ or $K^\up_{0,t}$, where $(K^+_{s,t})_{s\leq t}$ and $(K^\up_{s,t})_{s\leq t}$ are the versions of the Howitt-Warren flow defined in Theorem~\ref{T:HWconst}.
\ep
{\bf Proof.} In line with notation introduced in (\ref{SigT}), let $\Pi(\Sig_0):=\{\pi\in\Pi:\sig_\pi=0\}$. For each $t\geq 0$, define a continuous map $\psi_t:\Pi(\Sig_0)\to[-\infty,\infty]$ by $\psi_t(\pi):=\pi(t)$. Recall that if $E,F$ are Polish spaces, $\mu$ is a finite measure on $E$, $f:E\to F$ is a continuous function, and $f(\mu)$ denotes the image of $\mu$ under $f$, then $\supp\big(f(\mu)\big)=\ov{f\big(\supp(\mu)\big)}$. Then by Theorem~\ref{T:support},
\be\ba{l}
\supp(\rho_t)
=\supp\big(\psi_t\big(\int\rho_0(\di x)\Q^+_{(x,0)}\big)\big)
=\ov{\psi_t\big(\supp\big(\int\rho_0(\di x)\Q^+_{(x,0)}\big)\big)}\\[5pt]
\quad=\ov{\psi_t\big(\Ni_\infty(\supp(\rho_0)\times\{0\})\big)}
=\xi^{\supp(\rho_0)}_t\qquad(t\geq 0).
\ec
By the remarks above Theorem~\ref{T:support}, replacing $\Q^+_{(x,0)}$ with $\Q^\up_{(x,0)}$ makes no difference.\qed

\detail{The proof that $\supp\big(f(\mu)\big)=\ov{f\big(\supp(\mu)\big)}$ is easy if you use that $\supp(\mu)$ is the set of all $x$ such that each open ball around $x$ has positive measure.}

\noi
{\bf Proof of Theorems~\ref{T:speed} and \ref{T:supp}.} Immediate from Lemma~\ref{L:edge} and Proposition~\ref{P:rhosupp}.\qed

\section{Atomic or non-atomic}\label{S:atom}

In this section, we use our construction of the Howitt-Warren flows in
Theorems~\ref{T:HWconst} and \ref{T:HWconst2} to prove Theorem~\ref{T:atom} on
the atomicness/non-atomicness of the Howitt-Warren processes. Parts (a), (b)
and (c) are proved in Sections~\ref{S:detatom}, \ref{S:nonatom} and
\ref{S:erosion}, respectively.

\subsection{Atomicness at deterministic times}\label{S:detatom}

To prove Theorem~\ref{T:atom}~(a) on the atomicness of any Howitt-Warren process $(\rho_t)_{t\geq 0}$ at deterministic times, we need the following lemma.

\bl{\bf(Coincidence of points entered by a path)}\label{L:webatomic}
Let $(\Wi,\Wi')$ be a pair of sticky Brownian webs and for $t\in\R$, let $I(t):=\{\pi(t):\pi\in\Wi,\ \sig_\pi<t\}$ and let $I'(t)$ be defined similarly with $\Wi$ replaced by $\Wi'$. Then for each deterministic $t\in\R$, a.s.\ $I(t)=I'(t)$.

\el
{\bf Proof.}
By symmetry, it suffices to show that $I(t)\sub I'(t)$.
As in the proof of Theorem~\ref{T:webmod} (in Section~\ref{S:markproof}), we may without loss of generality assume that $\Wi$ and $\Wi'$ are embedded in a Brownian net $\Ni$ with associated left-right Brownian web $(\Wl,\Wr)$ and set of separation points $S$ in such a way that
\bc
\dis\Wi&=&\dis\big\{\pi\in\Ni:\sign_\pi(z)=\al_z
\ \forall z\in S\mbox{ s.t.\ }\pi\mbox{ enters }z\big\},\\[5pt]
\dis\Wi'&=&\dis\big\{\pi\in\Ni:\sign_\pi(z)=-\al_z
\ \forall z\in S\mbox{ s.t.\ }\pi\mbox{ enters }z\big\},
\ec
where conditional on $\Ni$, the $(\al_z)_{z\in S}$ are i.i.d.\ $\{-1,+1\}$-valued random variables. Then $I(t)\sub\{\pi(t):\pi\in\Ni:\sig_\pi<t\}$, and hence by the structure of special points in the Brownian net at deterministic times (see Proposition~\ref{P:detspec}), for each $z=(x,t)$ with $x\in I(t)$, there exist $l\in\Wl$ and $r\in\Wr$ such that $l\sim^z_{\rm in}r$. Since paths in $\Ni$ are contained between left-most and right-most paths \cite[Prop.~1.8]{SS08}, any path in $\Ni$ started between $l$ and $r$ must enter $z$. In particular,
it follows that there are paths in $\Wi'$ that enter $z$.\qed

\noi
{\bf Proof of Theorem~\ref{T:atom}~(a).} Let $(\Wi_0,\Mi,\Wi)$ be a reference
web, a set of marked points, and a sample web as
constructed in Theorem~\ref{T:HWconst}. For $t\in\R$, let
$I_0(t):=\{\pi(t):\pi\in\Wi_0,\ \sig_\pi<t\}$. By Theorems~\ref{T:HWconst} and
\ref{T:webmod}, $(\Wi_0,\Wi)$ is a pair of sticky Brownian webs, while by
Theorem~\ref{T:HWconst}, we have the representation
\be\label{rhoconst}
\rho_t=\int\rho_0(\di x)
\P\big[\pi^\up_{(x,0)}(t)\in\cdot\,\big|\,(\Wi_0,\Mi)\big].
\ee
Lemma \ref{L:webatomic} implies that a.s.\ $\{\pi^\up_{(x,0)}(t)\}_{x\in\R}\subset I_0(t)$, which is a countable
set. Therefore $\rho_t$ must be atomic with atoms in $I_0(t)$.\qed

The proof of Theorem~\ref{T:atom}~(a) shows that at deterministic times, $\rho_t$ is concentrated on the countable set $I_0(t)$ of points where there is an incoming path of the reference Brownian web. The following lemma, which will be needed in the proof of Theorem~\ref{T:erosion} below, tells us which points in $I_0(t)$ carry positive mass.

\bl{\bf(Position of atoms at deterministic times)}\label{L:posatom}
Let $\rho_0$ be a deterministic finite measure on $\R$, let $(\rho_t)_{t\geq 0}$ be a Howitt-Warren process constructed as in (\ref{rhoconst}) and let $I_0(t):=\{\pi(t):\pi\in\Wi_0,\ \sig_\pi<t\}$. Then for each deterministic $t>0$, one has
\be\label{posatom}
\big\{x\in\R:\rho_t(\{x\})>0\big\}
=I_0(t)\cap\supp(\rho_t)\qquad{\rm a.s.}
\ee
\el
{\bf Proof.} The inclusion $\subset$ follows from our proof of Theorem~\ref{T:atom}~(a). If both the left and right speeds of the Howitt-Warren flow are finite, then by Proposition~\ref{P:braco}~(d) and Theorem~\ref{T:supp}~(a), $\supp(\rho_t)$ consists of isolated points, hence $\big\{x\in\R:\rho_t(\{x\})>0\big\}=\supp(\rho_t)$, which by the inclusion $\subset$ implies that we have in fact equality in (\ref{posatom}).

If at least one of the speeds of the Howitt-Warren flow is infinite, then by Theorem~\ref{T:supp}~(b) and (c), either $\supp(\rho_t)=\R$ or $\supp(\rho_t)$ is a halfline. To prove the inclusion $\supset$ in (\ref{posatom}) in this case, let us fix a typical realization of the reference web and set of marked points $(\Wi_0, \Mi)$, which together with $\rho_0$ determines $(\rho_s)_{s\geq 0}$ by (\ref{rhoconst}). Let us take $x$ to be any point in $I_0(t)\cap\supp(\rho_t)$. We will treat the cases that $x$ lies in the interior or on the boundary of $\supp(\rho_t)$ separately.


Assume for the moment that $x$ lies in the interior of $\supp(\rho_t)$. We fix a regular version of the conditional law $\P[(\Wi,\hat\Wi)\in\cdot\,|\,(\Wi_0,\Mi)]$ of the sample Brownian web $\Wi$ and its dual $\hat\Wi$ given the marked reference web $(\Wi_0,\Mi)$. As discussed in the first remark below Theorem~\ref{T:webmod}, for a.e.\ $\om$ in our underlying probability space, the sets of paths $(\Wi,\hat\Wi)$ have the same a.s.\ properties under this conditional law as a double Brownian web and we can define probabilities for special paths such as $\pi^+_z,\pi^\up_z$ simultaneously for all $z\in\R^2$. Let $\hat\pi^-_z$ and $\hat\pi^+_z$ denote the left-most and right-most elements of $\hat\Wi(z)$. For each $y\in\R$, let $\tau_y:=\sup\{s<t:\hat\pi^-_{(y,t)}(s)=\hat\pi^+_{(y,t)}(s)\}$, and let $I(t):=\{\pi(t):\pi\in\Wi,\ \sig_\pi<t\}$. By the structure of special points of the Brownian web (Proposition~\ref{P:classweb}), a.s.\ $\tau_y<t$ for all $y\in I(t)$, which equals $I_0(t)$ by Lemma~\ref{L:webatomic}. In particular, for all $x\in I_0(t)$ that lie in the interior of $\supp(\rho_t)$ and for all $t_n\up t$, we must have
\be
\lim_{n\to\infty}\P\big[\tau_x<t_n\,\big|\,(\Wi_0,\Mi)\big]=1.
\ee
Since we have assumed that either the left or right speed of the Howitt-Warren flow is infinite, $\supp(\rho_s)$ is either a halfline for all $s\geq 0$ or the whole line for all $s\geq 0$. In the halfline case, Prop.~\ref{P:rhosupp} and Lemma~\ref{L:edge} further imply that the boundary of $\supp(\rho_t)$ is continuous in $t$. Therefore if $x$ lies in the interior of $\supp(\rho_t)$, then we can choose $t_n$ depending on $x$ such that
\be
\P\big[\tau<t_n,\ \big(\hat\pi^-_{(x,t)}(t_n),\hat\pi^+_{(x,t)}(t_n)\big)
\sub\supp(\rho_{t_n})\,\big|\,(\Wi_0,\Mi)\big]>0.
\ee
Note that every path in $\Wi$ started from $(\hat\pi_-(t_n),\hat\pi_+(t_n))\times\{t_n\}$ must pass through $x$. Therefore, by Prop.~\ref{P:regul}~(d),
\be\ba{l}
\dis\rho_t(\{x\})=\int_R\rho_{t_n}(\di y)K^\up_{t_n,t}(y,\{x\})
=\E\Big[\int_\R\rho_{t_n}(\di y)1_{\txt\{\pi^\up_{(y,t_n)}(t)=x\}}\,\Big|\,(\Wi_0,\Mi)\Big]\\[5pt]
\dis\quad\geq\E\big[\rho_{t_n}\big((\hat\pi^-_{(x,t)}(t_n),\hat\pi^+_{(x,t)}(t_n))\big)\,\big|\,(\Wi_0,\Mi)\big]>0.
\ec

If $x$ lies on the boundary of $\supp(\rho_t)$, then by symmetry, it suffices to consider the case that $\bet_-=-\infty$, $\bet_+<\infty$ and $\supp(\rho_t)=(-\infty,x]$. In this case, as in the proof of Lemma~\ref{L:Ngen}, we can without loss of generality assume that $(\Wi_0,\Mi,\Wi)$ are constructed as in Theorem~\ref{T:HWconst} with $\nu_{\rm l}=0$. In this case, $\Wi_0$ is the right Brownian web associated with the Brownian halfnet $\Ni_\infty$ from Lemma~\ref{L:Ngen}, hence by Proposition~\ref{P:rhosupp}, $\supp(\rho_t)=(-\infty,\pi^0(t)]$ $(t\geq 0)$ where $\pi^0$ is the a.s.\ unique path in $\Wi_0$ starting at time zero from $\sup(\supp(\rho_0))$. Set $\hat\pi^-:=\hat\pi^-_{(\pi_0(t),t)}$. Since $(\Wi,\Wi_0)$ is distributed as a left-right Brownian web when averaged over the law of $(\Wi_0,\Mi)$, by Proposition~\ref{P:detspec}, any path $\pi\in\Wi$ entering $x=\pi^0(t)$ must satisfy $\pi(s)=\pi^0(s)$ for some $\sig_\pi<s<t$, and hence, since paths in $\Wi$ cannot cross paths in $\Wi_0$ from left to right, $\tau':=\sup\{s<t:\pi^0(s)<\hat\pi^-(s)\}<t$ a.s. It follows that for deterministic $t_n\up t$,
\be
\lim_{n\to\infty}\P\big[\tau'<t_n\,\big|\,(\Wi_0,\Mi)\big]=1.
\ee
In particular, we can choose $n$ (depending on $(\Wi_0,\Mi)$) large enough such that this probability is positive. Again by Proposition~\ref{P:detspec} applied to $(\Wi,\Wi_0)$, we have $\hat\pi^-(t_n)<\pi^0(t_n)$ a.s.\ on the event $\tau'<t_n$. Now the argument proceeds as in the case when $x$ lies in the interior of $\supp(\rho_t)$, with the interval $(\hat\pi^-_{(x,t)}(t_n),\hat\pi^+_{(x,t)}(t_n))$ replaced by $(\hat\pi^-(t_n),\pi^0(t_n))$.\qed

\subsection{Non-atomicness at random times for non-erosion flows}\label{S:nonatom}

The proof of Theorem~\ref{T:atom}~(b) is a bit involved. When $\nu(0,1)>0$ and $\beta_+-\beta_-<\infty$, the Howitt-Warren process $(\rho_t)_{t\geq 0}$
can be constructed from a Brownain net as in Theorem~\ref{T:HWconst2}. In particular, at each separation point of $\Ni$, mass splits binarily. We can
then use the fact that separation points are dense in space and time to split each atom in $\rho_0$ infinitely often to reach a random time $t$ when
$\rho_t$ contains no atoms. This is closely related to, and in fact implies, the fact that the branching-coalescing point set $\xi_t$, whose properties are listed in
Prop.~\ref{P:braco}, admit random times when $\xi_t$ has no isolated points. When $\beta_+-\beta_-=\infty$, it turns out that as long as $\nu(0,1)>0$,
the picture of binary splitting of mass is still valid, although more work will be required. To avoid repetition, we will prove
Theorem~\ref{T:atom}~(b) directly under the general assumption $\nu(0,1)>0$, which we assume from now on. Who wants to see a similar, but simpler proof should read the proof that there exist random times when $\xi_t$ has no isolated points in \cite[Prop.~3.14]{SSS09}.

Let $(\Wi_0, \Mi, \Wi)$ be a reference web, a set of marked points, and a sample web with quenched law $\Q$ as in Theorem~\ref{T:HWconst}, where we
take $\nu_{\rm l}=\nu_{\rm r}=\nu$. We will use the version of Howitt-Warren flow
$(K^\up_{s,t})_{s<t}$ in Theorem~\ref{T:HWconst} to represent the Howitt-Warren process
\be\label{rhoflow}
\rho_t = \int\! \rho_s(\di x) K^\up_{s,t}(x,\cdot\,) =\int\! \rho_s(\di x) \P[\pi^\up_{(x,s)}(t)\in\cdot\,|(\Wi_0,\Mi)]= \int\! \rho_s(\di x) \Q[\pi^\up_{(x,s)}(t)\in\cdot\,],
\ee
where $K^\up$ satisfies Prop.~\ref{P:conpath}~(i)', which implies that (\ref{rhoflow}) holds a.s.\ for all $s<t$.

The heuristics outlined before require two ingredients. First we need to establish the existence of bottlenecks, i.e., points in space-time where
most of the mass in $\rho_0$ must enter. Secondly, at such bottlenecks, mass is split binarily. When the flow can be embedded in a Brownian net,
separation points provide the bottlenecks. In general, choose $\eta \in (0,1)$ such that $\nu([\eta, 1-\eta])>0$. Then
\be
A_{\rm l}^\eta := \{z: (z, \omega_z) \in \Mi, \omega_z\in [\eta, 1-\eta], z \mbox{ is of type } (1,2)_{\rm l} \mbox{ in } \Wi_0\}
\ee
will provide the bottlenecks we need, and we will show that an atom of size $1$ entering any $z\in A_{\rm l}^\eta$ will be split
into atoms of size no larger than $1-\eta/4$.

\bl\label{L:bneck}{\bf (Bottlenecks)} Let $\rho_0=\delta_0$, the delta mass at $0$. Then a.s.\ w.r.t.\ $(\Wi_0, \Mi)$, for any $\delta>0$,
we can find a point $z=(x,u)\in A_{\rm l}^\eta$ with $u\in (0,\delta)$ such that $\rho_u(\{x\}) \geq 1-\delta$.
\el
{\bf Proof.} By Theorems~\ref{T:HWconst} and \ref{T:webmod}, $(\Wi_0,\Wi)$ is a pair of sticky Brownian webs with drift $\beta$ and coupling parameter
$\kappa=2\nu([0,1])$. In particular, by Propositions~\ref{P:niden} and \ref{P:exchange}, if $\pi^0$ resp.\ $\pi$ denotes the a.s.\ unique path in $\Wi_0$ resp.\ $\Wi$ starting from the origin, then
$(\pi^0, \pi)$ solves the Howitt-Warren martingale problem, which is equivalent to (\ref{MP1b}) and (\ref{thetacouple}). Without loss of generality, we will assume $\beta=0$.

Our basic strategy is to first explore backward in time and find a $z=(x,u)\in A_{\rm l}^\eta$ that uses only information about $(\Wi_0,\Mi)$ above time $u$, and then get a lower bound on $\rho_u(\{x\})$ by using the fact that the increments of $(\pi^0,\pi)$ on the time interval $[0,u]$ are independent of $(\Wi_0, \Mi)$ restricted to the time interval $(u,\infty)$. The actual proof will consist of a direct analysis of $(\pi^0, \pi)$, using additional information from $(\Wi_0, \Mi)$ when necessary.

Fix any $\eps>0$. For each $0\leq s\leq t$, consider the following subset of $(1,2)_{\rm l}$ points along $\pi^0$:
$$
I^{\eps}_{[s,t]}:=\{(y,v) : s\leq v\leq t, \, y=\pi^0(v)=\hat\pi^0_{(\pi^0(u)+\eps,u)}(v) \mbox{ for some } u\in (v,t]  \},
$$
where $\hat\pi^0_{(\pi^0(u)+\eps,u)}$ is any path in the dual reference web $\hat\Wi_0$ starting from $(\pi^0(u)+\eps, u)$. It was shown in
\cite[Lemma 7.2]{NRS10} that $I^{\eps}_{[s,t]}$ contains exactly the points of intersection on the time interval $[s,t]$ between $\pi^0$ and a
backward Brownian motion $\hat \pi^0$ starting at $(\pi^0(t)+\eps, t)$, which is Skorohod reflected between $(\pi^0(r))_{0\leq r\leq t}$ and
$(\pi^0(r)+\eps)_{0\leq r\leq t}$. In particular, after reversing time and centering, the triple
$$
(\vec Z_r)_{0\leq r\leq t}:=\big(\pi^0(t-r)-\pi^0(t),\, \hat\pi^0(t-r)-\pi^0(t),\, \ell_{\rm l}(I^\eps_{[t-r,t]})\big)_{0\leq r\leq t}
$$
is a strong Markov process, where $\ell_{\rm l}(I^{\eps}_{[t-r,t]})$ is the intersection local time measure of the set $I^{\eps}_{[t-r,t]}$ and is a finite continuous increasing process. Also note that $(\pi^0(t-r)-\pi^0(r))_{0\leq r\leq t}$ is distributed as a standard Brownian motion.

By the definition of $A_{\rm l}^\eta$ and our construction of $\Mi$ in Theorem~\ref{T:HWconst}, conditional on $\Wi_0$, $A^\eta_{\rm l} \cap I^{\eps}_{[0,t]}$ is a Poisson subset of $I^{\eps}_{[0,t]}$ with intensity measure $21_{\{z\in I^{\eps}_{[0,t]}\}}\ell_{\rm l}({\rm d}z) \int_{[\eta, 1-\eta]} q^{-1}\nu({\rm d}q)$. In particular, $A^\eta_{\rm l} \cap I^{\eps}_{[0,t]}$ is a.s.\ a finite set. If $A^\eta_{\rm l} \cap I^{\eps}_{[0,t]}\neq\emptyset$, then let
$$
\tau^{\eps,t} := \sup\{r\in [0,t]: (\pi^0(r),r)\in A^\eta_{\rm l} \cap I^{\eps}_{[0,t]}\},
$$
and denote the corresponding point in $A^\eta_{\rm l}\cap I^{\eps}_{[0,t]}$ by $z^{\eps,t}=(x^{\eps, t}, \tau^{\eps, t})$. If $A^\eta_{\rm l} \cap I^{\eps}_{[0,t]}=\emptyset$, then we set $\tau^{\eps,t}=0$ and $z^{\eps, t}=(0,0)$. Therefore
$$
\P[\tau^{\eps, t} >0 | \Wi_0] = 1- e^{-c_{\eta,\nu} \ell_{\rm l}(I^{\eps}_{[0,t]})}, \quad \mbox{where}\ \ c_{\eta, \nu}:=2\int_{[\eta, 1-\eta]} q^{-1}\nu({\rm d}q).
$$
Now observe that we can construct $\tau^{\eps,t}$ by setting $t-\tau^{\eps,t}:=t\wedge \inf\{r\geq 0: c_{\eta,\nu}\ell_{\rm l}(I^{\eps}_{[t-r,t]})\geq L\}$,
where $L$ is an independent mean $1$ exponential random variable. In particular, conditional on $L$, $t-\tau^{\eps, t}$ is a stopping time for the process $\vec Z$. Therefore, conditional on $\tau^{\eps, t}$ and $(\vec Z_s)_{0\leq s\leq t-\tau^{\eps,t}}$, the law of $(\pi^0(s), \pi(s))_{0\leq s\leq \tau^{\eps, t}}$ remains the
same as before because conditional on $\tau^{\eps,t}$, $(\pi^0_s)_{0\leq s\leq \tau^{\eps, t}}$ is a Brownian motion independent of $(\vec Z_s)_{0\leq s\leq t-\tau^{\eps,t}}$ , and $\pi$ is constructed by switching among the truncated paths $\{(\gamma(s))_{s_0\leq s\leq \tau^{\eps, t}} : \gamma \in \Wi_0, \sigma_\gamma =s_0 <\tau^{\eps,t}\}$ using independent Poisson processes, which are all independent of $(\vec Z_s)_{0\leq s\leq t-\tau^{\eps,t}}$. Therefore
\begin{eqnarray*}
\P\big[\tau^{\eps, t}>0, \pi^0(\tau^{\eps,t}) = \pi(\tau^{\eps,t})\big]
&=& \E\big[ 1_{\{\tau^{\eps,t}>0\}}\,\P[\pi^0(\tau^{\eps,t}) = \pi(\tau^{\eps,t})\, |\, \tau^{\eps,t}, (\vec Z_s)_{0\leq s\leq t-\tau^{\eps,t}}] \big] \\
&=& \E[1_{\{\tau^{\eps,t}>0\}}\, \phi(\tau^{\eps, t})],
\end{eqnarray*}
where $\phi(r):=\P[\pi^0(r)=\pi(r)]$. Note that $|\pi^0(s)-\pi(s)|$ is distributed as a Brownian motion starting at $0$ and sticky reflected at $0$. It was shown in \cite[Lemma A.2]{SSS09} that $\phi(r)=\P[|\pi^0(r)-\pi(r)|=0]$ is monotone in $r$ and increases to $1$ as $r\downarrow 0$. Therefore
\be\label{bneck}
\P\big[\tau^{\eps, t}>0, \pi^0(\tau^{\eps,t}) = \pi(\tau^{\eps,t})\big] \geq \phi(t) (1- \E[e^{-c_{\eta,\nu} \ell_{\rm l}(I^\eps_{[0,t]})}]).
\ee
Note that almost surely,
\[
\ell_{\rm l}(I^\eps_{[0,t]})\astoo{\eps}{0}
\ell_{\rm l}(\{(\pi^0(s),s):s\in[0,t]\}),
\]
where the right-hand side is easily seen to be a.s.\ infinite using the fact that $\ell_{\rm l}(I^{\eps}_{[0,t]})\geq \sum_{i=1}^{\lfloor t\eps^{-2}\rfloor}\ell_{\rm l}(I^{\eps}_{[(i-1)\eps^2, i\eps^2]})$, and $\eps^{-1}\ell_{\rm l}(I^{\eps}_{[(i-1)\eps^2, i\eps^2]})$, $i\in\N$, are i.i.d.\ with the same distribution as
$\ell_{\rm l}(I^1_{[0,1]})$. Therefore (\ref{bneck}) implies
\[
\lim_{t\downarrow 0} \lim_{\eps\downarrow 0} \P\big[\tau^{\eps, t}>0, \pi^0(\tau^{\eps,t}) = \pi(\tau^{\eps,t})\big] =1.
\]
Since
\[
\begin{aligned}
\P[\tau^{\eps, t}>0,\ \pi^0(\tau^{\eps,t})=\pi(\tau^{\eps,t})] & = \E[1_{\{\tau^{\eps,t}>0\}}
\P[\pi^0(\tau^{\eps,t})=\pi(\tau^{\eps,t})|(\Wi_0,\Mi)]] \\
& = \E[1_{\{\tau^{\eps,t}>0\}}\rho_{\tau^{\eps,t}}(\{\pi^0(\tau^{\eps,t})\})],
\end{aligned}
\]
and $\rho_{\tau^{\eps,t}}(\{\pi^0(\tau^{\eps,t})\})\in [0,1]$, it follows that for any $\delta \in (0,1)$, we must have
\[
\lim_{t\downarrow 0} \lim_{\eps\downarrow 0} \P\big[\tau^{\eps, t}>0, \rho_{\tau^{\eps,t}}(\{\pi^0(\tau^{\eps,t})\})\geq 1-\delta\big] =1.
\]
Therefore for any sequence $(\eps_n, t_n)\to (0,0)$, a.s.\ with respect to $(\Wi_0,\Mi)$, the event $\{\tau^{\eps_n,t_n}>0,\ \rho_{\tau^{\eps_n,t_n}}(\{\pi^0(\tau^{\eps_n,t_n})\})\geq 1-\de\}$ must occur infinitely often by Borel-Cantelli. The lemma then follows.
\qed
\bigskip

Next we show that mass entering any $z\in A^\eta_{\rm l}$ must be split into smaller atoms.
\bl\label{L:split}{\bf (Splitting of mass)} Almost surely w.r.t.\ $(\Wi_0,\Mi)$, for each $z=(x,u)\in A^\eta_{\rm l}$ and for any $\eps, \delta>0$,
we can find $h>0$ such that if $(\rho_t)_{t\geq u}$ is defined as in (\ref{rhoflow}) with $\rho_u=\delta_x$, the delta mass at $x$, then for all $t\in (u, u+h)$,
we have $\rho_t([x-\eps, x+\eps])\geq 1-\delta$ and $|\rho_t|_\infty:=\sup_{y\in\R} \rho_t(\{y\}) \leq 1-\eta/4$.
\el
{\bf Proof.} As we will see, the bound $\rho_t([x-\eps, x+\eps])\geq 1-\delta$ will follow from the weak continuity of $\rho_t$ in $t\geq u$. To show the splitting of mass, $|\rho_t|_\infty<1-\eta/4$, the key is to show that for a typical realization of $(\Wi_0, \Mi)$, $z$ is of type $(1,2)_{\rm r}$ (resp.\ $(1,2)_{\rm l}$) in the sample web $\Wi$ with probability $\omega_z$ (resp.\ $1-\omega_z$). To prove this, we will embed the reference and sample webs $\Wi_0$ and $\Wi$ in a Brownian net, whose separation points contain all of $A^\eta_{\rm l}$.

Recall from Theorem~\ref{T:HWconst} that $\Wi$ is constructed from $(\Wi_0,\Mi)$ by switching paths in $\Wi_0$ at a countable collection of $(1,2)$ points
$A\cup B$, where $A$ is a random subset of the marked points in $\Mi$, with $z$ included in $A$ independently for each $(z,\omega_z)\in\Mi$ with probability
$\omega_z$ if $z$ is of type $(1,2)_{\rm l}$ in $\Wi_0$, and with probability $1-\omega_z$ if $z$ is of type $(1,2)_{\rm r}$; while $B$ is an independent Poisson
point set on $(1,2)$ points with intensity measure $2\nu(\{0\})\ell_{\rm l}+2\nu(\{1\})\ell_{\rm r}$. Note that conditional on $\Wi_0$,
$A\cup B\cup A^\eta_{\rm l}$ is a Poisson point process on $(1,2)$ points of $\Wi_0$ with intensity measure
$$
2\nu([0,1])\ell_{\rm r}({\rm d}z) + 2\Big(\nu([0,\eta))+\int_{[\eta, 1-\eta]}q^{-1}\nu({\rm d}q) + \nu((1-\eta,1])\Big)\ell_{\rm l}({\rm d}z).
$$
By Theorem~\ref{T:marknet}, if we allow hopping at all points in $A\cup B\cup A^\eta_{\rm l}$, then we obtain a Brownian net $\Ni$ with left and right speeds
$$
\beta_-=\beta - 2\nu([0,1]) \quad \mbox{and} \quad \beta_+= \beta + 2\Big(\nu([0,\eta))+\int_{[\eta, 1-\eta]}q^{-1}\nu({\rm d}q) + \nu((1-\eta,1])\Big).
$$
Clearly $\Wi_0$ and $\Wi$ are subsets of $\Ni$, and all points in $A^\eta_{\rm l}$ are separation points of $\Ni$. In particular,
a.s.\ each point in $A^\eta_{\rm l}$ is of type $(1,2)$ in $\Wi$. Given $z=(x,u)\in A^\eta_{\rm l}$, let $\pi^-_z$ resp.\ $\pi^+_z$ denote
the left resp.\ right of the two outgoing paths in $\Wi$ at $z$, and recall from after Prop.~\ref{P:classweb} that $\pi^\up_z$ picks from $\pi^+_z$
and $\pi^-_z$ the natural continuation of any paths in $\Wi$ entering $z$. By (\ref{rhoflow}), given $\rho_u=\delta_x$, we have
$\rho_t(\cdot) = \P[\pi^\up_z(t) \in \cdot\,|(\Wi_0,\Mi)]$ for all $t>u$. Since $\pi_z^\up$ is a.s.\ a continuous path starting at $z$, for
any $\eps,\delta>0$, we can choose $h>0$ sufficiently small such that $\rho_t([x-\eps,x+\eps])=\P[\pi^\up_z(t)\in [x-\eps, x+\eps]|(\Wi_0,\Mi)]\geq 1-\delta$ for all
$t\in [u,u+h]$. This establishes the first part of the lemma.

Since $z\in A^\eta_{\rm l}$, we have $(z, \omega_z)\in \Mi$ for some $\omega_z\in [\eta, 1-\eta]$. By our construction of $\Wi$ in Theorem~\ref{T:HWconst}
and the natural coupling between sticky Brownian webs and Brownian nets given in Lemma~\ref{L:stickynet}, we note that conditional on $(\Wi_0,\Mi)$ and $A\cup B\backslash\{z\}$, $\Wi$ is uniquely determined except for the orientation of paths in $\Wi$ entering $z$, which is then resolved by an
independent random variable $\alpha_z$ with $\P[\alpha_z=1]=\omega_z$ and $\P[\alpha_z=-1]=1-\omega_z$. More precisely, we set ${\rm sign}_{\Wi}(z)=\alpha_z$,
and $\pi^\up_z=\pi^+_z$ when $\alpha_z=1$ and $\pi^\up_z=\pi^-_z$ when $\alpha_z=-1$. Therefore with $\rho_u=\delta_x$, we have for all $t>u$,
\be\label{split1}
\rho_t(\cdot) = \omega_z\P[\pi^+_z(t) \in \cdot|(\Wi_0,\Mi)] + (1-\omega_z)\P[\pi^-_z(t)\in \cdot\,|(\Wi_0,\Mi)].
\ee
Almost surely, $\pi^-_z < \pi^+_z$ on $(u, u+h)$ for some $h>0$. Therefore for $h>0$ sufficiently small,
\be\label{split2}
\P\big[\pi^-_z(t)<\pi^+_z(t)\ \mbox{ for all } t\in (u, u+h) \,\big|\, (\Wi_0,\Mi)\big]>\frac{1}{2}.
\ee
We claim that this implies $\sup_{y\in\R} \rho_t(\{y\}) \leq 1-\eta/4$ for all $t\in (u,u+h)$. Otherwise if $\rho_t(\{y\})>1-\eta/4$ for some $t\in (u,u+h)$
and $y\in\R$, then by (\ref{split1}), we must have
$$
\P[\pi^+_z(t) \neq y |(\Wi_0,\Mi)] \leq  \frac{\eta}{4\omega_z} \leq \frac{1}{4} \qquad \mbox{and} \qquad
\P[\pi^-_z(t) \neq y |(\Wi_0,\Mi)] \leq  \frac{\eta}{4(1-\omega_z)} \leq \frac{1}{4},
$$
since $\omega_z\in [\eta, 1-\eta]$. This implies that
$$
\P\big[\pi^-_z(t)<\pi^+_z(t)\,\big|\, (\Wi_0,\Mi)\big] \leq \P[\pi^+_z(t) \neq y |(\Wi_0,\Mi)] + \P[\pi^-_z(t) \neq y |(\Wi_0,\Mi)]\leq \frac{1}{2},
$$
contradicting (\ref{split2}). This completes the proof of Lemma \ref{L:split}.
\qed
\medskip

Lemmas \ref{L:bneck} and \ref{L:split} immediately imply the following.
\bl\label{L:masssplit2}
Let $\rho_0 = \delta_x$ for some $x\in\R$, and let $\eta\in (0,1)$ satisfy $\nu([\eta, 1-\eta])>0$. Then for any
$\eps,\delta, h>0$, a.s.\ there exists $(u,v)\subset (0,h)$ such that $\rho_t([x-\eps, x+\eps])\geq 1-\delta$ for all $t\in [0,v]$
and $|\rho_t|_\infty:= \sup_{y\in\R} \rho_t(\{y\}) \leq 1-\eta/5$ for all $t\in [u,v]$.
\el

Now we can prove the existence of random times when $\rho_t$ is non-atomic, if $\nu(0,1)>0$.
\bigskip

\noindent{\bf Proof of Theorem~\ref{T:atom}~(b).} It suffices to show that for each interval $(u,v)$ with $0<u<v\in\Q$, a.s.\ there exists
$t\in (u,v)$ such that $\rho_t$ has no atoms. Denote $\lambda:= \sup_{y\in\R}\rho_u(\{y\})$, and let $\rho_t$ be
defined using $(K^\up_{s,t})_{s<t}$ as in (\ref{rhoflow}).

By Theorem \ref{T:atom}~(a), which we have already established, $\rho_u$ is a.s.\ atomic. Since $\rho_0$ is assumed to be a finite measure
(for infinite $\rho_0$, see the remark after Theorem~\ref{T:infmass}), for any $\eps>0$, we can find a finite set
of atoms of $\rho_u$ at $\{x_1, \cdots, x_k\}$ with $\lambda=\rho_u(\{x_1\})\geq \rho_u(\{x_2\})\geq \cdots$ such that
$\rho_u(\R\backslash \{x_1,\cdots, x_k\})\leq \eps\lambda$. For $t\geq u$, let
$$
\rho^{\rm a}_t:=\int_{x\notin\{x_1,\cdots, x_k\}}K^\up_{u,t}(x,\cdot)\rho_u({\rm d}x), \quad \mbox{and}\quad \rho^{(i)}_t:=\rho_u(\{x_i\})K^\up_{u,t}(x_i,\cdot)\quad \mbox{for } 1\leq i\leq k.
$$
We can find $\delta_1,\cdots, \delta_k>0$ such that
$$
[x_i-\delta_i, x_i+\delta_i] \cap [x_j-\delta_j, x_j+\delta_j]=\emptyset \quad\mbox{for all }\  i\neq j.
$$
By Lemma \ref{L:masssplit2}, we can choose $(u_1, v_1)\subset (u, v)$ with $u_1, v_1\in\Q$ such that
\be\label{escmass}
\sum_{i=1}^k \rho^{(i)}_t(\R\backslash [x_i-\delta_i, x_i+\delta_i]) \leq \eps \lambda \quad \mbox{for all  } t\in [u_1, v_1],
\ee
and
$$
\sup_{y\in \R}\rho^{(1)}_t(\{y\})\leq (1-\eta/5)\lambda  \quad \mbox{ for all  } t\in [u_1, v_1].
$$
Our construction guarantees that uniformly in $t\in [u_1, v_1]$, apart from the mass from $\rho^{\rm a}_t$ and from (\ref{escmass}) (with a total mass of at most $2\eps\lambda$) which we do not attempt to control, the atom of $\rho_u$ at $x_1$ is split into atoms contained in $[x_1-\delta_1, x_1+\delta_1]$ with size no larger than $(1-\eta/5+2\eps)\lambda$, the atom of $\rho_u$ at $x_i$, for $2\leq i\leq k$, stays within $[x_i-\delta_i, x_i+\delta_i]$ and can only gain a mass of at most $2\eps\lambda$, and there is no merging and formation of new atoms with size larger than $2\eps\lambda$. By choosing $\eps>0$ sufficiently small and iterating a finite number of times (say $m$ times) the same argument, we can split all the atoms of $\rho_u$ at $\{x_1,\cdots, x_k\}$ and find $u_m, v_m\in\Q$ such that $\sup_{y\in\R}\rho_t(\{y\})\leq (1-\eta/6)\lambda$
for all $t\in [u_m,v_m]$. By repeating the whole argument above, for each $n\in\N$, we can inductively find $u^{(n)},v^{(n)}\in\Q$ with
$u<u^{(n-1)}<u^{(n)}<v^{(n)}<v^{(n-1)}<v$, such that $\sup_{y\in\R}\rho_t(\{y\})\leq 1/n$ for all $t\in [u^{(n)}, v^{(n)}]$.
Any $t\in \bigcap_{n\in\N} [u^{(n)}, v^{(n)}]\neq \emptyset$ then gives a time when $\rho_t$ contains no atoms. In fact, $\bigcap_{n\in\N} [u^{(n)}, v^{(n)}]$ contains a single point since $\rho_t$ is a.s.\ atomic at deterministic times by Theorem~\ref{T:atom}~(a).
\qed

\subsection{Atomicness at all times for erosion flows}\label{S:erosion}

In this subsection, we prove Theorem~\ref{T:atom}~(c). In fact, we will prove the following stronger result. Note that below, when we apply Theorem~\ref{T:HWconst}, we deviate from the canonical choice $\nu_{\rm l}=\nu_{\rm r}=\nu$.

\bt{\bf(Atomicness of erosion flows)}\label{T:erosion}
Let $\rho=(\rho_t)_{t\geq 0}$ be a Howitt-Warren process with drift $\bet\in\R$ and characteristic measure of the form $\nu=c_0\de_0+c_1\de_1$, with $c_0,c_1\geq 0$, started in some deterministic, finite nonzero measure $\rho_0$ on $\R$. Let $\rho$ be constructed as $\rho_t=\int\rho_0(\di x)K^\up_{0,t}(x,\,\cdot\,)$, where $(K^\up_{s,t})_{s\leq t}$ is the Howitt-Warren flow constructed as in Theorem~\ref{T:HWconst} using a reference Brownian web $\Wi_0$ with drift $\bet_0=\bet-2c_0+2c_1$, $\nu_{\rm l}=c_0\de_0$ and $\nu_{\rm r}=c_1\de_1$. Then a.s., $\rho_t$ is purely atomic at each $t>0$ and
\be\ba{l}\label{atomerosion}
\dis\big\{(x,t)\in\R^2:t>0,\ \rho_t(\{x\})>0\big\}\\[5pt]
\dis\quad=\big\{(x,t)\in\R^2:t>0,\ x\in\supp(\rho_t),\ \exists\,\pi\in\Wi_0
\mbox{ s.t.\ }\sig_\pi<t \text{ and } \pi(t)=x\big\}.
\ec
\et
{\bf Remark.} Note that if $c_0$ and $c_1$ are both strictly positive, then by Theorem~\ref{T:supp}~(c), $\supp(\rho_t)=\R$, and hence the right-hand side of (\ref{atomerosion}) is just the set
\be\label{Wi0interior}
\big\{(\pi(t),t):t>0,\ \pi\in\Wi_0,\ \sig_\pi<t\big\}.
\ee
We state as an open problem whether in this case, $\Wi_0\big|_0^\infty$, the restriction of $\Wi_0$ to the time interval $[0,\infty]$, can a.s.\ be uniquely reconstructed from $(\rho_t)_{t\geq 0}$. In fact, it seems likely that $\Wi_0\big|_0^\infty$ consists of all paths $\pi\in\Pi$ starting at $\sig_\pi\geq 0$ such that $\rho_t(\{\pi(t)\})$ is locally uniformly bounded away from zero on $(\sig_\pi,\infty)$. Note that $\Wi_0\big|_0^\infty$ cannot be reconstructed from the set in (\ref{Wi0interior}), since switching the orientation of finitely many points of type $(1,2)$ does not change this set.\med

Recall that for an erosion flow, $\Q^\up_z:=\P[\pi^\up_z\in\,\cdot\,\big|\,\Wi_0\big]$ as defined in (\ref{HWquenz}) is the quenched law of a Markov process in a random environment. Theorem~\ref{T:erosion} says that conditional on $\Wi_0$, this Markov process has the property that at all fixed times $t$ (that may depend on $\Wi_0$ but not on $\Wi$), the motion is located in the countable set $I_0(t)$ of points where there is an incoming path from $\Wi_0$. This type of behavior is reminiscent of the FIN diffusion defined in \cite{FIN02}, which is concentrated on a random countable set at each deterministic time.

Our proof of Theorem~\ref{T:erosion} is based on the following lemma, which controls the speed of mass loss of an erosion flow along a path $\pi^0\in\Wi_0$. The proof of this lemma is somewhat involved. The intuitive idea behind it is that, for erosion flows, mass must dissipate continuously, which is contrary to the case $\nu(0,1)>0$ where mass undergoes binary splitting. However, a crude estimate on the loss of mass from $\pi^0$ will show that all mass is lost instantly. Indeed, we conjecture that $(\rho_t(\{\pi^0(t)\}))_{t\geq 0}$ as a function of time has locally unbounded variation, which means that in each positive time interval an infinite amount of mass leaves and rejoins $\pi^0$. Nevertheless, formula (\ref{massholder}) below shows that the decrease of this process is, in a sense, H\"older continuous for any exponent $\ga<1/2$.

\bl\label{L:atomonpath}{\bf (Atomic mass along a reference Brownian web path)}
In the set-up of Theorem~\ref{T:erosion}, assume that $\rho_0=\delta_0$ and let $\pi^0$ be the a.s.\ unique path starting from the origin in the reference web $\Wi_0$. Then a.s.\ with respect to $\Wi_0$, for each $\gamma\in (0, 1/2)$, there exists a constant $0<C_{\gamma,\Wi_0}<\infty$ depending on $\gamma$ and $\Wi_0$, such that
\be\label{massholder}
\rho_0(\{\pi^0(0)\})- \rho_t(\{\pi^0(t)\})\leq C_{\gamma, \Wi_0} t^\gamma \qquad \mbox{for all }\ 0\leq t\leq 1.
\ee
Moreover, for any $\delta>0$, there exists a deterministic constant $0<C_{\delta,\gamma}<\infty$, such that the random constant $C_{\gamma, \Wi_0}$ satisfies
\be\label{holderconst}
\P[C_{\gamma,\Wi_0}>u] \leq \frac{C_{\delta, \gamma}}{u^\delta} \wedge 1 \qquad \mbox{for all } u>0.
\ee
\el
{\bf Proof.} To prove (\ref{massholder}), we will apply a one-sided version of Kolmogorov's moment criterion, Theorem~\ref{T:moment}, to the process $X_t:=\rho_t(\{\pi^0(t)\})$. Note that we cannot expect the gain of mass at $\pi^0$ to be continuous due to the merging of atoms. Therefore the standard
version of Kolmogorov's moment criterion is not applicable.

Let $\pi$ denote the a.s.\ unique path starting from the origin in the sample web $\Wi$. Then by construction,
$$
X_s = \Q[\pi(s)=\pi^0(s)] = \P[\pi(s)=\pi^0(s)|\Wi_0] \qquad \mbox{a.s. for all } s\geq 0.
$$
For any $0\leq s\leq t$,
\begin{eqnarray*}
X_s-X_t &=& \P[\pi(s)=\pi^0(s)|\Wi_0]-\P[\pi(t)=\pi^0(t)|\Wi_0] \\
&\leq& \P[\pi(s)=\pi^0(s), \pi(t)\neq \pi^0(t)| \Wi_0].
\end{eqnarray*}
Let $(\Wi_i)_{i\in\N}$ be i.i.d.\ copies of the sample web $\Wi$ conditional on $\Wi_0$, with $\Wi_i(0,0)=\{\pi^i\}$.
Then for any $k\in\N$,
\be\label{kthmoment}
\E\big[\big((X_s-X_t)^+\big)^k\big] \leq \P[\pi^i(s)=\pi^0(s), \pi^i(t)\neq \pi^0(t) \mbox{ for all } 1\leq i\leq k].
\ee
Since $(\Wi^i)_{1\leq i\leq k}$ are constructed from $\Wi_0$ by independent Poisson marking and switching of $(1,2)$ points of $\Wi_0$
with intensity measure $2c_0\ell_{\rm l}+2c_1\ell_{\rm r}$, if we allow paths in $\Wi_0$ to hop at all such marked points, then
by Theorem~\ref{T:marknet}, we obtain a Brownian net $\Ni=\Ni(k)$ with left and right speeds $\beta_-=\beta_-(k)=\beta_0-2kc_1$, $\beta_+=\beta_+(k)=\beta_0+2kc_0$.
By Proposition~\ref{P:mark}, $\Wi_0,\cdots, \Wi_k$ are all contained in $\Ni$. Since the Poisson marked $(1,2)$ points of $\Wi_0$ are
a.s.\ distinct for different $\Wi_i$, at each separation point $z$ of $\Ni$, we must have ${\rm sign}_z(\Wi_0)={\rm sign}_z(\Wi_i)$
for all but one $i\in\{1,\cdots,k\}$. Note that at this place, we make essential use of the fact that we have an erosion flow; in fact, this is the only place in the proof of Theorem~\ref{T:erosion} where we will use this. Let $N_{s,t}=N_{s,t}(k)$ be the number of $s,t$-relevant separation points along $\pi^0$ on the time interval $(s,t)$. Then the event in the RHS of (\ref{kthmoment}) can only occur if $N_{s,t}\geq k$.

We will next bound $\P[N_{s,t}\geq k]$ and show that
$$
\E\big[\big((X_s-X_t)^+\big)^k\big] \leq \P[N_{s,t}\geq k] \leq C_k(t-s)^{\frac{k}{2}} \qquad (0\leq s\leq t\leq 1)
$$
for some $C_k$ depending only on $k$. We can then apply the one-sided Kolmogorov's moment criterion, Theorem~\ref{T:moment}, to deduce (\ref{massholder}).

Let $(\Wl,\Wr)$ be the left-right Brownian web associated with $\Ni$, and let $(\hat\Wl,\hat\Wr)$ be the corresponding dual left-right Brownian web. Note that $(\Wl, \Wi_0)$ and $(\Wi_0,\Wr)$ each form a left-right Brownian web. For any deterministic $0\leq s<t\leq 1$, the $s,t$-relevant separation points along $\pi^0$ which are of type $(1,2)_{\rm l}$ in $\Wi_0$ can be constructed by first following $\hat r_1\in \hat\Wr(\pi^0(t),t)$, starting on the right of $\pi^0$, until the first time $\tau_1$ when
$\hat r_1$ crosses $\pi^0$ from right to left. This gives the first $s,t$-relevant separation point $(\pi^0(\tau_1),\tau_1)$ of type
$(1,2)_{\rm l}$ along $\pi^0$. We then repeat the above procedure by following $\hat r_2\in \hat\Wr(\pi^0(\tau_1),\tau_1)$ until the first time $\tau_2$ when $\hat r_2$ crosses $\pi^0$ from right to left. Iterating this procedure until time $s$ exhausts all $s,t$-relevant separation points of type $(1,2)_{\rm l}$ along $\pi^0$, the total number of which will be denoted by $N^{\rm l}_{s,t}$ and is a.s.\ finite by Proposition~\ref{P:relev}~(b).
By \cite[Lemma~2.2]{SSS09} and its proof, if we define $\hat r=\hat r_1$ on $[\tau_1,t]$, $\hat r=\hat r_2$ on $[\tau_2, \tau_1]$, \ldots, then $\hat r$ is distributed as a Brownian motion
$B^{\rm l}$ with drift $-\beta_+$ starting from $(\pi^0(t),t)$ running backward in time and Skorohod reflected to the right of $\pi^0$. More precisely, $(\hat r(t-v))_{v\geq 0}$ solves the Skorohod equation
\be\label{skoro}
\begin{aligned}
{\rm d}\hat r(t-v) & = {\rm d}B^{\rm l}(v) + {\rm d}\Delta^{\rm l}(v), \qquad & 0\leq v\leq t-s, \\
{\rm d}\hat r(t-v) & = {\rm d}B^{\rm l}(v), & t-s\leq v,
\end{aligned}
\ee
where $\Delta^{\rm l}(v)$ is an increasing process with $\int_0^{t-s} 1_{\{\hat r(t-v)\neq \pi^0(t-v)\}}{\rm d}\Delta^{\rm l}(v) =0$, and $\hat r$ is subject to the constraint $\hat r(v)\geq \pi^0(v)$ for all $0\leq v\leq t-s$. Furthermore, by \cite[Lemma~2.2]{SSS09}, conditional on $\pi^0$, the set of $s,t$-relevant separation points of type $(0,1)_{\rm l}$ along $\pi^0$ is distributed as a Poisson point process along $\pi^0$
with intensity measure $2kc_0{\rm d}\Delta^{\rm l}(v)$ on the projected time interval $[s,t]$, where $2kc_0$ is the difference between the drifts of $\Wi_0$ and $\Wi_r$, and its appearance in the intensity measure can be deduced from the fact that Brownian nets of different left-right speeds are related by changing the drift and performing diffusive rescaling. In particular, $N^{\rm l}_{s,t}$ is distributed as a Poisson random variable with mean $2kc_0\Delta^{\rm l}(t-s)$. Therefore for any $k_{\rm l}\geq 0$,
\be\label{Nlst}
\P[N^{\rm l}_{s,t} \geq k_{\rm l}] = \sum_{j=k_{\rm l}}^\infty \frac{1}{j!}\E\big[e^{-2kc_0\Delta^{\rm l}(t-s)} (2kc_0\Delta^{\rm l}(t-s))^j\big]
\leq \sum_{j=k_{\rm l}}^\infty \frac{1}{j!}\E\big[(2kc_0\Delta^{\rm l}(t-s))^j\big].
\ee
The Skorohod equation (\ref{skoro}) admits a pathwise unique solution (see e.g.~\cite[Sec.~3.6.C]{KS91}) with $\Delta^{\rm l}(t-s)=-\inf_{0\leq v\leq t-s}(B^{\rm l}(v)-\pi^0(t-v))$. By the independence of $B^{\rm l}$ and $\pi^0$, $B^{\rm l}(v)-\pi^0(t-v)$ is distributed as $-\sqrt{2}\tilde B^{\rm l}(v)-2kc_0 v$ for a standard Brownian motion $\tilde B^{\rm l}$ starting from $0$. Therefore
$$
\Delta^{\rm l}(t-s) = \sqrt{t-s} \sup_{0\leq v\leq 1} (\sqrt{2}\tilde B^{\rm l}(v) + 2kc_0v\sqrt{t-s})\leq \sqrt{t-s}\big(\sqrt{2}\sup_{0\leq v\leq 1}\tilde B^{\rm l}(v) + 2kc_0\sqrt{t-s}\big).
$$
Letting $\phi_{k_{\rm l}}$ denote the function
$\phi_{k_{\rm l}}(z):=\sum_{j=k_{\rm l}}\frac{1}{j!}z^j$, we observe that
$\phi_{k_{\rm l}}(\la z)\leq\la^{k_{\rm l}}\phi_{k_{\rm l}}(z)$ for all
$z\geq 0$ and $0\leq\la\leq 1$. Substituting these bounds into (\ref{Nlst}),
we see that we may estimate
\[
\P[N^{\rm l}_{s,t} \geq k_{\rm l}]
\leq(t-s)^{k_{\rm l}/2}\E\big[\phi_{k_{\rm l}}
\big(\sqrt{2}\sup_{0\leq v\leq 1}\tilde B^{\rm l}(v)+2kc_0\big)\big]
\qquad \mbox{for all } 0\leq  t-s\leq 1.
\]
Since $\sup_{0\leq v\leq 1}\tilde B^{\rm l}(v)$ is equally distributed with
$|\tilde B^{\rm l}(1)|$, which has finite exponential moments, it follows that
$$
\P[N^{\rm l}_{s,t} \geq k_{\rm l}] \leq C_{k,k_{\rm l}} (t-s)^{k_{\rm l}/2} \qquad \mbox{for all } 0\leq  t-s\leq 1
$$
for some $C_{k,k_{\rm l}}$ depending only on $k$ and $k_{\rm l}$ and $c_0$. If we let $N^{\rm r}_{s,t}$ denote the number of $s,t$-relevant separation points along $\pi^0$ which are of type $(1,2)_{\rm r}$ in $\Wi_0$, then similarly, $N^{\rm r}_{s,t}$ is distributed as a Poisson random variable with mean $2kc_1\Delta^{\rm r}(t-s)$, where
$\Delta^{\rm r}(t-s)=\sup_{0\leq v\leq t-s}(B^{\rm r}(v)-\pi^0(t-v))$ for an independent Brownian motion $B^{\rm r}$ with drift $-\beta_-$ starting from $\pi^0(t)$. Also, for any $k_{\rm r}\geq 0$,
$$
\P[N^{\rm r}_{s,t} \geq k_{\rm r}] \leq C_{k,k_{\rm r}} (t-s)^{k_{\rm r}/2} \qquad \mbox{for all } 0\leq  t-s\leq 1
$$
for some $C_{k,k_{\rm r}}$ depending only on $k$ and $k_{\rm r}$ and $c_1$.

Observe that if we impose the partial order $\prec$ on ${\cal C}_0([0,t-s])$, the space of continuous functions with value $0$ at $0$, where $f,g\in {\cal C}_0([0,t-s])$ satisfies $f\prec g$ if $f(v)\leq g(v)$ for all $v\in [0,t-s]$, then conditional on $B^{\rm l}$
and $B^{\rm r}$,  $\Delta^{\rm l}(t-s)$ is increasing in $(\pi^0(t-v)-\pi^0(t))_{v\in [0, t-s]}$, while $\Delta^{\rm r}(t-s)$ is decreasing in $(\pi^0(t-v)-\pi^0(t))_{v\in [0, t-s]}$. This implies that for any $k_{\rm l}, k_{\rm r}\geq 0$, $\P[N^{\rm l}_{s,t}\geq k_{\rm l}|\pi^0, B^{\rm l}]$ is increasing in $(\pi^0(t-v)-\pi^0(t))_{v\in [0, t-s]}$, while $\P[N^{\rm r}_{s,t}\geq k_{\rm r}|\pi^0, B^{\rm r}]$ is decreasing in $(\pi^0(t-v)-\pi^0(t))_{v\in [0, t-s]}$, and hence the same holds for $\P[N^{\rm l}_{s,t}\geq k_{\rm l}|\pi^0]$ and $\P[N^{\rm r}_{s,t}\geq k_{\rm r}|\pi^0]$. Since the Brownian motion $(\pi^0(t-v)-\pi^0(t))_{v\in [0, t-s]}$ satisfies the FKG inequality w.r.t.\ the partial order $\prec$ (see e.g.~\cite{Bar05}), and the events $\{N^{\rm l}_{s,t}\geq k_{\rm l}\}$ and $\{N^{\rm r}_{s,t}\geq k_{\rm r}\}$ are independent conditional on $\pi^0$, we have
\begin{eqnarray*}
\P[N^{\rm l}_{s,t}\geq k_{\rm l}, N^{\rm r}_{s,t}\geq k_{\rm r}] &=& \E\big[\P[N^{\rm l}_{s,t}\geq k_{\rm l}|\pi^0]\P[N^{\rm r}_{s,t}\geq k_{\rm r}|\pi^0]\big] \\
&\leq& \P[N^{\rm l}_{s,t}\geq k_{\rm l}] \P[N^{\rm r}_{s,t}\geq k_{\rm r}] \leq C_{k, k_{\rm l}}C_{k, k_{\rm r}} t^{\frac{k_{\rm l}+k_{\rm r}}{2}}.
\end{eqnarray*}
Since $N_{s,t}=N^{\rm l}_{s,t}+N^{\rm r}_{s,t}$, substituting this bound into (\ref{kthmoment}) then gives
\be\label{kthmombd}
\E\big[\big((X_s-X_t)^+\big)^k\big] \leq \P[N_{s,t}\geq k] \leq \sum_{i=0}^k \P[N^{\rm l}_{s,t}\geq i, N^{\rm r}_{s,t}\geq k-i]
\leq C_k (t-s)^{\frac{k}{2}}
\ee
for some $C_k$ depending only on $k$.

We can now apply the one-sided Kolmogorov's moment criterion, Theorem~\ref{T:moment}, which implies that for any $\gamma \in (0, 1/2)$ and $\delta>0$, if we choose $k$ sufficiently large in (\ref{kthmombd}) such that $(k/2-1)/k>\gamma$ and $k/2-1-\gamma k>\delta$, then
there exists a random constant $C_{\gamma, \Wi_0}$ such that
\be\label{holderbd2}
X_0-X_t=\rho_0(\{\pi^0(0)\})-\rho_t(\{\pi^0(t)\}) \leq C_{\gamma,\Wi_0} t^\gamma \qquad \mbox{ for all }\ t\in \Q_2\cap [0,1],
\ee
where the distribution of $C_{\gamma,\Wi_0}$ satisfies (\ref{holderconst}).

To extend (\ref{holderbd2}) to all $t\in [0,1]$, we note that $\rho_t(\cdot)$ is a.s.\ weakly continuous in $t$, which implies that $\tilde\rho_t(\cdot):=\rho_t(\pi^0(t)+\cdot)$ is weakly continuous in $t$ because $\pi^0$ is a.s.\ continuous. Since $\{0\}$ is a closed set, $\tilde\rho_t(\{0\})=\rho_t(\{\pi^0(t)\})$ is upper semi-continuous in $t$. Approximating any $t\in [0,1]$ by a sequence $t_n\in \Q_2\cap [0,1]$, we can then
use (\ref{holderbd2}) to deduce (\ref{massholder}).\qed
\bigskip

\noi
{\bf Proof of Theorem~\ref{T:erosion}.} By the third remark after Theorem~\ref{T:infmass}, it suffices to consider the case $\rho_0$ is a probability measure. Let $I_0(t):=\{\pi(t):\pi\in\Wi_0,\ \sig_\pi<t\}$. We start by proving that
\be\label{inI0}
\P[\rho_t(I_0(t))=1\ \forall t>0]=1.
\ee
The proof of Theorem~\ref{T:atom}~(a) in Section~\ref{S:detatom} showed that
\be
\P[\rho_t(I_0(t))=1]=1\qquad(t>0).
\ee
Fix some deterministic $s>0$. For each $x\in\R$, let $\rho^{(x)}$ denote the Howitt-Warren process started at time $s$ from $\rho^{(x)}_s:=\de_x$, constructed from the same reference web $\Wi_0$ as $\rho$. Then $\rho_t=\sum_{x\in I_0(s)}a_x\rho^{(x)}_s$ $(t\geq s)$, where $a_x:=\rho_s(\{x\})$. If $x$ is deterministic or if $x\in I_0(s)$, then there is an a.s.\ unique path in $\Wi_0((x,s))$; let $\pi^0_{(x,s)}$ denote this path. Then, for each deterministic $\eps\in(0,1)$ and $u>s$,
\[\ba{l}
\dis\P\Big[\inf_{t\in[s,u]}\rho_t(I_0(t))\leq 1-\eps\Big]
\leq\P\Big[\inf_{t\in[s,u]}\sum_{x\in I_0(s)}\rho_t(\{\pi^0_{(x,s)}(t)\})
\leq 1-\eps\Big]\\[5pt]
\dis\quad\leq\P\Big[\inf_{t\in[s,u]}\sum_{x\in I_0(s)}a_x
\rho_t^{(x)}(\{\pi^0_{(x,s)}(t)\})\leq 1-\eps\Big]\\[5pt]
\dis\quad\leq\P\Big[\sum_{x\in I_0(s)}a_x\inf_{t\in[s,u]}
\rho_t^{(x)}(\{\pi^0_{(x,s)}(t)\})\leq 1-\eps\Big]\\[5pt]
\dis\quad=\P\Big[\sum_{x\in I_0(s)}a_x\sup_{t\in[s,u]}
\big(1-\rho_t^{(x)}(\{\pi^0_{(x,s)}(t)\})\big)\geq\eps\Big]\\[5pt]
\dis\quad\leq\eps^{-6}\E\Big[\Big(\sum_{x\in I_0(s)}a_x\sup_{t\in[s,u]}
\big(1-\rho_t^{(x)}(\{\pi^0_{(x,s)}(t)\})\big)\Big)^6\Big]\\[5pt]
\dis\quad\leq\eps^{-6}\E\Big[\sum_{x\in I_0(s)}a_x\big(\sup_{t\in[s,u]}
\big(1-\rho_t^{(x)}(\{\pi^0_{(x,s)}(t)\})\big)\big)^6\Big]\\[5pt]
\dis\quad=\eps^{-6}\E\Big[\big(\sup_{t\in[s,u]}
\big(1-\rho_t^{(0)}(\{\pi^0_{(0,s)}(t)\})\big)\big)^6\Big],
\ea\]
where in the last inequality we applied the H\"older inequality w.r.t.\ the
probability law given by the $(a_x)_{x\in I_0(s)}$, and in the last equality
we used the spatial translation invariance of $\Wi_0$.

By Lemma~\ref{L:atomonpath} with $\gamma=1/3$,
\be\label{gamma3}
\E\Big[\sup_{t\in [s,u]}\big(1-\rho^{(0)}_t(\{\pi^0_{(0,s)}(t)\})\big)^6\Big]
\leq \E\big[C_{\gamma, \Wi_0}^6 (u-s)^2\big]=C (u-s)^2
\ee
for some finite $C>0$, since $C_{\gamma,\Wi_0}$ has finite moments of all orders by (\ref{holderconst}). Therefore, by our previous calculation, uniformly for all deterministic $0<s<u$,
\be\label{masslossbd}
\P\Big[\inf_{t\in[s,u]}\rho_t(I_0(t))\leq 1-\eps\Big]\leq C\eps^{-6}(u-s)^2.
\ee
It follows that for each $n\geq 1$,
\be
\P\Big[\inf_{t\in[2^{-n},1]}\rho_t(I_0(t))\leq 1-\eps\Big]
\leq\sum_{k=2}^{2^n}\P\Big[\inf_{t\in[(k-1)2^{-n},k2^{-n}]}
\rho_t(I_0(t))\leq 1-\eps\Big]
\leq C\eps^{-6} 2^{-n}.
\ee
Letting first $n\to\infty$ and then $\eps\to 0$ shows that that a.s.\ $\rho_t(I_0(t))=1$ for all $t\in (0,1]$, and similarly for all $t>0$. This completes the proof of (\ref{inI0}). In particular, this shows that almost surely, $\rho_t$ is atomic for all $t>0$.

To complete the proof, we need to show that almost surely, for all $t>0$,
\be\label{Isupp}
\big\{x\in\R:\rho_t(\{x\})>0\big\}=I_0(t)\cap\supp(\rho_t).
\ee
The inclusion $\sub$ follows from (\ref{inI0}). The inclusion $\supset$ is trivial for the Arratia flow, for which the
characteristic measure $\nu=0$. To prove it for erosion flows with $\nu=c_0\delta_0 +c_1\delta_1$ ($c_0+c_1>0$), by Lemma~\ref{L:edge} and Proposition~\ref{P:rhosupp},
we only need to consider the cases that either $\supp(\rho_t)=\R$ for all $t>0$, or $\supp(\rho_t)$ is a halfline whose moving boundary is a path in $\Wi_0$. Let $\Ti\sub(0,\infty)$ be a deterministic countable dense set. Then for each $t>0$ and $x\in I_0(t)\cap\supp(\rho_t)$, we can find a time $s\in\Ti$ and $y\in I_0(s)\cap\supp(\rho_s)$ such that $\pi^0_{(y,s)}(t)=x$, where $\pi^0_{(y,s)}$ denotes the a.s.\ unique path in $\Wi_0$ starting from $(y,s)$. Since by Lemma~\ref{L:posatom}, we have $\big\{x\in\R:\rho_s(\{x\})>0\big\}\supset I_0(s)\cap\supp(\rho_s)$, it suffices to prove that if $\pi^0$ is a path in $\Wi_0$ and $\rho_s(\{\pi^0(s)\})>0$ for some $s\in\Ti$, then $\rho_t(\{\pi^0(t)\})>0$ for all $t\geq s$.

By translation invariance, and the fact that $I_0(s)$ is independent of the restriction of $\Wi_0$ to $[s,\infty]$, it suffices to prove that almost surely, if $\pi^0$ is the unique path in $\Wi_0$ starting from the origin and $\rho_0=\de_0$, then $\rho_t(\{\pi^0(t)\})>0$ for all $t\geq 0$. To see this, set
\be
\psi_{s,t}:=K^\up_{s,t}\big(\pi^0(s),\{\pi^0(t)\}\big)\qquad(0\leq s\leq t)
\ee
and observe that $\psi_{s,t}\psi_{t,u}\leq\psi_{s,u}$ $(0\leq s\leq t\leq u)$. It follows from (\ref{gamma3}) that
\be
\P[\inf_{t\in[s,u]}\psi_{s,t}\leq 1/2]\leq C(u-s)^2
\ee
for some finite $C>0$, and therefore
\be
\P\big[\inf_{t\in[0,1]}\psi_{0,t}=0\big]
\leq\sum_{k=1}^{2^n}\P\big[
\inf_{t \in [(k-1)2^{-n},k2^{-n}]} \psi_{(k-1)2^{-n}, t}\leq 1/2\big]
\leq C2^{-n}.
\ee
Letting $n\to\infty$ shows that a.s.\ $\rho_t(\{\pi^0(t)\})=\psi_{0,t}>0$ for all $t\in[0,1]$ and similarly for all $t>0$.\qed
\bigskip

\noi
{\bf Proof of Theorem~\ref{T:atom}~(c).} Immediate from Theorem~\ref{T:erosion}.\qed

\section{Infinite starting mass and discrete approximation}\label{S:inf}

In this section, we prove Theorems \ref{T:infmass}--\ref{T:disflow}, which will be based on our construction
of the Howitt-Warren flows in Theorems~\ref{T:HWconst} and \ref{T:HWconst2}.

\subsection{Proof of Theorem~\ref{T:infmass}}

\label{S:Feller}
We first prove part (b), assuming $\beta_+-\beta_-<\infty$. By Theorem~\ref{T:HWconst2}, we have the following
representation for a Howitt-Warren process with drift $\beta$ and characteristic measure $\nu$:
\be\label{rhotrep}
\rho_t = \int \rho_s(\di x) K^+_{s,t}(x, \cdot) = \int \rho_s(\di x) \P[\pi^+_{(x,s)}(t)\in \cdot|(\Ni, \omega)] \qquad\mbox{for all } s<t,
\ee
where $\Ni$ is the Brownian net with left and right speeds $\beta_-$ resp.\ $\beta_+$, $\omega:=(\omega_z)_{z\in S}$ are
i.i.d.\ marks attached to the separation points $S$ of $\Ni$, and $\pi^+_{(x,s)}$ is the rightmost path starting from $(x,s)$ in
the sample web $\Wi$ conditional on the random environment $(\Ni, \omega)$. By construction, $\Wi\subset \Ni$ a.s., therefore
for any $L>0$, $K^+_{0,t}(x, [-L,L])=0$ if $|x|$ is sufficiently large. This implies that if $\rho_0\in\Mi_{\rm loc}$, then a.s.\
$\rho_t\in \Mi_{\rm loc}$ for all $t\geq 0$.

Let $s_n$, $t_n$, $t$, $\rho^{\langle n\rangle}_{s_n}$ and $\rho_0$ be as in Theorem~\ref{T:infmass}. By the representation (\ref{rhotrep}),
proving (\ref{rhonconv}) with vague convergence in $\Mi_{\rm loc}(\R)$ amounts to showing that, for all $f\in C_c(\R)$,
\be\label{Feller2}
\int \rho^{\langle n\rangle}_{s_n}(\di x)\, \E\big[f(\pi^+_{(x,s_n)}(t_n)) \,|\, (\Ni, \om)\big]
\asto{n} \int \rho_0(\di x)\, \E\big[f(\pi^+_{(x,0)}(t)) \,|\, (\Ni, \om)\big] \qquad a.s.
\ee
Note that the sample web $\Wi$ constructed in Theorem~\ref{T:HWconst2} and from which we draw $\pi^+_{(\cdot,\cdot)}$, is distributed
as a Brownian web with drift $\beta$. Therefore for each $(x,s)\in \R^2$, a.s.\ $(x,s)$ is of type $(0,1)$ in $\Wi$, and hence
$\pi^+_{(x_n,s_n)}\to \pi^+_{(x,s)}$ for any $(x_n, s_n)\to (x,s)$. By Fubini, a.s.\ w.r.t.\ $(\Ni, \om)$, there exists
a set $A_{(\Ni, \om)} \subset \R$ with full $\rho_0$ measure such that for all $x\in A_{(\Ni, \om)}$, $(x,0)$ is of type $(0,1)$
in $\Wi$ a.s.\ w.r.t.\ the quenched law $\Q:=\P(\Wi\in \cdot |(\Ni,\omega)$. This implies that
$$
\phi_{(\Ni, \om)}(x,s; t):= \E\big[f(\pi^+_{(x,s)}(t)) \,|\, (\Ni, \om)\big]
$$
is continuous at all $x\in A_{(\Ni, \om)}$, $s=0$ and $t\geq 0$. Since $\Wi\subset \Ni$ a.s.\ and $f$ has compact support,
there exists $L>0$ such that
\be\label{Feller3}
 \phi_{(\Ni, \om)}(x,s_n; t_n) = \E\big[f(\pi^+_{(x,s_n)}(t_n)) \,|\, (\Ni, \om)\big] =  0 \qquad \mbox{for all } |x|> L, n\in\N.
\ee
Choose $-L_1<-L$ and $L_2>L$ to be points of continuity of $\rho_0$. Note that restricted to $(-L_1, L_2)$, $\rho^{\langle n\rangle}_{s_n}$ converges
weakly to $\rho_0$. The convergence in (\ref{Feller2}) then follows from the continuous mapping theorem for weak convergence. The almost sure
path continuity of $(\rho_t)_{t\geq 0}$ follows from (\ref{Feller2}) by setting $s_n=0$ and $\rho^{\langle n\rangle}_{s_n}=\rho_0$. This proves
part (b).

We now prove part (a). When $\beta_+-\beta_-=\infty$, (\ref{Feller3}) may fail and we need to control in any finite region the inflow of measure
from arbitrarily far away. By Theorem~\ref{T:HWconst}, we have the representation:
\be\label{rhotrep2}
\rho_t = \int \rho_s(\di x) K^+_{s,t}(x, \cdot) = \int \rho_s(\di x) \P[\pi^+_{(x,s)}(t)\in \cdot |(\Wi_0, \Mi)] \qquad\mbox{for all } s<t.
\ee
Here $(\Wi_0, \Mi)$ are the reference web and its associated set of marked points as in Theorem~\ref{T:HWconst}, where for definitiveness
we can choose $\nu_{\rm l}=\nu_{\rm r}=\nu$, so that $\Wi_0$ and the sample web $\Wi$ are both Brownian webs with drift $\beta$. In particular,
$\pi^+_{(x,0)}$ is distributed as a Brownian motion with drift $\beta$ starting from $x$ at time $0$. Without loss of generality, assume $\beta=0$.
Since $\rho_0\in \Mi_g(\R)$, it is then easy to check that for any $c>0$,
$$
\E\Big[\int e^{-cy^2}\rho_t(\di y)\Big] = \frac{1}{\sqrt{2\pi t}} \int\rho_0(\di x) \int e^{-cy^2} e^{-\frac{(y-x)^2}{2t}}\di y \leq C_1\int \rho_0(\di x)e^{-C_2x^2} <\infty
$$
for some $C_1, C_2>0$, which implies (\ref{EK}) and that $\rho_t\in \Mi_g(\R)$ almost surely.

By the representation (\ref{rhotrep2}), proving (\ref{rhonconv}) with convergence in $\Mi_g(\R)$ amounts to showing that,
for all $c>0$ and all bounded continuous function $f:\R\to\R$, we have
\be\label{Feller4}
\begin{aligned}
& \int e^{-\eps x^2} \rho^{\langle n\rangle}_{s_n}(\di x)\  e^{\eps x^2} \E\big[f(\pi^+_{(x,s_n)}(t_n))e^{-c\pi^+_{(x,s_n)}(t_n)^2} \,|\, (\Wi_0, \Mi)\big] \\
& \qquad \qquad \asto{n} \qquad \int e^{-\eps x^2} \rho_0(\di x)\  e^{\eps x^2} \E\big[f(\pi^+_{(x,0)}(t))e^{-c\pi^+_{(x,0)}(t)^2} \,|\, (\Wi_0, \Mi)\big] \qquad a.s.,
\end{aligned}
\ee
where $\eps>0$ is chosen small. Denote
\be\label{phiwm}
\phi_{(\Wi_0, \Mi)}(x,s; t):= e^{\eps x^2} \E\big[f(\pi^+_{(x,s)}(t))e^{-c\pi^+_{(x,s)}(t)^2} \,|\, (\Wi_0, \Mi)\big].
\ee
As before, a.s.\ w.r.t.\ $(\Wi_0, \Mi)$, there exists $A_{(\Wi_0,\Mi)}\subset \R$ with full $\rho_0$ measure such that
$\phi_{(\Wi_0, \Mi)}(x,s; t)$ is continuous at all $x\in A_{(\Wi_0,\Mi)}$, $s=0$ and $t\geq 0$. Our assumption $\rho^{\langle n\rangle}_{s_n}\to \rho_0$ in
$\Mi_{\rm g}(\R)$ implies that $e^{-\eps x^2}\rho^{\langle n\rangle}_{s_n}(\di x)$ converges weakly to $e^{-\eps x^2}\rho_0(\di x)$. Therefore
(\ref{Feller4}) follows from the continuous mapping theorem for weak convergence, provided we show that a.s.\ w.r.t.\ $(\Wi_0,\Mi)$,
\be\label{Feller5}
\sup_{y\in \R, 0\leq u\leq v\leq t} |\phi_{(\Wi_0,\Mi)}(y,u;v)| <\infty \qquad \mbox{for all } t>0,
\ee
so that we can apply the bounded convergence theorem. We verify (\ref{Feller5}) by Borel-Cantelli.

Without loss of generality, assume $|f|_\infty=1$. Fix $t>0$. For each $m\in\Z$, we have
\be\label{Feller6}
\E\Big[\sup_{y\in [m, m+1] \atop 0\leq u\leq v\leq t} |\phi_{(\Wi_0, \Mi)}(y,u; v)|\Big]
\leq  C e^{3\eps m^2} \E\Big[\sup_{y\in [m,m+1] \atop 0\leq u\leq v\leq t} e^{-c\pi^+_{(y,u)}(v)^2} \Big] <\infty.
\ee
Now consider $m\geq L$ for some fixed large $L$. By the coalescing property of paths in the sample web $\Wi$, if $\pi^+_{(3m/4,0)}\in\Wi$ starting
from $(3m/4,0)$ stays within $[m/2, m]$ on the time interval $[0,t]$, then $\inf_{y\in [m,m+1], 0\leq u\leq v\leq t} \pi^+_{(y,u)}(v) \geq m/2$.
Therefore
$$
\E\Big[\sup_{y\in [m,m+1] \atop 0\leq u\leq v\leq t} e^{-c\pi^+_{(y,u)}(v)^2} \Big]
\leq \P\Big( \sup_{0\leq s\leq t}\big|\pi^+_{(3m/4,0)}(s)-\frac{3m}{4}\big|\geq \frac{m}{4}\Big) + e^{-\frac{cm^2}{4}}
\leq C e^{-\frac{m^2}{32t}} + e^{-\frac{cm^2}{4}},
$$
and hence
\be\label{Feller7}
\E\Big[\sup_{y\in [m, m+1] \atop 0\leq u\leq v\leq t} |\phi_{(\Wi_0, \Mi)}(y,u; v)|\Big]
\leq C e^{-(\frac{1}{32t}-3\eps)m^2} + C e^{-(\frac{c}{4}-3\eps)m^2} \leq C_1e^{-C_2m^2}
\ee
for some $C_1, C_2>0$ depending only on $t$ and $c$ if we choose $\eps>0$ sufficiently small. Thus
$$
\P\Big(\sup_{y\in [m, m+1] \atop 0\leq u\leq v\leq t} |\phi_{(\Wi_0, \Mi)}(y,u; v)| > 1 \Big) \leq C_1e^{-C_2m^2}.
$$
By Borel-Cantelli, a.s.\ w.r.t.\ $(\Wi_0,\Mi)$, there exists a random $N_+>L$ sufficiently large such that
$\sup_{y\geq N_+, 0\leq u\leq v\leq t} |\phi_{(\Wi_0, \Mi)}(y,u; v)| \leq 1$. Similarly, a.s.\ there exists $N_-<-L$ such that
$\sup_{y\leq N_-, 0\leq u\leq v\leq t} |\phi_{(\Wi_0, \Mi)}(y,u; v)| \leq 1$. Combined with (\ref{Feller6}), this implies
(\ref{Feller5}), and hence (\ref{Feller4}). The almost sure path continuity of $(\rho_t)_{t\geq 0}$ in $\Mi_{\rm g}(\R)$
follows from (\ref{Feller4}) by setting $s_n=0$ and $\rho^{\langle n\rangle}_{s_n}=\rho_0$.
\qed

\subsection{Proof of Theorem \ref{T:disflow}}\label{S:disflow}

The proof is similar to that of Theorem \ref{T:infmass}. The complication lies again with infinite $\bar\rho^{\langle k\rangle}_t$
and $\rho_t$. Without loss of generality, assume that the Howitt-Warren process $(\rho_t)_{t\geq 0}$ has drift $\beta=0$.
First we note that there exists a countable family of bounded continuous functions $\{f_n\}_{n\in\N}$ such that a sequence of
finite measures $\xi_k\in\Mi(\R)$ converges weakly to $\xi\in\Mi(\R)$ if and only if $\int f_n(x)\xi_k(\di x) \to\int f_n(x)\xi(\di x)$ for all $n\in\N$,
see e.g.\ \cite[Proof of Prop.~3.17]{Res87}. Since $\xi_k\to\xi$ in $\Mi_{\rm g}(\R)$ is equivalent to weak convergence of
$e^{-cx^2}\xi_k(\di x)$ to $e^{-cx^2}\xi(\di x)$ for all $c>0$, to prove the weak convergence in (\ref{disflow}) on path
space with uniform topology, it suffices to show that for any finite sets $K\subset (0,\infty)$ and $\Lambda\subset \N$, we have
\be\label{df1}
\Big(F^{\langle k\rangle}_{c,n}(t) := \int  e^{-cx^2}f_n(x) \bar\rho^{\langle k\rangle}_t(\di x)\Big)_{c\in K,n\in \Lambda}\ \Asto{n}\ \Big(F_{c,n}(t):= \int e^{-cx^2} f_n(x)\rho_t(\di x)\Big)_{c\in K, n\in\Lambda},
\ee
where $\Rightarrow$ denotes weak convergence of $\Ci([0,T],\R)^{|K|+|\Lambda|}$-valued random variables.

For any $c>0$ and $n\in\N$, $F_{c,n}\in \Ci([0,T],\R)$ a.s.\ by Theorem \ref{T:infmass}. By similar reasoning, $\bar\rho^{\langle k\rangle}_t$ has a.s.\
continuous sample path in $\Mi_{\rm g}(\R)$ and hence
$F^{\langle k\rangle}_{c,n}\in \Ci([0,T],\R)$. Since $\bar\rho^{\langle k\rangle}_0$ converges weakly to $\rho_0$ as $\Mi_{\rm g}(\R)$-valued random variables, the Skorohod
representation theorem for weak convergence (see e.g.~\cite[Theorem~6.7]{Bil99}) allows a coupling between $(\bar\rho^{\langle k\rangle}_0)_{k\in\N}$ and
$\rho_0$ such that $\bar\rho^{\langle k\rangle}_0\to \rho_0$ in $\Mi_{\rm g}(\R)$ almost surely. Therefore we may assume that
$(\bar\rho^{\langle k\rangle}_0)_{n\in\N}$ and $\rho_0$ are deterministic and $\bar\rho^{\langle k\rangle}_0\to \rho_0$.

Recall from (\ref{quench}) and Section \ref{S:approx} the discrete quenched law $\Qdis_{\langle k\rangle}$ associated with a discrete Howitt-Warren flow with
characteristic measure $\mu_k$, and recall from (\ref{Seps2}) and Theorem~\ref{T:quench} the diffusive scaling map $S_\eps$ and its action on a quenched law $\Qdis$.
We have the representation
$$
\bar\rho^{\langle k\rangle}_t = \int \bar\rho_0(\di x) S_{\eps_k}(\Qdis_{\langle k\rangle})[\pi^{\langle k\rangle}_{(x,0)}(t)\in \cdot],
$$
where for $(x,s) \in S_{\eps_k}(\Z^2_{\rm even})$, $\pi^{\langle k\rangle}_{(x,s)}$ is the unique path starting from $(x,s)$ in a discrete sample web $\Wi^{\langle k\rangle}$ with quenched law
$S_{\eps_k}(\Qdis_{\langle k\rangle})$. Similarly, by Theorem~\ref{T:HWconst},
$$
\rho_t =\int\rho_0(\di x) \Q[\pi^+_{(x,0)}(t)\in \cdot],
$$
with the quenched law $\Q$ defined as in (\ref{HWquen}). For any $L>0$, we can then write
$$
\begin{aligned}
F^{\langle k\rangle}_{c,n}(t) & = F^{\langle k\rangle, [-L,L]}_{c,n}(t) + F^{\langle k\rangle, [-L,L]^c}_{c,n}(t)  \\
F_{c,n}(t) & = F^{[-L,L]}_{c,n}(t) + F^{[-L,L]^c}_{c,n}(t),
\end{aligned}
$$
where for any $I\subset\R$,
$$
\begin{aligned}
F^{\langle k\rangle, I}_{c,n}(t) &= \int_I \bar\rho_0(\di x) S_{\eps_k}(\Qdis_{\langle k\rangle})\Big[f_n(\pi^{\langle k\rangle}_{(x,0)}(t)) e^{-c \pi^{\langle k\rangle}_{(x,0)}(t)^2}\Big], \\
F^I_{c,n}(t) &= \int_I \rho_0(\di x) \Q\Big[f_n(\pi^+_{(x,0)}(t)) e^{-c \pi^+_{(x,0)}(t)^2}\Big].
\end{aligned}
$$
To prove (\ref{df1}), it suffices to show that for any $\eps>0$, $c>0$ and $n\in\N$, we can choose $L$ large such that
\begin{gather}
\limsup_{k\to\infty} \E\big[\,|F^{\langle k\rangle, [-L,L]^c}_{c,n}|_\infty \big] \leq \eps,  \label{df2}\\
\E\big[\,|F^{[-L,L]^c}_{c,n}|_\infty \big] \leq \eps,   \label{df3}
\end{gather}
where $|\cdot|_\infty$ denotes the supremum norm on $\Ci([0,T], \R)$, and furthermore,
\be
(F^{\langle k\rangle, [-L,L]}_{c,n})_{c\in K, n\in\Lambda} \Asto{k} (F^{[-L,L]}_{c,n})_{c\in K, n\in\Lambda} \label{df4}
\ee
with $\Rightarrow$ denoting weak convergence of $\Ci([0,T],\R)^{|K|+|\Lambda|}$-valued random variables.

Fix $0<\eps < \inf K$ and define
\be\label{phiq}
\begin{aligned}
\phi^{\langle k\rangle,c,n}_{\Qdis_{\langle k\rangle}}(x,t) &:= e^{\eps x^2} (S_{\eps_k}\Qdis_{\langle k\rangle})\Big[f_n(\pi^{\langle k\rangle}_{(x,0)}(t)) e^{-c \pi^{\langle k \rangle}_{(x,0)}(t)^2}\Big], \\
\phi^{c,n}_\Q(x,t) &:=  e^{\eps x^2} \Q\Big[f_n(\pi^+_{(x,0)}(t)) e^{-c \pi^+_{(x,0)}(t)^2}\Big].
\end{aligned}
\ee
Then
$$
| F^{[-L,L]^c}_{c,n}(\cdot)|_\infty = \Big|\int_{[-L,L]^c} e^{-\eps x^2} \rho_0(\di x) \phi^{c,n}_{\Q}(x,\cdot)\Big|_\infty
\leq  \int_{[-L,L]^c} e^{-\eps x^2} \rho_0(\di x) |\phi^{c,n}_\Q(x,\cdot)|_\infty .
$$
Note that $\phi^{c,n}_\Q(x,t)$ is exactly $\phi_{(\Wi_0,\Mi)}(x,s;t)$ in (\ref{phiwm}) if we set $s=0$ and $f=f_n$. Since $e^{-\eps x^2}\rho_0(\di x)$
is a finite measure, (\ref{df3}) then follows from (\ref{Feller7}). Note that (\ref{Feller7}) is based on Brownian motion estimates, and
analogues of (\ref{Feller7}) for $\phi^{\langle k\rangle, c,n}_{\Qdis_{\langle k\rangle}}$ can be established using random walk estimates, which are
furthermore uniform in $k\in\N$. Such a uniform estimate and the fact that $e^{-\eps x^2}\bar\rho_0(\di x)$ converges weakly to $e^{-\eps x^2}\rho_0(\di x)$ then
imply (\ref{df2}). We omit the details here.

Lastly, we prove (\ref{df4}), where we may assume that $-L$ and $L$ are points of continuity of $\rho_0(\di x)$ so that
restricted to $[-L,L]$, $e^{-\eps x^2}\bar\rho_0(\di x)$ converges weakly to $e^{-\eps x^2}\rho_0(\di x)$. By Skorohod representation, the
weak convergence of $S_{\eps_k}\Qdis_{\langle k\rangle}$ to $\Q$ in Theorem \ref{T:quench} can be turned into a.s.\ convergence in $\Mi_1(\Ki(\Pi))$
via a suitable coupling, which we now assume.  Then the discrete sample web $\Wi^{\langle k\rangle}$ with law $S_{\eps_k}\Qdis_{\langle k\rangle}$
converges weakly to the sample web $\Wi$ with law $\Q$, where the convergence can again assumed to be a.s.\ in $\Ki(\Pi)$ by Skorohod representation.
Then for each $c>0$ and $k, n\in\N$,
$$
|F^{\langle k\rangle, [-L,L]}_{c,n} - F^{[-L,L]}_{c,n}|_\infty =  \Big|\int_{-L}^L  e^{-\eps x^2}\bar\rho^{\langle k\rangle}_0(\di x) \phi^{\langle k\rangle, c,n}_{\Qdis_{\langle k\rangle}}(x,\cdot)
- \int_{-L}^L e^{-\eps x^2} \rho_0(\di x) \phi^{c,n}_\Q(x,\cdot)\Big|_\infty.
$$
First we claim that a.s.\ w.r.t.\ $\Qdis_{\langle k\rangle}$ and $\Q$, for a.e.\ $x\in [-L,L]$ w.r.t.\ the measure $e^{-\eps x^2} \rho_0(\di x)$, if $x_k\to x$ for some sequence $x_k\in \eps_k\Z_{\rm even}$, then for each $c>0$ and $n\in\N$,
\be\label{supphi}
|\phi^{\langle k\rangle,c,n}_{\Qdis_{\langle k\rangle}}(x_k,\cdot) - \phi^{c,n}_\Q(x,\cdot)|_\infty \asto{k} 0.
\ee
Indeed, since the law of $\Wi$ averaged over the randomness of $\Q$ is that of a Brownian web, each deterministic $(x,0)$ is a.s.\ of type
$(0,1)$ in $\Wi$, and by Fubini, a.s.\ w.r.t.\ $\Q$, the same is true for $e^{-\eps x^2} \rho_0$ a.e.\ $x$. Therefore $\Wi^{\langle k\rangle }\to\Wi$ a.s.\ in $\Ki(\Pi)$
implies that $\pi^{\langle k\rangle}_{(x_k,0)}\to \pi^+_{(x,0)}$ in $\Ci([0,T],\R)$, which when plugged into the definitions in (\ref{phiq}) then implies
(\ref{supphi}). If we regard $\phi^{\langle k\rangle, c,n}_{\Qdis_{\langle k\rangle}}(x,\cdot)$ and $\phi^{c,n}_\Q(x,\cdot)$ as mappings from $\R$ to
$\Ci([0,T],\R)$ and note that $|\phi^{\langle k\rangle,c,n}_{\Qdis_{\langle k\rangle}}(x,\cdot)|_\infty$ and $|\phi^{c,n}_\Q(x,\cdot)|_\infty$ are bounded
uniformly in $\Qdis_{\langle k\rangle}$, $\Q$ and $x\in [-L,L]$, then the continuous mapping theorem for weak convergence implies that
$|F^{\langle k\rangle, [-L,L]}_{c,n} - F^{[-L,L]}_{c,n}|_\infty\to 0$ a.s.\ for each $c\in K$ and $n\in\Lambda$, which then implies (\ref{df4}).
\qed

\section{Ergodic properties}\label{S:ergproof}

In this section, we prove Theorem~\ref{T:HIL}--\ref{T:invsupp} on homogeneous invariant laws for Howitt-Warren processes. By the
observation that $\rho_t$ in (\ref{rho}) depends linearly on the initial condition $\rho_0$, the Howitt-Warren process falls in
the class of linear systems, the theory of which for processes on $\Z^d$ was developed by Liggett and Spitzer, see e.g.\ \cite{LS81} and
\cite[Chap.~IX]{Lig05}. \index{linear system} We will adapt the theory of linear systems to our continuum setting.
The main tools are duality, second moment calculations, and coupling, which will be developed in successive
subsections. Duality is used to give a simple construction of the family of ergodic homogeneous invariant laws. Second moment
calculations determine spatial correlations for the homogeneous invariant laws, and are used to prove the uniform integrability of the
Howitt-Warren process $(\rho_t)_{t\geq 0}$ over time, as well as to show that certain spatial ergodic properties of the initial measure
$\rho_0$ are preserved by the dynamics even in the limit $t\to\infty$. The last point will be crucial for proving convergence to
homogeneous invariant laws, which is based on coupling arguments. Most of our arguments are adapted from \cite{LS81} and
\cite[Chap.~IX]{Lig05}, to which we will refer many of the details. The main difference in our proof lies in the second moment
calculations of Lemma \ref{L:2pt}, for which we need to devise a different and perhaps more robust approach than the one
used in \cite{LS81, Lig05}.

\subsection{Dual smoothing process}

\index{smoothing process}
Similar to the linear systems on $\Z^d$ studied in \cite{LS81, Lig05}, the Howitt-Warren process $(\rho_t)_{t\geq 0}$ is dual to
a function-valued smoothing process with random kernels. Analogous to the construction of the Howitt-Warren process from the
Howitt-Warren flow $(K^+_{s,t})_{s<t}$ as in (\ref{rho}), we can define a function-valued dual process
$(\zeta_t)_{t\geq 0}$ by
\be\label{zetat}
\zeta_t(x) = \int K^+_{-t,0}(x, \di y) \zeta_0(y) = \int \Q[\zeta_0(\pi^+_{(x,-t)}(0))]  \qquad \mbox{for all } x\in\R,
\ee
where $\Q$ is the quenched law of a sample web $\Wi$ defined as in Theorem~\ref{T:HWconst}, and $\pi^+_{(x,s)}$ is the a.s.\ unique
rightmost path in $\Wi$ starting from $(x,s)$. A natural state space for $(\zeta_t)_{t\geq 0}$ is $D_b(\R)$, the space of bounded c\'adl\'ag
functions on $\R$. Note that $\pi^+_{(x,-t)}(0)$ is c\'adl\'ag in $x$. With this
observation, it is easy to see that if $\zeta_0\in D_b(\R)$, then $\zeta_t\in D_b(\R)$ for all $t>0$.

We have the following duality relation between $(\rho_t)_{t\geq 0}$ and $(\zeta_t)_{t\geq 0}$.
\bl{\bf (Duality)}\label{L:duality}
Let $\rho_0\in \Mi_{\rm loc}(\R)$, or $\rho_0\in \Mi_{\rm g}(\R)$ if $\beta_+-\beta_-=\infty$ in (\ref{speeds}). Let $\zeta_0\in D_b(\R)$. Assume
that either $\rho_0$ is a finite measure or $\zeta_0$ has bounded support. Then for all $t\geq 0$,
\be
\int \zeta_0(x) \rho_t(\di x) \stackrel{\rm dist}{=} \int \zeta_t(x) \rho_0(\di x).
\ee
\el
{\bf Proof.} Follows from the definition of $\rho_t$ and $\zeta_t$, and the equality in law between $K^+_{0,t}$ and $K^+_{-t,0}$.
\qed

The advantage of working with the smoothing process $(\zeta_t)_{t\geq 0}$ is that there is a natural martingale associated with it.
\bl{\bf (Extinction vs uniform integrability)}\label{L:masslim} Let $\zeta_0\in D_b(\R)$ have bounded support. Then
$[\zeta]_t :=\int \zeta_t(x) \di x$ is a martingale which a.s.\ has a limit $[ \zeta]_\infty$ as
$t\to\infty$. Furthermore, either $[ \zeta]_\infty=0$ a.s.\ for all $\zeta_0\in D_b(\R)$ with bounded support,
or $([ \zeta]_t)_{t\geq 0}$ is uniformly integrable for all $\zeta_0\in D_b(\R)$ with bounded support. We say the
finite smoothing process $\zeta$ dies out in the first case, and survives in the second case.
\el
{\bf Proof.} By separating $\zeta_0$ into its positive and negative parts and by the linear dependence of $\zeta_t$ on $\zeta_0$, we
may assume $\zeta_0\geq 0$. Note that for $0\leq s<t$,
\be\label{dicho}
\begin{aligned}
{[ \zeta]}_t & = \int \zeta_t(x) \di x = \iint K^+_{-t,0}(x,\di y) \zeta_0(y) \di x \\
& = \iiint K^+_{-t,-s}(x, \di z) K^+_{-s,0}(z, \di y) \zeta_0(y) \di x
= \iint \zeta_s(z)\, K^+_{-t,-s}(x, \di z)  \di x.
\end{aligned}
\ee
By the independence of $K^+_{-t,-s}$ and $(K^+_{-u,0})_{0\leq u\leq s}$, $K^+_{-t,-s}$ is independent of $(\zeta_u)_{0\leq u\leq s}$.
By the translation invariance in law of $K^+_{-t,-s}(x,\cdot)$ in $x$, we note that $\int \E[K^+_{-t,-s}(x,\cdot)] \di x$ is simply
the Lebesgue measure. Therefore
$$
\E\big[[ \zeta] _t \,\big|\, ([ \zeta]_u)_{0\leq u\leq s} \big] = \iint \zeta_s(z) \E[K^+_{-t,-s}(x,\di z)] \di x
=\int \zeta_s(z) \di z = [ \zeta]_s,
$$
which proves the martingale property of $[ \zeta]_t$. Since $[ \zeta]_t$ is furthermore non-negative, a.s.\ it has a
limit $[ \zeta]_\infty$.

The dichotomy between a.s.\ extinction and uniform integrability of $([ \zeta]_t)_{t\geq 0}$ follows from a similar argument as in
the proof of \cite[Thm.~IX.2.4.(a)]{Lig05}. Given $\zeta_0(x)=1_{[0,1]}(x)$, let $\lambda:=\E[[\zeta]_\infty]$. First we claim that for any
$\zeta_0\in D_b(\R)$ with bounded support, we have $\E[[\zeta]_\infty] = \lambda [ \zeta]_0$. By the linear dependence
of $\zeta_\infty$ on $\zeta_0$ and the translation invariance of the kernels $(K^+_{s,t})_{s\leq t}$, the claim holds for all $\zeta_0$ which
are linear combinations of characteristic functions of finite intervals. Since all $\zeta_0\in D_b(\R)$ with bounded support can be approximated
from above and below by such functions, and $[\zeta]_\infty$ depends monotonically on $\zeta_0$, the claim holds for all
$\zeta_0\in D_b(\R)$ with bounded support. The dichotomy amounts to showing either $\lambda=0$ or $\lambda=1$.

Note that the RHS of (\ref{dicho}) can be interpreted as $[ \ti \zeta]_{t-s}:= \int \ti\zeta_{t-s}(x)\di x$ for a smoothing process
$\ti\zeta$ defined from the time-shifted kernels $(K^+_{-r-s, -s})_{r\geq 0}$ with initial condition $\ti\zeta_0=\zeta_s$. In particular,
a.s.\ $[ \ti\zeta]_{t-s}$ tends to a limit $[ \ti\zeta]_\infty$ as $t\to\infty$. Letting $t\to\infty$ in (\ref{dicho})
then gives $[ \zeta]_\infty = [ \ti\zeta]_\infty$. Therefore by Jensen's inequality,
$$
\E[e^{-[ \zeta]_\infty}] = \E[e^{-[ \ti\zeta]_\infty}] \geq \E\big[e^{-\E[[ \ti\zeta]_\infty \,|\, \ti\zeta_0=\zeta_s]}\big]
= \E[e^{-\lambda[ \zeta]_s}],
$$
where we now take the limit $s\to\infty$ and obtain by the bounded convergence theorem that
\be\label{dichot2}
\E[e^{-[ \zeta]_\infty}] \geq \E[e^{-\lambda[ \zeta]_\infty}].
\ee
Since $\zeta_0\geq 0$ by assumption in (\ref{dicho}), we have $[\zeta]_\infty\geq 0$. Assume further that $[ \zeta]_0>0$. If
$[ \zeta]_\infty=0$ a.s., then $\lambda=0$. If $[ \zeta]_\infty>0$ with positive probability, then because $\lambda\in [0,1]$, (\ref{dichot2})
can only hold if $\lambda=1$.
\qed

Lemmas \ref{L:duality} and \ref{L:masslim} imply the weak convergence of the Howitt-Warren process $\rho_t$ with initial condition
$\rho_0(\di x)=c\,\di x$ to a homogeneous invariant law. Recall the set of invariant laws $\Ii$ and $\Ti$ from Theorem \ref{T:HIL}.

\bl{\bf (Construction of homogeneous invariant laws)}\label{L:existlim}
Assume that $\rho_0(\di x)=c\,\di x$ for some $c\geq 0$. Then there exists $\Lambda_c \in \Ii\cap\Ti$ such that
\be
\Li(\rho_t) \Asto{t} \Lambda_c,
\ee
where $\Li(\cdot)$ denotes law and $\Rightarrow$ denotes weak convergence of probability laws on $\Mi_{\rm loc}(\R)$. If the finite smoothing process $\zeta$ survives,
then $\int \rho([0,1]) \Lambda_c(\di \rho) = c$; otherwise $\Lambda_c= \delta_0$, the delta measure concentrated in the zero measure on $\R$. Furthermore,
$\Lambda_c(\di (c\rho))= \Lambda_1(\di \rho)$.
\el
{\bf Proof.} Since $\rho_0(\di x) = c\,\di x$, by the translation invariance of $(K^+_{0,t}(x,\cdot))_{x\in\R}$ in space, we have
$\E[\rho_t(\di x)]=c\,\di x$ for all $t\geq 0$. In particular, for any bounded interval $I\subset\R$, $(\rho_t(I))_{t\geq 0}$ is a tight
family of random variables, which implies that $(\rho_t)_{t\geq 0}$ is a tight family of $\Mi_{\rm loc}(\R)$-valued random variables (see e.g.\
\cite[Lemma 3.20]{Res87}). In fact $(\rho_t)_{t\geq 0}$ is also a tight family of $\Mi_{\rm g}(\R)$-valued random variables (recall (\ref{MiG})).
This follows from the additional observation that for any $a>0$, $\big(\int e^{-ax^2}\rho_t(\di x)\big)_{t\geq 0}$ is a tight family of
real-valued random variables, because
$$
\sup_{t\geq 0} \E\Big[\int e^{-ax^2} \rho_t(\di x)\Big] = c \int e^{-ax^2} \di x <\infty.
$$

Let $\zeta_t$ be the dual smoothing process with initial condition $\zeta_0\in C_c(\R)$, the space of continuous functions with compact support.
If $\rho_t$ converges weakly to a $\Mi_{\rm g}(\R)$-valued
random variable $\rho^*$ along a subsequence $t_n\uparrow \infty$, then by Lemmas \ref{L:duality} and \ref{L:masslim}, we must have equality
in distribution between $\int \zeta_0(x) \rho^*(\di x)$ and $c[ \zeta]_\infty$. Since the law of $[ \zeta]_\infty$
does not depend on $t_n\uparrow \infty$, and $\zeta_0$ can be any function in $C_c(\R)$, the law of $\rho^*$ is uniquely determined (see e.g.\
\cite[Prop.~3.19]{Res87}). Together with tightness, this implies that $\rho_t \Asto{t} \rho^*$ as $\Mi_{\rm g}(\R)$-valued random variables,
and we denote $\Lambda_c:=\Li(\rho^*)$. The fact that $\Lambda_c \in \Ii$ then follows from the Feller property of $(\rho_t)_{t\geq 0}$ implied by
Theorem \ref{T:infmass}, and clearly $\Lambda_c \in \Ti$. Since $\E[\rho^*([0,1])] = c\,\E[[\zeta]_\infty]$ with $\zeta_0=1_{[0,1]}$,
the dichotomy between $\int \rho([0,1]) \Lambda_c(\di \rho)=c$ and $\Lambda_c=\delta_0$ follows from Lemma \ref{L:masslim}. The scaling
relation between $\Lambda_1$ and $\Lambda_c$ is trivial.
\qed

When the characteristic measure $\nu$ for the Howitt-Warren flow is not zero so that the flow is not purely coalescing, the possibility
of $\Lambda_c =\delta_0$ in Lemma \ref{L:existlim} can be ruled out by showing the uniform integrability of $\rho_t([0,1])$ in $t\geq 0$.
This can be accomplished by the second moment calculation in Lemma \ref{L:2pt} below. In any event, we can deduce the extremality of
$\Lambda_c$ in $\Ii\cap\Ti$ using Lemmas \ref{L:duality} and \ref{L:existlim}.

\bl{\bf (Extremality of $\Lambda_c$)}\label{L:extremal}
For all $c\geq 0$, we have $\Lambda_c \in (\Ii\cap\Ti)_{\rm e}$.
\el
{\bf Proof.}  The proof is the same as that of \cite[Lemma IX.2.9]{Lig05}. We include it here for the reader's
convenience. Assume that $\Lambda_c = \alpha \mu_1 + (1-\alpha)\mu_2$ for some $\alpha\in (0,1)$ and
$\mu_1, \mu_2\in \Ii\cap\Ti$ with $\int \rho([0,1])\mu_i(\di \rho) = c_i$, where $c=\alpha c_1+(1-\alpha)c_2$.
Then for any $\zeta_0\in \Ci_{\rm c}(\R)$ and $i=1,2$,
\begin{eqnarray}
 \int e^{ -\int\zeta_0(x)\rho_0(\di x)} \mu_i(\di \rho_0) &=&  \int \E[e^{ -\int\zeta_0(x) \rho_t(\di x)}] \mu_i(\di \rho_0)
= \int  \E[e^{ -\int \zeta_t(x) \rho_0(\di x)}] \mu_i(\di \rho_0) \nonumber \\
&\geq&  \E[e^{ -\iint \zeta_t(x) \rho_0(\di x) \mu_i(\di \rho_0)}] = \E[e^{ -c_i\int\zeta_t(x)\di x}]  \label{extre} \\
&=& \E\big[e^{-\frac{c_i}{c}\int \zeta_0(x) \rho_t(\di x) } \big| \rho_0 \equiv c\big]
\underset{t\to\infty}{\longrightarrow} \int e^{-\frac{c_i}{c}\int \zeta_0(x) \rho(\di x)} \Lambda_c(\di \rho), \nonumber
\end{eqnarray}
where we used the fact that $\mu_i\in \Ii\cap\Ti$, duality, Jensen's inequality, and Lemma \ref{L:existlim}.
Denote $\phi(a) = \int e^{-a\int\zeta_0(x) \rho(\di x)} \Lambda_c(\di \rho)$. Since $\Lambda_c=\alpha \mu_1 + (1-\alpha) \mu_2$,
(\ref{extre}) implies that
\be\label{extre2}
\phi(1) \geq \alpha \phi\big(\frac{c_1}{c}\big) + (1-\alpha)\phi\big(\frac{c_2}{c}\big).
\ee
If $\Lambda_c=\delta_0$, then the extremality of $\Lambda_c$ is trivial; otherwise we can find $\zeta_0\in \Ci_{\rm c}(\R)$ such that $\phi$ is
strictly convex, which implies equality in (\ref{extre2}) and hence $c_1=c_2=c$. Therefore we have equality in (\ref{extre2}) for all
$\zeta_0\in \Ci_{\rm c}(\R)$, and we can then deduce from (\ref{extre}) that
$$
\int e^{ -\int \zeta_0(x) \rho_0(\di x)} \mu_i(\di \rho_0)=\int e^{ -\int\zeta_0(x) \rho(\di x)} \Lambda_c(\di \rho),
$$
which implies that $\mu_i = \Lambda_c$.
\qed
\medskip

We remark that Lemmas \ref{L:existlim} and \ref{L:extremal} can also be deduced from the convergence to invariant laws proved below using
coupling. However, the proof by duality illustrates a useful tool.

\subsection{Second moment calculations}

Following \cite{LS81} and \cite{Lig05}, we first introduce for each $c>0$ the subset of translation invariant probability laws $\Ti_c \subset \Ti$,
where $\Gamma \in \Ti$ is in $\Ti_c$ if and only if $\int \rho([0,1]) \Gamma(\di \rho)=c$, $\int \rho([0,1])^2 \Gamma(\di \rho) < \infty$,
and
\be \label{Tic}
\int \Big(\int \frac{1}{\sqrt{2\pi t}} e^{-\frac{x^2}{2t}} \rho(\di x) -c\Big)^2 \Gamma(\di \rho) \asto{t} 0.
\ee
For a Howitt-Warren process with initial law $\Li(\rho_0)\in \Ti_c$, we can perform second moment calculations for $\rho_t$ as $t\to\infty$
(see Lemma \ref{L:2pt}). Furthermore, if $\Li(\rho_0)\in \Ti_c$, then any weak limit of $\rho_t$ as $t\to\infty$ is also in $\Ti_c$ (see Corollary
\ref{C:limTi1}), which will be crucial for proving the convergence of $\Li(\rho_t)$ to the homogeneous ergodic law $\Lambda_c$. First we note that

\bl\label{L:ergTic}
If $\Gamma \in \Ti_{\rm e}$, $\int \rho([0,1])\Gamma(\di \rho) =c>0$ and $\int \rho([0,1])^2\Gamma(\di \rho) <\infty$, then $\Gamma\in \Ti_c$.
Conversely, any $\Gamma\in \Ti_c$ is a mixture of laws in $\Ti_{\rm e}$ satisfying the conditions above.
\el
{\bf Proof.} Our assumption implies that
\be \label{rhoL2to0}
\int \Big(\frac{\rho([-L,L])}{2L}-c\Big)^2 \Gamma(\di \rho) \asto{L} 0
\ee
by the $L_2$ ergodic theorem. By the layercake representation,
\begin{eqnarray*}
\int_\R \frac{1}{\sqrt{2\pi t}} e^{-\frac{x^2}{2t}} \rho(\di x) -c
&=& \frac{1}{\sqrt{2\pi}} \int_\R e^{-\frac{y^2}{2}}\Big(\frac{\rho(\sqrt{t}\,\di y)}{\sqrt t} -c\,\di  y\Big) \\
&=& \frac{1}{\sqrt{2\pi}} \int_\R \int_0^1 1_{\{z<e^{-\frac{y^2}{2}}\}} \di z \Big(\frac{\rho(\sqrt{t}\,\di y)}{\sqrt t} -c\,\di  y\Big) \\
&=& \frac{1}{\sqrt{2\pi}} \int_0^1  \Big(\frac{\rho([-\sqrt{-2t\ln z}, \sqrt{-2t\ln z}])}{2\sqrt{-2t\ln z}}-c\Big)
2\sqrt{-2\ln z}\ \di z,
\end{eqnarray*}
where we note that $\frac{2\sqrt{-2\ln z}}{\sqrt{2\pi}}\, \di z$ is a probability distribution on $[0,1]$ independent of $t$. Substituting
this representation into the left hand side of (\ref{Tic}), applying the H\"older inequality with respect to
$\frac{2\sqrt{-2\ln z}}{\sqrt{2\pi}}\, \di z$, and applying (\ref{rhoL2to0}) then proves (\ref{Tic}). A more
general argument using Bochner's theorem can be found in the proof of \cite[Theorem 5.6]{Lig73} or \cite[Corollary II.8.20]{Lig05}. The second statement
in Lemma \ref{L:ergTic} follows from the ergodic decomposition of $\Gamma\in \Ti_c$.
\qed

\bl{\bf (Second moment calculation)}\label{L:2pt}
Let $(\rho_t)_{t\geq 0}$ be a Howitt-Warren process with drift $\beta\in\R$ and characteristic measure $\nu\neq 0$. If $\Li(\rho_0)\in \Ti_1$,
then for all $\phi, \psi\in C_c(\R)$, we have
\be\label{2pt}
\lim_{t\to\infty} \E\Big[\int \phi(x) \rho_t(\di x) \int \psi(y) \rho_t(\di y)\Big]
= \int \phi(x) \di x \int \psi(y) \di y + \frac{\int \phi(x)\psi(x) \di x}{2\nu([0,1])}.
\ee
\el
{\bf Proof.} We may assume $\phi, \psi\in C_c(\R)$ are non-negative. Since such functions can be approximated from above and below by finite
linear combinations of indicator functions of finite intervals, it suffices to prove (\ref{2pt}) for $\phi=1_{I_1}$ and $\psi=1_{I_2}$ for
some finite intervals $I_1$ and $I_2$. Since $2\phi\psi=(\phi+\psi)^2-\phi^2-\psi^2$, it suffices to consider only $I_1=I_2$, and we may
even take $I_1=[0,1]$, so that (\ref{2pt}) reduces to showing
\be\label{2pt1}
\lim_{t\to\infty}\E[\rho_t([0,1])^2] = 1 + \frac{1}{2\nu([0,1])}.
\ee
By Theorem~\ref{T:HWconst}, we have the representation
\be\label{2pt2}
\E[\rho_t([0,1])^2] = \E\Big[\iint \Q^{\otimes 2}\Big(\big(\pi^{1,+}_{(x,0)}(t), \pi^{2,+}_{(y,0)}(t)\big)\in [0,1]^2\Big) \rho_0(\di x)\rho_0(\di y)\Big],
\ee
where $\Q^{\otimes 2}$ denotes the $2$-fold product measure, and $\pi^{1,+}_{(\cdot,\cdot)}$ resp.\ $\pi^{2,+}_{(\cdot,\cdot)}$ are rightmost elements
in two independent sample webs $\Wi^1$ resp.\ $\Wi^2$, both with quenched law $\Q$. With respect to $\E[\Q^{\otimes 2}]$,
$(\pi^{1,+}{(x,0)}, \pi^{2,+}_{(y,0)})$ is the two-point motion of the Howitt-Warren flow with drift $\beta$ and characteristic measure $\nu$, and hence
solves the Howitt-Warren martingale problem under conditions (\ref{MP1b}) and (\ref{thetacouple}). In particular, $R_t:= \pi^{2,+}_{(y,0)}(t)-\pi^{1,+}_{(x,0)}(t)$
is an autonomous Brownian motion with stickiness at the origin, and conditional on $(R_t)_{t\geq 0}$, $S_t := \pi^{1,+}_{(x,0)}(t)+\pi^{2,+}_{(y,0)}(t)$ is
distributed as a time change of an independent Brownian motion with drift $2\beta$. We leave the verification of this statement as an exercise to
the reader. A similar statement for a pair of Brownian motions satisfying (\ref{lrsde}) can be found in \cite[Lemma~2.2]{SS08}.

Let $\bar \rho_0^{\otimes 2}:= \E[\rho_0^{\otimes 2}]$, and let $\hat\rho_0^{\otimes2}$ denote the image measure of $\bar\rho_0^{\otimes 2}$ under the change
of coordinates $(x,y) \to (r,s):=(y-x, x+y)$. Then by Lemma \ref{L:ergod} below,
$\hat\rho_0^{\otimes2}(\di r\,\di s) = \alpha(\di r)\,\di s$. Therefore we can rewrite (\ref{2pt2}) as
\begin{eqnarray}
\E[\rho_t([0,1])^2] &=& \iint \P_{(r,s)} \big(|R_t|\leq 1, |R_t|\leq S_t \leq 2-|R_t|\big)\, \alpha(\di r)\, \di s \nonumber \\
&=& 2\int \E_r\big[(1-|R_t|)1_{\{|R_t|\leq 1\}}\big]\, \alpha(\di r) \nonumber \\
&=& 2\int_0^1 \int_\R \P_r(|R_t|\leq a)\, \alpha(\di r)\ \di a,                 \label{2pt3}
\end{eqnarray}
where $\P_{(r,s)}$ denotes probability for $(R_t, S_t)$ starting at $(r,s)$, and we have used the fact that conditioned on $(R_t)_{t\geq 0}$,
$(S_t)_{t\geq 0}$ with differential initial conditions can be coupled together simply by translation.

In (\ref{2pt3}), let $f_{t,a}(r):=\P_r(|R_t|\leq a)$, which is even, and strictly decreasing on $[0,\infty)$. The latter follows from the
fact that $|R_t|$ is a reflected Brownian motion with stickiness at the origin, and there is a natural coupling through coalescence for $|R_t|$ starting at
differential initial conditions on $[0,\infty)$. By the layercake representation, we may rewrite (\ref{2pt3}) as
\begin{eqnarray}
\E[\rho_t([0,1])^2] &=& 2 \int_0^1  \int_\R \int_0^\infty 1_{\{y\leq f_{t,a}(r)\}} \di y\, \alpha(\di r)\, \di a  \nonumber \\
&=& 2\int_0^1 \int_0^\infty \alpha([-f^{-1}_{t,a}(y), f^{-1}_{t,a}(y)]) \di y\, \di a \nonumber \\
&=& 2\int_0^1 \int_0^\infty -\frac{\alpha([-r,r])}{r}\, r\, \di f_{t,a}(r) \, \di a.    \label{2pt4}
\end{eqnarray}
Note that $-r \di f_{t,a}(r)$ is a finite measure on $(0,\infty)$ with total mass $\int_0^\infty f_{t,a}(r)\di r$, and for
any $u>0$, integrating by parts gives
$$
\int_0^u -r \di f_{t,a}(r) = -u f_{t,a}(u) + \int_0^u f_{t,a}(r)\di r \asto{t} 0,
$$
since $f_{t,a}(r)$ is decreasing on $[0,\infty)$ and $f_{t,a}(0)\to 0$ as $t\to\infty$ by basic properties of $|R_t|$. Therefore the sequence of
measures $-r \di f_{t,a}(r)$ shifts its mass to $\infty$ as $t\to\infty$. Since $\lim_{r\to\infty} \frac{\alpha([-r,r])}{r} = 1$ by
Lemma \ref{L:ergod}  below, we deduce from (\ref{2pt4}) that
\be\label{2pt5}
\lim_{t\to\infty} \E[\rho_t([0,1])^2] = 2\int_0^1 \lim_{t\to\infty} \int_0^\infty f_{t,a}(r)\di r \ \di a.
\ee
Since $\di r + \frac{1}{4\nu([0,1])}\delta_0(r)$ is an invariant measure for $|R_t|$ by Lemma \ref{L:2ptinv}, for any $t\geq 0$,
$$
\frac{1}{4\nu([0,1])}f_{t,a}(0) + \int_0^\infty f_{t,a}(r)\,\di r = \frac{1}{4\nu([0,1])}\P_0(|R_t|\leq a) + \int_0^\infty \P_r(|R_t|\leq a)\,\di r = \frac{1}{4\nu([0,1])} + a.
$$
Since $f_{t,a}(0)\to 0$ as $t\to\infty$, we obtain $\lim_{t\to\infty} \int_0^\infty f_{t,a}(r)\,\di r = \frac{1}{4\nu([0,1])} + a$.
Substituting this into (\ref{2pt5}) then gives (\ref{2pt1}).
\qed

\bl{\bf (Ergodicity of the second moment measure)}\label{L:ergod} Let $\Li(\rho)\in \Ti_1$. Let $\bar\rho^{\otimes2}:=\E[\rho^{\otimes2}]$ denote
the second moment measure of $\rho$, and let $\hat \rho^{\otimes2}$ denote the image measure of $\bar\rho^{\otimes2}$ under the change of
coordinates $(x,y) \to (r,s) := (y-x, x+y)$. Then $\hat\rho^{\otimes2}(\di r\,\di s) = \alpha(\di r)\,\di s$,
where $\alpha(A) = \hat\rho^{\otimes2}(A\times [0,1])$ for all $A\in \Bi(\R)$, and $\lim_{L\to\infty} \frac{1}{L}\alpha([-L,L])=1$.
\el
{\bf Proof.} By translation invariance in law of $\rho$, $\bar\rho^{\otimes2}(\cdot) = \bar\rho^{\otimes2}(\cdot+(a,a))$ for all $a\in\R$. Therefore
$\hat\rho^{\otimes2}(\di r\,\di s)$ is translation invariant in $s$, and hence $\hat\rho^{\otimes2}$ admits the desired
factorization. Therefore
$$
\alpha([-L,L]) = \int 1_{[-L,L]}(r) 1_{[0,1]}(s) \hat\rho^{\otimes2}(\di r\, \di s) = \int 1_{[-L,L]}(r) 1_{[0,1]}(r+s) \hat\rho^{\otimes2}(\di r\, \di s).
$$
Transforming back into the variables $(x,y)$ and the measure $\bar\rho^{\otimes2}$ then gives
$$
\alpha([-L,L]) = \bar\rho^{\otimes2}\big(\{(x,y) : 0\leq y\leq 1/2, |x-y|\leq L \}\big),
$$
which is bounded between $\bar\rho^\otimes([-L+1, L-1]\times [0,1/2])$ and $\bar\rho^\otimes([-L-1, L+1]\times [0,1/2])$. Note that
\begin{eqnarray}
\frac{1}{L}\bar\rho^\otimes([-L+1, L-1]\times [0,1/2]) \!\!\!\!\! &=&\!\!\!\!\!
\frac{1}{L}\E\big[\rho([0,1/2]) \rho([-L+1, L-1])\big] \nonumber \\
&\asto{L}& 2\, \E[\rho([0,1/2])] = 1, \label{2pt6}
\end{eqnarray}
provided that $\frac{1}{2L} \rho([-L+1, L-1])\to 1$ in $L_2$. Indeed, by Lemma \ref{L:ergTic}, there exists a probability measure
$\gamma(\di \Lambda)$ on $\Ti_1\cap \Ti_{\rm e}$ such that $\Li(\rho) = \int_{\Ti_1\cap\Ti_{\rm e}} \Lambda \gamma(\di \Lambda)$. Therefore
\begin{eqnarray*}
&& \E\Big[\Big(\frac{1}{2L} \rho([-L+1, L-1])-1\Big)^2\Big] \\
&=& \int_{\Ti_1\cap\Ti_{\rm e}}\int \Big(\frac{1}{2L} \rho([-L+1, L-1])-1\Big)^2\Lambda(\di \rho) \gamma(\di \Lambda) \asto{L} 0,
\end{eqnarray*}
since the integrand w.r.t.\ $\gamma(\di \Lambda)$ tends to $0$ $\gamma$ a.s.\ by the $L_2$ ergodic theorem applied to $\Lambda$, and is
dominated uniformly in $L$ by $2+2\int \rho([0,1])^2\Lambda(\di \rho)$, which is integrable by the assumption that
$\E[\rho([0,1])^2]<\infty$. This proves (\ref{2pt6}), and the same can be proved for $\bar\rho^\otimes([-L-1,L+1]\times [0,1/2])$,
the upper bound on $\alpha[-L,L]$. Therefore $\lim_{L\to\infty} \frac{1}{L}\alpha([-L,L])=1$.
\qed

\bl\label{L:2ptinv}{\bf (Invariant measure for the two point motion)} Let $(X_t,Y_t)$ be two coupled Brownian motions solving the Howitt-Warren
martingale problem under conditions (\ref{MP1b}) and (\ref{thetacouple}). Let $R_t=Y_t-X_t$. Then $\di r + \frac{1}{4\nu([0,1])}\delta_0(r)$ on
$[0,\infty)$ is an invariant measure for $|R_t|$, and $\di x\,\di y + \frac{1}{2\nu([0,1])} \delta_x(y)\di x$ on $\R^2$ is an invariant measure
for $(X_t, Y_t)$.
\el
{\bf Proof.} Note that $|R_t|$ is uniquely characterized in law by the following two properties: (1)
$|R_t|- 4\nu([0,1])\int_0^t 1_{\{|R_s|=0\}}\di s$ is a martingale; (2) $\langle |R_t|, |R_t|\rangle = 2\int_0^t 1_{\{|R_s|\neq 0\}}\di s$.
These two properties are clearly satisfied by the solution of the following SDE
\be\label{SDER}
\di Z_t = 1_{\{Z_t\neq 0\}} \sqrt{2} \di B_t + 1_{\{Z_t=0\}} 4\nu([0,1]) \di t,
\ee
where $B_t$ is a standard Brownian motion, and $Z_t$ is constrained to be non-negative. For the existence and uniqueness of a weak solution
to this SDE, see e.g.\ \cite[Prop.~2.1]{SS08}. The solution of (\ref{SDER}) generates a Feller semigroup $(S_t)_{t\geq 0}$ on the Banach
space $C_0([0,\infty))$, the space of continuous functions on $[0,\infty)$ which vanish at $\infty$ and equipped with the supremum norm.
By It\^o's formula, the generator $L$ for $S_t$ is given by
\be\label{SDEgen}
Lf(x) = 1_{\{x\neq 0\}} f''(x) + 1_{\{x=0\}} 4\nu([0,1]) f'(x).
\ee
Let $\Di\subset C_0([0,\infty))$ denote the domain of $L$. If $f\in \Di$, then $Lf\in C_0([0,\infty))$. In particular, we must have $f''(0)=4\nu([0,1]) f'(0)$
and $f''\in C_0([0,\infty))$. Together with $f\in C_0([0,\infty))$, this also implies that $f'\in C_0([0,\infty))$. Conversely, if
$f\in D:=\{f\in C_0([0,\infty)) : f',f''\in C_0([0,\infty)) \ \mbox{and}\ f''(0)=4\nu([0,1]) f'(0)\}$, then it is not difficult to see from It\^o's formula
that $f\in \Di$. Therefore $\Di=D$. If we denote $\mu(\di x)= \di x+ \frac{1}{4\nu([0,1])}\delta_0(x)$, then we have
\be\label{invgen}
\int_0^\infty Lf(x) \mu(\di x) = \int_0^\infty  f''(x) \di x + f'(0) = 0 \qquad \mbox{ for all } f\in \Di.
\ee
From (\ref{invgen}), we can deduce that $\mu$ is an invariant measure for $(|R_t|)_{t\geq 0}$. Indeed, for any $f\in \Di\cap C_c([0,\infty))$ with
compact support and for any $t>0$, we have
\begin{eqnarray}
&& \int_0^\infty S_tf(x) \mu(\di x) - \int_0^\infty f(x) \mu(\di x)
= \int_0^\infty (S_tf(x)-f(x)) \mu(\di x)   \nonumber \\
&=& \int_0^\infty \int_0^t \frac{\rm d\ }{\di s}S_sf(x) \di s\, \mu(\di x)
= \int_0^\infty \int_0^t LS_sf(x) \di s\, \mu(\di x) \nonumber \\
&=& \int_0^t \int_0^\infty LS_sf(x) \mu(\di x)\, \di s = 0, \label{laplace}
\end{eqnarray}
where we have used Fubini based on the fact that $LS_sf(x)=S_sLf(x)$ is decaying super-exponentially in $x$ because $Lf$ has compact support and
$|R_t|$ is distributed as a Brownian motion on $(0,\infty)$; and in (\ref{laplace}), we have  applied (\ref{invgen}) using the fact that $f\in \Di$
implies $S_tf\in \Di$. Since (\ref{laplace}) holds for all $f\in \Di\cap C_c([0,\infty))$, which is a measure determining class, $\mu$ is an invariant
measure.

The symmetry of $R_t$ implies that $\ti\mu(\di r):=\di r + \frac{1}{2\nu([0,1])} \delta_0(r)$ on $\R$ is an invariant measure for $R_t$, and
the translation invariance of $(X_t, Y_t)$ along the diagonal implies that $\ti\mu(\di r)\,\di s$ is an invariant measure for
$(R_t, S_t)$ with $R_t= Y_t-X_t$ and $S_t=Y_t+X_t$. A change of coordinates then verifies that
$\di x\,\di y + \frac{1}{2\nu([0,1])} \delta_x(y)\di x$ is an invariant measure for $(X_t, Y_t)$.
\qed
\medskip

\noi
{\bf Remark.} In \cite{LS81} and \cite{Lig05}, the analogue of Lemma \ref{L:2pt} is proved by treating the two-point motion as a perturbation of two independent
one-point motions. This approach requires exact calculations involving the two-point motion and is not clear how to implement in the continuous space
setting. Our approach reduces the task to first identifying the invariant measure for the two-point motion, which when integrated over the test function
$\phi(x)\psi(y)$ gives the RHS of (\ref{2pt}), and then using qualitative properties of the two-point motion together with the ergodicity of the initial
condition to remove the dependence on the initial condition.
\medskip

The following corollary of Lemma \ref{L:2pt} is the analogue of \cite[Lemma (5.3)(b)]{LS81} for our model, and will be crucial in proving convergence
to the homogeneous invariant laws.
\bcor\label{C:limTi1}{\bf (Preservation of $\Ti_1$)}
Assume the same conditions as in Lemma \ref{L:2pt}. Then the law of any subsequential weak limit of $(\rho_t)_{t\geq 0}$ is also in $\Ti_1$.
\ecor
{\bf Proof.} Let $\rho^*$ be the weak limit of $\rho_{t_n}$ along a subsequence $t_n\uparrow \infty$. Clearly $\Li(\rho^*)\in \Ti$. Lemma \ref{L:2pt}
implies the uniform integrability of $(\rho_{t_n}([0,1]))_{n\in\N}$, and hence $\E[\rho^*([0,1])]=1$ since $\E[\rho_t([0,1])]=1$ for all $t\geq 0$.
By Fatou's lemma, Lemma \ref{L:2pt} also implies that for all non-negative $\phi,\psi\in C_c(\R)$,
\be
\E\Big[\int \phi(x) \rho^*(\di x) \int \psi(y) \rho^*(\di y)\Big]
\leq \int \phi(x) \di x \int \psi(y) \di y + \frac{\int \phi(x)\psi(x) \di x}{2\nu([0,1])}.
\ee
Approximating $\frac{1}{\sqrt{2\pi t}}e^{-\frac{x^2}{2t}}$ from below by functions with compact support, we then have
$$
\E\Big[ \Big(\int \frac{1}{\sqrt{2\pi t}} e^{-\frac{x^2}{2t}} \rho^*(\di x)\Big)^2\Big] \leq 1+ \frac{1}{4\nu([0,1])\sqrt{\pi t}}.
$$
Together with $\E\big[\int \frac{1}{\sqrt{2\pi t}} e^{-\frac{x^2}{2t}} \rho^*(\di x)\big]=1$, this implies (\ref{Tic}) with $c=1$ and $\Gamma=\Li(\rho^*)$.
Therefore $\Li(\rho^*)\in \Ti_1$.
\qed

\subsection{Coupling and convergence}

The definition of the Howitt-Warren process $(\rho_t)_{t\geq 0}$ from the kernels $(K^+_{s,t})_{s<t}$ constructed in Theorem~\ref{T:HWconst} gives
a natural coupling between $(\rho_t)_{t\geq 0}$ with different initial conditions. This coupling is monotone in the sense that if $\rho^1_0(A)\geq \rho^2_0(A)$
for all $A\in \Bi(\R)$, which we denote by $\rho^1_0 \succ \rho^2_0$, then $\rho^1_t \succ \rho^2_t$ a.s.\ for all $t>0$. Through this coupling, we will prove the weak convergence of $\rho_t$ to a mixture of homogeneous invariant laws under suitable assumption on $\Li(\rho_0)$. The first observation is the following.

\bl\label{L:coupling} {\bf (Coupled Howitt-Warren processs)} Let $(\rho^1_t)_{t\geq 0}$ and $(\rho^2_t)_{t\geq 0}$ be Howitt-Warren
processes with drift $\beta$ and characteristic measure $\nu$, defined from the same Howiit-Warren flow $(K^+_{s,t})_{s<t}$ as in (\ref{rho}). Assume that
$\Li(\rho^1_0), \Li(\rho^2_0)\in \Ti$ and $\E[\rho^1_0([0,1])]<\infty$, $\E[\rho^2_0([0,1])]<\infty$. Then any weak limit point
$(\rho^{1*}, \rho^{2*})$ of $(\rho^1_t, \rho^2_t)_{t\geq 0}$ as $t\to\infty$ satisfies $\P(\rho^{1*} \succ \rho^{2*} \mbox{ or } \rho^{2*} \succ \rho^{1*})=1$.
\el
{\bf Proof.} Let $(\rho^1_t-\rho^2_t) = (\rho^1_t-\rho^2_t)^+ - (\rho^1_t-\rho^2_t)^-$ denote the Jordan decomposition of $\rho^1_t-\rho^2_t$, and let $|\rho^1_t-\rho^2_t|:=(\rho^1_t-\rho^2_t)^+ + (\rho^1_t-\rho^2_t)^-$ denote the total variation measure of $\rho^1_t-\rho^2_t$. Recall the quenched law $\Q$
from (\ref{HWquen}). For any $0<s<t$, almost surely we have
$$
|\rho^1_t-\rho^2_t| =  \Big|\int (\rho^1_s-\rho^2_s)(\di x)\, \Q[\pi^+_{(x,s)}(t) \in\cdot ]\Big| \prec \int |\rho^1_s-\rho^2_s|(\di x)\, \Q[\pi^+_{(x,s)}(t) \in \cdot].
$$
Therefore
\begin{eqnarray*}
\E\big[|\rho^1_t-\rho^2_t|([0,1])\big] &\leq& \E\Big[\int |\rho^1_s-\rho^2_s|(\di x)\, \Q\big[\pi^+_{(x,s)}(t) \in [0,1]\big] \Big] \\
&=& \int \E[|\rho^1_s-\rho^2_s|](\di x)\, \E\big[\Q\big[\pi^+_{(x,s)}(t) \in [0,1]\big]\big]
= \E[|\rho^1_s-\rho^2_s|([0,1])],
\end{eqnarray*}
where we have used the independence between $\rho^1_s,\rho^2_s$ and $\big(\Q\big[\pi^+_{(x,s)}(t) \in [0,1]\big]\big)_{x\in\R}$, the translation
invariance of $\E|\rho^1_s-\rho^2_s|$, and the fact that $\pi^+_{(x,s)}$ is distributed as a standard Brownian motion starting at $x$ at time $s$
under the law $\E\Q$. Therefore $\E\big[|\rho^1_t-\rho^2_t|([0,1])\big]$ decreases monotonically to a non-negative limit as $t\uparrow \infty$.

Note that the lemma follows once we show that for any $\eps>0$ and any $\phi,\psi\in C_c(\R)$ with $0\leq \phi,\psi\leq 1$, we have
\be\label{discre1}
\P\Big( \int \phi(x)\rho^1_t(\di x)-\! \int \phi(x)\rho^2_t(\di x)>\eps \mbox{ and } \int \psi(x)\rho^2_t(\di x)-\! \int \psi(x)\rho^1_t(\di x)>\eps\Big)
\asto{t} 0.
\ee
Suppose that (\ref{discre1}) fails so that for some $\eps>0$ and $\phi, \psi\in C_c(\R)$ with $0\leq \phi,\psi\leq 1$, the probability in
(\ref{discre1}) is bounded uniformly from below by $\delta>0$ along a sequence $t_i\uparrow\infty$. Choose $L>0$ large such that $\phi$ and $\psi$
vanish outside $[-L,L]$. Given $\rho^1_{t_i}$ and $\rho^2_{t_i}$ satisfying the conditions in the probability in (\ref{discre1}), we have
$$
\eps < \int \phi(x) (\rho^1_{t_i}-\rho^2_{t_i})(\di x) \leq \int \phi(x) (\rho^1_{t_i}-\rho^2_{t_i})^+(\di x) \leq (\rho^1_{t_i}-\rho^2_{t_i})^+([-L,L]),
$$
and similarly $(\rho^1_{t_i}-\rho^2_{t_i})^-([-L,L])>\eps$. For such a realization of $\rho^1_{t_i}$ and $\rho^2_{t_i}$,
\begin{eqnarray*}
&& |\rho^1_{t_i+1}-\rho^2_{t_i+1}|([0,1])  \\
&\leq&\! \int |\rho^1_{t_i}-\rho^2_{t_i}|(\di x) \Q\big[\pi^+_{(x,t_i)}(t_i+1) \in [0,1]\big]
 - 2\eps \Q\big[\pi^+_{(-L,t_i)}(t_i+1)
=\pi^+_{(L,t_i)}(t_i+1)\in [0,1]\big],
\end{eqnarray*}
where we observed that the mass assigned by $(\rho^1_{t_i}-\rho^2_{t_i})^+$ and $(\rho^1_{t_i}-\rho^2_{t_i})^-$ to $[-L,L]$ are carried by
$(\pi^+_{(x,t_i)}(t_i+1))_{x\in \R}$ to the same point in $[0,1]$ when $\pi^+_{(-L,t_i)}(t_i+1)=\pi^+_{(L,t_i)}(t_i+1)\in [0,1]$. Recall that the event in
(\ref{discre1}) is assumed to have probability at least $\delta$ along $(t_i)_{i\in\N}$, we thus have
\be\label{discre2}
\E\big[|\rho^1_{t_i+1}-\rho^2_{t_i+1}|([0,1])\big] \leq \E\big[|\rho^1_{t_i}-\rho^2_{t_i}|([0,1])\big] - 2\eps \delta h,
\ee
where $h= \E\Q\big[\pi^+_{(-L,t_i)}(t_i+1)=\pi^+_{(L,t_i)}(t_i+1)\in [0,1]\big] >0$ is independent of $t_i$. Since (\ref{discre2}) holds for all
$t_i$, this contradicts the fact that $\E\big[|\rho^1_t-\rho^2_t|([0,1])\big]$ decreases monotonically to a non-negative limit as $t\uparrow \infty$.
\qed

\bl{\bf (Convergence to $\Lambda_c$)}\label{L:conv}
Let $(\rho_t)_{t\geq 0}$ be a Howitt-Warren process with drift $\beta$ and characteristic measure $\nu\neq 0$. If $\Li(\rho_0)\in \Ti_{\rm e}$ and
$\E[\rho_0([0,1])]=c<\infty$, then $\rho_t$ converges weakly to $\Lambda_c$, which was defined in Lemma \ref{L:existlim}. If $\Li(\rho_0)\in \Ti_{\rm e}$ and
$\E[\rho_0([0,1])]=\infty$, then $\rho_t$ has no weak limit which is supported on $\Mi_{\rm loc}(\R)$.
\el
{\bf Proof.} Without loss of generality, assume $c=1$. Let $(\rho^2_t)_{t\geq 0}$ be a Howitt-Warren process with initial condition $\rho^2_0$
such that $\Li(\rho^2_0) = \Lambda_1$, and let $\rho_t$ and $\rho^2_t$ be defined from the same Howitt-Warren flow $(K^+_{s,t})_{s<t}$. As in
the proof of Lemma \ref{L:existlim}, we note that $(\rho_t)_{t\geq 0}$ (resp.\ $(\rho_t, \rho^2_t)_{t\geq 0}$) is a tight family of $\Mi_{\rm g}(\R)$
(resp.\ $\Mi_{\rm g}(\R)^2$) valued random variable. Therefore any subsequential weak limit $\Gamma\in\Mi_1(\Mi_{\rm g}(\R))$ of $\Li(\rho_t)$ as
$t\to\infty$ can be realized as the marginal law of the first component of a random couple $(\rho^*, \rho^{2*})\in \Mi_{\rm g}(\R)^2$, which
arises as a subsequential weak limit of $(\rho_t, \rho^2_t)_{t\geq 0}$. By Lemma \ref{L:coupling},
$\P(\rho^* \succ \rho^{2*} \mbox{ or } \rho^{2*} \succ \rho^*)=1$.  Therefore for any rational $a<b$, by the translation invariance in law of
$(\rho^*-\rho^{2*}) 1_{\{\rho^* \succ \rho^{2*}\}}$ and $(\rho^{2*}-\rho^*) 1_{\{\rho^{2*} \succ \rho^*\}}$, we have
\begin{eqnarray}
&& \E|\rho^*([a,b))- \rho^{2*}([a,b))| = \E\Big|\int \frac{b-a}{\sqrt{2\pi t}} e^{-\frac{x^2}{2t}} (\rho^*-\rho^{2*})(\di x)\Big| \nonumber\\
&\leq& (b-a) \E\Big|\int\frac{1}{\sqrt{2\pi t}} e^{-\frac{x^2}{2t}}\rho^*(\di x) -1\Big|
+ (b-a)\E\Big|\int\frac{1}{\sqrt{2\pi t}} e^{-\frac{x^2}{2t}}\rho^{2*}(\di x) -1\Big|. \label{conv1}
\end{eqnarray}
If we first restrict ourselves to the case $\E[\rho_0([0,1])^2]<\infty$, then $\Li(\rho_0)\in \Ti_1$ by Lemma \ref{L:ergTic}, and by Corollary
\ref{C:limTi1}, $\Li(\rho^*)\in \Ti_1$ and $\Li(\rho^{2*})=\Lambda_1 \in \Ti_1$. Definition of $\Ti_1$ implies that both terms in
(\ref{conv1}) vanish as $t\to\infty$, and hence $\rho^*([a,b))=\rho^{2*}([a,b))$ a.s.\ for all rational $a<b$. Since $\{[a,b)\}_{a<b\in\Q}$
is measure determining, we have $\rho^*=\rho^{2*}$ a.s., and hence $\Li(\rho_t)$ converges weakly to $\Lambda_1$.

If $\E[\rho_0([0,1])^2]=\infty$, then we can approximate $\rho_0$ by $(\rho^n_0)_{n\in\N}$ with $\Li(\rho^n_0)\in\Ti_{\rm e}$ such that
$\E[\rho^n_0([0,1])^2]<\infty$ and $\rho^n_0$ increases monotonically to $\rho_0$ almost surely. For instance, given $\rho_0$, we can
sample a uniform random variable $U$ on $[0,1]$ and then define $\rho_0^n$ on $[U+k, U+k+1)$ for each $k\in\Z$ by $\rho_0^n=\rho_0$ on $[U+k, U+k+1)$
if $\rho_0([U+k, U+k+1))\leq n$, and set $\rho_0^n=0$ on $[U+k, U+k+1)$ otherwise. Then $\E[\rho^n_0([0,1])]=1-\eps_n$ for some
$\eps_n\downarrow 0$. Our argument above shows that $\Li(\rho^n_t)$ converges weakly to $\Lambda_{1-\eps_n}$. Since $\rho_0\succ \rho^n_0$
a.s.\ for all $n\in\N$, any weak limit point $\Gamma$ of $\Li(\rho_t)_{t\geq 0}$ stochastically dominates $\Lambda_c$ for all $c<1$. Since
$\int\rho([0,1])\Gamma(\di \rho)\leq 1$ by Fatou, we must have $\Gamma=\Lambda_1$.

If $\Li(\rho_0)\in \Ti_{\rm e}$ and $\E[\rho_0([0,1])]=\infty$, then by the same argument as above, any weak limit point of $\rho_t$ stochastically
dominates $\Lambda_c$ for all $c>0$, which is not possible for an $\Mi_{\rm loc}(\R)$-valued random variable since $\Lambda_1$ is not concentrated on the
zero measure by our assumption $\nu\neq 0$.
\qed

From Lemma \ref{L:conv}, we can deduce that
\bl{\bf (Extremal measures in $\Ii\cap \Ti$)}\label{L:extremal2} For the Howitt-Warren processs with drift $\beta$ and characteristic measure
$\nu\neq 0$, we have $(\Ii\cap\Ti)_{\rm e}=\{\Lambda_c : c\geq 0\}$.
\el
{\bf Proof.} If $\Li(\rho_0)\in (\Ii\cap\Ti)_{\rm e}$, then $\Li(\rho_0)$ can be decomposed into measures in $\Ti_{\rm e}$ with different mean densities, which
by Lemma \ref{L:conv} converges to mixtures of $(\Lambda_c)_{c\geq 0}$. Therefore by the extremality of $\rho_0$, we must have $\Li(\rho_0)=\Lambda_c$ for some
$0\leq c<\infty$, and hence $(\Ii\cap\Ti)_{\rm e}\subset\{\Lambda_c : c\geq 0\}$. The converse $\{\Lambda_c : c\geq 0\}\subset(\Ii\cap\Ti)_{\rm e}$ has been established in Lemma
\ref{L:extremal}.
\qed

\subsection{Proof of Theorems \ref{T:HIL}--\ref{T:invsupp}}\label{S:HIL}

{\bf Proof of Theorem \ref{T:HIL}.} Part (a) follows from Lemma \ref{L:extremal2}, where the scaling relation
$\Lambda_c(\di (c\rho))=\Lambda_1(\di \rho)$ is trivial, while (\ref{1stmom}) and (\ref{2ndmom}) follow from Lemma \ref{L:existlim} and
Lemma \ref{L:2pt} applied to $\rho_0(\di x)=\di x$. Parts (b) and (c) follow from Lemma \ref{L:conv} and \ref{L:2pt}, while part
(d) follows from spatial ergodic decomposition and Lemma \ref{L:conv}.
\qed
\bigskip

\noi
{\bf Proof of Theorem \ref{T:invsupp}.} Part (a) follows from Theorem~\ref{T:supp}~(a) and Proposition~\ref{P:braco}~(c). Part (b) follows from
Theorem~\ref{T:atom}~(a).
\qed

\appendix

\section{The Howitt-Warren martingale problem}\label{A:HWMP}

Howitt and Warren \cite[Thm~2.1]{HW09a} formulated a martingale problem for a class of sticky Brownian motions on $\R$, for which they showed
that for each deterministic initial state $\vec x\in\R^n$, there exists a unique solution in distribution to their martingale
problem. Moreover, they showed that the family of all solutions to their martingale problem forms a consistent Feller
family \cite[Prop.~8.1]{HW09a}, which defines a stochastic flow of kernels we call a Howitt-Warren flow.
In this appendix, we show that our formulation of the Howitt-Warren martingale problem in Definition~\ref{D:HWMP2} is equivalent to Howitt and Warren's
original formulation in \cite{HW09a}. The advantage of our formulation is that we use a much simpler
set of test functions, which somewhat simplies the proof of the convergence of the $n$-point
motions of discrete Howitt-Warren flows to their continuous counterparts. This convergence result is formulated in
Proposition~\ref{P:nconv}, and is used to verify that the flows we construct in Theorem~\ref{T:HWconst} are indeed
Howitt-Warren flows. A similar convergence result for the $n$-point motions of a continuous time version of the discrete Howitt-Warren flows
was established previously in \cite{HW09a}. We will also give some new parametrizations of Howitt-Warren martingale problems in Lemma \ref{L:para}.

\index{Howitt-Warren!martingale problem}

\subsection{Different formulations}

Let us first recall the original formulation of the Howitt-Warren martingale problem from \cite{HW09a}, and then state two
lemmas that show how one can go from their formulation to ours in Definition~\ref{D:HWMP2} and vice versa. The proof of
these lemmas will be given in the next subsection.

Recall that if $Y$ is a continuous semimartingale, then there exists a unique
continuous process $Y^\Ci$ with bounded variation such that $Y-Y^\Ci$ is a
martingale. The process $Y^\Ci$ is called the {\em compensator} of $Y$. Now if
$Y_1$ and $Y_2$ are continuous, square integrable semimartingales, then by
definition, the {\em covariance process} $\li Y_1,Y_2\re$ of $Y_1$ and $Y_2$ is
the compensator of $(Y_1-Y^\Ci_1)(Y_2-Y^\Ci_2)$, i.e., $\li Y_1,Y_2\re$ is the
unique continuous process of bounded variation such that
\be
t\mapsto\big(Y_1(t)-Y^\Ci_1(t)\big)\big(Y_2(t)-Y^\Ci_2(t)\big)-\li Y_1,Y_2\re(t)
\ee
is a martingale. We generalize our definition of the Howitt-Warren martingale
problem as follows.
\begin{remark}{\bf(Initial states with infinite second moments)}\label{R:infvar}
The solutions to a Howitt-Warren martingale problem (for given $\bet$, $\nu$
and $n$) form a Feller process. Therefore, if $\P_{\vec x}$
denotes the law of the solution of the Howitt-Warren martingale problem with
initial state $\vec x$, and $\rho$ is any probability law on $\R^n$, then
$\int\rho(\di\vec x)\,\P_{\vec x}$ is the law of some Markov process in
$\R^n$. Generalizing Definition~\ref{D:HWMP2}, we may call such a process $\vec
X$ the solution to the Howitt-Warren martingale problem with initial law
$\rho$, even though $\vec X$ is not square integrable if $\rho$
does not have a finite second moment.
\end{remark}

We now turn our attention to the original formulation of Howitt and Warren's
martingale problem in \cite{HW09a}. Recall that our formulation of the
Howitt-Warren martingale in Definition~\ref{D:HWMP2} is based on the constants $(\bet_+(m))_{m\geq 1}$
defined in (\ref{betplusm}). Instead, Howitt and Warren's formulation of their
martingale problem is based on real constants $(\tet(k,l))_{k,l\geq 0}$
satisfying
\be\ba{rll}\label{tetprop}
{\rm(i)}&\tet(k,l)\geq 0\qquad&(k,l\geq 1),\\[5pt]
{\rm(ii)}&\tet(k,l)=\tet(k+1,l)+\tet(k,l+1)\qquad&(k,l\geq 0).
\ec
The $\tet(k,l)$'s are related to the $\bet_+(m)$'s by
\be\label{tetbet}
\bet_+(m)=\tet(0,0)-2\tet(0,m)\qquad(m\geq 1),
\ee
while their relation to the constant $\bet$ and measure $\nu$ is described by
\be\ba{l}\label{tetnu}
\dis\tet(k,l)=\int\nu(\di q)\,q^{k-1}(1-q)^{l-1}\qquad(k,l\geq 1),\\[5pt]
\tet(1,0)-\tet(0,1)=\bet.
\ec
Note that we have now three ways to parametrize Howitt-Warren martingale
problems: we may use the pair $(\bet,\nu)$, the constants $(\bet_+(m))_{m\geq
  1}$, or the constants $(\tet(k,l))_{k,l\geq 0}$. The next lemma shows how to
go from one parametrization to another.
\bl{\bf(Different parametrizations)}\label{L:para}\newline
{\bf(a)} Let $(\tet(k,l))_{k,l\geq 0}$ be real constants satisfying (\ref{tetprop}). Then there exists a unique $\bet\in\R$ and a finite
measure $\nu$ on $[0,1]$ such that (\ref{tetnu}) holds.\med

\noi
{\bf(b)} Let $\bet\in\R$ and let $\nu$ be a finite measure on
$[0,1]$. Then there exists a function $\tet:\N^2\to\R$ satisfying
(\ref{tetprop}) such that (\ref{tetnu}) holds. Any other $\tet'$ satisfies (\ref{tetprop}) and (\ref{tetnu})
if and only if
\be\label{teteq}
\tet'(k,l)=\tet(k,l)+c\big(1_{\{k=0\}}+1_{\{l=0\}}\big)\qquad(k,l\geq 0)
\ee
for some $c\in\R$, and we say that $\theta$ and $\theta'$ are equivalent.
\med

\noi
{\bf(c)} Let $\bet\in\R$ and let $\nu$ be a finite measure on
$[0,1]$. Let $(\tet(k,l))_{k,l\geq 0}$ be real constants satisfying
(\ref{tetprop}), and let $(\bet_+(m))_{m\geq 1}$ be real constants. Then of the
relations (\ref{betplusm}), (\ref{tetbet}), and (\ref{tetnu}), any two imply the
third one.
\el

By definition, a {\em weak total order} on $\{1,\ldots,n\}$ is a relation
$\prec$ such that
\be\ba{rl}
{\rm(i)}&\dis i\prec i,\\
{\rm(ii)}&\dis i\prec j\prec k\mbox{ implies }i\prec k,\\
{\rm(iii)}&\dis\mbox{there exist no }i,j\mbox{ with }i\not\prec j
\mbox{ and }j\not\prec i.
\ec
Each weak total order $\prec$ on $\{1,\ldots,n\}$ defines a nonempty cell
$C_\prec\sub\R^n$ by
\be\label{precel}
C_\prec:=\{\vec x\in\R^n:x_i\leq x_j\mbox{ if and only if }i\prec j\}.
\ee
We note that cells defined by different weak total orders are disjoint, and
that the union of all such cells is $\R^n$. For example:
\be
\{\vec x:x_1<x_3<x_2\},\quad\{\vec x:x_2=x_3<x_1\},
\quad\mbox{and}\quad\{\vec x:x_1=x_2=x_3\}
\ee
are three of the thirteen cells that make up $\R^3$.
Let $L_n$ be the linear space consisting of all continuous real functions
on $\R^n$ that are piecewise linear on each cell $C_\prec$, i.e.,
\be\ba{r@{\,}l}
\dis L_n:=\big\{f:&\dis f\mbox{ is a continuous function }f:\R^n\to\R
\mbox{ such that for each weak total}\\[5pt]
&\dis\mbox{order $\prec$ there exists a linear function }
l:\R^n\to\R\mbox{ with }f=l\mbox{ on }C_\prec\big\}.
\ec
For each $\vec x\in\R^n$, let us define
\be\label{Ranx}
{\rm Ran}(\vec x):=\bigcup_{i=1}^n\{x_i\},
\ee
and for each $x\in{\rm Ran}(\vec x)$, let us write
\be\label{Jx}
J_x=J_x(\vec x):=\big\{i\in\{1,\ldots,n\}:x_i=x\big\}.
\ee
For disjoint $I,J\sub\{1,\ldots,n\}$ let us define a vector
$\vec v_{I,J}\in\R^n$ by
\be\label{vIJ}
v_{I,J}(i):=\left\{\ba{rl}
1\quad&\mbox{if }i\in I,\\
-1\quad&\mbox{if }i\in J,\\
0\quad&\mbox{otherwise.}\ea\right.
\ee
For any $\vec v\in\R^n$, let $\nab_{\vec v}$ denote the one-sided derivative
\be\label{nabv}
\nab_{\vec v}f(\vec x):=\lim_{\eps\down 0}\eps^{-1}
\big(f(\vec x+\eps\vec v)-f(\vec x)\big).
\ee
Let $(\tet(k,l))_{k,l\geq 0}$ be real constants satisfying
(\ref{tetprop}). Then, by definition, $\Ai^\tet_n$ is the linear operator
acting on functions in $L_n$ defined by
\be
\Ai^\tet_nf(\vec x):=\sum_{x\in{\rm Ran}(\vec x)}\sum_{I\sub J_x}
\tet(|I|,|J_x\beh I|)\nab_{\vec v_{I,J_x\beh I}}f(\vec x).
\ee
The original formulation of the Howitt-Warren martingale problem in
\cite{HW09a} differs from our formulation in that formula (\ref{MP2b}) is
replaced by the requirement that for each $f\in L_n$
\be\label{MP2t}
f\big(\vec X(t)\big)-\int_0^t\Ai^\tet_n f(\vec X(s))\di s,
\ee
is a martingale with respect to the filtration generated by $\vec X$.
To see that this is equivalent to the formulation in Definition~\ref{D:HWMP2},
we need the following lemma, the proof of which is not entirely trivial.
\bl{\bf(Action of operator on basis vectors)}\label{L:basis}
Let $f_\De,g_\De$ be defined as in (\ref{fgDe}). Then:\med

\noi
{\bf(a)} The functions
\be
\big\{f_\De:\emptyset\neq\De\sub\{1,\ldots,n\}\big\}
\ee
form a basis for the space $L_n$.\med

\noi
{\bf(b)} Let $(\tet(k,l))_{k,l\geq 0}$ be real constants satisfying
(\ref{tetprop}) and let $(\bet_+(m))_{m\geq 1}$ be given by (\ref{tetbet}). Then
for each nonempty $\De\sub\{1,\ldots,n\}$, one has
\be\label{AfDe}
\Ai^\tet_nf_\De(\vec x)=\bet_+(g_\De(\vec x))\qquad(\vec x\in\R^n).
\ee

\noi
{\bf(c)} If $\tet$ and $\tet'$ satisfy (\ref{tetprop}) and are equivalent in the sense of (\ref{teteq}), then
$\Ai^\tet_n=\Ai^{\tet'}_n$.
\el

\subsection{Proof of the equivalence of formulations}

To prepare for the proof of Lemma~\ref{L:para}, we start with the following
lemma.
\bl{\bf(Moments defining a measure)}\label{L:phi}
Let $(\phi(k,l))_{k,l\geq 0}$ be real constants such that
\be\ba{rl}\label{phiprop}
{\rm(i)}&\phi(k,l)\geq 0,\\[5pt]
{\rm(ii)}&\phi(k,l)=\phi(k+1,l)+\phi(k,l+1)
\ec
for all $k,l\geq 0$. Then there exists a unique finite measure $\nu$ on
$[0,1]$ such that
\be\label{phirep}
\phi(k,l)=\int\nu(\di q)\,q^k(1-q)^l\qquad(k,l\geq 0).
\ee
\el
{\bf Proof.} Let $-\De$ be the operator, acting on sequences of real constants
$(a_k)_{k\geq 0}$ as $((-\De)a)_k:=a_k-a_{k+1}$. Setting $a_k:=\phi(k,0)$, we
observe that $((-\De)a)_k=\phi(k,1)$ $(k\geq 0)$ and more generally
$((-\De)^la)_k=\phi(k,l)\geq 0$ $(k,l\geq 0)$. This qualifies $(a_n)_{n\in\N}$
as a completely monotone sequence, which by \cite[Theorem VII.3.2]{Fel66} can
be represented as $a_k=\int\nu(\di q)q^k$ for some finite measure $\nu$ on
$[0,1]$. Using (\ref{phiprop})~(ii), this implies (\ref{phirep}).\qed

\noi
{\bf Proof of Lemma~\ref{L:para}} Part~(a) is a straightforward consequence of
Lemma~\ref{L:phi}. To prove part~(b), note that by (\ref{tetnu}), $\nu$ uniquely
determines $\theta(k,l)$ for $k,l\geq 1$, which is easily seen to satisfy (\ref{tetprop})
for $k,l\geq 1$. Once $\theta(1,0)$ and $\theta(0,1)$ are chosen, $\theta(k,0)$
and $\theta(0,l)$ for $k,l\geq 0$ are uniquely determined from the recursion relation
(\ref{tetprop})~(ii). Since $\theta(1,0)-\theta(0,1)=\beta$, it follows that $\theta$
is uniquely determined up to the equivalence defined in (\ref{teteq}).

To prove part~(c), we observe that (\ref{betplusm}) and (\ref{tetnu}), together
with (\ref{tetprop})~(ii), imply that
\be
\bet_+(1)=\bet=\tet(1,0)-\tet(0,1)=\tet(0,0)-2\tet(0,1)
\ee
and
\begin{eqnarray*}
\dis\bet_+(m)&=& \dis\bet+2\int\nu(\di q)\sum_{k=1}^{m-1}(1-q)^{k-1}
=\dis\tet(1,0)-\tet(0,1)+2\sum_{k=1}^{m-1}\tet(1,k)\\
&=&\dis\tet(1,0)-\tet(0,1)+2\big(\tet(0,1)-\tet(0,m)\big)= \dis\tet(0,0)-2\tet(0,m)\qquad(m\geq 2).
\end{eqnarray*}
This shows that (\ref{betplusm}) and (\ref{tetnu}) imply (\ref{tetbet}). Running
the argument backward, we also see that (\ref{tetnu}) and (\ref{tetbet}) imply
(\ref{betplusm}). Finally, (\ref{betplusm}) and (\ref{tetbet}) imply that
\bc\label{vartetnu}
\dis\tet(0,0)-2\tet(0,1)&=&\dis\bet,\\[5pt]
\dis\tet(0,0)-2\tet(0,m)
&=&\dis\bet+2\int q^{-1}\big(1-(1-q)^{m-1}\big)\nu(\di q)
\qquad(m\geq 2),
\ec
from which it is not hard to derive (\ref{tetnu}) using (\ref{tetprop})~(ii).
\qed

\detail{Formula (\ref{vartetnu}) implies
\be
\tet(1,0)-\tet(0,1)=\tet(1,0)+\tet(0,1)-2\tet(0,1)=\bet
\ee
and
\be
\tet(0,1)-\tet(0,m)=\int q^{-1}\big(1-(1-q)^{m-1}\big)\nu(\di q)
\qquad(m\geq 2).
\ee
It follows that
\be
\tet(1,1)=\tet(0,1)-\tet(0,2)=\int q^{-1}\big(1-(1-q)^{2-1}\big)\nu(\di q)
=\int\nu(\di q)
\ee
and, for $m\geq 2$,
\be\ba{l}
\dis\tet(1,m)=\tet(0,m)-\tet(0,m+1)
=\int\nu(\di q)\,\Big(q^{-1}\big(1-(1-q)^m\big)
-q^{-1}\big(1-(1-q)^{m-1}\big)\Big)\\[5pt]
\dis\qquad=\int\nu(\di q)q^{-1}\,\big((1-q)^{m-1}-(1-q)^m\big)\\[5pt]
\dis\qquad=\int\nu(\di q)q^{-1}(1-q)^{m-1}\big(1-(1-q)\big)\\[5pt]
\dis\qquad=\int\nu(\di q)(1-q)^{m-1},
\ec
which by (\ref{tetprop})~(ii) implies (\ref{tetnu}).}

\noi
{\bf Proof of Lemma~\ref{L:basis}} As a first step towards proving part~(a),
we start by proving that the functions
$\big\{f_\De:\emptyset\neq\De\sub\{1,\ldots,n\}\big\}$ are linearly
independent. Consider the set $\{0,1\}^n\sub\R^n$. For each
$A\sub\{1,\ldots,n\}$, define $g_A:\{0,1\}^n\to\R$ by
\be
g_A(\vec x):=\left\{\ba{ll}
1-f_A(\vec x)\quad&\mbox{if }A\neq\emptyset,\\
1\quad&\mbox{if }A=\emptyset.
\ea\right.
\ee
It is not hard to see that $g_Ag_B=g_{A\cup B}$ and that the functions
$\big\{g_A:A\sub\{1,\ldots,n\}\big\}$ separate points. Therefore, by the
Stone-Weierstrass theorem, they span the space of all real functions on
$\{0,1\}^n$. Since this space has dimension $2^n$ and since
$\big\{g_A:A\sub\{1,\ldots,n\}\big\}$ has $2^n$ elements, we conclude that the
$g_A$'s are linearly independent and hence the same is true for the $f_\De$'s.

We next prove that the $f_\De$'s span $L_n$. Obviously $f_\De\in L_n$ for each
$\emptyset\neq\De\sub\{1,\ldots,n\}$. Therefore, since the $f_\De$'s are
linearly independent and since
$\big\{f_\De:\emptyset\neq\De\sub\{1,\ldots,n\}\big\}$ has $2^n-1$ elements,
it suffices to show that ${\rm dim}(L_n)\leq 2^n-1$. We proceed by
induction. It is easy to check that $L_1$ is the space of all linear functions
from $\R$ to $\R$, which has dimension one. Now assume that ${\rm
  dim}(L_n)\leq 2^n-1$. We claim that ${\rm dim}(L_{n+1})\leq 2^{n+1}-1$. Each
function $f\in L_{n+1}$ can be uniquely written as
\be
f(\vec x)=\sum_{i=1}^{n+1}c_i(\vec x)x_i,
\ee
where the functions $c_1,\ldots,c_{n+1}$ are piecewise constant on each cell
$C_\prec$. In fact, since functions in $L_n$ are continuous, we must have that
the function $c_{n+1}$ depends only on the relative order of $x_{n+1}$ with
respect to the first $n$ coordinates and does not change if we interchange the
order of two other coordinates $x_j,x_k$ with $j,k\leq n$. More precisly, for
each $A\sub\{1,\ldots,n\}$, if we set
\be
U_A:=\big\{\vec x\in\R^{n+1}:
x_i<x_{n+1}\ \forall i\in A,\ x_i>x_{n+1}
\ \forall i\in\{1,\ldots,n\}\beh A\big\},
\ee
then
\be
c_{n+1}(\vec x)=l_A\qquad(x\in U_A)
\ee
for some constant $l_A\in\R$. Let $l$ be the linear map defined by
\[
l(f):=\big(l_A(f))_{A\sub\{1,\ldots,n\}}\qquad(f\in L_n).
\]
Then ${\rm Ker}(l)$ consists of all functions in $L_{n+1}$ that do not depend
on the variable $x_{n+1}$, hence ${\rm Ker}(l)\sub L_n$. It follows that
\[
{\rm dim}(L_{n+1})={\rm dim}({\rm Ker}(l))+{\rm dim}({\rm Ran}(l))
\leq (2^n-1)+2^n=2^{n+1}-1,
\]
as claimed.

To prove part~(b) of the lemma, we need to calculate
\be\label{Aform1}
\Ai^\tet_nf_\De(\vec x)
=\sum_{x\in{\rm Ran}(\vec x)}\sum_{I\sub J_x}
\tet(|I|,|J_x\beh I|)\nab_{\vec v_{I,J_x\beh I}}f_\De(\vec x).
\ee
Let us define
\bc
\dis H(\vec x)&:=&\dis J_{f_\De(\vec x)}
=\big\{i\in\{1,\ldots,n\}:x_i=f_\De(\vec x)\big\},\\[5pt]
\dis G(\vec x)&:=&\dis H(\vec x)\cap\De
=\big\{i\in\De:x_i=f_\De(\vec x)\big\}.
\ec
Recalling (\ref{vIJ}) and (\ref{nabv}), we see that
\be\label{nabDe}
\nab_{\vec v_{I,J}}f_\De(\vec x)=\left\{\ba{rl}
+1\quad&\mbox{if }I\cap G(\vec x)\neq\emptyset,\\
-1\quad&\mbox{if }J\supset G(\vec x),\\
0\quad&\mbox{otherwise.}\ea\right.
\ee
Inserting this into (\ref{Aform1}) we see that
\be\label{Aform2}
\Ai^\tet_nf_\De(\vec x)
=\sum_{I\sub H(\vec x)}
\tet(|I|,|H(\vec x)\beh I|)\big(1_{\{I\cap G(\vec x)\neq\emptyset\}}
-1_{\{I\cap G(\vec x)=\emptyset\}}\big),
\ee
where we have used that for $I\sub H(\vec x)$, one has $(H(\vec x)\beh
I)\supset G(\vec x)$ if and only if $I\cap G(\vec x)=\emptyset$.

We claim that (\ref{Aform2}) can be rewritten as
\be\label{Aform3}
\Ai^\tet_nf_\De(\vec x)
=\sum_{I\sub G(\vec x)}
\tet(|I|,|G(\vec x)\beh I|)\big(1_{\{I\cap G(\vec x)\neq\emptyset\}}
-1_{\{I\cap G(\vec x)=\emptyset\}}\big).
\ee
To see this, note that if $H'$ is a set such that $G(\vec x)\sub H'\sub H(\vec
x)$ and $H'$ contains one element less than $H(\vec x)$, then since it does
not make a difference for the sign of a term in (\ref{Aform2}) whether we
include this element in $I$ or in $H(\vec x)\beh I$, we have
\bc\label{Aform2a}
\dis\Ai^\tet_nf_\De(\vec x)
&=&\dis\sum_{I\sub H'}
\big(\tet(|I|+1,|H'\beh I|)+\tet(|I|,|H'\beh I|+1)\big)
\big(1_{\{I\cap G(\vec x)\neq\emptyset\}}
-1_{\{I\cap G(\vec x)=\emptyset\}}\big)\\[5pt]
&=&\dis\sum_{I\sub H'}\tet(|I|,|H'\beh I|)
\big(1_{\{I\cap G(\vec x)\neq\emptyset\}}
-1_{\{I\cap G(\vec x)=\emptyset\}}\big),
\ec
where we have used (\ref{tetprop})~(ii). Continuing this process of removing
points from $H(\vec x)$ we arrive at (\ref{Aform3}).

We may rewrite (\ref{Aform3}) as
\bc\label{Aform3a}
\dis\Ai^\tet_nf_\De(\vec x)
&=&\dis\sum_{I\sub G(\vec x)}
\tet(|I|,|G(\vec x)\beh I|)
\big(1-2\cdot 1_{\{I\cap G(\vec x)=\emptyset\}}\big)\\[5pt]
&=&\dis\Big(\sum_{I\sub G(\vec x)}\tet(|I|,|G(\vec x)\beh I|)\Big)
-2\tet(0,|G(\vec x)|).
\ec
The same sort of argument as in (\ref{Aform2a}) shows that
the first term on the right-hand side of (\ref{Aform3a}) equals $\tet(0,0)$
and hence, recalling (\ref{tetbet}) and the fact that $|G(\vec x)|=g_\De(\vec
x)$ (see (\ref{fgDe})), we arrive at (\ref{AfDe}).

Part~(c) is a trivial consequence of parts~(a) and (b) and the fact that if
$\tet$ and $\tet'$ are equivalent in the sense of (\ref{teteq}), then they
define the same $(\bet_+(m))_{m\geq 1}$ through (\ref{tetbet}).\qed

\subsection{Convergence of discrete n-point motions}

In this section we prove that if $\mu_k$ are probability measures on $[0,1]$
satisfying (\ref{mucon}), then the diffusively rescaled discrete $n$-point
motions associated with the $\mu_k$ converge in law to the Markov process
defined by the Howitt-Warren martingale problem with drift $\bet$ and
characteristic measure $\nu$. To formulate this precisely, fix $\mu_k$
satisfying (\ref{mucon}), let $\vec X^{\langle k\rangle}$ be discrete
$n$-point motions associated with the $\mu_k$, started in deterministic
initial states $\vec x^{\langle k\rangle}$ and linearly interpolated between
integer times, and let $\vec Y^{\langle k\rangle}$ defined by
\be\label{vecY}
Y^{\langle k\rangle}_i(t):=\eps_kX^{\langle k\rangle}_i(t/\eps_k^2)
\qquad(i=1,\ldots,n,\ t\geq 0)
\ee
denote the process $\vec X^{\langle k\rangle}$, diffusively rescaled with
$\eps_k$. Let $\Ci_{\R^n}\half$ denote the space of continuous functions from
$\half$ to $\R^n$, equipped with the topology of local uniform
convergence. Then, in analogy with \cite[Thm.~8.1]{HW09a}, we have the following
result.
\bp{\bf(Convergence of the $n$-point motions)}\label{P:nconv}
Assume that the initial states satisfy
\be
\eps_k\vec x^{\langle k\rangle}\asto{k}\vec x
\ee
for some $\vec x\in\R^n$. Then
\be\label{nconv}
\P\big[\big(\vec Y^{\langle k\rangle}(t)\big)_{t\geq 0}\in\cdot\,\big]
\Asto{k}
\P\big[\big(\vec X(t)\big)_{t\geq 0}\in\cdot\,\big],
\ee
where $\Rightarrow$ denotes weak convergence of probability laws on
$\Ci_{\R^n}\half$ and $\vec X$ is the unique solution of the Howitt-Warren
martingale problem with drift $\bet$ and characteristic measure $\nu$, started
in the initial state $\vec X_0=\vec x$.
\ep
{\bf Remark.} There is an analoguous statement for random initial states, see
Remark~\ref{R:infvar}.\med

\noi
We will actually prove a somewhat stronger statement than the convergence in
(\ref{nconv}), since we will show that the intersection times of the
rescaled discrete process also converge to those of the limiting process. For
technical reasons, it will be convenient to interpolate in a piecewise
constant, rather than in a linear way. Therefore, we set (compare
(\ref{vecY}))
\be\label{vecY2}
\Yr^{\langle k\rangle}_i(t):=\eps_kX^{\langle k\rangle}_i(\lfloor t/\eps_k^2\rfloor)
\qquad(i=1,\ldots,n,\ t\geq 0).
\ee
We view $\vec\Yr^{\langle k\rangle}$ as a process with paths in
$\Di_{\R^n}\half$, the space of c\'adl\'ag functions from $\half$ to $\R^n$,
equipped with the Skorohod topology. Letting $\vec Y^{\langle
  k\rangle},\vec\Yr^{\langle k\rangle}$ denote the linearly interpolated and
piecewise constant processes, respectively, we have
\be
\sup_{t\geq 0}|\Yr^{\langle k\rangle}_i(t)-Y^{\langle k\rangle}_i(t)|=\eps_k\asto{k}0.
\ee
{F}rom this, it is easy to see that Proposition~\ref{P:nconv} is implied by the
following, somewhat stronger result.
\bp{\bf(Convergence including intersection times)}\label{P:nconv2}
Let $\vec X^{\langle k\rangle}$ be discrete $n$-point motions associated with
probability measures $\mu_k$ satisfying (\ref{mucon}), started from initial
states $\vec x^{\langle k\rangle}$, and let $\vec\Yr^{\langle k\rangle}$ denote
$\vec X^{\langle k\rangle}$ diffusively rescaled as in (\ref{vecY2}).
Let $\vec X$ be the unique solution of the Howitt-Warren
martingale problem with drift $\bet$ and characteristic measure $\nu$, started
in $\vec X_0=\vec x$. Define $n\times n$ matrix valued processes
$Z^{\langle k\rangle}$ and $Z$ by
\be\ba{rr@{\,}c@{\,}l}\label{ZZ}
{\rm(i)}&\dis Z^{\langle k\rangle}_{ij}(t)
&:=&\dis\int_0^t1_{\{\Yr^{\langle k\rangle}_i(s)
=\Yr^{\langle k\rangle}_j(s)\}}\di s,\\[5pt]
{\rm(ii)}&\dis Z_{ij}(t)
&:=&\dis\int_0^t1_{\{X_i(s)=X_j(s)\}}\di s.
\ec
Then, assuming that the initial states satisfy
\be
\eps_k\vec x^{\langle k\rangle}\asto{k}\vec x,
\ee
one has
\be\label{nconv2}
\P\big[\big(\vec\Yr^{\langle k\rangle}(t),
Z^{\langle k\rangle}(t)\big)_{t\geq 0}\in\cdot\,\big]
\Asto{k}
\P\big[\big(\vec X(t),Z(t)\big)_{t\geq 0}\in\cdot\,\big],
\ee
where $\Rightarrow$ denotes weak convergence of probability laws on path space.
\ep
{\bf Proof.} When $\vec X^{\langle k\rangle}$ is the $n$-point motion of a continuous
time version of the discrete Howitt-Warren flow, the same result has been proved by Howitt and Warren in
\cite[Prop.~6.3]{HW09a} (for tightness in their case, see the remarks above their
formula (6.13).) Our proof copies their proof in many places, except that
we use a different argument to get convergence of the compensators of
$f_\De(\vec\Yr^{\langle k\rangle})$ and we have also simplified their proof
somewhat due to our reformulation of their martingale problem.

Let $P^{\langle k\rangle}$ be the transition kernel from $\Z^n$ to
$\Z^n$ defined by
\be
P^{\langle k\rangle}(\vec x,\vec y):=\prod_{x\in{\rm Ran}(\vec x)}
\int\mu_k(\di q)\prod_{i\in J_x}
\big(1_{\{y_i=x_i+1\}}q+1_{\{y_i=x_i-1\}}(1-q)\big)\qquad({\vec x},{\vec y}\in\Z^n),
\ee
where ${\rm Ran}(\vec x)$ and $J_x$ are defined in (\ref{Ranx}) and
(\ref{Jx}). We adopt the notation
\be
P^{\langle k\rangle}f(x):=\sum_{y\in\Z^n}P^{\langle k\rangle}(x,y)f(y)
\qquad(x\in\Z^n,\ f:\Z^n\to\R),
\ee
whenever the infinite sum is well-defined.

We observe that $\vec X^{\langle k\rangle}$ is a Markov chain with transition
kernel $P^{\langle k\rangle}$. Since we start $\vec X^{\langle k\rangle}$ in
an initial state $\vec X^{\langle k\rangle}(0)=x^{\langle k\rangle}\in(\Z_{\rm
  even})^n$, because of the nature of the transition mechanism, we have $\vec
X^{\langle k\rangle}(t)\in(\Z_{\rm even})^n$ at even times and $\vec
X^{\langle k\rangle}(t)\in(\Z_{\rm odd})^n$ at odd times.

For $\emptyset\neq\De\sub\{1,\ldots,n\}$, let $f_\De,g_\De$ be the functions
defined in (\ref{fgDe}). By standard theory, for each
$\emptyset\neq\De\sub\{1,\ldots,n\}$, the discrete-time
process
\be
f_\De(\vec X^{\langle k\rangle}(t))
-\sum_{s=0}^{t-1}\big(P^{\langle k\rangle}f_\De(\vec X^{\langle k\rangle}(s))
-f_\De(\vec X^{\langle k\rangle}(s))\big)
\ee
is a martingale with respect to the filtration generated by $\vec X^{\langle
  k\rangle}$. We observe that if either $x\in(\Z_{\rm even})^n$ or $x\in(\Z_{\rm
  odd})^n$, then under the transition kernel $P^{\langle k\rangle}$ the
maximum $f_\De(x)=\max_{i\in\De}x_i$ moves down by one with probability
$\mu_k(\di q)(1-q)^{g_\De(x)}$ and up by one with the remaining probability,
hence
\bc
\dis P^{\langle k\rangle}f_\De(x)-f_\De(x)
&=&\dis\int\mu_k(\di q)\big(1-2(1-q)^{g_\De(x)}\big)\\[5pt]
&=&\dis\bet_k(g_\De(x))
\qquad\big(x\in(\Z_{\rm even})^n\ \mbox{or}\ (\Z_{\rm even})^n\big),
\ec
where we have introduced the notation
\be
\bet_k(m):=\int\mu_k(\di q)\big(1-2(1-q)^m\big)\qquad(m\geq 1).
\ee
Setting $\bet_k:=\bet_k(1)=\int\mu_k(\di q)(2q-1)$, by standard theory, one
may moreover check that
\be
\big(X^{\langle k\rangle}_i(t)-\bet_kt\big)
\big(X^{\langle k\rangle}_j(t)-\bet_kt\big)
-\sum_{s=0}^{t-1}\Ga^{\langle k\rangle}_{ij}(\vec X^{\langle k\rangle}(s))
\ee
is a martingale, where
\bc
\dis\Ga^{\langle k\rangle}_{ij}(\vec x)
&:=&\dis\sum_{{\vec y}\in\Z^n}P^{\langle k\rangle}({\vec x},{\vec y})(y_i-x_i)(y_j-x_j) - \beta_k^2\\[15pt]
&=&\dis\left\{\ba{cl}
1-\beta_k^2\quad&\mbox{if }i=j,\\
\int\mu_k(\di q)\big(1-4q(1-q)\big) -\beta_k^2
\quad&\mbox{if }i\neq j,\ x_i=x_j,\\
0\quad&\mbox{otherwise}
\ea\right.\\[20pt]
&=&\dis\big(1-\beta_k^2-2\big(\bet_k(2)-\bet_k(1)\big)1_{\{i\neq j\}}\big)1_{\{x_i=x_j\}}.
\ec
For the process $\vec\Yr^{\langle k\rangle}$ defined in (\ref{vecY2}), our
arguments so far show that for each $\emptyset\neq\De\sub\{1,\ldots,n\}$,
\be\label{kMP}
f_\De(\vec\Yr^{\langle k\rangle}(t))-\eps_k^{-1}
\int_0^{\lfloor t\rfloor_k}\bet_k(g_\De(\vec\Yr^{\langle k\rangle}(s)))\di s
\ee
is a martingale, where $\lfloor t\rfloor_k:=\eps_k^2\lfloor t/\eps_k^2\rfloor$
denotes the time $t$ rounded downwards to the next time in $\eps_k^2\N$.
Moreover, for each $1\leq i,j\leq n$, the process
\be\ba{l}\label{qvark}
\dis\big(\Yr^{\langle k\rangle}_i(t)-\eps_k^{-1}\bet_kt\big)
\big(\Yr^{\langle k\rangle}_j(t)-\eps_k^{-1}\bet_kt\big)\\[5pt]
\dis\qquad-\big(1-\beta_k^2-2\big(\bet_k(2)-\bet_k(1)\big)1_{\{i\neq j\}}\big)
\int_0^{\lfloor t\rfloor_k}
1_{\{\Yr^{\langle k\rangle}_i(s)=\Yr^{\langle k\rangle}_i(s)\}}\di s
\ec
is a martingale with respect to the filtration generated by $\vec\Yr^{\langle
  k\rangle}$. It follows from our assumption (\ref{mucon}) (see also
(\ref{mucon2a})) that
\be\label{conscon}
\lim_{k\to\infty}\eps_k^{-1}\bet_k(m)=\bet_+(m)\qquad(m\geq 1).
\ee
Standard results (Donsker's invariance principle) tell us that for
$k\to\infty$, each component $\Yr^{\langle k\rangle}_i$ of the rescaled
process converges weakly in law, on the space $\Di_\R\half$, to a Brownian
motion with drift $\bet=\bet_+(1)$. This implies that the laws of the processes
$\vec\Yr^{\langle k\rangle}$ (viewed as probability laws on $\Di_{\R^n}\half$)
are tight. Let $Z^{\langle k\rangle}$ be the matrix valued processes defined
in (\ref{ZZ})~(i). Since the slope of each $Z^{\langle k\rangle}_{ij}$ is
between zero and one, tightness for these processes is immediate.

By going to a subsequence if necessary, we may assume that the joint processes
$(\vec\Yr^{\langle k\rangle},Z^{\langle k\rangle})$ converges weakly in law,
and by Skorohod's representation theorem (see
e.g.\ \cite[Theorem~6.7]{Bil99}), we can couple the $(\vec\Yr^{\langle
  k\rangle},Z^{\langle k\rangle})$'s such that the convergence is almost
sure. Let $(\vec X,\ti Z)$ denote the limiting process. Then, taking the limit
in (\ref{qvark}), using (\ref{conscon}), we see that
\be
\big(X_i(t)-\bet t\big)\big(X_j(t)-\bet t\big)-\ti Z_{ij}(t)
\ee
is a martingale, hence
\be\label{qvar}
\li X_j,X_j\re(t)=\ti Z_{ij}(t)=\lim_{k\to\infty}
\int_0^t1_{\{\Yr^{\langle k\rangle}_i(s)=\Yr^{\langle k\rangle}_j(s)\}}\di s
\qquad\forall t\geq 0,\ 1\leq i,j\leq n\quad{\rm a.s.}
\ee
Since, for given $t>0$, the function $w\mapsto\int_0^t1_{\{w_i(s)=w_j(s)\}}\di
s$ is upper semicontinuous with respect to the topology on $\Di_{\R^n}\half$,
formula (\ref{qvar}) implies that
\be\label{Zleq}
\li X_i,X_j\re(t)\leq\int_0^t1_{\{X_i(s)=X_j(s)\}}\di s
\qquad(t\geq 0,\ 1\leq i,j\leq n).
\ee
To prove also the other inequality in (\ref{Zleq}), we use an argument due to
Howitt and Warren (see the proof of formula (6.9) in \cite{HW09a}). For any real
square integrable semimartingale $W$, one can define a `local time' $L(x,t)$
such that
\be\label{loctime}
\int_0^t f(W(s))\di\li W,W\re(s)=\int_\R f(x)L(x,t)\di x.
\ee
(See \cite[formula~(3)]{BY81}.) Applying this to the semimartingale
$X_i-X_j$ and the function $f=1_{\{0\}}$, we find that
\be\label{octi}
\int_0^t1_{\{X_i(s)=X_j(s)\}}\di\li X_i-X_j,X_i-X_j\re(s)= \int_\R 1_{\{0\}}(x) L(x,t)\di x = 0.
\ee
Since $X_i,X_j$ are Brownian motions, we have
\bc
\dis\li X_i-X_j,X_i-X_j\re(t)
&=&\dis\li X_i,X_i\re(t)+\li X_j,X_j\re(t)-2\li X_i,X_j\re(t)\\[5pt]
&=&\dis 2t-2\li X_i,X_j\re(t).
\ec
Inserting this into (\ref{octi}) yields
\be\label{difplace}
\int_0^t1_{\{X_i(s)=X_j(s)\}}\di s
=\int_0^t1_{\{X_i(s)=X_j(s)\}}\di\li X_i,X_j\re(s).
\ee
On the other hand, (\ref{Zleq}) implies that
\be\label{Xonly}
\int_0^t1_{\{X_i(s)\neq X_j(s)\}}\di\li X_i,X_j\re(s)=0.
\ee
Combining this with (\ref{difplace}) yields
\be\label{toge}
\int_0^t1_{\{X_i(s)=X_j(s)\}}\di s=\li X_i,X_j\re(t)
\qquad(t\geq 0,\ 1\leq i,j\leq n),
\ee
as claimed.

We now show that ${\vec X}$ solves the Howitt-Warren martingale problem.
By (\ref{qvar}) and (\ref{toge}), we conclude that
\be\label{intime}
\int_0^t1_{\{\Yr^{\langle k\rangle}_i(s)=\Yr^{\langle k\rangle}_j(s)\}}\di s
\asto{k}\int_0^t1_{\{X_i(s)=X_j(s)\}}\di s
\qquad\forall t\geq 0,\ 1\leq i,j\leq n\quad{\rm a.s.}
\ee
The lower semicontinuity of the map
$w\mapsto\int_0^t1_{\{w_i(s)<w_j(s)\}}\di s$
implies that
\be\label{leslim}
\liminf_{k\to\infty}
\int_0^t1_{\{\Yr^{\langle k\rangle}_i(s)<\Yr^{\langle k\rangle}_j(s)\}}\di s
\geq\int_0^t1_{\{X_i(s)<X_j(s)\}}\di s
\qquad\forall t\geq 0,\ 1\leq i,j\leq n\quad{\rm a.s.}
\ee
Combining this with (\ref{intime}) we see that a.s., for all $t\geq 0$ and
$1\leq i,j\leq n$,
\be\ba{l}
\dis\limsup_{k\to\infty}
\int_0^t1_{\{\Yr^{\langle k\rangle}_i(s)>\Yr^{\langle k\rangle}_j(s)\}}\di s\\[5pt]
\dis\quad\leq 1-\lim_{k\to\infty}
\int_0^t1_{\{\Yr^{\langle k\rangle}_i(s)=\Yr^{\langle k\rangle}_j(s)\}}\di s
-\liminf_{k\to\infty}
\int_0^t1_{\{\Yr^{\langle k\rangle}_i(s)<\Yr^{\langle k\rangle}_j(s)\}}\di s\\[5pt]
\dis\quad \leq 1 - \int_0^t1_{\{X_i(s)=X_j(s)\}}\di s - \int_0^t1_{\{X_i(s)<X_j(s)\}}\di s
=\int_0^t1_{\{X_i(s)>X_j(s)\}}\di s,
\ec
which together with (\ref{leslim}) shows that
\be\label{lestime}
\int_0^t1_{\{\Yr^{\langle k\rangle}_i(s)<\Yr^{\langle k\rangle}_j(s)\}}\di s
\asto{k}\int_0^t1_{\{X_i(s)<X_j(s)\}}\di s
\qquad\forall t\geq 0,\ 1\leq i,j\leq n\quad{\rm a.s.}
\ee
By Lemma~\ref{L:intcon} below, this implies that
\be
\int_0^t\big|1_{\{\Yr^{\langle k\rangle}_i(s)<\Yr^{\langle k\rangle}_j(s)\}}
-1_{\{X_i(s)<X_j(s)\}}\big|\di s\asto{k}0
\qquad\forall t\geq 0,\ 1\leq i,j\leq n\quad{\rm a.s.},
\ee
which in turn implies that a.s., for each $t\geq 0$ and weak total order
$\prec$ on $\{1,\ldots,n\}$, one has
\be\label{precon}
\int_0^t\big|1_{\{\Yr^{\langle k\rangle}(s)\in C_\prec\}}
-1_{\{X(s)\in C_\prec\}}\big|\di s\asto{k}0,
\ee
where $C_\prec$ is the cell defined in (\ref{precel}). Since $g_\De(x)$
depends only on the relative order of the coordinates $x_1,\ldots,x_n$,
formulas (\ref{precon}) and (\ref{conscon}) imply that for each
$\emptyset\neq\De\sub\{1,\ldots,n\}$,
\be\label{compcon}
\eps_k^{-1}\int_0^{\lfloor t\rfloor_k}
\bet_k(g_\De(\vec\Yr^{\langle k\rangle}(s)))\di s
\asto{k}\int_0^t\bet_+(g_\De(\vec X(s)))\di s
\qquad\forall t\geq 0\quad{\rm a.s.}
\ee
Taking the limit $k\to\infty$ in (\ref{kMP}) using (\ref{compcon}) and the
fact that $g_\De$ is a bounded function (indeed, $1\leq g_\De(x)\leq|\De|$),
we find that for each $\emptyset\neq\De\sub\{1,\ldots,n\}$,
\be
f_\De\big(\vec X(t)\big)-\int_0^t \bet_+\big(g_\De(\vec X(s))\big)\di s,
\ee
is a martingale with respect to the filtration generated by $\vec X$.
Together with (\ref{toge}) this shows that $\vec X$ solves the
Howitt-Warren martingale problem, completing our proof.\qed

\bl{\bf(Convergence of integrals)}\label{L:intcon}
Let $T>0$ and let $\rho$ be a finite measure on $[0,T]$.\med

\noi
{\bf(a)} Let $f_k,f$ be Borel measurable real functions on $[0,T]$ such that
$\sup_k\|f_k\|<\infty$, where $\|\,\cdot\,\|$ denotes the supremum norm.
Assume that
\be
\int_{[0,t]}\rho(\di s)f_k(s)\asto{k}\int_{[0,t]}\rho(\di s)f(s)
\qquad(0\leq t\leq T).
\ee
Then
\be\label{Acon}
\int_A\rho(\di s)f_k(s)\asto{k}\int_A\rho(\di s)f(s)
\ee
for each Borel measurable $A\sub[0,T]$.\med

\noi
{\bf(b)} Let $A_k,A\sub[0,T]$ be Borel measurable. Assume that
\be
\int_{[0,t]}\rho(\di s)1_{A_k}(s)\asto{k}\int_{[0,t]}\rho(\di s)1_A(s)
\qquad(0\leq t\leq T).
\ee
Then
\be
\int_{[0,T]}\rho(\di s)\big|1_{A_k}(s)-1_A(s)\big|\asto{k}0.
\ee
\el
{\bf Proof.} To prove part~(a), let $\Gi$ be the set of Borel measurable
subsets $A\sub[0,T]$ for which (\ref{Acon}) holds. It is clear that
$A,B\in\Gi$, $A\supset B$ implies that $A\beh B\in\Gi$. We claim that
moreover, if $A_n\in\Gi$ satisfy $A_n\up A$ for some $A\sub[0,T]$, then
$A\in\Gi$. To see this, write
\be\ba{l}
\dis\Big|\int_A\rho(\di s)f_k(s)-\int_A\rho(\di s)f(s)\Big|\\[5pt]
\dis\quad\leq\Big|\int_A\rho(\di s)f_k(s)-\int_{A_n}\rho(\di s)f_k(s)\Big|
+\Big|\int_{A_n}\rho(\di s)f(s)-\int_A\rho(\di s)f(s)\Big|\\[5pt]
\dis\quad\phantom{\leq}
+\Big|\int_{A_n}\rho(\di s)f_k(s)-\int_{A_n}\rho(\di s)f(s)\Big|\\[5pt]
\dis\quad\leq2\rho(A\beh A_n)\sup_m\|f_m\|
+\Big|\int_{A_n}\rho(\di s)f_k(s)-\int_{A_n}\rho(\di s)f(s)\Big|.
\ec
By choosing $n$ large enough, we see that
\be
\limsup_{k\to\infty}\Big|\int_A\rho(\di s)f_k(s)-\int_A\rho(\di s)f(s)\Big|
\leq\eps
\ee
for all $\eps>0$, proving our claim. Since the set $\Hi:=\{[0,t]:0\leq t\leq
T\}$ is closed under intersections and contained in $\Gi$, Sierpi\'nski's
$\pi/\la$-theorem \cite[Theorem~1.1]{Kal02} tells us that $\Gi$ contains the
\si-field generated by $\Hi$, completing our proof.

To prove part~(b), we note that
\bc
\dis\int_{[0,T]}\rho(\di s)\big|1_{A_k}(s)-1_A(s)\big|
&=&\dis\Big(\int_A\rho(\di s)1_A(s)-\int_A\rho(\di s)1_{A_k}(s)\Big)\\[5pt]
&&\dis+\Big(\int_{[0,T]\beh A}\rho(\di s)1_{A_k}(s)
-\int_{[0,T]\beh A}\rho(\di s)1_A(s)\Big),
\ec
which tends to zero by part~(a).\qed

\section{The Hausdorff topology}\label{A:Haus}

Let $(E,d)$ be a metric space, let $\Ki(E)$ be the space of all compact
subsets of $E$ and set $\Ki_+(E):=\{K\in\Ki(E):K\neq\emptyset\}$. Then the
{\em Hausdorff metric} $d_{\rm H}$ on $\Ki_+(E)$ is defined as
\bc\label{Hmet2}
\dis d_{\rm H}(K_1,K_2)&:=&\dis\sup_{x_1\in K_1}\inf_{x_2\in K_2}d(x_1,x_2)
\vee\sup_{x_2\in K_2}\inf_{x_1\in K_1}d(x_1,x_2)\\[12pt]
&=&\dis\sup_{x_1\in K_1}d(x_1,K_2)
\vee\sup_{x_2\in K_2}d(x_2,K_1),
\ec
where $d(x,A):=\inf_{y\in A}d(x,y)$ denotes the distance between a point $x\in E$ and a set $A\sub E$. The corresponding topology is called the {\em Hausdorff topology}. We extend this topology to $\Ki(E)$ by adding $\emptyset$ as an isolated point. The next lemma shows that the Hausdorff topology depends only on the topology on $E$, and not on the choice of the metric.

\bl{\bf(Convergence criterion)}\label{L:Hauconv}
Let $K_n,K\in\Ki_+(E)$ $(n\geq 1)$. Then $K_n\to K$ in the Hausdorff topology if and only if there exists a $C\in\Ki_+(E)$ such that $K_n\sub C$ for all $n\geq 1$ and
\bc\label{Haulim}
K&=&\dis\{x\in E:\exists x_n\in K_n\mbox{ s.t.\ }x_n\to x\}\\[5pt]
&=&\dis\{x\in E:\exists x_n\in K_n
\mbox{ s.t.\ $x$ is a cluster point of } (x_n)_{n\in\N}\}.
\ec
\el

The following lemma shows that $\Ki(E)$ is Polish if $E$ is.

\bl{\bf(Properties of the Hausdorff metric)}\label{L:Hauprop}\ \\[5pt]
{\bf(a)} If $(E,d)$ is separable, then so is $(\Ki_+(E),d_{\rm H})$.\med

\noi
{\bf(b)} If $(E,d)$ is complete, then so is $(\Ki_+(E),d_{\rm H})$.
\el

Recall that a subset $A$ of a metric space is {\em precompact} if its closure is compact. This is equivalent to the statement that each sequence of points $x_n\in A$ has a convergent subsequence.

\bl{\bf(Compactness in the Hausdorff topology)}\label{L:Haucomp}
A set $\Ai\sub\Ki(E)$ is precompact if and only if there exists a $C\in\Ki(E)$ such that $K\sub C$ for each $K\in\Ai$.
\el

The following lemma is useful when proving convergence of $\Ki(E)$-valued random variables.

\bl{\bf(Tightness criterion)}\label{L:Hautight}
Assume that $E$ is a Polish space and let $K_n$ $(n\geq 1)$ be $\Ki(E)$-valued random variables. Then the collection of laws $\{\P[K_n\in\cdot\,]:n\geq 1\}$ is tight if and only if for each $\eps>0$ there exists a compact $C\sub E$ such that $\P[K_n\sub C]\geq 1-\eps$ uniformly in $n\in\N$.
\el

If $E$ is compact, then the Hausdorff topology on $\Ki(E)$ coincides with the Fell topology defined in \cite[Thm.~A.2.5]{Kal02}. The Hausdorff metric may more generally be defined on the space of nonempty bounded closed subsets of $(E,d)$. In particular, if $d$ is bounded, then $d_{\rm H}(A_1,A_2)$ can be defined for any nonempty closed $A_1,A_2$. In this more general set-up, Lemma~\ref{L:Hauprop}~(b) and the `if' part of Lemma~\ref{L:Haucomp} remain true, as well as the `if' part of Lemma~\ref{L:Hautot} below. This is Excercise~7 (with some hints for a possible solution) in \cite[\S~45]{Mun00}. A detailed solution of this excercise can be found in \cite{Hen99}. We are not aware of any reference for the other statements in Lemmas~\ref{L:Hauconv}--\ref{L:Hautight}, although they appear to be well-known. For completeness, we provide self-contained proofs of all these lemmas. We start with some preparations.

Recall that for any metric space $(E,d)$, a set $A\sub E$ is {\em totally
  bounded} if for every $\eps>0$ there exists a finite collection of points
$x_1,\ldots,x_n\in E$ such that $A\sub\bigcup_{i=1}^nB_\eps(x_i)$, where
$B_\eps(x)$ denotes the open ball of radius $\eps$ around $x$. This is
equivalent to the statement that every sequence $x_n\in A$ has a Cauchy
subsequence. As a consequence, a set $A\sub E$ is compact if and only if it is
complete and totally bounded.

\bl{\bf(Totally bounded sets in the Hausdorff metric)}\label{L:Hautot}
A set $\Ai\sub\Ki_+(E)$ is totally bounded in $(\Ki_+(E),d_{\rm H})$ if and only if the set $A:=\{x\in E:\exists K\in\Ai\mbox{ s.t.\ }x\in K\}$ is totally bounded in $(E,d)$.
\el
{\bf Proof.} Assume that $A$ is totally bounded. Let $\eps>0$ and let $\De\sub E$ be a
finite set such that $A=\bigcup_{x\in\De}B_\eps(x)$. Let $K\in\Ki_+(E)$ and set
$\De':=\{x\in\De:B_\eps(x)\cap K\neq\emptyset\}$. Then for all $y\in K$ there
is an $x\in\De'$ such that $d(x,y)<\eps$ and for all $x\in\De'$ there is a
$y\in K$ such that $d(x,y)<\eps$ proving that $d_{\rm H}(\De',K)<\eps$. This
shows that $\Ai$ is covered, in the Hausdorff metric, by the collection of
open balls of radius $\eps$ centered around finite subsets of $\De$. Since
$\eps$ is general, we conclude that $\Ai$ is totally bounded.

Conversely, if $\Ai$ is totally bounded, then for each $\eps>0$ we can find $K_1,\ldots,K_n\in\Ki_+(E)$ such that $\Ai\sub\bigcup_{k=1}^n\Bi_{\eps/2}(K_n)$, where $\Bi_\eps(K)$ denotes the open ball in the Hausdorff metric of radius $\eps$ centered around a compact set $K$. Since each $K_k$ is compact, there exist $x_{k,1},\ldots,x_{k,m_k}$ such that $K_k\sub\bigcup_{j=1}^{m_k}B_{\eps/2}(x_{k,j})$, hence $A\sub\bigcup_{k=1}^n\bigcup_{j=1}^{m_k}B_\eps(x_{k,j})$.\qed

\bl{\bf(Cauchy sequences in the Hausdorff metric)}\label{L:HauCau}
Let $K_n\in\Ki_+(E)$ be a Cauchy sequence in $(\Ki_+(E),d_{\rm H})$. Then there exists a closed set $K$ such that (\ref{Haulim}) holds.
\el
{\bf Proof.} If the sets on the first and second line of the right-hand side of (\ref{Haulim}) are not equal, then there exists some $x\in E$ such that $x$ is a cluster point of some $x_n\in K_n$ but there do not exist $x'_n\in K_n$ such that $x'_n\to x$. It follows that there is some $\eps>0$ such that for each $k\geq 1$ we can find $n,m\geq k$ such that $K_n\cap B_\eps(x)\neq\emptyset$ and $K_m\cap B_{2\eps}(x)=\emptyset$, hence $d_{\rm H}(K_n,K_m)\geq\eps$, contradicting the asumption that the $K_n$ form a Cauchy sequence.

To see that $K$ is closed, assume that $x_n\in K$ satisfy $x_n\to x$ for some $x\in E$. Since $d_{\rm H}(K_n,K)\to 0$ we can choose $x'_n\in K_n$ such that $d(x'_n,x_n)\to 0$. It follows that $d(x'_n,x)\leq d(x'_n,x_n)+d(x_n,x)\to 0$ and hence $x\in K$.\qed

\bl{\bf(Sufficient conditions for convergence)}\label{L:Hausuff}
The conditions for convergence in the Hausdorff topology given in Lemma~\ref{L:Hauconv} are sufficient.
\el
{\bf Proof.} Our assumptions imply that $d(x,K_n)\to 0$ for each $x\in K$. We wish to show that in fact $\sup_{x\in K}d(x,K_n)\to 0$. If this is not the case, then by going to a subsequence if necessary we may assume that there exist $x_n\in K$ and $\eps>0$ such that $\liminf_{n\to\infty}d(x_n,K_n)\geq\eps$. Since $K$ is compact, by going to a further subsequence if necessary, we may assume that $x_n\to x\in K$. But then $\liminf_{n\to\infty}d(x,K_n)\geq\liminf_{n\to\infty}(d(x_n,K_n)-d(x,x_n))\geq\eps$ for this subsequence, contradicting the fact that for the original sequence, $d(x,K_n)\to 0$ for each $x\in K$.

The proof that $\sup_{x\in K_n}d(x,K)\to 0$ is similar. If this is not true, then we can go to a subsequence of the $K_n$ and then find $x_n\in K_n$ such that $d(x_n,K)\geq\eps$ for all $n$, for some $\eps>0$. Using the compactness of $C$, we can select a further subsequence such that $x_n\to x\in C$. Now $x$ is a cluster point of some $x_n\in K_n$ but $d(x,K)\geq\eps$, contradicting the fact that the two sets on the right-hand side of (\ref{Haulim}) are equal.\qed

\noi
{\bf Proof of Lemma~\ref{L:Hauprop}.}  To prove part~(a), it suffices to show
that if $\Di$ is a countable dense subset of $(E,d)$, then the collection of
finite subsets of $\Di$ is a countable dense subset of $(\Ki_+(E),d_{\rm
  H})$. Since a compact set $K\sub E$ is totally bounded, for each $\eps>0$,
we can find a finitely many points $x_1,\ldots,x_n\in E$ such that
$K\sub\bigcup_{i=1}^nB_{\eps/2}(x_i)$. Since $\Di$ is dense, we can choose
$x'_i\in\Di\cap B_{\eps/2}(x_i)$. Then $d_{\rm
  H}(K,\{x'_1,\ldots,x'_n\})\leq\eps$, proving our claim.

To prove part~(b), let $K_n\in\Ki_+(E)$ be a Cauchy sequence. Then, by
Lemma~\ref{L:HauCau}, there exists a closed set $K$ such that (\ref{Haulim})
holds. Since each sequence in the set $\{K_n:n\geq 1\}$ contains a Cauchy
subsequence, the set $\{K_n:n\geq 1\}$ is totally bounded, hence by
Lemma~\ref{L:Hautot}, there exists some totally bounded set containing all of
the $K_n$. Let $C$ denote its closure. Then $C$ is compact since $E$ is
complete, hence also $K\sub C$ is compact and Lemma~\ref{L:Hausuff} implies
that $K_n\to K$.\qed

\noi
{\bf Proof of Lemma~\ref{L:Haucomp}.} It suffices to prove the statement for $\Ai\sub\Ki_+(E)$. Let $\ov\Ai$ be the closure of $\Ai$ and set $C:=\{x\in E:\exists K\in\ov\Ai\mbox{ s.t.\ }x\in K\}$. By Lemma~\ref{L:Hautot}, $\ov\Ai$ is totally bounded if and only if $A$ is. Moreover, by Lemma~\ref{L:Hauprop}~(b), if $A$ is complete then so is $\{K\in\Ki_+(E):K\sub C\}$ and hence the same is true for $\ov\Ai$, being a closed subset of the former. Therefore, since compactness is equivalent to total boundedness and completeness, it suffices to show that compactness of $\ov\Ai$ implies completeness of $C$. Assume that $\ov\Ai$ is compact and that $x_n\in C$ is a Cauchy sequence. We need to show that the sequence $x_n$ has a cluster point $x\in C$. Choose $K_n\in\ov\Ai$ such that $x_n\in K_n$. Since $\ov\Ai$ is compact, by going to a subsequence if necessary, we may assume that $K_n\to K$ for some $K\in\ov\Ai$. Choose $x'_n\in K$ such that $d(x_n,x'_n)\to 0$. Since $K$ is compact, by going to a further subsequence if necessary, we may assume that $x'_n\to x$ for some $x\in K$. Since $d(x_n,x)\leq d(x_n,x'_n)+d(x'_n,x)\to 0$ this proves that the sequence $x_n$ has a cluster point $x\in K\sub C$.\qed

\noi
{\bf Proof of Lemma~\ref{L:Hautight}.} Immediate from Lemma~\ref{L:Haucomp} and the definition of tightness.\qed

\noi
{\bf Proof of Lemma~\ref{L:Hauconv}.} By Lemma~\ref{L:Hausuff}, we only need to prove that if $K_n\in\Ki_+(E)$ converge to a limit $K$, then there exists a $C\in\Ki_+(E)$ such that $K_n\sub C$ for all $n$ and (\ref{Haulim}) holds. If $K_n\to K$ then the set $\{K_n:n\geq 1\}$ is precompact, hence by Lemma~\ref{L:Hautight} there exists a $C\in\Ki_+(E)$ such that $K_n\sub C$ for all $n$. Formula (\ref{Haulim}) follows from the facts that if $x\in K$, then $d(x,K_n)\to 0$ hence there exist $K_n\ni x_n\to x$, while if $x\not\in K$, then $B_\eps(x)\cap K_n=\emptyset$ for all $n$ large enough such that $\sup_{x'\in K}d(x',K_n)<\eps$, hence $x$ is not a cluster point of some $x_n\in K_n$.\qed

\section{Some measurability issues}

Let $E,F$ be Polish spaces. By definition, the {\em pointwise closure} of a
set $\Fi$ of functions $f:E\to F$ is the smallest set containing $\Fi$ that is
closed under taking of pointwise limits, i.e., it is the intersection of all
sets $\Gi$ of functions from $E$ to $F$, such that $\Gi\supset\Fi$ and
$f_n\in\Gi$, $\lim_{n\to\infty}f_n(x)=f(x)$ $(x\in E)$ imply $f\in\Gi$.

\bl{\bf(Pointwise closure of functions to the unit interval)}\label{L:pclosunit}
Let $E$ be a Polish space and let $\Ci_{[0,1]}(E)$ be the set of all
continuous functions $f:E\to[0,1]$. Then the pointwise closure of
$\Ci_{[0,1]}(E)$ is the set $B_{[0,1]}(E)$ of all Borel measurable functions
$f:E\to[0,1]$. If $E$ is locally compact, then the same conclusion holds with
$\Ci_{[0,1]}(E)$ replaced by the space of continuous and compactly supported
functions $f:E\to[0,1]$.
\el
{\bf Proof.} By definition, one says that a sequence $f_n$ of real functions
on $E$ converges in a bounded pointwise way to a limit $f$ if $f_n(x)\to f(x)$
for each $x\in E$ and there exists some constant $C>0$ such that $|f_n|\leq C$
for all $n\geq 0$. The {\em bp-closure} of a set $\Fi$ of real functions on
$E$ is the smallest set containing $\Fi$ that is closed under taking of
bounded pointwise limits. By copying the proof of \cite[Lemma~3.4.1]{EK86}, we
see that the bp-closure of a convex set is convex. Let $\Bi$ be the set of all
subsets $A\sub E$ such that $1_A$ is in the bp-closure of
$\Ci_{[0,1]}(E)$. Then $\Bi$ is a Dynkin class containing all open sets, hence
by the Dynkin class theorem \cite[Thm.~A.4.2]{EK86} (resp.\ the
$\pi/\la$-theorem \cite[Theorem~1.1]{Kal02}), $\Bi$ contains all Borel
measurable subsets of $E$. Since indicator functions are the extremal elements
of the convex set consisting of all simple functions in $B_{[0,1]}(E)$, it is
easy to see that every simple function can be written as a convex combinations
of indicator functions. Since every function in $B_{[0,1]}(E)$ is an
increasing limit of simple functions in $B_{[0,1]}(E)$, the first claim
follows. In case $E$ is locally compact, it is easy to see that each
continuous function $f:E\to[0,1]$ is the pointwise limit of compactly
supported continuous functions $f:E\to[0,1]$, proving the second claim.\qed

\noi
We will need the following generalization of Lemma~\ref{L:pclosunit}. Below,
$[0,1]^\N$ denotes the space of all functions $x:\N\to[0,1]$, equipped with
the product topology. Note that the statement of Lemma~\ref{L:pclosure} is
false if we replace $[0,1]^\N$ by a general compact metrizable space
$F$. E.g., it is already wrong if $F$ consists of two isolated points, since
in this case all continuous functions are constant but there are lots of
measurable functions.

\bl{\bf(Pointwise closure)}\label{L:pclosure}
Let $E$ be a Polish space and let $\Ci_{[0,1]^\N}(E)$ be the set of all
continuous functions $f:E\to[0,1]^\N$. Then the pointwise closure of
$\Ci_{[0,1]^\N}(E)$ is the set $B_{[0,1]^\N}(E)$ of all Borel measurable
functions $f:E\to[0,1]^\N$.
\el
{\bf Proof.} Let $E,F,G$ be Polish spaces and let $\Fi,\Gi$ be sets of
functions $f:E\to F$ and $g:E\to G$, respectively. We claim that ${\rm
  pclos}(\Fi\times\Gi)\supset{\rm pclos}(\Fi)\times{\rm pclos}(\Gi)$, where
${\rm pclos}(\,\cdot\,)$ denotes the pointwise closure of a set and we regard
a pair of functions $(f,g)$ as a function from $E$ to $F\times G$ (equipped
with the product topology). To prove our claim, for any $f\in{\rm
  pclos}(\Fi)$, let $\Gi_f$ be the space of functions $g\in{\rm pclos}(\Gi)$
such that $(f,g)\in{\rm pclos}(\Fi\times\Gi)$. Then $\Gi_f$ is closed under
pointwise limits since ${\rm pclos}(\Fi\times\Gi)$ is. If $f\in\Fi$, then
moreover $\Gi_f$ contains $\Gi$ so $\Gi_f=\Gi$. Next, let $\hat\Fi$ be the
space of functions $f\in{\rm pclos}(\Fi)$ such that $(f,g)\in{\rm
  pclos}(\Fi\times\Gi)$ for all $g\in{\rm pclos}(\Gi)$. Then $\hat\Fi$ is
closed under pointwise limits since ${\rm pclos}(\Fi\times\Gi)$ is and
$\hat\Fi$ contains $\Fi$ by what we have just proved, so $\hat\Fi={\rm
  pclos}(\Fi)$, proving our claim.

Applying our clain inductively to $\Ci_{[0,1]^\N}(E)=(\Ci_{[0,1]}(E))^\N$,
using Lemma~\ref{L:pclosunit}, we see that $(f_1,\ldots,f_n,0,\ldots)$ lies in
the pointwise closure of $\Ci_{[0,1]^\N}(E)$ for each $f_1,\ldots,f_n\in
B_{[0,1]}(E)$ and $n\geq 1$. By taking pointwise limits, we see that each
infinite sequence $(f_1,f_2,\ldots)$ of Borel measurable functions
$f_i:E\to[0,1]$ lies in the pointwise closure of $\Ci_{[0,1]^\N}(E)$.\qed

\bl{\bf(Measurability of image measure map)}\label{L:imeas}
Let $E,F,G$ be Polish spaces and let $\Mi(E),\Mi(F)$ be the spaces of finite
measures on $E$ and $F$, respectively, equipped with the topology of weak
convergence and the associated Borel \si-field. Then, for any measurable map
$E\times G\ni(x,z)\mapsto f_z(x)\in F$, setting $\psi^f_z(\mu):=\mu\circ
f_z^{-1}$ defines a measurable map $\Mi_1(E)\times
G\ni(\mu,z)\mapsto\psi^f_z(\mu)\in\Mi_1(F)$.
\el
{\bf Proof.} We first prove the statement if $E,G$ are compact and
$F=[0,1]^\N$. In this case, we claim that if $E\times G\ni(x,z)\mapsto
f_z(x)\in F$ is continuous, then also $\Mi_1(E)\times
G\ni(\mu,z)\mapsto\psi^f_z(\mu)\in\Mi_1(F)$ is continuous. To see this, it
suffices to observe that $\mu_n\Rightarrow\mu$ and $z_n\to z$ imply that for
any continuous $h:F\to\R$,
\be\ba{l}
\dis\Big|\int\psi^f_z(\mu)(\di y)h(y)
-\int\psi^f_{z_n}(\mu_n)(\di y)h(y)\Big|
=\Big|\int\mu(\di x)h(f_z(x))-\int\mu_n(\di x)h(f_{z_n}(x))\Big|\\[5pt]
\dis\quad=\Big|\int\mu(\di x)h(f_z(x))-\int\mu_n(\di x)h(f_z(x))\Big|
+\Big|\int\mu_n(\di x)h(f_z(x))-\int\mu_n(\di x)h(f_{z_n}(x))\Big|.
\ec
Here the first term on the right-hand side converges to zero by our assumption
that $\mu_n$ converges weakly to $\mu$, while the second term can be bounded
by $\|h\circ f_{z_n}-h\circ f_z\|_\infty$, which tends to zero since $E\times
G\ni(x,z)\mapsto h\circ f_z(x)\in\R$ is continuous and $E,G$ are compact
spaces.

We next claim that if $f^n\to f$ pointwise, then also $\psi^{f^n}\to\psi^f$
pointwise. Indeed, if $f^n_z(x)\to f_z(x)$ for all $x,z$, then, for
any continuous (and hence bounded) $h:F\to\R$,
\be
\int\psi^{f^n}_z(\mu)(\di y)h(y)
=\int\mu(\di x)h(f^n_z(x))\asto{n}\int\mu(\di x)h(f_z(x))
=\int\psi^f_z(\mu)(\di y)h(y),
\ee
showing that $\psi^{f^n}_z(\mu)\Asto{n}\psi^f_z(\mu)$ for all $\mu,z$. It
follows that the set $\Gi$ of all $E\times G\ni(x,z)\mapsto f_z(x)\in F$ such
that $\Mi_1(E)\times G\ni(\mu,z)\mapsto\psi^f_z(\mu)\in\Mi_1(F)$ is measurable
is closed under pointwise limits and contains all continuous functions
$(x,z)\mapsto f_x(x)$. By Lemma~\ref{L:pclosure}, it follows that $\Gi$
contains all measurable $(x,z)\mapsto f_z(x)$.

To treat the general case, where $E,G$ need not be compact and $F$ may be
different from $[0,1]^\N$, we will use a compactification argument. We need
the following three facts: 1.\ Each separable metrizable space is isomorphic
to a subset of $[0,1]^\N$. 2.\ A subset of a Polish space is Polish in the
induced topology if and only if it is a $G_\de$-set, i.e., a countable
intersection of open sets \cite[\S 6 No.~1, Thm.~1]{Bou58}. 3.\ If $E_1\sub
E_2$ are Polish spaces and $\Mi_1(E_i)$ is the space of probability measures
on $E_i$ $(i=1,2)$, equipped with the topology of weak convergence, then
$\Mi_1(E_1)$ is isomorphic to the set $\{\mu\in\Mi_1(E_2):\mu(E_1)=1\}$. (The
fact that the topology on $\Mi_1(E_1)$ coincides with the one induced by the
embedding in $\Mi_1(E_2)$ follows, for example, from Skorohod's representation
theorem \cite[Theorem~6.7]{Bil99}.) Note that facts~2 and 3 and the fact that
$\Mi_1(E_i)$ $(i=1,2)$ are Polish spaces imply that $\Mi_1(E_1)$ is a
$G_\de$-subset of $\Mi_1(E_2)$.

In view of facts 1 and 2 above, we may without loss of generality assume that
$E,G$ are $G_\de$-subsets of some compact metrizable spaces $\ov E,\ov G$ and
that $F$ is a $G_\de$-subset of $[0,1]^\N$. Then each measurable function
$E\times G\ni(x,z)\mapsto f_z(x)\in F$ may be extended to a measurable
function from $\ov E\times\ov G$ to $[0,1]^\N$ by setting $f_z(x)$ equal to
some constant if $(x,z)\not\in E\times G$. By what we have already proved, the
associated map $\Mi_1(\ov E)\times\ov
G\ni(\mu,z)\mapsto\psi^f_z(\mu)\in\Mi_1([0,1]^\N)$ is measurable. Since
$\Mi_1(E)$ and $G$ are measurable subsets of $\Mi_1(\ov E)$ and $\ov G$,
respectively, the restriction of the map $(\mu,z)\mapsto\psi^f_z(\mu)$ to
$\Mi_1(E)\times G$ yields a measurable map from $\Mi_1(E)\times G$ to $F$.\qed

\section{Thinning and Poissonization}\label{A:thin}

Let $E$ be a Polish space and let $\Mi(E)$ be the space of finite measures on
$E$ equipped with the topology of weak convergence, under which it is
Polish. We let $\Mi_{\rm count}(E)$ denote the space of finite counting
measures on $E$, i.e., measures of the form $\sum_{i=1}^n\de_{x_i}$ with
$n\geq 0$ and $x_1,\ldots,x_n\in E$. Since $\Mi_{\rm count}(E)$ is a closed
subset of $\Mi(E)$, it is also Polish (under the topology of weak
convergence). We let $B_+(E)$ denote the space of measurable functions
$f:E\to\half$ and write $B_{[0,1]}(E)$ for the space of measurable functions
$f:E\to[0,1]$. For any $f\in B_{[0,1]}(E)$ and $\nu\in\Mi_{\rm count}(E)$, we
introduce the notation
\be
f^\nu:=\prod_{i=1}^nf(x_i)\quad\mbox{where}\quad\nu=\sum_{i=1}^n\de_{x_i},
\ee
with the convention that $f^0:=1$. Let $\mu\in\Mi(E)$. By definition, a {\em
  Poisson point measure} with {\em intensity} $\mu$ is an $\Mi_{\rm
  count}(E)$-valued random variable $\nu$ such that
\be\label{Poisdef}
\E[(1-f)^\nu]=\ex{-\int\!f\,\di\mu}\qquad\big(f\in B_{[0,1]}(E)\big).
\ee
An explicit way to construct such a Poisson point measure is to write $\mu=\la\mu'$ where $\la\geq 0$ and $\mu'$ is a probability measure, and to put $\nu=\sum_{i=1}^N\de_{X_i}$ where $(X_i)_{i\geq 1}$ are i.i.d.\ with law $\mu'$ and $N$ is an independent Poisson distributed random variable with mean $\la$. By \cite[Prop.~3.5]{Res87}, the law of $\nu$ is uniquely characterized by (\ref{Poisdef}). The proof there is stated for locally compact spaces only, which in the present paper is actually all we need, but the statement holds more generally for Polish spaces. If $\mu$ is non-atomic, then $\nu$ a.s.\ contains no double points, i.e.,
\be
\nu=\sum_{x\in\supp(\nu)}\de_x\quad{\rm a.s.},
\ee
see \cite[Prop.~10.4]{Kal02}. In this case, we call $\supp(\nu)$ a {\em Poisson point set} with {\em intensity} $\mu$.

If $\nu\in\Mi_{\rm count}(E)$ is a
(deterministic) finite counting measure and $g\in B_{[0,1]}(E)$, then by
definition a {\em $g$-thinning}\index{thinning} of $\nu$ is an $\Mi_{\rm count}(E)$-valued
random variable $\nu$ such that
\be\label{thindef}
\E[(1-f)^{\nu'}]=(1-gf)^\nu\qquad\big(f\in B_{[0,1]}(E)\big).
\ee
An explicit way to construct such a {\em $g$-thinning}, when
$\nu=\sum_{i=1}^n\de_{x_i}$, is to construct independent $\{0,1\}$-valued
random variables $\chi_1,\ldots,\chi_n$ with $\P[\chi_i=1]=g(x_i)$ and to put
$\nu':=\sum_{i=1}^n\chi_i\de_{x_i}$. By \cite[Prop.~3.5]{Res87}, the law of $\nu'$ is uniquely
characterized by (\ref{thindef}).

It is easy to see that the class of functions $f:E\to[0,1]$ for which
(\ref{Poisdef}) or (\ref{thindef}) hold is closed under taking of pointwise
limits. Therefore, by Lemma~\ref{L:pclosunit}, in order to check
(\ref{Poisdef}) or (\ref{thindef}), it suffices to verify the relation for all
continuous functions $f:E\to[0,1]$, and in case $E$ is locally compact, even
the continuous functions with compact support suffice.

We also need Poisson point sets with \si-finite, but in general locally
infinite intensities. To this aim, let $\Count(E)$ be the space of all
countable subsets of $E$. We equip $\Count(E)$ with the \si-field generated by
all mappings $A\mapsto1_{\{A\cap B=\emptyset\}}$ where $B\sub E$ is Borel
measurable.

\bl{\bf(Poisson point sets with \si-finite intensity)}\label{L:siPois}
For each $f\in B_{[0,1]}(E)$, the map $\Count(E)\ni A\mapsto\prod_{x\in A}\big(1-f(x)\big)\in[0,1]$ is measurable. Moreover, for each \si-finite non-atomic measure $\mu$ on $E$, there exists a $\Count(E)$-valued random variable $C$, unique in law, such that
\be\label{Poisdef2}
\E\big[\prod_{x\in C}\big(1-f(x)\big)\big]=\ex{-\int\!f\,\di\mu}
\qquad\big(f\in B_{[0,1]}(E)\big),
\ee
where $e^{-\infty}:=0$.
\el
{\bf Proof.} We claim that for all (Borel) measurable $B\sub E$, the function $\Count(E)\ni A\mapsto|A\cap B|\in\{0,1,\ldots\}\cup\{\infty\}$ is measurable. To see this, let $\Di\sub E$ be countable and dense and let $\Oi:=\{B_{1/k}(x):x\in\Di,\ k\geq 1\}$, where $B_\eps(x)$ denotes the open ball of radius $\eps$ around $x$. Then
\be
\{|A\cap B|\geq n\}=\big\{\exists U_1,\ldots,U_n\in\Oi\mbox{ disjoint, s.t.\ }A\cap B\cap U_i\neq\emptyset\ \forall i=1,\ldots,n\big\}
\ee
is a countable union of finite intersections of measurable sets, and hence itself measurable. It follows that $A\mapsto\sum_{x\in A}f(x)$ is measurable for each $f$ of the form $f=\sum_{i=1}^nb_i1_{B_i}$ with $B_i$ (Borel) measurable and $b_i\in[0,\infty)$. By taking increasing limits it follows that $A\mapsto\sum_{x\in A}f(x)$ is measurable for each measurable $f:E\to[0,\infty]$. Since $\prod_{x\in A}\big(1-f(x)\big)=\exp\{\sum_{x\in A}\log(1-f(x))\}$, we conclude that $A\mapsto\prod_{x\in A}\big(1-f(x)\big)$ is measurable for each $f\in B_{[0,1]}(E)$.

Since $\mu$ is \si-finite, there exist disjoint measurable $B_i\sub E$ such that $\mu(B_i)<\infty$ $(i\geq 1)$. Let $C_i$ be independent Poisson point sets with intensity $\mu_i:=\mu(B_i\cap\,\cdot\,)$ $(i\geq 1)$ and set $C:=\bigcup_{i\geq 1}C_i$. Then $\{C\cap B=\emptyset\}=\bigcup_{i\geq 1}\{C_i\cap B=\emptyset\}$ is measurable for all measurable $B\sub E$, hence $C$ is a measurable $\Count(E)$-valued random variable. Since the $C_i$ are disjoint and independent and the $\mu_i$ are non-atomic, we have
\be
\E\big[\prod_{x\in C}\big(1-f(x)\big)\big]
=\prod_{i\geq 1}\E\big[\prod_{x\in C_i}\big(1-f(x)\big)\big]
=\prod_{i\geq 1}\ex{-\int\!f\,\di\mu_i}
=\ex{-\int\!f\,\di\mu}
\qquad\big(f\in B_{[0,1]}(E)\big).
\ee
In particular, setting $f=1_B$ we see that $\P[C\cap B=\emptyset]=e^{-\mu(B)}$ for all measurable $B\sub E$. Set $\Ai_B:=\{A:A\cap B=\emptyset\}$. Then $\Ai_B\cap\Ai_{B'}=\Ai_{B\cup B'}$, $\Ai_\emptyset=\Om$, and the class of all $\Ai_B$ with $B\sub E$ measurable generates the \si-field on $\Count(E)$, hence by the $\pi/\la$-theorem \cite[Theorem~1.1]{Kal02},
(\ref{Poisdef2}) uniquely determines the law of $C$.\qed

\section{A one-sided version of Kolmogorov's moment criterion}

\newcommand{\Nomeg}{N}

We prove a variant of Kolmogorov's moment criteria (see e.g.\ \cite[Chap.~7, Theorem (1.5)]{Dur96}) for the H\"older continuity of a stochastic process, with bounds on the distribution of the H\"older constant. We assume a one-sided moment condition, which in turn gives one-sided H\"older continuity at deterministic times.

\bt\label{T:moment} Let $(X_t)_{t\in [0,T]}$ be a real-valued stochastic process. If for all $0\leq s<t\leq T$,
\be\label{moment}
\E\big[((X_s-X_t)^+)^\beta\big] \leq K (t-s)^{1+\alpha}
\ee
for some $\alpha, \beta>0$ and $K<\infty$, then for any $0<\gamma <\frac{\alpha}{\beta}$, there exists a random constant $C\in (0,\infty)$ such that a.s.\
\be\label{holder}
(X_r-X_q)^+ \leq C(q-r)^\gamma \qquad \mbox{ for all }  r,q \in \Q_2 \cap [0,T] \mbox{ with } r<q,
\ee
where $\Q_2=\{m 2^{-n}: m,n \geq 0\}$ is the set of dyadic rationals. Furthermore, for any $0<\delta<\alpha-\beta\gamma$, there exists a deterministic constant
$C_{\delta,\gamma}$ depending only on $\gamma,\delta, K,\alpha$ and $\beta$, such that
\be\label{Comegatail}
\P[C \geq u] \leq \frac{C_{\delta,\gamma}}{u^{\delta}} \wedge 1 \qquad \mbox{for all } u>0.
\ee
The same results hold if we replace $(\cdot)^+$ by $(\cdot)^-:=-(\cdot \wedge 0)$ or $|\cdot|$.
\et
{\bf Proof.} The proof is essentially the same as that for the standard version of Kolmogorov's moment criterion.
Fix $0<\gamma <\frac{\alpha}{\beta}$. Without loss of generality, assume $T=1$ and let $D_n:= \{i 2^{-n}: 0\leq i\leq 2^n\}$.
For any $s:=i2^{-n}< t:=j2^{-n}\in D_n$, by (\ref{moment}) and the Chebychev inequality,
\be\label{hol2}
\P[(X_s-X_t)^+ > (t-s)^\gamma ] \leq K (t-s)^{1+\alpha-\beta\gamma} = K (j-i)^{1+\alpha-\beta\gamma} 2^{-n(1+\alpha-\beta\gamma)}.
\ee
If we let $G_n:=\{ (X_{i2^{-n}}-X_{j2^{-n}})^+ \leq (j-i)^{\gamma}2^{-n\gamma} \ \mbox{for all } 0\leq i\leq j\leq 2^n, j-i\leq 2^{n\eta}\}$ for
some fixed $\eta\in (0,1)$, then by (\ref{hol2}),
\be
\P[G_n^c] \leq \sum_{0\leq i< j\leq 2^n \atop j-i\leq 2^{n\eta}} K (j-i)^{1+\alpha-\beta\gamma} 2^{-n(1+\alpha-\beta\gamma)}
\leq K 2^{-[(1-\eta)(\alpha-\beta\gamma)-2\eta]n}=K 2^{-\delta n},
\ee
where we have chosen $\eta>0$ such that $(1-\eta)(\alpha-\beta\gamma)-2\eta=\delta\in (0,\alpha-\beta\gamma)$. Then by Borel-Cantelli, a.s.\ $\Nomeg := \inf\{ n\in\N: \cap_{i\geq \Nomeg} G_i \mbox{ occurs}\}<\infty$, and furthermore,
for any $L\in\N$,
\be\label{Nomegatail}
\P[\Nomeg > L] \leq \sum_{n=L}^\infty \P[G_n^c]\leq K \sum_{n=L}^\infty 2^{-\delta n} = \frac{K2^{-\delta L}}{1-2^{-\delta}}.
\ee

Note that
\be\label{triangle}
(X_u-X_w)^+ \leq (X_v-X_w)^+ +(X_u-X_v)^+ \qquad \mbox{for any}\quad u<v<w.
\ee
We will use this triangle inequality to deduce (\ref{holder}) on the event $\cap_{n\geq \Nomeg} G_n$.  First assume that
$r<q\in \Q_2\cap [0,1]$ and $q-r < 2^{-\Nomeg(1-\eta)}$. We can find an $m\geq \Nomeg$ such that
\be\label{hol3}
2^{-(m+1)(1-\eta)} \leq q-r < 2^{-m(1-\eta)}.
\ee
By binary expansion for $q$ and $r$, we can write
$$
\begin{aligned}
q & = & j 2^{-m} + 2^{-q_1}+\cdots 2^{-q_k} \\
r & = & i 2^{-m} - 2^{-r_1} -\cdots 2^{-r_l},
\end{aligned}
$$
where $m<q_1<\cdots <q_k$ and $m<r_1<\cdots <r_l$. By (\ref{hol3}),
$$
2^{(m+1)\eta-1} \leq (q-r)2^m \leq j-i+2.
$$
Since $m\geq \Nomeg$, if we replace $\Nomeg$ with $\Nomeg \vee 2/\eta$, then we are guaranteed that $j\geq i$.
Since $q-r\geq (j-i)2^{-m}$, again by (\ref{hol3}), we have $j-i\leq 2^{m\eta}$. Since the event $\cap_{n\geq \Nomeg} G_n$ occurs
by definition, we have
\be\label{hol4}
(X_{i2^{-m}}-X_{j2^{-m}})^+ \leq (j-i)^\gamma 2^{-m\gamma} \leq 2^{-m(1-\eta)\gamma}.
\ee
By (\ref{triangle}),
\be
(X_{j2^{-m}}-X_q)^+ \leq \sum_{\sigma=1}^k 2^{-q_\sigma\gamma} \leq \sum_{\sigma> m} 2^{-\sigma \gamma} \leq \frac{2^{-m\gamma}}{2^{\gamma}-1}.
\ee
Similarly, the same bound also holds for $(X_r-X_{i2^{-m}})^+$. Combining the above estimates and applying (\ref{triangle}) once more,
we get
$$
(X_r-X_q)^+ \leq 2^{-m(1-\eta)\gamma} + \frac{2^{1-m\gamma}}{2^\gamma-1} = \big(1+\frac{2^{1-m\eta\gamma}}{2^\gamma-1}\big)2^{(1-\eta)\gamma} 2^{-(m+1)(1-\eta)\gamma} \leq C_\gamma (q-r)^\gamma
$$
for $C_\gamma = 2^\gamma(1+\frac{2}{2^\gamma-1})$. This verifies (\ref{holder}) for $r<q\in Q_2\cap [0,1]$ with $q-r<2^{-(\Nomeg\vee 2/\eta)(1-\eta)}$. For general $r<q\in \Q_2\cap [0,1]$, we can apply the triangle inequality (\ref{triangle}) at most $2^{(\Nomeg\vee 2/\eta)(1-\eta)}$ times to obtain (\ref{holder}) with $C=C_\gamma 2^{(\Nomeg\vee 2/\eta)(1-\eta)}$. The distributional tail bound (\ref{Comegatail}) then follows from (\ref{Nomegatail}).

When we replace $(\cdot)^+$ by $(\cdot)^-$ or $|\cdot|$, the proof is the same since analogues of the triangle inequality
(\ref{triangle}) still hold.
\qed
\bigskip

\noi
{\bf Acknowledgement} We thank Jon Warren for helpful discussions on the
Howitt-Warren martingale problem and an unknown referee for taking on the
daunting task of refereeing such a long paper, and for providing helpful
comments.

\vspace{2cm}

\vfill

\parbox[t]{4.7cm}{\small
Emmanuel Schertzer\\
LPMA, Universit\'e  Pierre \\
et Marie Curie, France \\
email: \\
emmanuel.schertzer@gmail.com}
\hspace{.3cm}
\parbox[t]{4.7cm}{\small
Rongfeng~Sun\\
Department of Mathematics\\
National University of\\ Singapore\\
10 Lower Kent Ridge Road\\
119076, Singapore\\
e-mail: matsr@nus.edu.sg}
\hspace{.3cm}
\parbox[t]{4.7cm}{\small
Jan M.~Swart\\
Institute of Information\\ Theory and Automation\\ of
the ASCR (\' UTIA)\\
Pod vod\'arenskou v\v e\v z\' i 4\\
18208 Praha 8\\
Czech Republic\\
e-mail: swart@utia.cas.cz}

\newpage

\noi
{\bf \Large Notation List}
\bigskip

\noi
{\bf General notation:}
\med

\begin{tabular}{@{\hspace{-10pt}}r@{\ }l}
$\Z^2_{\rm even}$\;: & The even sublattice of $\Z^2$, $\{(x,t)\in\Z^2: x+t \text{ is even}\}$. \hspace{4cm}\\
$\Li(\,\cdot\,)$\;:
& the law of a random variable.\\
$S_{\eps}$\;: & the diffusive scaling map, applied to subsets of $\R^2$, paths, and sets of\\
& paths, quenched laws, etc. See (\ref{Seps2}). \\
$(\eps_k)_{k\in\N}$\;: & a sequence of constants decreasing to $0$, acting as scaling parameters. \\
\hspace{-20pt}
$(K_{s,t}(x,\cdot))_{s\leq t, x\in E}$\;: & a stochastic flow of kernels on the space $E$. \\
$\Mi_1(E)$\;: & the space of probability measures on the space $E$. \\
$\Mi(\R), \Mi_{\rm loc}(\R)$\;: & the space of finite and locally finite measures on $\R$. \\
$\Mi_{\rm g}(\R)$\;: & the subset of $\Mi_{\rm loc}(\R)$ satisfying the growth constraint (\ref{MiG}). \\
${\rm supp}(\cdot)$\;: & support of a measure.
\end{tabular}
\bigskip

\noi
{\bf Paths, Space of Paths:}
\med

\begin{tabular}{@{\hspace{-10pt}}r@{\ }l}
\hspace{60pt}
$\Rc$\;: & the compactification of $\R^2$, see Figure \ref{fig:compac}.\\
$z=(x,t)$\;: & a point in $\Rc$, with
  position $x$ and time $t$.\\
$s,t,u,S,T,U$\;: & times.\\
$\Ki(\Rc)$\;: & the space of compact subsets of $\Rc$. \\
$(\Pi,d)$\;: & the space of continuous paths in $\Rc$, with metric $d$, see (\ref{dPi}).\\
$(\Ki(\Pi), d_{\rm H})$\;: & the space of compact subsets of $\Pi$, with Hausdorff metric $d_{\rm H}$, see (\ref{Hmet}). \\
$\Pi(A),\Pi(z)$\;:
& the set of paths in $\Pi$ starting from a set $A\sub\Rc$
resp.\ a point $z\in\Rc$.\\
& The same notatation applies to any subset of $\Pi$ such as
$\Wi,\Ni$, etc.\\
$\pi$\;: & a path in $\Pi$.\\
$\sigma_\pi$\;: & the starting time of the path $\pi$\\
$\pi(t)$\;: & the position of $\pi$ at time $t\geq\sigma_\pi$.\\
$(\hat\Pi,\hat d)$\;: & the space of continuous backward paths in $\Rc$
  with metric $\hat d$.\\
$\hat\Pi(A),\hat\Pi(z)$\;:
& the set of backward paths in $\hat\Pi$ starting from
$A\sub\Rc$ resp.\ $z\in\Rc$.\\
$\hat\pi$\;: & a path in $\hat\Pi$.\\
$\hat\sigma_{\hat\pi}$\;: & the starting time of the backward path $\hat\pi$.\\
$\sim^z_{\rm in}, \sim^z_{\rm out}$\;: & equivalence of paths entering,
resp.\ leaving $z\in\R^2$, see Definition \ref{D:equiv}.\\
$=^z_{\rm in}, =^z_{\rm out}$\;: & strong equivalence of paths entering,
resp.\ leaving $z\in\R^2$, see \\
& Definition \ref{D:inout}.
\end{tabular}

\newpage

\noindent
{\bf Discrete environments, paths, webs, and flows:}
\med

\begin{tabular}{@{\hspace{-23pt}}r@{\ }l}
\hspace{4pt}
$\omega:=(\omega_z)_{z\in\Z^2_{\rm even}}$\;: & an i.i.d.\ environment for random walks on $\Z^2_{\rm even}$. \\
$\mu$\;: & the law of $\omega_o\in[0,1]$. \\
$K^\omega_{(s,t)}(x,y)$\;: & transition probability of a random walk from $(x,s)$ to $(y,t)\in\Z^2_{\rm even}$ in \\
& the environment $\omega$.\\
$\Qdis^\omega_z$\;: & the law of a random walk on $\Z^2_{\rm even}$ starting from $z$ in the environment $\omega$. \\
$\Qdis^\omega$\;: & discrete Howitt-Warren quenched law.\\
$p_z$\;: & a discrete path on $\Z^2_{\rm even}$ starting from $z\in\Z^2_{\rm even}$. \\
$\hat p_z$\;: & a dual discrete path on $\Z^2_{\rm odd}:=\Z^2\backslash\Z^2_{\rm even}$ starting from $z\in\Z^2_{\rm odd}$.\\
$\alpha:=(\alpha_z)_{z\in\Z^2_{\rm even}}$\;: & a family of independent $\pm1$-valued random variables. \\
$(\Ui, \hat\Ui)$\;: &  a discrete web and its dual. \\
$(\Ui^{\rm l}, \Ui^{\rm r}, \hat\Ui^{\rm l}, \hat\Ui^{\rm r})$\;: & a discrete left-right web and its dual. \\
$(\Vi, \hat\Vi)$\;: & a discrete net and its dual. \\
$(\omega^{\langle k\rangle})_{k\in\N}$\;: & a sequence of i.i.d.\ environments on $\Z^2_{\rm even}$, with $\Li(\omega^{\langle k\rangle}_o)=\mu_k$.\\
$(\mu_k)_{k\in\N}$\;: & a sequence of probability laws on $[0,1]$, satisfying (\ref{mucon}). \\
$\Qdis_{\langle k\rangle}$\;: & the discrete Howitt-Warren quenched law associated with $\omega^{\langle k\rangle}$.
\end{tabular}
\bigskip

\noindent
{\bf Brownian webs:}
\med

\begin{tabular}{@{\hspace{-12pt}}r@{\ }l}
$(\Wi,\hat\Wi)$\;: & a double Brownian web consisting of a Brownian web and its dual.\\
$\pi_z,\hat\pi_z$\;: & the a.s.\ unique path in $\Wi$ resp.\ $\hat\Wi$ starting from a deterministic \\
& $z\in\R^2$.\\
$\pi^-_z, \pi^+_z$\;: & the leftmost, resp.\ rightmost path starting from $z\in\R^2$ in the \\
& Brownian web $\Wi$. \\
$\pi^\up_z$\;: & same as $\pi^+_z$, except when there is an incoming path in $\Wi$ at $z$, then $\pi^\up_z$ \\
& is defined to be the continuation of the incoming path. \\
$\Wi_{\rm in}(z), \Wi_{\rm out}(z)$\;: & the set of paths in $\Wi$ entering, resp.\ leaving $z$.\\
${\rm sign}_\Wi(z)$\;: & the orientation of a $(1,2)$ point $z\in\R^2$ in $\Wi$. See (\ref{signz}).\\
${\rm switch}_z(\Wi)$\;: & a modification of $\Wi$ by switching the orientation of all paths in $\Wi$ \\
& entering $z$. See (\ref{switch}).\\
${\rm hop}_z(\Wi)$\;: & $\Wi\cup {\rm switch}_z(\Wi)$.\\
$\ell, \ell_{\rm l}, \ell_{\rm r}$\;: & the intersection local time measure between $\Wi$ and $\hat\Wi$, and its \\
& restriction to the set of $(1,2)_{\rm l}$, resp.\ $(1,2)_{\rm r}$ points. See Proposition~\ref{P:refloc}.\\
$(\Wl,\Wr, \hat\Wl, \hat\Wr)$\;: & the left-right Brownian web and its dual. \\
$l,r,\hat l,\hat r$\;: & a path in $\Wl$, resp. $\Wr$, $\hat\Wl$, $\hat\Wr$.\\
$l_z,r_z,\hat l_z,\hat r_z$\;: & the a.s.\ unique
path in $\Wl(z)$, resp. $\Wr(z)$, $\hat\Wl(z)$, $\hat\Wr(z)$, starting \\
& from a deterministic $z\in\R^2$. \\
$W(\hat r, \hat l)$\;: & a wedge defined by the paths $\hat r\in\hat\Wr$ and
$\hat l\in\hat\Wl$, see (\ref{wedge}).\\
$M(r, l)$\;: & a mesh defined by the paths $r\in\Wr$ and $l\in\Wl$, see
(\ref{mesh}).
\end{tabular}

\newpage

\noi
{\bf Brownian net:}
\med

\begin{tabular}{@{\hspace{-10pt}}r@{\ }l}
\hspace{20pt}
$(\Ni,\hat\Ni)$\;: & the Brownian net and the dual Brownian net. See Theorem \ref{T:net}.\\
$\xi_t$\;: & a branching-coalescing point set. See (\ref{braco}) and Proposition~\ref{P:braco}.\\
$S$\;: & the set of separation points of $\Ni$.\\
${\rm sign}_\pi(z)$\;: & the orientation of a path $\pi\in\Ni$ at the separation point $z$. See (\ref{signznet}).\\
$\Hi_-$\;: & a Brownian half-net with infinite left speed and finite right speed. \\
$\Hi_+$\;: & a Brownian half-net with finite left speed and infinite right speed. \\
$R_{S,U}$\;: & the set of $S,U$-relevant separation points in $\Ni$. \\
\end{tabular}
\bigskip

\noindent
{\bf Howitt-Warren flows and processes:}
\med

\begin{tabular}{@{\hspace{-10pt}}r@{\ }l}
$\beta, \nu$\;: & the drift and characteristic measure of a Howitt-Warren flow. See \\
& Definition~\ref{D:HWMP2}.\\
$\beta_-, \beta_+$\;: & the left and right speeds of a Howitt-Warren flow. See (\ref{speeds}). \\
$(\Wi_0, \Mi, \Wi)$\;: & the reference web, the set of marked points, and the sample web. See \\
& Section~\ref{S:flowcons}.\\
$(\Wi_i)_{i\in\N}$\;: & i.i.d.\ copies of the sample web $\Wi$.\\
$\beta_0, \beta$\;: & the drift of the reference, resp.\ sample web.\\
$\nu_{\rm l}, \nu_{\rm r}$\;: & a decomposition of $\nu$ via (\ref{nuform}). \\
$\pi^+_z$\;: & the rightmost path starting from $z\in\R^2$ in the sample web $\Wi$. \\
$\pi^\up_z$\;: & same as $\pi^+_z$, except when there is an incoming path in $\Wi$ at $z$, then $\pi^\up_z$ \\
& is defined to be the continuation of the incoming path. \\
$(K^+_{s,t})_{s\leq t}$\;: & the version of the Howitt-Warren flow constructed using $\pi^+$. See (\ref{HWconst}). \\
$(K^\uparrow_{(s,t)})_{s\leq t}$\;: & the version of the Howitt-Warren flow constructed using $\pi^\up$.\\
$\Q$\;: & the Howitt-Warren quenched law of $\Wi$ conditional on $(\Wi_0,\Mi)$. \\
$\rho_t$\;: & the Howitt-Warren process defined from either $K^+$ or $K^\up$. See (\ref{rho}). \\
$(\Lambda_c)_{c\geq 0}$\;: & ergodic homogeneous invariant laws for the Howitt-Warren process. See
\\ & Theorem~\ref{T:HIL}. \\
$\zeta_t$\;: & the smoothing process dual to the Howitt-Warren process $\rho_t$. See (\ref{zetat}).
\end{tabular}
\med

\printindex


\bigskip

\begin{thebibliography}{FINRS07}

\bibitem[Arr79]{Arr79}
R.~Arratia.
Coalescing Brownian motions on the line.
Ph.D. Thesis, University of Wisconsin, Madison, 1979.

\bibitem[Arr81]{Arr81}
R.~Arratia.
Coalescing Brownian motions and the voter model on $\Z$.
Unpublished partial manuscript. Available from rarratia@math.usc.edu.

\bibitem[Bar05]{Bar05}
D.\ Barbato.
FKG inequality for Brownian motion and stochastic differential equations.
{\em Electron.\ Comm.\ Probab.}~10, 7--16, 2005.


\bibitem[Bil99]{Bil99}
P.~Billingsley.
{\em Convergence of probability measures}, 2nd edition.
John Wiley \& Sons, 1999.

\bibitem[Bou58]{Bou58}
N.~Bourbaki.
{\em \'El\'ements de Math\'ematique. VIII. Part. 1: Les Structures Fondamentales de l'Analyse. Livre III: Topologie G\'en\'erale. Chap. 9: Utilisation des Nombres R\'eels en Topologie G\'en\'erale. 2i\'eme \'ed.}
Actualit\'es Scientifiques et Industrielles~1045. Hermann \& Cie, Paris, 1958.

\bibitem[BY81]{BY81}
M.T.~Barlow and M.~Yor.
(Semi-) martingale inequalities and local times.
{\em Z.\ Wahrscheinlichkeitstheor.\ Verw.\ Geb.}~55, 237--254, 1981.


\bibitem[Daw91]{Daw91}
D.A.~Dawson.
{\em Measure-valued Markov processes.}
Lecture Notes in Math.~1541, 1--260, Springer, Berlin,  1993.

\bibitem[Dur96]{Dur96}
R.~Durrett.
{\em Probability: Theory and Examples}, 2nd edition, Duxbury Press, 1996.

\bibitem[EK86]{EK86}
S.N.~Ethier and T.G.~Kurtz.
{\em Markov Processes: Characterization and Convergence.}
John Wiley \& Sons, New York, 1986.

\bibitem[EMS13]{EMS13}
S.N.~Evans, B.~Morris and A.~Sen.
Coalescing systems of non-Brownian particles.
{\em Probab.\ Theory Related Fields}~156, 307--342, 2013

\bibitem[Fel66]{Fel66}
W.~Feller.
{\em An introduction to probability theory and its applications. Vol.~II.}
John Wiley \& Sons, Inc.,
New York-London-Sydney, 1966.


\bibitem[FIN02]{FIN02}
L.R.G.~Fontes, M.~Isopi, and C.M.~Newman.
Random walks with strongly inhomogeneous rates and singular diffusions:
Convergence, localization and aging in one dimension.
{\em Ann.\ Probab.}~30(2), 579--604, 2002.

\bibitem[FINR02]{FINR02}
L.R.G.~Fontes, M.~Isopi, C.M.~Newman, and K.~Ravishankar.
The Brownian web.
{\em Proc.\ Natl.\ Acad.\ Sci.\ USA}~99, no.~25, 15888--15893, 2002.

\bibitem[FINR04]{FINR04}
L.R.G.~Fontes, M.~Isopi, C.M.~Newman, and K.~Ravishankar.
The Brownian web: characterization and convergence.
{\em Ann.\ Probab.}~32(4), 2857--2883, 2004.

\bibitem[FINR06]{FINR06}
L.R.G.~Fontes, M.~Isopi, C.M.~Newman, and K.~Ravishankar.
Coarsening, nucleation, and the marked Brownian web.
{\em Ann. Inst. H. Poincar\'e Probab. Statist.} 42, 37--60, 2006.





\bibitem[Hen99]{Hen99}
J.~Henrikson.
Completeness and total boundedness of the Hausdorff metric.
{\em MIT Undergraduate Journal of Mathematics}~1 Number 1, 69--79, 1999.

\bibitem[Hoe63]{Hoe63}
W.~Hoeffding.
Probability inequalities for sums of bounded random variables.
{\em Journal of the American Statistical Association}~58, 13--30, 1963.

\bibitem[HW09a]{HW09a}
C.~Howitt and J.~Warren.
Consistent families of Brownian motions and stochastic flows of kernels.
{\em Ann.\ Probab.}~37, 1237--1272, 2009.

\bibitem[HW09b]{HW09b}
C.~Howitt and J.~Warren.
Dynamics for the Brownian web and the erosion flow.
{\em Stochastic Processes Appl.}~119, 2028--2051, 2009.

\bibitem[JM11]{JM11}
A.St.~John and H.~Mathur.
Correlations and critical behavior of the $q$ model.
{\em Phys.\ Rev.\ E}~84, 051303, 2011.

\bibitem[Kal02]{Kal02}
O.\ Kallenberg.
{\em Foundations of modern probability.} 2nd ed.
Springer, New York, 2002.

\bibitem[KS91]{KS91}
I.~Karatzas and S.E.~Shreve.
{\em Brownian Motion and Stochastic Calculus},
2nd edition, Springer-Verlag, New York, 1991.


\bibitem[Kur98]{Kur98}
T.G.~Kurtz.
Martingale problems for conditional distributions of Markov processes.
{\em Electronic J.\ Probab.}~3, Paper no.~9, 1--29, 1998.

\bibitem[LL04]{LL04}
Y.~Le Jan and S.~Lemaire.
Products of beta matrices and sticky flows.
{\em Probab.\ Th.\ Relat.\ Fields}~130, 109--134, 2004.

\bibitem[LR04a]{LR04AOP}
Y.~Le Jan and O.~Raimond.
Flows, Coalecence and Noise.
{\em Annals of Probab.}~32, 1247--1315, 2004.

\bibitem[LR04b]{LR04PTRF}
Y.~Le Jan and O.~Raimond.
Sticky flows on the circle and their noises.
{\em Probab.\ Th.\ Relat.\ Fields}~129, 63--82, 2004.

\bibitem[LMY01]{LMY01}
M.~Lewandowska, H.~Mathur, and Y.-K.~Yu.
Dynamics and critical behavior of the $q$ model.
{\em Phys.\ Rev.\ E}~64, 026107, 2001.

\bibitem[Lig73]{Lig73}
T.M.~Liggett.
A characterization of the invariant measures for an infinite particle system with interactions.
{\em Trans.\ Amer.\ Math.\ Soc.}~179, 433--453, 1973.

\bibitem[LS81]{LS81}
T.M.~Liggett and F.~Spitzer.
Ergodic theorems for coupled random walks and other systems with locally interacting components.
{\em Z.\ Wahrsch.\ Verw.\ Gebiete}~56, 443--468, 1981.

\bibitem[Lig05]{Lig05}
T.M.~Liggett.
{\em Interacting particle systems}.
Reprint of the 1985 original. Classics in Mathematics. Springer-Verlag, Berlin, 2005.


\bibitem[MKM78]{MKM78}
K.~Matthes, J.~Kerstan, and J.~Mecke.
{\em Infinitely Divisible Point Processes.}
Wiley, Chichester, 1978.

\bibitem[Mun00]{Mun00}
J.R.~Munkres.
{\em Topology, 2nd ed.}
Prentice Hall, Upper Saddle River, 2000.


\bibitem[NRS10]{NRS10}
C.M.~Newman, K.~Ravishankar, and E.~Schertzer.
Marking $(1,2)$ points of the Brownian web and applications.
{\em Ann.\ Inst.\ Henri Poincar\'e Probab.\ Statist.}~46, 537--574, 2010.



\bibitem[Res87]{Res87}
S.I.~Resnick.
{\em Extreme values, regular variation, and point processes}.
Springer-Verlag, New York, 1987.

\bibitem[RP81]{RP81}
L.C.G.~Rogers and J.W.~Pitman.
Markov functions.
{\em Ann.\ Probab.}~9(4), 573--582, 1981.


\bibitem[SS08]{SS08}
R.~Sun and J.M.~Swart.
The Brownian net.
{\em Ann.\ Probab.}~36(3), 1153-1208, 2008.

\bibitem[SSS09]{SSS09}
E.~Schertzer, R.~Sun, and J.M.~Swart.
Special points of the Brownian net.
{\em Electron.\ J.\ Prob.}~14, Paper 30, 805--864, 2009.

\bibitem[Sto67]{Sto67}
C.~Stone.
On Local and Ratio Limit Theorems,
{\em Proc.~Fifth Berkeley Sympos.~Math.~Statist.~and Probability},
Vol.~II, 217--224,
Univ.~California Press, Berkeley, CA, 1967.

\bibitem[STW00]{STW00}
F.~Soucaliuc, B.~T\'oth, and W.~Werner.
Reflection and coalescence between independent one-dimensional Brownian paths.
{\em Ann.\ Inst.\ Henri Poincar\'e Probab.\ Statist.}~36, 509--536, 2000.

\bibitem[TW98]{TW98}
B.~T\'oth and W.~Werner.
The true self-repelling motion.
{\em Probab.\ Theory Related Fields}~111, 375--452, 1998.




\end{thebibliography}
\end{document}